\documentclass[11pt, reqno]{amsart}

\usepackage[usenames,dvipsnames]{color}

\usepackage[OT2, T1]{fontenc}
\usepackage{url}
\usepackage{amsmath}
\usepackage{array}
\usepackage{graphicx}
\usepackage{amssymb}
\usepackage{amsthm}
\definecolor{link}{RGB}{11,0,128}
\usepackage[colorlinks=true,citecolor=NavyBlue,linkcolor=OliveGreen,urlcolor=link]{hyperref}
\usepackage{hyperref}
\usepackage{paralist}
\usepackage{colonequals}			
\usepackage{color}
\usepackage{mathabx}			
\usepackage{amscd}
\usepackage[all,cmtip]{xy}			
\usepackage{sseq}				
\usepackage{verbatim}
\usepackage{parskip}
\usepackage{microtype}
\usepackage[in]{fullpage}
\usepackage{mathrsfs} 			
\usepackage{cleveref}
\usepackage{enumitem}
\usepackage{moreenum}			
\usepackage[alphabetic,lite]{amsrefs} 	
\usepackage{stmaryrd}			
\usepackage{subfiles}
\usepackage{ragged2e}
\usepackage[bbgreekl]{mathbbol}     
\DeclareSymbolFontAlphabet{\mathbb}{AMSb}   
\usepackage{stackengine}

\newcommand{\gA}{\alpha}
\newcommand{\gB}{\beta}

\newcommand{\bA}{\mathbb{A}}
\newcommand{\bC}{\mathbb{C}}
\newcommand{\bD}{\mathbb{D}}
\newcommand{\bE}{\mathbb{E}}
\newcommand{\bF}{\mathbb{F}}
\newcommand{\bG}{\mathbb{G}}

\newcommand{\bL}{\mathbb{L}}
\newcommand{\bM}{\mathbb{M}}

\newcommand{\bP}{\mathbb{P}}
\newcommand{\bQ}{\mathbb{Q}}
\newcommand{\bR}{\mathbb{R}}

\newcommand{\bZ}{\mathbb{Z}}

\newcommand{\bbB}{\mathbf{B}}

\newcommand{\cD}{\mathcal{D}}

\newcommand{\cH}{\mathcal{H}}

\newcommand{\cO}{\mathcal{O}}
\newcommand{\cP}{\mathcal{P}}

\newcommand{\ff}{\mathfrak{f}}

\newcommand{\fm}{\mathfrak{m}}

\newcommand{\fp}{\mathfrak{p}}
\newcommand{\fq}{\mathfrak{q}}

\newcommand{\fU}{\mathfrak{U}}

\newcommand{\fX}{\mathfrak{X}}

\newcommand{\sA}{\mathscr{A}}
\newcommand{\sB}{\mathscr{B}}
\newcommand{\sC}{\mathscr{C}}
\newcommand{\sD}{\mathscr{D}}
\newcommand{\sE}{\mathscr{E}}
\newcommand{\sF}{\mathscr{F}}

\newcommand{\sH}{\mathscr{H}}

\newcommand{\sL}{\mathscr{L}}

\newcommand{\sO}{\mathscr{O}}

\newcommand{\sS}{\mathscr{S}}
\newcommand{\sT}{\mathscr{T}}
\newcommand{\sU}{\mathscr{U}}

\DeclareMathOperator{\Ani}{Ani}		
\DeclareMathOperator{\Br}{Br}			
\DeclareMathOperator{\Char}{char}		
\DeclareMathOperator{\Coker}{Coker}	
\DeclareMathOperator{\depth}{depth}	
\DeclareMathOperator{\Ext}{Ext}		
\DeclareMathOperator{\Fib}{Fib}		
\DeclareMathOperator{\Frac}{Frac}		
\DeclareMathOperator{\Fun}{Fun}		
\DeclareMathOperator{\Hom}{Hom}		
\DeclareMathOperator{\im}{Im}			
\DeclareMathOperator{\Ind}{Ind}		
\DeclareMathOperator{\Ker}{Ker}		
\DeclareMathOperator{\Mod}{Mod}		
\DeclareMathOperator{\PGL}{PGL}		
\DeclareMathOperator{\Pic}{Pic}		
\DeclareMathOperator{\Proj}{Proj}		
\DeclareSymbolFont{cyrletters}{OT2}{wncyr}{m}{n}
\DeclareMathSymbol{\Sha}{\mathalpha}{cyrletters}{"58}	
\DeclareMathOperator{\Spa}{Spa}		
\DeclareMathOperator{\Spec}{Spec}		
\DeclareMathOperator{\Tor}{Tor}		
\DeclareMathOperator{\vdim}{vdim}		
\DeclareMathOperator{\Vect}{Vect}		

\providecommand{\Ab}{\mathrm{Ab}}
\newcommand{\ad}{\mathrm{ad}}			

\newcommand{\ce}{\colonequals}
\newcommand{\cris}{{\mathrm{cris}}}		

\newcommand{\eps}{\epsilon}
\newcommand{\et}{\mathrm{\acute{e}t}}	
\newcommand{\fet}{\mathrm{f\acute{e}t}}	
\newcommand{\fppf}{\mathrm{fppf}}		
\newcommand{\Frob}{\mathrm{Frob}}			
\providecommand{\Gp}{\mathrm{Gp}}
\newcommand{\hra}{\hookrightarrow}
\DeclareMathOperator{\Idem}{\mathrm{Idem}}	
\renewcommand{\i}{^{-1}}
\newcommand{\id}{\mathrm{id}}			
\newcommand{\Inf}{\mathrm{inf}}		

\newcommand{\invlim}{\varprojlim}
\newcommand{\isomto}{\overset{\sim}{\longrightarrow}}

\newcommand{\llb}{\llbracket}			
\newcommand{\llp}{(\!(}			
\newcommand{\op}{{\mathrm{op}}}			
\newcommand{\ov}{\overline}
\providecommand{\p}[1]{\left(#1\right)}
\newcommand{\perf}{\mathrm{perf}}			

\DeclareSymbolFontAlphabet{\mathbbl}{bbold}
\newcommand{\Prism}{{ \mathbbl{\Delta}}}
\newcommand{\profppf}{\mathrm{pro\x{-}fppf}}		

\newcommand{\ra}{\rightarrow}
\newcommand{\Ra}{\Rightarrow}
\newcommand{\red}{{\mathrm{red}}}		
\providecommand{\Ring}{\mathrm{Ring}}
\newcommand{\rrb}{\rrbracket}			
\newcommand{\rrp}{)\!)}			
\providecommand{\Set}{\mathrm{Set}}
\newcommand{\sfp}{\mathrm{sfp}}		
\newcommand{\sh}{\mathrm{sh}}		

\newcommand{\surjects}{\twoheadrightarrow}
\newcommand{\tensor}{\otimes} 			
\newcommand{\tors}{\mathrm{tors}}		
\newcommand{\un}{\underline}
\newcommand{\wh}{\widehat}
\newcommand{\wt}{\widetilde}
\newcommand{\xra}{\xrightarrow}

\newcommand{\csub[1]}{\centering\subsection{\text{#1} \nopunct} \hfill \justify}

\providecommand{\abs}[1]{\left\lvert#1\right\rvert}

\providecommand{\up}[1]{{\upshape(}#1{\upshape)}}
\providecommand{\resp}[1]{{\upshape(}resp.,~#1{\upshape)}}
\providecommand{\uref}[1]{{\upshape\ref{#1}}}
\providecommand{\uS}{{\upshape\S}}
\providecommand{\f}[2]{\frac{#1}{#2}}

\renewcommand{\b}{\textbf}
\providecommand{\ucolon}{{\upshape:} }
\providecommand{\uscolon}{{\upshape;} }
\providecommand{\uc}{{\upshape:} }

\newcommand{\brems}{\begin{rems} \hfill \begin{enumerate}[label=\b{\thenumberingbase.},ref=\thenumberingbase]}
\newcommand{\remi}{\addtocounter{numberingbase}{1} \item}
\newcommand{\erems}{\end{enumerate} \end{rems}}
\newcommand{\begs}{\begin{egs} \hfill \begin{enumerate}[label=\b{\thenumberingbase.},ref=\thenumberingbase]}

\newcommand{\eegs}{\end{enumerate} \end{egs}}
\newcommand{\m}{\item}
\newcommand{\bsm}{\begin{smallmatrix}}
\newcommand{\esm}{\end{smallmatrix}}
\newcommand{\blem}{\begin{lemma}}
\newcommand{\elem}{\end{lemma}}

\newcommand{\bconj}{\begin{conj}}
\newcommand{\econj}{\end{conj}}
\newcommand{\bprob}{\begin{Problem}}
\newcommand{\eprob}{\end{Problem}}

\newcommand{\bq}{\begin{Q}}
\newcommand{\eq}{\end{Q}}
\newcommand{\benum}{\begin{enumerate}[label={{\upshape(\alph*)}}]}
\newcommand{\benuma}{\begin{enumerate}[label={{\upshape(\arabic*)}}]}
\newcommand{\benumr}{\begin{enumerate}[label={{\upshape(\roman*)}}]}
\newcommand{\eenum}{\end{enumerate}}
\newcommand{\bitem}{\begin{itemize}}
\newcommand{\eitem}{\end{itemize}}
\newcommand{\bc}{}
\newcommand{\bd}{\begin{defn}}
\newcommand{\ed}{\end{defn}}

\newcommand{\beg}{\begin{eg}}
\newcommand{\eeg}{\end{eg}}

\newcommand{\bcl}{\begin{claim}}
\newcommand{\ecl}{\end{claim}}
\newcommand{\lab}{\label}

\newcommand{\x}{\text}

\providecommand{\qxq}[1]{\quad\text{#1}\quad}
\providecommand{\qx}[1]{\quad\text{#1}}
\providecommand{\xq}[1]{\text{#1}\quad}
\newcommand{\q}{\quad}
\newcommand{\qq}{\quad\quad}
\newcommand{\qqq}{\quad\quad\quad}

\newcommand{\tst}{\textstyle}
\newcommand{\ba}{\begin{aligned}}
\newcommand{\ea}{\end{aligned}}
\newcommand{\be}{\begin{equation}}
\newcommand{\ee}{\end{equation}}
\newcommand{\bpf}{\begin{proof}}
\newcommand{\epf}{\end{proof}}
\newcommand{\bthm}{\begin{thm}}
\newcommand{\ethm}{\end{thm}}
\newcommand{\bprop}{\begin{prop}}
\newcommand{\eprop}{\end{prop}}
\newcommand{\bcor}{\begin{cor}}
\newcommand{\ecor}{\end{cor}}
\newcommand{\brem}{\begin{rem}}
\newcommand{\erem}{\end{rem}}
\newcommand*{\QED}{\hfill\ensuremath{\qed}}
\newcommand*{\QEDD}{\hfill\ensuremath{\qed\qed}}

\usepackage{aliascnt}
\newaliascnt{numberingbase}{subsubsection}
\numberwithin{equation}{numberingbase}

\newtheoremstyle{thms}{0pt}{0pt}{\itshape}{}{\bfseries}{.}{ }{}
\theoremstyle{thms}
\newtheorem{conj}[numberingbase]{Conjecture}
\newtheorem{corollary}[numberingbase]{Corollary}
\newtheorem{cor}[numberingbase]{Corollary}
\newtheorem{lemma}[numberingbase]{Lemma}

\newtheorem{prop}[numberingbase]{Proposition}
\newtheorem{proposition}[numberingbase]{Proposition}
\newtheorem{Q}[numberingbase]{Question}
\newtheorem{thm}[numberingbase]{Theorem}
\newtheorem{theorem}[numberingbase]{Theorem}

\newtheoremstyle{claims}{0pt}{0pt}{}{}{\itshape}{.}{ }{}
\theoremstyle{claims}
\newtheorem{claim}[equation]{Claim}

\newtheoremstyle{defs}{0pt}{0pt}{}{}{\bfseries}{.}{ }{}
\theoremstyle{defs}
\newtheorem{defn}[numberingbase]{Definition}

\newtheorem{eg}[numberingbase]{Example}
\newtheorem{example}[numberingbase]{Example}
\newtheorem*{egs}{Examples}
\newtheorem{rem}[numberingbase]{Remark}

\newtheorem*{rems}{Remarks}

\Crefname{claim}{Claim}{Claims}
\Crefname{sublemma}{Lemma}{Lemmas}
\Crefname{conj}{Conjecture}{Conjectures}
\Crefname{cor}{Corollary}{Corollaries}
\Crefname{defn}{Definition}{Definitions}
\Crefname{eg}{Example}{Examples}
\Crefname{prop}{Proposition}{Propositions} 
\Crefname{Q}{Question}{Questions}
\Crefname{rem}{Remark}{Remarks}
\Crefname{thm}{Theorem}{Theorems}
\Crefname{variant}{Variant}{Variants}

\theoremstyle{thms}
\newtheorem{thm-tweak}[subsection]{Theorem}
\Crefname{thm-tweak}{Theorem}{Theorems}
\newtheorem{lemma-tweak}[subsection]{Lemma}
\Crefname{lemma-tweak}{Lemma}{Lemmas}
\newtheorem{cor-tweak}[subsection]{Corollary}
\Crefname{cor-tweak}{Corollary}{Corollaries}
\newtheorem{prop-tweak}[subsection]{Proposition}
\Crefname{prop-tweak}{Proposition}{Propositions} 
\newtheorem{conj-tweak}[subsection]{Conjecture}
\Crefname{conj-tweak}{Conjecture}{Conjectures} 

\theoremstyle{defs}
\newtheorem{defn-tweak}[subsection]{Definition}
\Crefname{defn-tweak}{Definition}{Definitions}
\newtheorem{eg-tweak}[subsection]{Example}
\Crefname{eg-tweak}{Example}{Examples}
\newtheorem*{rems-tweak}{Remarks}
\newtheorem{rem-tweak}[subsection]{Remark}
\Crefname{rem-tweak}{Remark}{Remarks}

\newtheoremstyle{subsection-tweak}
   {0pt}
   {0pt}%
   {}
   {}%
   {\bfseries}
   {}%
   {.5em}
   {\thmnumber{\@{#1}{}\@{#2}.}%
    \thmnote{~{\bfseries#3.}}}    
    
\theoremstyle{subsection-tweak}
\newtheorem{pp}[numberingbase]{}
\newcommand{\bpp}{\begin{pp}}
\newcommand{\epp}{\end{pp}}

\theoremstyle{subsection-tweak}

\newtheorem*{pp-tweak}{}

\makeatletter
\def\@tocline#1#2#3#4#5#6#7{
    \begingroup 
    \@ifempty{#4}{}{}

    \parindent\z@ \leftskip#3\relax \advance\leftskip\@tempdima\relax
    #5\hskip-\@tempdima
      \ifcase #1
       \or\or \hskip 2em \or \hskip 1em \else \hskip 3em \fi%
      #6\nobreak\relax
    \dotfill\hbox to\@pnumwidth{\@tocpagenum{#7}}\par
    \nobreak
    \endgroup
 }
 \def\l@section{\@tocline{1}{0pt}{1pc}{}{}}

\renewcommand{\tocsection}[3]{%
  \indentlabel{\@ifnotempty{#2}{\makebox[1.3em][l]{%
    \ignorespaces#1 \bfseries{#2}.\hfill}}}\bfseries{#3}
    \vspace{-5pt}}

\renewcommand{\tocsubsection}[3]{%
  \indentlabel{\@ifnotempty{#2}{\hspace*{-0.5em}\makebox[2.1em][l]{%
    \ignorespaces#1#2.\hfill}}}#3
    \vspace{-5pt}}
   
\makeatother 

\makeatletter 
\newcommand\appendix@section[1]{%
  \refstepcounter{section}%
  \orig@section*{Appendix \@Alph\c@section. #1}%
}
\let\orig@section\section
\g@addto@macro\appendix{\let\section\appendix@section}
\makeatother

\author{K\k{e}stutis \v{C}esnavi\v{c}ius}
\address{CNRS, UMR 8628, Laboratoire de Math\'{e}matiques d'Orsay, Universit\'{e} Paris-Saclay, 91405 Orsay, France}
\email{kestutis@math.u-psud.fr}

\date{\today}

\newcommand{\bb}{\boldsymbol}		

\usepackage{epigraph}

\begin{document}

\title{Purity for flat cohomology}
\author{Peter Scholze}
\address{Mathematisches Institut, Universit\"{a}t Bonn, Endenicher Allee 60, D-53115, Bonn, Germany}
\email{scholze@math.uni-bonn.de}

\date{\today}


\begin{abstract}\vspace{-10pt}
 We establish the flat cohomology version of the Gabber--Thomason purity for \'{e}tale cohomology: for a complete intersection Noetherian local ring $(R, \fm)$ and a commutative, finite, flat $R$-group $G$, the flat cohomology $H^i_\fm(R, G)$ vanishes for $i < \dim(R)$. For small $i$, this settles conjectures of Gabber that extend the Grothendieck--Lefschetz theorem and give purity for the Brauer group for schemes with complete intersection singularities. For the proof, we reduce to a flat purity statement for perfectoid rings, establish $p$-complete arc descent for flat cohomology of perfectoids, and then relate to coherent cohomology of $\bA_\Inf$ via prismatic Dieudonn\'{e} theory. We also present an algebraic version of tilting for \'{e}tale cohomology, use it to reprove the Gabber--Thomason purity, and exhibit general properties of fppf cohomology of (animated) rings with finite, locally free group scheme coefficients, such as excision, agreement with fpqc cohomology, and continuity.  \end{abstract}


\maketitle

\vspace{-45pt}

\hypersetup{
    linktoc=page,     
}

\renewcommand*\contentsname{}
\tableofcontents


\newpage

\section{Absolute cohomological purity for flat cohomology}

\numberwithin{equation}{subsection}

\csub[Purity theorems]

Purity in algebraic and arithmetic geometry is the phenomenon of various invariants of schemes being insensitive to removing closed subsets of large enough codimension, perhaps the most basic instance being the Hartogs' extension principle in complex geometry. Our main goal is to exhibit purity in the context of flat cohomology, more precisely, to show that on Noetherian schemes with complete intersection singularities, flat cohomology classes with coefficients in commutative, finite, flat group schemes extend uniquely over closed subsets of sufficiently large codimension. In its key local case, this  amounts to the following vanishing (see \Cref{P-d-pf} for a general global statement).

\bthm[Theorems~\ref{main-pf} and \ref{thm:main-pf-regular}] \label{main}
For a Noetherian local ring $(R, \fm)$ that is a complete intersection\footnote{Recall that $(R, \fm)$ is a \emph{complete intersection} if its 
completion is a quotient of a regular 
ring by a regular sequence.} and a commutative, finite, flat $R$-group scheme $G$,
\[
H^i_\fm(R, G) \cong 0 \q \text{for} \q \begin{cases} i < \dim(R); \\ i \le \dim(R), \qx{if $R$ is regular and not a field.} \end{cases}
\]
\ethm

\Cref{main} is the flat cohomology version of absolute cohomological purity\footnote{\label{foot:semipurity}In the terminology of \cite{SGA2new}*{expos\'{e} XIV, th\'{e}or\`{e}me 1.10} or \cite{SGA4.5}*{expos\'{e} Cycle, rappel 2.2.8}, vanishing of cohomology with supports in low cohomological degrees as in \Cref{main,abs-coh-pur-ann} goes by the name \emph{semipurity}, as opposed to \emph{purity} that would also include high cohomological degrees. In this article, for the sake of brevity, we do not make this distinction.} for \'{e}tale cohomology that had been conjectured by Grothendieck. The latter, stated in \Cref{abs-coh-pur-ann}, was proved by Gabber: first in \cite{Fuj02} by building on the $K$-theoretic approach of Thomason \cite{Tho84}, and then again in \cite{ILO14}*{expos\'{e} XVI} in the framework of general structural results on \'{e}tale cohomology of Noetherian schemes. We give a third proof that uses perfectoid techniques to reduce to positive characteristic.\footnote{A proof that uses perfectoids was also discovered by Fujiwara.}

\bthm[\Cref{abs-coh-pur}] \label{abs-coh-pur-ann}
For a regular local ring $(R, \fm)$ and a commutative, finite, \'{e}tale $R$-group $G$ whose order is invertible in $R$, 
\[
H^i_\fm(R, G) \cong 0 \qxq{for} i < 2\dim(R). 
\]
\ethm


Gabber used \Cref{abs-coh-pur-ann} to deduce the case of \Cref{main} when the order of $G$ is invertible in $R$ in \cite{Gab04}*{Theorem 3}. We review one such deduction based on the Lefschetz hyperplane theorem in local \'{e}tale cohomology in \S\ref{et-depth-geom-depth}. Since the Lefschetz isomorphism range is roughly in degrees $< \dim(R)$, the weaker condition $i < 2\dim(R)$ is specific to \Cref{abs-coh-pur-ann}. 

\Cref{main} for regular $R$ was conjectured in \cite{Pop19}*{Conjecture A.1} and desired reductions quickly lead to including complete intersections (see \Cref{finite-cover}).
In unpublished work, Gabber obtained it 
for $G = \bZ/p\bZ$ and  also for $G = \mu_p$ with $i \le 3$ by building on the combination of perfectoid techniques that were used to settle the complete intersection case of the weight monodromy conjecture of Deligne in \cite{Sch12} and the purity for the Brauer group conjecture of Grothendieck in \cite{brauer-purity}.

The following corollary of \Cref{main} settles conjectures of Gabber \cite{Gab04}*{Conjectures~2~and~3}.

\bthm[\Cref{deg-2-case,deg-3-case}] \label{thm:Gabber-conjectures}
Let $(R, \fm)$ be a Noetherian local ring  that is a complete intersection and let $U_R \ce \Spec(R) \setminus \{ \fm\}$ be its punctured spectrum.
\benum
\item \lab{item:GC-a}
If $\dim(R) \ge 3$, then 
\[
\q\Pic(U_R)_\tors \cong 0
\]
\up{recall that if $\dim(R) \ge 4$, then even $\Pic(U_R) \cong 0$}.

\item \lab{item:GC-b}
If $\dim(R) \ge 4$ or if both $R$ is regular and $\dim(R) \ge 2$, then 
\[
\qq\Br(R) \isomto \Br(U_R).
\]
\eenum
\ethm

The parenthetical aspect of \ref{item:GC-a} is the Grothendieck--Lefschetz theorem \cite{SGA2new}*{expos\'{e} XI, th\'{e}or\`{e}me 3.13~(ii)}. Although the statement of \Cref{thm:Gabber-conjectures}~\ref{item:GC-a} is relatively basic, we do not know how to argue it without ideas that go into proving \Cref{main}. Nevertheless, 
\begin{itemize}
\m
its case when $\wh{R}$ is a quotient of a regular local ring by a \emph{principal} ideal (the hypersurface case) was settled by Dao \cite{Dao12}*{Corollary 3.5}, who even found a version for vector bundles; 

\m
its case when $R$ is an $\bF_p$-algebra was settled in \cite{Gab04}*{Theorem~5~(1)} (also in \cite{DLM10}*{Corollary 2.10});~and 

\m
its case for torsion of order invertible in $R$ was settled in \cite{Rob76} (also in 
\cite{SGA2new}*{expos\'{e}~X, th\'{e}or\`{e}me~3.4}). 
\end{itemize}

\Cref{thm:Gabber-conjectures}~\ref{item:GC-b} reproves the purity for the Brauer group from \cite{brauer-purity} and extends it to schemes with complete intersection singularities. In the cases when $R$ is an $\bF_p$-algebra or when $\dim(R) \ge 5$, this extension was obtained by Gabber in \cite{Gab04}*{Theorem~5}. For regular $R$, even though  the proof is more complex than the one in \cite{brauer-purity}, it does not require treating the case $\dim(R) = 3$ separately (this case was settled in \cite{Gab81}*{Chapter I, Theorem~2$'$} and used in \cite{brauer-purity} as an input).

The global version of \Cref{thm:Gabber-conjectures}~\ref{item:GC-b} may be formulated as follows.

\bthm[\Cref{Brauer-purity}] 
For a Noetherian scheme $X$ and a closed subset $Z \subset X$ such that each $\sO_{X,\, z}$ with $z \in Z$ is either a complete intersection of dimension $\ge 4$ or regular of dimension $\ge 2$, 
\[
H^2(X, \bG_m)_\tors \isomto H^2(X \setminus Z, \bG_m)_\tors.
\]
\ethm


As for \Cref{main} itself, except for its assertion about the cohomological degree $i = \dim(R)$ that requires further arguments, we exploit a suitable version of Andr\'{e}'s lemma to eventually reduce the key case when $G$ is of $p$-power order with $p = \Char(R/\fm) > 0$ to the following purity for flat cohomology in an (integral) perfectoid setting (for a basic review of perfectoid rings, see \S\ref{perfectoid-def}).

\bthm[\Cref{thm:perfectoid-purity}] \label{thm:perfectoid-purity-ann}
For a perfectoid $\bZ_p$-algebra $A$, a commutative, finite, locally free $A$-group $G$ of $p$-power order, and a closed subset $Z \subset \Spec (A/pA)$ such that $\depth_Z(A) \ge d$ in the sense that there is an $A$-regular sequence $a_1, \dotsc, a_d \in A$ that vanishes on $Z$, we have
\[
H^i_Z(A, G) \cong 0 \qxq{for} i < d.
\]
\ethm

For instance, a basic case is when $A$ is a perfect $\bF_p$-algebra. Then, by results of Berthelot \cite{Ber80}, Gabber (unpublished), and Lau \cite{Lau13}, such $A$-groups $G$ are classified by their crystalline Dieudonn\'{e} modules $\bM(G)$, which are $p$-power torsion, finitely presented $W(A)$-modules (that is, $\bA_\Inf(A)$-modules) of projective dimension $\le 1$ equipped with semilinear Frobenius and Verschiebung endomorphisms $F$ and $V$. We use this classification to describe the flat cohomology with coefficients in $G$ in terms of the quasi-coherent cohomology with coefficients in $\bM(G)$: we show in \Cref{KT-input} that
\be\tag{$1.1.6$} \label{eqn:key-formula}
R\Gamma_Z(A, G) \cong R\Gamma_Z(\bA_\Inf(A), \bM(G))^{V - 1}.
\ee
Since $p$ is a nonzerodivisor in $\bA_\Inf(A)$, the sequence $p, a_1, \dotsc, a_d$ is regular in $\bA_\Inf(A)$ and vanishes on $Z$. By expressing $\bM(G)$ as the cokernel of a map between finite projective $\bA_\Inf(A)$-modules, we may then deduce the vanishing of the right side of \eqref{eqn:key-formula} in the desired range of degrees from the fact that ``enough depth'' implies the vanishing of quasi-coherent cohomology with supports.

An analogous argument proves \Cref{thm:perfectoid-purity-ann} in general, except that now the key formula \eqref{eqn:key-formula} lies significantly deeper. To make sense of it, we replace crystalline Dieudonn\'{e} theory used to define $\bM(G)$ by its prismatic generalization developed in \cite{ALB20}, which built on the classification of $G$ in terms of $\bM(G)$ over perfectoid rings due to Lau and the second named author \cite{Lau13}, \cite{SW19}*{Appendix to Lecture 17}. Our strategy for settling \eqref{eqn:key-formula} in general is to first show that both of its sides satisfy hyperdescent in the $p$-complete arc topology of Bhatt--Mathew \cite{BM20} (reviewed in \uS\uref{pp:arc-topology}) and to then use the resulting ability to replace $A$ by a $p$-complete arc cover to reduce to the case when $A$ is a product of perfectoid valuation rings with algebraically closed fraction fields, a case that admits a reasonably direct argument. 
 With \eqref{eqn:key-formula} in place, the regular sequence $p, a_1, \dotsc, a_d$ gets replaced by $\xi, a_1, \dotsc, a_d$, where $\xi$ is a generator of $\Ker(\theta\colon \bA_\Inf(A) \surjects A)$ (that is, $\xi$ is an orientation of the perfect prism that corresponds to $A$), and the same depth argument gives~\Cref{thm:perfectoid-purity-ann}.

Overall, a critical point of the proof of \Cref{main} is the implication
\[
\depth_Z(A) \ge d \q \Ra \q \depth_Z(\bA_\Inf(A)) \ge d + 1,
\]
where we understand $\depth_Z$ na\"{i}vely, that is, in terms of regular sequences. Indeed, in the end it seems critical to work over $\bA_\Inf(A)$---direct reductions of \Cref{main} to positive characteristic, in cases in which they are available, seem to always produce a ``one off'' cohomological degree problem, and hence to not give optimal statements. For instance, under the weaker assumption $\dim(R) \ge 5$, Gabber proved \Cref{thm:Gabber-conjectures}~\ref{item:GC-b} in \cite{Gab04}*{Theorem~5}  by first reducing to $p$-torsion free $R$ and then further to the complete intersection $\bF_p$-algebra $R/p$ of dimension $\ge 4$.

\csub[Flat cohomology of animated rings]

The $p$-complete arc hyperdescent for the flat cohomology side of the key formula \eqref{eqn:key-formula}  is a major portion of the overall argument of \Cref{main}, a portion for which we resort to flat cohomology in the more flexible setting of derived algebraic geometry (as defined in \S\ref{pp:simplicial-coho-def}). For the latter, we use simplicial rings, for which 
we decided to use 
different terminology than the usual one because we believe it to be confusing to continue calling  the objects of the resulting $\infty$-category ``simplicial rings''---certainly, we do not think of them as simplicial objects in the category of rings.\footnote{The following standard example explains why we do not like to think in terms of simplicial rings: if $A_\bullet$ is a simplicial ring, then any scheme $X$ gives rise to the simplicial set $X(A_\bullet)$; however, for general $X$ this functor does not preserve weak equivalences. There is another functor $A\mapsto X(A)$ whose input is a simplicial ring up to weak equivalence and whose output is a simplicial set up to weak equivalence. This functor is slightly tricky to define in the simplicial language, but it is the one that will be relevant to us.}

We refer to the $\infty$-category obtained from simplicial rings (resp.,~simplicial abelian groups; resp.,~simplicial sets, etc.)~by inverting weak equivalences as the $\infty$-category of \emph{animated rings} (resp.,~\emph{animated abelian groups}; resp.,~\emph{animated sets}, etc.). In the background there is a general ``free generation by sifted colimits'' procedure described in \S\ref{animation} that from any reasonable category $\sC$ produces an $\infty$-category $\mathrm{Ani}(\sC)$, the \emph{animation} of $\sC$, that contains $\sC$ as full subcategory: $\mathrm{Ani}(\sC)$ is nothing else than a ``nonabelian derived category'' in the sense of Quillen, compare with \cite[Section 5.5.8]{HTT}. The inclusion $\sC\hookrightarrow \mathrm{Ani}(\sC)$ has a left adjoint $\pi_0\colon \mathrm{Ani}(\sC)\to \sC$.

For example, the $\infty$-category of animated sets (the case $\sC = \Set$) is exactly the $\infty$-category of ``spaces'' in the sense of Lurie. We prefer the term ``animated set,'' or ``anima'' for brevity, suggested by the general naming convention: we believe the term ``space,''  whose origins seem to be historical, to be highly nondescriptive---it is arguable whether something as combinatorial as a simplicial set should count as a space, and also note that ``spaces'' in the sense of Lurie do not have an underlying set of points. 
Philosophically, ``anima''  means something like ``soul''---and, indeed, the functor from topological spaces to their homotopy category extracts something like the soul of a space: it only remembers data independent of any worldly representation in terms of physical points.

The animation procedure is quite powerful: for example, the $\infty$-category of pairs consisting of an animated ring $A$ and an animated $A$-module (also known as a connective $A$-module) is obtained by animating the usual category of rings equipped with a module. In particular, by passing to the fibre over any given animated ring $A$, we obtain the $\infty$-category of animated $A$-modules. Derived tensor products (of animated modules or animated rings) are obtained by animating the usual functors.\footnote{Taking up the previous footnote: in this language, if $X=\Spec(R)$ is an affine scheme and $A$ is an animated ring, then $X(A)$ refers to the anima of maps $R\to A$ of animated rings (which is now the only option that suggests itself). One can extend to nonaffine schemes by Zariski sheafification.}

We show the following properties of fppf cohomology of animated rings with coefficients in commutative, finite, locally free group schemes. These properties are new already for usual commutative rings but their proofs greatly benefit from the flexibility of the more general setting: intermediate steps, such as passage to derived $p$-adic completions or derived base changes, leave the realm of usual rings. 

\bthm[\Cref{cor:excision} with \Cref{comp-explain}] \label{thm:excision-ann}
For a ring $R$, a map $f\colon A \ra A'$ of animated $R$-algebras, and a finitely generated ideal 
\[
I = (a_1, \dotsc, a_r) \subset \pi_0(A)
\] 
such that $f$ induces an isomorphism after iteratively forming derived $a_i$-adic completions for $i = 1, \dotsc, r$, 
\[
R\Gamma_I(A, G) \isomto R\Gamma_I(A', G)
\]
for every commutative, finite, locally free $R$-group $G$.
\ethm

For instance, this excision result allows us to replace $R$ in \Cref{main} by its completion $\wh{R}$.

\bthm[\Cref{thm:strongdescent}] \label{thm:fppf-fpqc-ann}
For a ring $R$ and a commutative, finite, locally free $R$-group $G$, the functor $A \mapsto R\Gamma_\fppf(A, G)$ satisfies hyperdescent in the fpqc topology on animated $R$-algebras $A$.
\ethm

The following result is the $p$-complete arc hyperdescent for the left side of the key formula \eqref{eqn:key-formula}. An important input to its proof is the analogous $p$-complete arc descent for the structure (pre)sheaf functor $A \mapsto A$ on perfectoids that was exhibited in \cite{BS19}*{Proposition~8.10}.


\bthm[\Cref{thm:strongdescent-perfectoid}] \label{thm:descent-ann}
For a $p$-complete arc hypercover $A \ra A^\bullet$  
of perfectoid $\bZ_p$-algebras, a closed $Z \subset \Spec(A/pA)$, and a commutative, finite, locally free $A$-group $G$ of $p$-power order,
\[
\tst R\Gamma_Z(A, G) \isomto R\lim_{\Delta} \p{R\Gamma_Z(A^\bullet, G)}\!, \qx{where $\Delta$ is the simplex category.}
\]
\ethm

The following continuity formula, among other things, computes the flat cohomology of complete Noetherian local rings with commutative, finite, flat group coefficients and has consequences for invariance of flat cohomology under Henselian pairs, see \Cref{eg:continuity} and \Cref{inv-Hens-pair}.

\bthm[\Cref{adic-continuity}] \label{adic-continuity-ann}
For a ring $R$, an animated $R$-algebra $A$, elements $a_1, \dotsc, a_r \in A$ such that $A$ agrees with its iterated derived $a_i$-adic completion for $i = 1, \dotsc, r$, and a commutative, finite, locally free $R$-group $G$,
\[
\tst R\Gamma(A, G) \isomto R\lim_{n > 0} (R\Gamma(A/^\bL (a_1^n, \dotsc, a_r^n), G)). 
\]
\ethm


For the derived quotient notation used above, see \uS\uref{ani-rings}. 
Roughly speaking, we deduce Theorems~\ref{thm:excision-ann}--\ref{adic-continuity-ann} from the positive characteristic case of the key formula \eqref{eqn:key-formula}, that is, from crystalline Dieudonn\'{e} theory. More precisely, for $G$ of $p$-power order we first analyze $R\Gamma((-)[\f 1p], G)$ by identifying with \'{e}tale cohomology and using arc descent results of \cite{BM20} recalled in \Cref{arc-descent} and \Cref{BM-hens-comp} (animated aspects disappear in this step because the $\infty$-category of \'{e}tale $A$-algebras is equivalent to that of \'{e}tale $\pi_0(A)$-algebras, see \Cref{rem:no-new-etale}). We may then work along $\{p = 0\}$ to assume that $A$ is $p$-Henselian and consequently reduce to $\bF_p$-algebras by combining animated deformation theory with the following general $p$-adic continuity result that we establish by a more or less direct attack.

\bthm[\Cref{thm:padiclimit}] \label{thm:p-adic-continuity-ann}
For a prime $p$, a ring $R$, a commutative, finite, locally free $R$-group $G$ of $p$-power order, and an animated $R$-algebra $A$ for which the ring $\pi_0(A)$ is $p$-Henselian,
\[
\tst R\Gamma_\fppf(A, G) \isomto R\lim_{n > 0}\p{R\Gamma_\fppf(A/^\bL p^n, G)}\!. 
\]
\ethm

For instance, $A$ in \Cref{thm:p-adic-continuity-ann} could be a $p$-adically complete (usual) ring, although even in this case, unless $A$ is $p$-torsion free, the derived reductions appearing in the limit are animated rings.

\csub[An overview of the proof of purity for flat cohomology]

In summary, the overall proof of the purity for flat cohomology of \Cref{main} proceeds as follows.

\benuma
\m  
Use a Lefschetz hyperplane theorem in local \'{e}tale cohomology to deduce the prime to the residue characteristic aspects from the purity for \'{e}tale cohomology of \Cref{abs-coh-pur-ann} 
(see~\S\ref{et-depth-geom-depth}).

\m
Use crystalline Dieudonn\'{e} theory to establish the positive characteristic case of the key formula \eqref{eqn:key-formula} 
(see \S\ref{section:Dieudonne-positive-characteristic}); this already mostly settles \Cref{main} when $R$ is an $\bF_p$-algebra.

\m \label{step-3}
Use the positive characteristic case of the key formula \eqref{eqn:key-formula} and animated deformation theory to show the new properties of fppf cohomology stated in Theorems \ref{thm:excision-ann}--\ref{thm:descent-ann} (see \S\S\ref{sec:defthy}--\ref{sec:descent}).

\m
Combine $p$-complete arc descent of \Cref{thm:descent-ann} with prismatic Dieudonn\'{e} theory to establish the key formula \eqref{eqn:key-formula} in general; deduce the perfectoid purity \Cref{thm:perfectoid-purity-ann} (see \S\ref{perfectoid-version}).

\m \label{step-5}
Combine excision obtained in \Cref{thm:excision-ann}, a version of Andr\'{e}'s lemma (see \S\ref{And-lem}), and deformation theory to reduce \Cref{main} to the perfectoid purity \Cref{thm:perfectoid-purity-ann} (see \S\ref{grand-finale}).
\eenum

Andr\'{e}'s lemma says that every element of a perfectoid ring attains compatible $p$-power roots after passing to a flat modulo powers of $p$ perfectoid cover. We build on ideas of Gabber--Ramero to generalize it: in \Cref{Andre-lem} below, the cover is flat, and even ind-syntomic, before reducing modulo powers of $p$. This is well suited for us, although the version of \cite{BS19}*{Theorem~7.14, Remark~7.15} combined with \Cref{thm:p-adic-continuity-ann} suffices as well. 
We apply Andr\'{e}'s lemma to elements $f_i$ that cut out our complete intersection  inside a regular ring:  only regular rings have flat covers by perfectoids (see \cite{BIM19}), so, by \Cref{recognize}~\ref{R-b} below, we need to kill all the $f_i^{1/p^{\infty}}$ to reach a perfectoid.

Deformation theory used in Step \ref{step-5} is where the complete intersection assumption manifests itself. Namely, on flat cohomology the difference between killing $f_i$ and, say, $f_i^{\f{p - 1}{p}}$ amounts to quasi-coherent cohomology, and if the $f_i$ form a regular sequence, then the intervening square-zero ideals are module-free and so of large enough depth (see \Cref{mod-out-by-more} for this argument). In general, \Cref{main} fails for Cohen--Macaulay $R$ (even over $\bC$), see \Cref{no-free-lunch}. 

In contrast, if $G$ is \'{e}tale, then the complete intersection assumption is a red herring: as the following refinement of \Cref{main} shows, then the purity of cohomology is controlled by the \emph{virtual dimension} $\vdim(R)$ of the Noetherian local ring $(R, \fm)$. This numerical invariant is defined in terms of the number of equations that cut $\wh{R}$ out in a regular ring (see \eqref{eqn:vdim-def}) and satisfies $\vdim(R) \le \dim(R)$ with equality precisely for complete intersection $R$.

\bthm[\Cref{ell-semipurity-p}] \label{ell-semipurity-ann}
For a Noetherian local ring $(R, \fm)$ and a commutative, finite, \'{e}tale $R$-group $G$,
\[
H^i_\fm(R, G) \cong 0 \qxq{for} i < \vdim(R).
\]
\ethm

Informally, this result says that the ``\'{e}tale depth'' of $R$ is at least $\vdim(R)$ (the former was defined in \cite{SGA2new}*{expos\'{e} XIV, d\'{e}finitions 1.2 et 1.7}). In \S\ref{nonabelian}, we exhibit a nonabelian version: by \Cref{non-ab-et}, the purity for the \'{e}tale fundamental group proved in \cite{SGA2new}*{expos\'{e} X, th\'{e}or\`{e}me~3.4} for complete intersections of dimension $\ge 3$ continues to hold for arbitrary Noetherian local rings of virtual dimension $\ge 3$.

\csub[Notation and conventions] \label{conv}
All our rings are commutative and unital. We use the definition \cite{EGAIV1}*{chapitre 0, d\'{e}finition~15.1.7, paragraphe 15.2.2} of a regular sequence (so there is no condition on the quotients being nonzero).	 A regular local ring $(R, \fm)$ is \emph{unramified} if it is of mixed characteristic $(0, p)$ and $p \not \in \fm^2$. By the Cohen structure theorem \cite{EGAIV1}*{chapitre 0, th\'{e}or\`{e}me~19.8.8~(i)}, the completion of a Noetherian local ring $(R, \fm)$ is a quotient of a regular ring $\wt{R}$ that may be chosen unramified if  $\Char(R/\fm) = p > 0$. Such an $R$ is a \emph{complete intersection} if the ideal that cuts $\wh{R}$ out in $\wt{R}$ is generated by a regular sequence. We recall that every ideal that cuts out a complete intersection in a regular ring is generated by a regular sequence (see \cite{SP}*{Lemma~\href{https://stacks.math.columbia.edu/tag/09Q1}{09Q1}}).  

For a module $M$ over a ring $A$, we write $M\langle a \rangle$ for the kernel of the scaling by $a \in A$, and we set 
\[
\tst M\langle a^\infty \rangle \ce \bigcup_{n \ge 0} M\langle a^n\rangle;
\]
we say that $M$ has \emph{bounded $a^\infty$-torsion} if $M\langle a^\infty\rangle = M\langle a^N\rangle $ for some  $N > 0$. 
An $\bF_p$-algebra is \emph{perfect} (resp.,~\emph{semiperfect}) if its absolute Frobenius endomorphism $a \mapsto a^p$ is bijective (resp.,~surjective); these conditions ascend along \'{e}tale maps (see \cite{SGA5}*{expos\'{e}~XV, proposition~2~c)~2)} and \cite{SP}*{Lemma~\href{https://stacks.math.columbia.edu/tag/04D1}{04D1}}). For an implicit prime $p$, we let $W(-)$ denote the $p$-typical Witt vectors and  indicate Teichm\"{u}ller lifts by $[-]$. We use the (somewhat nonstandard) notation 
\be \lab{nst-terminology}
\tst A\llb x_1^{1/p^\infty},\ldots,x_n^{1/p^\infty}\rrb \ce \varinjlim_m \p{A\llb x_1^{1/p^m},\ldots,x_n^{1/p^m}\rrb}
\ee
(we do not form an additional $(x_1,\ldots,x_n)$-adic completion). We use the derived quotient notation
\[
A/^\bL a \ce A\otimes^\bL_{\mathbb Z[X]} \mathbb Z,
\]
where
\[
\bZ[X] \ra A \ \ \x{via}\ \ X \mapsto a \qxq{and} \bZ[X] \ra \bZ \ \ \x{via}\ \ X \mapsto 0.
\]
We let $(-)^*$ indicate the dual of a vector bundle, or of a $p$-divisible group, or of a commutative, finite, locally free group scheme $G$.
  We often use the \emph{B\'{e}gueri resolution} of the latter by commutative, smooth, affine $S$-group schemes (see \cite{Beg80}*{proposition 2.2.1} and \cite{SP}*{Lemma~\href{http://stacks.math.columbia.edu/tag/01ZT}{01ZT}}):
\be \label{Begueri-0}
0\to G\to \mathrm{Res}_{G^\ast/S}( \mathbb G_m)\to Q\to 0.
\ee
Unless indicated otherwise, we form cohomology in the fppf topology and make the identifications with \'{e}tale (smooth coefficients) or Zariski cohomology (quasi-coherent coefficients) implicitly. 

 We let $\Delta$ be the simplex category, whose opposite indexes simplicial objects. We write $\cD(\bZ)$ for the derived $\infty$-category of~$\bZ$.
 We say that a functor $F$ defined on some subcategory of rings and, for the sake of concreteness, valued in $\cD(\bZ)$ \emph{satisfies descent} (resp.,~\emph{satisfies hyperdescent}) for a Grothendieck topology $\sT$ if for every $\sT$-cover $A \ra A'$ with its \v{C}ech nerve $A'^\bullet$ (resp.,~for every $\sT$-hypercover $A \ra A'^\bullet$), we have 
\[
\tst F(A) \isomto R\lim_{\Delta} (F(A'^\bullet)).
\]
We say that an ($\infty$-)category $\sC$ is \emph{complete} (resp.,~\emph{cocomplete}) if it has all small limits (resp.,~colimits). A \emph{strong limit cardinal of uncountable cofinality} is a limit cardinal $\kappa$ such that for every sequence $\kappa_0, \kappa_1, \dotsc$ of cardinals $< \kappa$ we have $2^{\kappa_0} < \kappa$ and $\sup_{n \ge 0} \kappa_n < \kappa$; there exist arbitrarily large such $\kappa$, see \cite{Sch17}*{Lemma 4.1 and its proof}. We use such cardinals $\kappa$ to avoid set-theoretic problems when working with large sites (such as the fpqc site); of course, in such cases we check along the way that our assertions do not depend on the choice of  $\kappa$. We say that a scheme $S$ is of size $<\kappa$ if the cardinality of its underlying topological space is $< \kappa$ and $\abs{\Gamma(U, \sO_S)} < \kappa$ for every affine open $U \subset S$.

\begin{pp-tweak}[\!\!Acknowledgements]
We thank the referee for helpful comments and suggestions. We thank Johannes Ansch\"{u}tz, Bhargav Bhatt, Alexis Bouthier, Hailong Dao, Aise Johan de Jong, Ofer Gabber, Benjamin Hennion, Luc Illusie, Kazuhiro Ito, Teruhisa Koshikawa, Arnab Kundu, Arthur-C\'{e}sar Le Bras,  Bernard Le Stum, Shang Li, Jacob Lurie, Linquan Ma, Akhil Mathew, Kei Nakazato, Martin Olsson, Burt Totaro, Zijian Yao, and Yifei Zhao for helpful conversations or correspondence. In particular, we thank Jacob Lurie for explaining to us the proof of \Cref{cor:deftheorygroup}. Moreover, we thank Dustin Clausen for the term ``animation.'' We thank the MSRI for support during a part of the preparation period of this article during the Spring semester of 2019 supported by the National Science Foundation under Grant No.~1440140. The first-named author thanks the University of Bonn and Vi\stackon[-9.3pt]{\^{e}$\mkern0.5mu$}{\d{}}n To\'{a}n H\d{o}c for hospitality during visits. This project has received funding from the European Research Council (ERC) under the European Union's Horizon 2020 research and innovation programme (grant agreement No.~851146).
\end{pp-tweak}

\numberwithin{equation}{numberingbase}



\section{The geometry of integral perfectoid rings} \lab{global}

We begin with generalities about perfectoids that will be important in multiple steps of the overall argument of purity for flat cohomology. We review the definitions and expose some basic properties in \S\ref{sec:perfectoid-structure}. We then present an algebraic approach to controlling cohomology under tilting in \S\ref{sec:tilting-cohomology} that avoids adic spaces and the almost purity theorem in favor of arc descent. Finally, in \S\ref{And-lem}, we generalize Andr\'{e}'s lemma: we improve its flatness aspect to ind-syntomicity and we avoid completions.

\bpp[The implicitly fixed prime] \label{fixed-prime}
Throughout \S\S\ref{sec:perfectoid-structure}--\ref{And-lem}, to discuss perfectoids, we fix a prime~$p$.
\epp


\csub[Structural properties of perfectoid rings] \label{sec:perfectoid-structure}

Perfectoids play a central role in our approach to purity, so we summarize their most relevant basic properties in this section. Our perfectoids are what some authors call ``integral perfectoids'': these appear to be the ones most directly related to commutative algebra. Perfectoid rings generalize perfect $\bF_p$-algebras beyond the setting of positive characteristic, and their definition is intimately related to the properties of the following tilting adjunction.

\bpp[The tilting adjunction] \label{pp:tilt-adj}
In positive characteristic, the inclusion of perfect $\bF_p$-algebras into all $\bF_p$-algebras admits a right adjoint given by the inverse limit perfection $B \mapsto B^\flat \ce \varprojlim_{b \mapsto b^p} B$. We define the \emph{tilt} of general ring $A$ as the inverse limit perfection $A^\flat \ce (A/pA)^\flat$ of $A/pA$. When restricted to $p$-adically complete $A$, the tilting functor is the right adjoint of $p$-typical Witt vectors: 
\[
\xymatrix{
\{ \x{perfect $\bF_p$-algebras}\} \ar@<0.5ex>[rr]^-{W(-)} && \ar@<0.5ex>[ll]^-{(-)^\flat} \{\x{$p$-adically complete $\bZ_p$-algebras}\},
}
\]
see \cite{SZ18}*{Proposition 3.12}. The Fontaine functor $\bA_\Inf(-)$ is the following composition of these adjoints:
\be \label{eqn:map-theta}
\bA_\Inf(A) \ce W(A^\flat),
\ee
so it comes with the counit of the adjunction $\theta \colon \bA_\Inf(A) \ra A$. 
The functors $(-)^\flat$ and $W(-)$ commute with limits of rings, so the same holds for $\bA_\Inf(-)$. By \cite{BMS18}*{Lemma~3.2~(i)}, for $p$-adically complete $A$, the reduction modulo $p$ map
\be \label{monoid-id}
\tst \varprojlim_{a \mapsto a^p} A \xra{\sim} \varprojlim_{a \mapsto a^p} A/pA  \cong A^\flat
\ee
is a multiplicative isomorphism, and we let $a \mapsto a^\sharp$ denote the resulting multiplicative projection $A^\flat \ra A$ onto the last coordinate. The counit map $\theta$ satisfies $\theta([a]) = a^\sharp$ for $a \in A^\flat$\footnote{The elements $\theta([a]), a^\sharp \in A$ agree modulo $p$, and the same holds for their $p^n$-th roots $\theta([a^{1/p^n}]), (a^{1/p^n})^\sharp \in A$. Thus, the $p$-adic completeness of $A$ gives the middle equality in 
\[
\tst \theta([a]) = \lim_{n \ra \infty} (\theta([a^{1/p^n}])^{p^n}) = \lim_{n \ra \infty} (((a^{1/p^n})^\sharp)^{p^n}) = a^\sharp.
\]} and is uniquely determined by this. By \emph{loc.~cit.},~if $A$ is $\varpi$-adically complete for a $\varpi \in A$ with $\varpi \mid p$, then $A^\flat$ may be defined with $\varpi$ in place of $p$: then
\be \label{monoid-id-2}
\tst A^\flat = \varprojlim_{a \mapsto a^p} A/pA \xra{\sim} \varprojlim_{a \mapsto a^p} A/\varpi A.  
\ee


 %

\epp

\bpp[Perfectoid rings]  \label{perfectoid-def}
As in \cite{BMS18}*{Definition~3.5}, we say that a ring $A$ is \emph{perfectoid} if  
\benumr
\item \label{PD-i}
there is a $\varpi \in A$ with $\varpi^p \mid p$ such that $A$ is $\varpi$-adically (so also $p$-adically) complete, and

\item \label{PD-ii}
the counit map $\theta\colon \bA_\Inf(A) \surjects A$ of \eqref{eqn:map-theta} is surjective and its kernel is principal.


\eenum

Since $A$ is $p$-adically complete, the surjectivity of $\theta$ is equivalent to the semiperfectness of $A/pA$. 

One may choose $\varpi$ to be the $p$-th root of a unit multiple of $p$, but it is useful to develop the theory relative to a general $\varpi$. More precisely, by \cite{BMS18}*{Lemma~3.9}, the condition \ref{PD-i} and the semiperfectness of $A/pA$ alone ensure that some unit multiples of $\varpi$ and $p$ have compatible $p$-power roots in $A$: there are $\varpi^\flat, p^\flat \in A^\flat$ such that $(\varpi^\flat)^\sharp, (p^\flat)^\sharp \in A$ are unit multiples of $\varpi$ and $p$, respectively. 

By \cite{Lau18}*{Remark~8.6} or \cite{BS19}*{Theorem~3.10}, the conditions \ref{PD-i}--\ref{PD-ii} may be synthesized: a ring $A$ is perfectoid if and only if there are a perfect $\bF_p$-algebra $B$ and a $\xi = (\xi_0, \xi_1,  \dotsc) \in W(B)$ such that 
\be \label{alt-def}
A  \simeq W(B)/(\xi) \qxq{and} B \qxq{is $\xi_0$-adically complete with} \xi_1 \in B^\times.
\ee 
We necessarily have $B \cong A^\flat$, the displayed identification is induced by $\theta$, and $\xi$ is a nonzerodivisor in $\bA_\Inf(A)$. In fact, $A^\flat$ is $\xi'_0$-adically complete for any $\xi' \in \Ker(\theta)$ with Witt vector coordinates $\xi' = (\xi'_0, \xi'_1, \dotsc)$ and, by \cite{BMS18}*{Remark~3.11}, such a $\xi'$ generates $\Ker(\theta)$ if and only if $\xi'_1 \in (A^\flat)^\times$, in which case $\xi'$ is a  nonzerodivisor in $\bA_\Inf(A)$.  In particular, $\xi$ continues to generate $\Ker(\theta)$ for any perfectoid $A$-algebra, and an $\bF_p$-algebra is perfectoid if and only if it is perfect (choose $\xi' = p$). Explicitly, there is an $x \in \bA_\Inf(A)$ such that $p + [p^\flat]x \in \Ker(\theta)$, and this element is a possible choice for $\xi$: the equality $x = [x \bmod p] + px'$ shows that the first Witt coordinate of $p + [p^\flat]x$ is a unit because subtracting the Teichm\"{u}ller of its zeroth Witt coordinate gives $p(1 + [p^\flat]x')$. More precisely, we may choose this  $x$ in such a way that $\theta(x)$ be a unit in $A$; then \eqref{monoid-id} shows that $x \bmod p$ is a unit in $A^\flat$, to the effect that we may adjust our choice of $p^\flat$ to arrange that even $p^\flat = \xi_0$.

By \ref{PD-ii} and the proof of \cite{BMS18}*{Lemma~3.10~(i)}, the map $a \mapsto a^\sharp$ induces isomorphisms
\be \label{same-mod-pi}
A^\flat/(\varpi^\flat)^p A^\flat \isomto 
A/\varpi^p A  \qxq{and}  A^\flat/p^\flat A^\flat \isomto 
A/pA.
\ee
Thus, $(\varpi^\flat)^p \mid p^\flat$, and \eqref{monoid-id-2} shows that $A^\flat$ is $\varpi^\flat$-adically complete (every $\varpi^\flat$-adic Cauchy sequence in $A^\flat$ stabilizes in each term of $\varprojlim_{a \mapsto a^p} A/\varpi A$). 
Although $\varpi^\flat$ and $p^\flat$ are noncanonical, \eqref{same-mod-pi} determines the ideals $(\varpi^\flat)$ and $(p^\flat)$ of $A^\flat$, so we will use $\varpi^\flat$ and $p^\flat$ when only $(\varpi^\flat)$ and $(p^\flat)$ matter. 





By \cite{BMS18}*{Lemma~3.10}, for a perfectoid $A$ that is $\varpi$-adically  complete for a $\varpi \in A$ with $\varpi^p \mid p$, the $p$-power map
\be \label{p-oid-mod-pi}
\xymatrix@C=30pt{A/\varpi A \ar[r]_-{\sim}^-{a\, \mapsto\, a^p}  & A/\varpi^pA} 
\ee
is an isomorphism, so if there is a $\varpi^{1/p^n} \in A$, then, by applying this to $\varpi^{1/p^j}$ with $0 \le j \le n$, we get
\be \label{p-oid-mod-pin}
\xymatrix@C=35pt{A/\varpi^{1/p^n}A \ar[r]_-{\sim}^-{a\, \mapsto\, a^p} & A/\varpi^{1/p^{n - 1}}A \ar[r]_-{\sim}^-{a\, \mapsto\, a^p} & \dotsc  \ar[r]_-{\sim}^-{a\, \mapsto\, a^p} & A/\varpi^pA.} 
\ee
Conversely, by \cite{BMS18}*{Lemmas~3.9 and~3.10}, if a ring $A$ is $\varpi$-adically complete for a \emph{nonzerodivisor} $\varpi \in A$ with $\varpi^p \mid p$ (the nonzerodivisor condition is automatic if $A$ is $p$-torsion free) and \eqref{p-oid-mod-pi} holds, then $A$ is perfectoid. This gives a very practical criterion for recognizing perfectoid rings.




For a perfectoid ring $A$, both $p$ and $\xi$ are nonzerodivisors in $\bA_\Inf(A)$, so the $\bA_\Inf(A)$-module
\[
\{(x, y) \in \bA_\Inf(A)^2 \, \vert\, \xi x = p y \}/ \{(pz, \xi z)\, \vert\, z \in \bA_\Inf(A)\}
\]
is isomorphic to both $A^\flat\langle p^\flat\rangle$ and $A\langle p\rangle$
via the maps $(x, y) \mapsto x \bmod p$ and $(x, y) \mapsto y \bmod \xi$, respectively. Consequently, for every $\varpi \in A$ with $\varpi^p \mid p$ such that $A$ is $\varpi$-adically complete, \eqref{same-mod-pi} supplies an $\bA_\Inf(A)$-module isomorphism 
\be \label{tor-free-transfer}
A^\flat\langle \varpi^\flat\rangle \cong A\langle\varpi\rangle, \qx{so $A^\flat$ is $\varpi^\flat$-torsion free iff $A$ is $\varpi$-torsion free.}
\ee
After harmlessly replacing $\varpi$ by $(\varpi^\flat)^\sharp$, we apply \eqref{tor-free-transfer} with $\varpi^{1/p^n}$ in place of $\varpi$ and conclude that $A\langle \varpi \rangle$ is killed by the $\varpi^{1/p^n}$, so that 
\[
A\langle \varpi^{1/p^\infty}\rangle = A\langle \varpi \rangle = A\langle \varpi^\infty \rangle
\]
(see also \eqref{no-big-tor}).

\epp

It is useful to decompose perfectoids as follows, in the style of \cite{GR18}*{Section 16.4.18, Remark~16.4.19} or \cite{Lau18}*{Remark~8.9}.

\bpp[Canonical decompositions of perfectoids] \label{perfectoid-structure}
Let $A$ be a perfectoid ring that is $\varpi$-adically complete for a $\varpi \in A$ with $\varpi^p \mid p$, and set 
\[
\ov{A} \ce A/A\langle\varpi^\infty\rangle \qxq{and} \ov{A^\flat} \ce A^\flat/A^\flat\langle(\varpi^\flat)^\infty\rangle.
\]
By the end of \S\ref{perfectoid-def}, we have $A\langle \varpi^\infty \rangle \cap (\varpi^{1/p^\infty}) = 0$ in $A$, and analogously for $A^\flat$, so
\[
A \isomto \ov{A} \times_{\ov{A}/(\varpi^{1/p^\infty})} A/(\varpi^{1/p^\infty}) \qxq{and} A^\flat \isomto \ov{A^\flat} \times_{\ov{A^\flat}/((\varpi^\flat)^{1/p^\infty})} A^\flat/((\varpi^\flat)^{1/p^\infty}).
\]
In particular, $pA \isomto p\ov{A}$, so the ring $\ov{A}$ is $p$-adically complete.\footnote{For a ring $B$ and an ideal $I \subset B$, the property that $B$ be $I$-adically complete only depends on $I$ as a nonunital ring: it amounts to the property that every $I$-adic Cauchy sequence with values in $I$ have a unique limit in $I$.} Moreover, \eqref{same-mod-pi} and \eqref{tor-free-transfer} imply that $\ov{A}/(\varpi^{1/p^\infty}) \cong \ov{A^\flat}/((\varpi^\flat)^{1/p^\infty})$, so this quotient is a perfect $\bF_p$-algebra.  Since the functor $\bA_\Inf(-)$ preserves limits (see \S\ref{pp:tilt-adj}), we obtain the decomposition
\[
\bA_\Inf(A) \isomto \bA_\Inf(\ov{A}) \times_{\bA_\Inf(\ov{A}/(\varpi^{1/p^\infty}))} \bA_\Inf(A/(\varpi^{1/p^\infty})).
\]
By considering the counit maps $\theta$ (see \S\ref{pp:tilt-adj}) and using \S\ref{perfectoid-def} and the snake lemma, we now conclude that the generator $\xi$ of $\Ker(\theta)$ for $A$ is a nonzerodivisor in $\bA_\Inf(\ov{A})$ and that $\bA_\Inf(\ov{A})/\xi\bA_\Inf(\ov{A}) \cong \ov{A}$. In particular, by \S\ref{perfectoid-def} again, $\ov{A}$ is a perfectoid ring and $\ov{A^\flat}$ is its tilt. Thus, by \eqref{same-mod-pi},
\[
A/(\varpi^{1/p^\infty}) 
\cong (A/(\varpi))^\red \qxq{and} \ov{A}/(\varpi^{1/p^\infty}) 
\cong (\ov{A}/(\varpi))^\red.
\]
In conclusion, $A$ is a glueing of the $\varpi$-torsion free perfectoid $\ov{A}$ and the $\varpi$-torsion one $(A/(\varpi))^\red$:
\be \lab{perf-str-eq}
A \isomto \ov{A} \times_{(\ov{A}/(\varpi))^\red} (A/(\varpi))^\red
\ee
and, compatibly,
\[
A^\flat \isomto \ov{A}^\flat \times_{(\ov{A}/(\varpi))^\red} (A/(\varpi))^\red.
\]
For instance, we may choose $\varpi^p$ to be a unit multiple of $p$, in which case $\ov{A} \cong A/A\langle p^\infty\rangle$. We deduce that every perfectoid ring $A$ is reduced: \eqref{perf-str-eq} allows us to pass to $p$-torsion free $A$, and then we iteratively apply \eqref{p-oid-mod-pi} to argue that the nilradical lies in $\bigcap_{n \ge 0} p^n A = 0$. Reducedness then implies that for every $a \in A$ that has compatible $p$-power roots $a^{1/p^n} \in A$, we have 
\be \label{no-big-tor}
A\langle a^{1/p^n}\rangle = A\langle a\rangle = A\langle a^\infty\rangle \qxq{for} n \ge 0.
\ee

\epp

The decomposition \eqref{perf-str-eq} admits the following converse that is useful for recognizing perfectoids. 

\bprop \label{lem:fib-prod}
For a surjective morphism $f\colon A \surjects A'$ of perfectoid $\bZ_p$-algebras, the map 
\[
\bA_\Inf(f)\colon \bA_\Inf(A) \surjects \bA_\Inf(A')
\]
induces a surjection on the kernels of the counit maps $\theta$, more generally, it induces surjections
\be \label{eqn:Ainf-surjs}
\bA_\Inf(A) \surjects \bA_\Inf(A')\times_{A'} A \qxq{and} [p^\flat]\bA_\Inf(A) \surjects [p^\flat]\bA_\Inf(A')\times_{pA'} pA.
\ee
Moreover, for any map $B \ra A'$ with $B$ perfectoid, $C \ce A \times_{A'} B$ is perfectoid with tilt $C^\flat \cong A^\flat \times_{A^{\prime\flat}} B^\flat$.
\eprop

\bpf
First of all, the map $\bA_\Inf(f) $ inherits its indicated surjectivity from $A \surjects A'$: by completeness, this may be checked modulo $p$ and then modulo $p^\flat$, where it follows from \eqref{same-mod-pi}. Thus, letting $\xi_1$ be a generator for the kernel of $\theta$ for $A$, we apply the snake lemma and use \S\ref{perfectoid-def} to conclude that 
\[
\Ker(\bA_\Inf(f))/\xi_1\Ker(\bA_\Inf(f)) \cong \Ker(f).
\] 
The snake lemma then also shows that $\bA_\Inf(f)$ induces a surjection on the kernels of $\theta$.   Since the map $\theta\colon \bA_\Inf(A) \surjects A$ is surjective, the first surjection in \eqref{eqn:Ainf-surjs} follows. The second surjection in \eqref{eqn:Ainf-surjs} follows from the first applied to the quotients of $A$ and $A'$ by their $p$-torsion (see \S\ref{perfectoid-structure}).

As for the ring $C$, it inherits $p$-adic completeness from the $A$, $A'$, and $B$ (the unique limit of a $p$-adic Cauchy sequence is computed componentwise), so \S\ref{pp:tilt-adj} identifies the tilt and shows that 
\[
\bA_\Inf(C) \isomto \bA_\Inf(A) \times_{\bA_\Inf(A')} \bA_\Inf(B).
\] 
Thus, by \S\ref{perfectoid-def} and the second surjection in \eqref{eqn:Ainf-surjs}, there is a $\xi \in \bA_\Inf(C)$ of the form $p + [p^\flat]x$ that maps to a generator of the kernel of $\theta$ for $A$, $A'$, and $B$. A final application of the snake lemma then shows that $\bA_\Inf(C)/(\xi) \cong C$, so that $C$ is indeed perfectoid by \eqref{alt-def}.
\epf

\beg \label{eg-glue}
Consider a perfectoid ring $A$ that is $\varpi$-adically complete for a $\varpi \in A$ with $\varpi \mid p$. \Cref{lem:fib-prod} implies that for any map $B \ra (A/(\varpi))^\red$ from a perfect  $\bF_p$-algebra $B$, the ring
\[
 A \times_{(A/(\varpi))^\red} B \qx{is perfectoid.}
\]
\eeg

\bcor \label{cor:ind-etale}
For a perfectoid ring $A$ that is $\varpi$-adically complete for a $\varpi \in A$ with $\varpi^p \mid p$, the $\varpi$-adic \up{for instance, the $p$-adic, see \uS\uref{perfectoid-def}} completion  of every ind-\'{e}tale $A$-algebra is perfectoid.
\ecor

\bpf
\Cref{eg-glue} and the decomposition $A \isomto \ov{A} \times_{(\ov{A}/(\varpi))^\red} (A/(\varpi))^\red$ supplied by \eqref{perf-str-eq} reduce to $A$ being either $\varpi$-torsion free or an $\bF_p$-algebra. The $\varpi$-torsion free case follows from the criterion \eqref{p-oid-mod-pi} mentioned at the end of \S\ref{perfectoid-def}. The $\bF_p$-algebra case follows from the fact that an ind-\'{e}tale algebra over a perfect $\bF_p$-algebra is again perfect (see \cite{SGA5}*{expos\'{e}~XV, proposition~2~c)~2)}).
\epf

The perfectoids $A$ as above are more general than those in the rigid analytic approach to the theory. For instance, in the $p$-torsion free case we do not build the integral closedness of $A$ in $A[\f 1p]$ into the definitions. As we now recall, the $p$-primary aspect of this closedness is nevertheless automatic.

\bpp[$p$-integral closedness of perfectoid rings] \label{p-int-cl}
For an inclusion of rings $A \subset A'$, we recall 
that $A$ is \emph{$p$-integrally closed} in $A'$ if every $a' \in A'$ with $a'^{p} \in A$ lies in $A$. In general, the \emph{$p$-integral closure} of $A$ in $A'$, constructed as $\bigcup_{n \ge 0} A_n$ where $A_0 \ce A$ and $A_{n + 1} \subset A'$ is the $A_n$-subalgebra generated by all the $a' \in A'$ with $a'^p \in A_n$, is the smallest $p$-integrally closed subring of $A'$ containing $A$. Evidently, the $p$-integral closure lies in the integral closure of $A$ in $A'$.

The relevance of $p$-integral closedness to perfectoids was pointed out by Andr\'{e} in \cite{And18a}*{section~2.3}. For instance, 
if $\varpi \in A$ is a nonzerodivisor with $\varpi^p \mid p$ in $A$, then the map 
\be \label{p-int-cl-lem}
\tst A/\varpi A \xra{a\, \mapsto\, a^p} A/\varpi^pA
\ee
is injective if and only if $A$ is $p$-integrally closed in $A[\frac 1\varpi]$.
Indeed, the `if' direction follows from the definition, and for the `only if' one notes that the injectivity of the map ensures that any $a'=\frac{a}{\varpi^n}\in A[\frac 1\varpi] \setminus A$ with $a'^p \in A$ has its numerator $a$ divisible by $\varpi$. 

Thus, by \eqref{p-oid-mod-pi} and \eqref{perf-str-eq}, for every perfectoid $A$ that is $\varpi$-adically complete for a 
$\varpi \in A$ with $\varpi^p \mid p$, the image $\ov{A} \subset A[\f 1\varpi]$ of $A$ is $p$-integrally closed, 
so we have compatible multiplicative identifications 
\be \label{monoid-id-pi}
\tst \varprojlim_{a \mapsto a^p} (A[\f 1{\varpi}]) \isomto A^\flat[\f 1{\varpi^\flat}] \qxq{and} \xymatrix@C=40pt{\varprojlim_{a \mapsto a^p} A \ar[r]^-{\sim}_-{\eqref{monoid-id}} &A^\flat.}
\ee
The $p$-integral closedness of perfectoids has the following converse that is a variant of \cite{GR18}*{Corollary~16.9.15}. 
\epp


\bprop \label{p-oid-rec}
Let $A$ be a ring and let $\varpi \in A$ be a nonzerodivisor with $\varpi^p \mid p$ that has compatible $p$-power roots $\varpi^{1/p^n} \in A$ and is such that the map 
\[
A/\varpi A \xra{a\, \mapsto\, a^p} A/\varpi^pA
\]
is surjective. The $\varpi$-adic completion of the $p$-integral closure $\wt{A}$ of $A$ in $A[\frac 1\varpi]$ is perfectoid.
\eprop

\bpf
As in \emph{loc.~cit.},~ 
for each $a \in A$ with $a^p \in \varpi^pA$, we choose a sequence $\{a_n\}_{n \ge 0}$ in $A$ such that
\[
a_0 \ce a \q \x{and} \q a_n^p \equiv a_{n - 1} \bmod \varpi^pA \q \x{for} \q n > 0. 
\]
By construction, $a_n^{p^{n + 1}} \in \varpi^pA$, so also $\f{a_n}{\varpi^{1/p^{n}}} \in \wt{A}$, 
where, as in the statement, $\wt{A}$ is the $p$-integral closure of $A$ in $A[\f1\varpi]$. We consider the $A$-subalgebra
\[
\tst A_1 \ce A[\f{a_n}{\varpi^{1/p^{n}}}\, \vert \, n \ge 0,\, a \in A \x{ with } a^p \in \varpi^pA] \subset \wt{A}.
\]
By construction, the map $A_1/\varpi A_1 \xra{a\,\mapsto\, a^p} A_1/\varpi^pA_1$ is surjective, so we may repeat the construction with $A_1$ in place of $A$ to likewise build an $A_1$-subalgebra $A_2 \subset A$. Proceeding in this way, we obtain an $A$-subalgebra 
\[
\tst A_\infty \ce \bigcup_{i \ge 1} A_i  \subset \wt{A}
\]
for which the map $A_\infty/\varpi \xra{x\, \mapsto\, x^p} A_\infty/\varpi^p$ is both surjective and, since every $x \in A_i$ with $x^p \in \varpi^pA_i$ is divisible by $\varpi$ in $A_{i + 1}$, also injective. Thus, the observation about \eqref{p-int-cl-lem} ensures that $A_\infty = \wt{A}$, 
and \eqref{p-oid-mod-pi} then ensures that the $\varpi$-adic completion of $\wt{A}$ is perfectoid.
\epf

We turn to categorical properties of tilting that are analogues of their counterparts in the adic theory. 

\bprop \label{tilt-summary}
For a perfectoid ring $A$, there is an equivalence of categories
\[
\left\{\x{perfectoid $A$-algebras $A'$} \right\} \isomto \left\{\x{$\xi_0$-adically complete perfect $A^\flat$-algebras $B$} \right\}\!,
\]
where $\xi = (\xi_0, \xi_1, \dotsc )$ is a generator of $\Ker(\theta \colon  W(A^\flat) \surjects A)$ and the pair of inverse functors are
\[
\ \ A' \mapsto A'^\flat \ \ \x{and}\ \  B \mapsto W(B)/(\xi).
\]
Moreover, $A'^\flat$ is $\varpi^\flat$-adically complete for a $\varpi^\flat \in A'^\flat$ with $\varpi^\flat \mid \xi_0$ if and only if $A'$ is $\varpi$-adically complete for $\varpi \ce (\varpi^\flat)^\sharp$, and 
 $A'$ is a valuation ring \up{resp.,~with an algebraically closed fraction field} if and only if so is $A'^\flat$, in which case the value groups~agree\ucolon
\[
\Frac(A'^\flat)^\times/(A'^\flat)^\times \isomto \Frac(A')^\times/(A')^\times \qxq{induced by the map} x \mapsto x^\sharp.
\]
\eprop

\bpf
By \S\ref{perfectoid-def}, the functors are well-defined, inverse, and 
map $\varpi$-adically complete $A'$ to $\varpi^\flat$-adically complete $A'^\flat$. We now assume that $A'^\flat$ is $\varpi^\flat$-adically complete, so that $W(A'^\flat)$ is $[\varpi^\flat]$-adically complete, and seek to show that $\xi W(A'^\flat)$  is closed in $W(A'^\flat)$ for the $[\varpi^\flat]$-adic topology: $A'$ will then be $\varpi$-adically complete by \cite{SP}*{Lemma~\href{https://stacks.math.columbia.edu/tag/031A}{031A}}. For this, it suffices to show~that 
\[
\xi W(A'^\flat) \cap ([\varpi^\flat]^n) = \xi(W(A'^\flat) \cap ([\varpi^\flat]^{n})) \qxq{for every} n \ge 1;
\]
indeed, then the limit of every $[\varpi^\flat]$-adic Cauchy series of $W(A'^\flat)$ with terms in $\xi W(A'^\flat)$ would lie in $\xi W(A'^\flat)$. We then fix a $w \ce (w_0, w_1, \dotsc ) \in W(A'^\flat)$ with $\xi w \in ([\varpi^\flat]^n)$ and seek to show that $w \in ([\varpi^\flat]^n)$, that is, that $w_m \in (\varpi^\flat)^{np^m} A'^\flat$ for $m \ge 0$. This is clear when $\varpi^\flat = 0$, 
so \eqref{perf-str-eq} (with $A$ there equal to our $A'^\flat$) allows us to replace $A'^\flat$ by $\ov{A'^\flat} \ce A'^\flat/A'^\flat\langle(\varpi^\flat)^\infty\rangle$. 
Then $\varpi^\flat$ becomes a nonzerodivisor and induction on $n$ reduces us to $n = 1$. We fix the smallest hypothetical $m$ with 
\[
w_m \not \in (\varpi^\flat)^{p^m} \ov{A'^\flat}
\]
and, by \cite{BouAC}*{chapitre~IX, section~1, num\'{e}ro~6, lemme~4}, may assume that $w_{m'} = 0$ for $m' < m$. Then 
\[
w = V^m((w_m, w_{m + 1}, \dotsc)),
\]
so
\[
\xi w = V^m((w_m\xi_0^{p^m}, w_{m + 1}\xi_0^{p^{m + 1}} + w_m^p \xi_1^{p^m}, \dotsc)).
\]
Since $\varpi^\flat \mid \xi_0$ and $\xi_1 \in (A^\flat)^\times$, the assumption $\xi w \in ([\varpi^\flat])$ then implies that $(\varpi^\flat)^{p^{m + 1}} \mid w_m^p$, so that, by the perfectness of $\ov{A'^\flat}$, also $(\varpi^\flat)^{p^{m}} \mid w_m$, which is a desired contradiction to the existence of $m$.

If $A'$ is a valuation ring, then \eqref{monoid-id} and \eqref{same-mod-pi} show that $A'^\flat$ is a local domain in which for $a, a' \in A'^\flat$ either $a \mid a'$ or $a' \mid a$, so  $A'^\flat$ is a valuation ring. Conversely, if $A'^\flat$ is a valuation ring, then, by \eqref{same-mod-pi}, \eqref{tor-free-transfer}, and \S\ref{perfectoid-structure}, the $p$-adically complete ring $A'$ is local, $p$-torsion free unless $p = 0$ in $A'$ (in which case $A' \cong A'^\flat$), and reduced. To conclude that $A'$ is a valuation ring and also settle the claim about the value groups, we now show that every $a \in A'$ is of the form $a = ub^\sharp$ for some $u \in A'^\times$ and $b \in A'^\flat$. For this, we follow \cite{GR18}*{Proposition~16.5.50}, namely, by dividing by a power of $(p^\flat)^\sharp$, we may assume that $a$ is nonzero in $A'/pA'$, so that, by \eqref{same-mod-pi}, we have $a = b^\sharp + (p^\flat)^\sharp c$ for a $b \in A'^\flat$ that is nonzero modulo $p^\flat$ and a $c \in A'$. Since $A'^\flat$ is a valuation ring, $b$ strictly divides $p^\flat$, so it remains to set $u \ce 1 + (\frac{p^\flat}b)^\sharp c$. In the case of valuation rings of dimension $\le 1$, the remaining parenthetical assertion follows from \cite{Sch12}*{Theorem~3.7~(ii)}. 
To then deduce it for any perfectoid valuation ring $A'$ of dimension $\ge 1$, we may assume that $A'$ is of mixed characteristic $(0, p)$ and it suffices to argue that the valuation ring $A'_\fp$ that is the localization of $A'$ at the height $1$ prime $\fp \subset A'$ (concretely, at the intersection of all the primes of $A'$ containing $p$) is still perfectoid and that its tilt is the localization of $A'^\flat$ at its height $1$ prime (concretely, at the intersection of all the primes of $A'^\flat$ containing $p^\flat$).

For this last claim, the $p$-adic (resp.,~$p^\flat$-adic) topology of $A'$ (resp.,~of $A'^\flat$) is the valuation topology, so, due to \eqref{same-mod-pi}, it suffices to argue that for any valuation ring $V$ that is $a$-adically complete for some $a \in V$ and any prime ideal $\fq \subset V$ containing $a$, the localization $V_\fq$  is also $a$-adically complete.\footnote{
It is also true that if a valuation ring $V$ is complete for its valuation topology, then its localization $V_\fq$ at any prime ideal $\fq \subset V$ is also complete for its valuation topology. To see this, first note that the valuation topology is characterized by every nonzero ideal of $V$ being open, alternatively, since $a^2 V \subset a\fq \subset a V$ for $a \in \fq$, by every principal ideal of the nonunital ring $\fq$ being open. Thus, by considering Cauchy nets, we see that $V$ is complete for its valuation topology if and only if the nonunital ring $\fq$ is complete for its topology in which the principal ideals are all open. It then remains to recall that $\fq \isomto \fq V_\fq$, to the effect that replacing $V$ by $V_\fq$ does not change the nonunital ring $\fq$.
} However, by the definition of a valuation ring, every element of $V \setminus \fq$ divides every element of $\fq$, so $\fq$ maps isomorphically to $\fq V_\fq$ and kills the quotient $V_\fq/V$. It follows that in the inverse system
\[
\{0 \ra V_\fq/V \ra V/(a^n) \ra V_\fq/(a^n) \ra V_\fq/V \ra 0\}_{n > 0}
\]
of exact sequences the transition maps at the term that is the left copy of $V_\fq/V$ all vanish because they are induced by multiplication by $a$. Thus, by forming the inverse limit and applying the snake lemma we see that the $a$-adic completeness of $V$ implies that of $V_\fq$.
\epf

\brem
As the proof shows, any localization of a perfectoid valuation ring is still a perfectoid valuation ring, granted that we exclude the $0$-dimensional localization in mixed characteristic.
\erem

We will use the following further compatibilities that concern tilting. They also complement \Cref{lem:fib-prod} with additional general stability properties of perfectoid rings. 

\bprop \label{recognize}
Let $A$ be a perfectoid ring that is $\varpi$-adically complete for a $\varpi \in A$ with $\varpi^p \mid p$, let $I$ be a set, let $\{A_i \}_{i \in I}$ be  $\varpi$-adically complete perfectoid $A$-algebras, and let $S \subset A$ be a subset.
\benum
\item \label{R-0}
The $\varpi$-adic completion of $A[X_i^{1/p^\infty}]_{i \in I}$ is perfectoid and its tilt is the $\varpi^\flat$-adic completion of $A^\flat[(X_i^\flat)^{1/p^\infty}]_{i \in I}$, where $X_i^\flat$ corresponds to the $p$-power compatible sequence $(X_i^{1/p^n})_{n \ge 0}$. 

\item \label{R-a} \up{See also \cite{GR18}*{Proposition~16.3.9}}. \!\!The $\varpi$-adically completed tensor product $\wh{\bigotimes}_{i \in I} A_i$ over $A$ is perfectoid and its tilt is the $\varpi^\flat$-adically completed tensor product $\wh{\bigotimes}_{i \in I} A_i^\flat$ over $A^\flat$.

\item \label{R-b}
Suppose that the ideal $(S \bmod \varpi^n) \subset A/(\varpi^n)$ is generated by the $p^n$-th powers of its elements for $n > 0$
\up{for instance, that each $s \in S$ has a root $s^{1/p^N} \in S$ with $N > 0$}. Then the $\varpi$-adic completion of $A/(S)$ is perfectoid and its tilt is the $\varpi^\flat$-adic completion of $A^\flat/(S^\flat)$ where 
\[
\tst\qq S^\flat \ce \varprojlim_{a\, \mapsto\, a^p} (S \bmod \varpi) \subset \varprojlim_{a\, \mapsto\, a^p} A/(\varpi) \cong A^\flat.
\]

\m \label{R-c}
A product $\prod_{i \in I} B_i$ of $\bZ_p$-algebras is perfectoid iff so is each $B_i$, and then 
\[
\tst\qq (\prod_{i \in I} B_i)^\flat \cong \prod_{i \in I} B_i^\flat.
\]

\m \label{R-e}\up{See also \cite{GR18}*{Theorem~16.3.76} and \cite{And20}*{Proposition~2.2.1}}.
For $a_1^\flat, \dotsc, a_r^\flat \in A^\flat$ and $a_j \ce (a_j^\flat)^\sharp \in A$, the $(a_1, \dotsc, a_r)$-adic completion 
of $A$ is perfectoid, agrees with the derived $(a_1, \dotsc, a_r)$-adic completion of $A$, and has the $(a_1^\flat, \dotsc, a_r^\flat)$-adic completion 
of $A^\flat$ as its tilt.
\eenum
\eprop

\bpf 
For \ref{R-0}, we first note that, by \Cref{tilt-summary}, the $A$-algebra 
\[
W((A^\flat[(X_i^\flat)^{1/p^\infty}]_{i \in I})\wh{\ }\, )/(\xi),
\]
where the completion is $\varpi^\flat$-adic, is perfectoid. It then remains to note that, since $p^n \in (\xi, [\varpi^\flat]^n)$, the map that sends each $X_i^{1/p^m}$ to $((X_i^\flat)^{\sharp})^{1/p^m}$ exhibits it as the $\varpi$-adic completion of $A[X_i^{1/p^\infty}]_{i \in I}$. 

For \ref{R-a}, a tensor product indexed by $I$ is defined as the direct limit of subproducts over the finite subsets of $I$  and is a categorical coproduct. With tensor products over $W_n(A^\flat)$ and $A^\flat$, 
\be \label{WW-trick-map}
\tst \bigotimes_{i \in I} W_n(A_i^\flat) \isomto W_n(\bigotimes_{i \in I} A_i^\flat)
\ee
because both the source and the target are initial among the $\bZ/p^n\bZ$-algebras whose reduction modulo $p$ is equipped with a map from $ \bigotimes_{i \in I} A_i^\flat$ (
see, for instance, \cite{SZ18}*{Proposition 3.12}). By reducing \eqref{WW-trick-map} modulo 
\[
(p^n, p^{n - 1}[\varpi^\flat]^{p}, \dotsc, p [\varpi^\flat]^{p^{n - 1}}, [\varpi^\flat]^{p^{n}}),
\]
we obtain
\[
\tst \bigotimes_{i \in I} W_n(A_i^\flat/([\varpi^\flat]^{p^{n}})) \isomto W_n(\bigotimes_{i \in I} (A_i^\flat/([\varpi^\flat]^{p^{n}}))),
\]
so, since $(p, [\varpi^\flat])^{p^n} \subset (p^n, p^{n - 1}[\varpi^\flat]^{p}, \dotsc, p [\varpi^\flat]^{p^{n - 1}}, [\varpi^\flat]^{p^{n}}) \subset (p, [\varpi^\flat])^n$, also 
\be \label{WW-map}
\tst  \wh{\bigotimes}_{i \in I} W(A_i^\flat) \isomto W(\wh{\bigotimes}_{i \in I} A_i^\flat)
\ee
where the first $\wh{\tensor}$ is $(p, [\varpi^\flat])$-adically completed and over $W(A^\flat)$. 
Since $\varpi^\flat \mid \xi_0$ in $A^\flat$ for a generator $\xi = (\xi_0, \xi_1, \dotsc)$ of $\Ker(\theta\colon W(A^\flat) \surjects A)$ (see \S\ref{perfectoid-def}), the perfect $A^\flat$-algebra $\wh{\bigotimes}_{i \in I} A_i^\flat$ is $\xi_0$-adically complete. Thus, by \Cref{tilt-summary}, 
the reduction of  \eqref{WW-map} modulo $\xi$ is a map of $\varpi$-adically complete perfectoids, so, since $(\xi, [\varpi^\flat]) = (p, [\varpi^\flat])$ in $W(A^\flat)$ (see \eqref{same-mod-pi}), it is 
the desired
\[
\tst \wh{\bigotimes}_{i \in I} A_i \isomto (W(\wh{\bigotimes}_{i \in I} A_i^\flat))/(\xi).
\]
For \ref{R-b}, 
by the assumption on $S$ and by construction, $S^\flat$ surjects onto 
\[
(S \bmod \varpi) \subset A/(\varpi) \cong A^\flat/(\varpi^\flat)
\]
and is stable under $p$-power roots. Moreover, the $p^n$-th power of an element of $A/(\varpi^{n + 1})$ depends only on its residue class modulo $\varpi$, so 
\[
(S \bmod \varpi^n) = ((S^\flat)^\sharp \bmod \varpi^n) \qxq{in} A/(\varpi^n).
\]
Thus we lose no generality by assuming that $S = (S^\flat)^\sharp$, in other words, that every $s \in S$ admits a $p$-th root $s^{1/p} \in S$. Both $A^\flat/(S^\flat)$ and its $\varpi^\flat$-adic completion $\wh{A^\flat/(S^\flat)}$ are perfect $A^\flat$-algebras, and 
\[
W_n(A^\flat)/([S^\flat]) \isomto W_n(A^\flat/(S^\flat)).
\]
Thus, by \Cref{tilt-summary}, the $A$-algebra $W(\wh{A^\flat/(S^\flat)})/(\xi)$ is a $\varpi$-adically complete perfectoid. In conclusion, since $p^n \in (\xi, [\varpi^\flat]^n)$ and $S = (S^\flat)^\sharp$, the following map exhibits its perfectoid target as the $\varpi$-adic completion of the source:
\[
A/(S) \cong (W(A^\flat)/(\xi))/(S) \cong (W(A^\flat)/([S^\flat]))/(\xi) 
\ra W(\wh{A^\flat/(S^\flat)})/(\xi).
\]
Part \ref{R-c} is immediate from the definition of \S\ref{perfectoid-def} because 
\[
\tst\bA_\Inf(\prod_{i \in I} B_i) \cong \prod_{i \in I} \bA_\Inf(B_i).
\] 
For \ref{R-e}, it suffices to argue that the derived $(a_1, \dotsc, a_r)$-adic completion $\wh{A}$ of $A$ is perfectoid (so, in particular, is a classical ring) and that its tilt is the derived $(a_1^\flat, \dotsc, a_r^\flat)$-adic completion $\wh{A^\flat}$ of $A^\flat$. Indeed, this will imply the claimed agreement with the usual $(a_1, \dotsc, a_r)$-adic completion (and likewise for $A^\flat$): perfectoid rings are reduced (see \S\ref{perfectoid-structure}), so \cite{SP}*{Lemma~\href{https://stacks.math.columbia.edu/tag/0G3I}{0G3I}} will ensure that $\wh{A}$ is $(a_1, \dotsc, a_r)$-adically separated, and hence, by \cite{SP}*{Proposition~\href{https://stacks.math.columbia.edu/tag/091T}{091T}}, even $(a_1, \dotsc, a_r)$-adically complete, so the map $A \ra \wh{A}$ will be initial among maps to $(a_1, \dotsc, a_r)$-adically complete $A$-algebras. 

The derived $(a_1, \dotsc, a_r)$-adic completion of $A$ agrees with the iterated derived $a_i$-adic completion for $i = 1, \dotsc, r$, so we lose no generality by assuming that $r = 1$ and renaming $a \ce a_1$ and $a^\flat \ce a_1^\flat$. Since $A^\flat$ is a perfect $\bF_p$-algebra, $A^\flat[a^\flat] = A^\flat[(a^\flat)^{1/p^\infty}]$, to the effect that the inverse system $\{ A^\flat[(a^\flat)^n]\}_{n > 0}$ is almost zero. Thus, the derived $a^\flat$-adic completion $\wh{A^\flat}$ of $A^\flat$ agrees with the classical $a^\flat$-adic completion of $A^\flat$. In particular, $\wh{A^\flat}$ is a perfect $\bF_p$-algebra that inherits derived $\xi_0$-adic completeness from $A^\flat$. Thus, $\wh{A^\flat}$ is reduced and we conclude as in the previous paragraph that it is $\xi_0$-adically complete. This already settles the positive characteristic case, in which $A = A^\flat$.

By arguing via Witt vector coordinates, we see that each $W_n(\wh{A^\flat})$ is $[a^\flat]$-adically complete, so that $W(\wh{A^\flat})$ is also $[a^\flat]$-adically complete. Moreover, the derived $[a^\flat]$-adic completion $\wh{W(A^\flat)}$ of $W(A^\flat)$ inherits derived $p$-adic completeness and its derived reduction modulo $p$ is the derived $[a^\flat]$-adic completion $\wh{A^\flat}$ of $A^\flat$. Thus, we may check on derived reductions modulo $p$ that 
\[
\wh{W(A^\flat)} \isomto W(\wh{A^\flat}).
\]
However, \S\ref{perfectoid-def} ensures that $\xi$ is a nonzerodivisor in $W(\wh{A^\flat})$, so this isomorphism shows that $W(\wh{A^\flat})/(\xi)$ is the derived $a$-adic completion $\wh{A}$ of $A$ and, simultaneously, that $\wh{A}$ is a classical ring. To then conclude that $\wh{A}$ is perfectoid with tilt $\wh{A^\flat}$ it remains to review \S\ref{perfectoid-def}.
\epf




The following proposition is sometimes useful for reducing to $p$-torsion free perfectoid rings. 


\bprop \lab{lem:perfectoid-quotient-of-torsion-free}
Every perfectoid ring $A$ that is $\varpi$-adically complete for a $\varpi \in A$ with $\varpi^p \mid p$ is a quotient of a  perfectoid ring $\wt{A}$ that is \emph{$\wt{\varpi}$-torsion free} and $\wt{\varpi}$-adically complete for a lift $\wt{\varpi} \in \wt{A}$ of $\varpi$ with $\wt{\varpi}^p \mid p$. In addition, every perfectoid ring is a quotient of a $p$-torsion free perfectoid ring. 
\eprop

\bpf
As in \eqref{alt-def}, we write $A \cong W(A^\flat)/(\xi)$ with $\xi = (\xi_0, \xi_1, \dotsc)$ in Witt coordinates such that $\xi_1 \in (A^\flat)^\times$ and $A^\flat$ is $\xi_0$-adically complete. In fact, $A^\flat$ is even $\varpi^\flat$-adically complete and we fix a choice of $\varpi^\flat \in A^\flat$, so that $\varpi = (\varpi^\flat)^\sharp u $ with $u \in A^\times$ (see \S\ref{perfectoid-def}). 
We consider the perfect $\bF_p$-algebra 
\[
\tst B_0 \ce \bF_p[X_a^{1/p^\infty}\, \vert\, a \in A^\flat] \qxq{and the surjection} B_0 \surjects A^\flat \qxq{given by} 
X_a\, \mapsto\, a.
\]
Since $(\varpi^\flat)^p \mid \xi_0$ in $A^\flat$ (see \S\ref{perfectoid-def}), we may lift the $\xi_i$ to $\wt{\xi}_i \in B_0$ with $(X_{\varpi^\flat})^p \mid \wt{\xi}_0$ 
and let $B$ be the $X_{\varpi^\flat}$-adic completion of $B_0[\f{1}{\wt{\xi}_1}]$. Certainly, $B$ is a $\wt{\xi_0}$-adically complete (see \cite{SP}*{Lemma~\href{https://stacks.math.columbia.edu/tag/090T}{090T}}), perfect $\bF_p$-algebra equipped with a surjection $B \surjects A^\flat$. Letting $\wt{\xi} \in W(B)$ be defined by its Witt coordinates $\wt{\xi}_i$ and using \eqref{alt-def}, we obtain a surjection of perfectoid rings 
\[
\wt{A}' \ce W(B)/(\wt{\xi}) \surjects W(A^\flat)/(\xi) \cong A, \qxq{and we set} \wt{\varpi}' \ce (X_{\varpi^\flat})^\sharp \in \wt{A}'. 
\]
Since $(X_{\varpi^\flat})^p \mid \wt{\xi}_0$ in $B$, we have $(\wt{\varpi}')^p \mid p$ in $\wt{A}'$ (see \S\ref{perfectoid-def}), so \Cref{tilt-summary} ensures that $\wt{A}'$ is $\wt{\varpi}'$-adically complete, and \eqref{tor-free-transfer}  then  ensures that $\wt{A}'$ is $\wt{\varpi}'$-torsion free. We lift $u \in A$ to a $\wt{u} \in \wt{A}'$, and we let $\wt{A}$ be the $\wt{\varpi}'$-adic completion  of $\wt{A}'[\f 1{\wt{u}}]$, so that we have the induced surjection $\wt{A} \surjects A$ and the lift $\wt{\varpi} \ce \wt{\varpi}' \wt{u} \in \wt{A}$ of $\varpi$. By \Cref{cor:ind-etale}, the ring $\wt{A}$ is perfectoid and, by construction, it is $\wt{\varpi}$-adically complete and $\wt{\varpi}$-torsion free. The proof of the $p$-torsion free variant 
is similar but simpler: it suffices to choose $\wt{\xi}_0 \ce X_{\xi_0}$, replace $X_{\varpi^\flat}$ by $X_{\xi_0}$ in the subsequent argument, and set $\wt{A} \ce \wt{A}'$.
\epf






\csub[Tilting \'{e}tale cohomology algebraically] \label{sec:tilting-cohomology}

Guided by the idea that comparing a perfectoid ring $A$ and its tilt $A^\flat$ is close in spirit to an Elkik-type comparison of a Henselian ring and its completion, in \Cref{baby-tilting} we exhibit ``algebraic'' incarnations of the paradigm that tilting preserves topological information, specifically, idempotents (that is, clopen subschemes) and \'etale cohomology. The idea of the proof is that the  idempotent case is pretty much immediate from \eqref{monoid-id} with \eqref{monoid-id-pi} and, by $p$-complete arc descent, it implies the assertion about the \'{e}tale cohomology. This style of argument bypasses any recourse to adic spaces, although, of course, the conclusion is not as strong as an equivalence of \'{e}tale sites.

\bpp[The $I$-complete arc-topology] \label{pp:arc-topology}\label{eg:fflat-is-arc}
We recall from \cite{BM20}*{Definition 1.2} that a ring map $A \ra A'$ is an \emph{arc cover} if any $A \ra V$ with $V$ a valuation ring of dimension $\le 1$ fits in a commutative diagram
\be\label{test-diagram}\ba
\xymatrix{
A \ar[r] \ar[d] & A' \ar@{-->}[d] \\
V  \ar@{-->}[r] & V' 
}
\ea\ee
in which $V'$ is a valuation ring of dimension $\le 1$ and $V \ra V'$ is faithfully flat (that is, an extension of valuation rings). For a fixed finitely generated ideal $I \subset A$ (example: $I = (p)$), if the same holds whenever $V$ is, in addition, $I$-adically complete, 
then $A \ra A'$ is an \emph{$I$-complete arc cover} (called a \emph{$\varpi$-complete arc cover} when $I = (\varpi)$ is principal). 
An arc cover is simply a $0$-complete arc cover, and 
an $I$-complete arc cover is an $I'$-complete arc cover whenever $I \subset I'$ (see \cite{SP}*{Lemma~\href{https://stacks.math.columbia.edu/tag/090T}{090T}}). In particular, for every $I$, an arc cover is an $I$-complete arc cover and the reduction modulo $I$ of an $I$-complete arc cover is an 
arc cover. 

In fact, there is no need to assume that $V'$ be of dimension $\le 1$: one can arrange $\dim(V) = \dim(V')$ \emph{a posteriori} by the argument of \cite{BM20}*{Proposition~2.1}.
In addition, by extending $V'$ (of dimension $\le 1$) to a valuation ring of dimension $\le 1$ on the algebraic closure of $\Frac(V')$ (see \cite{BouAC}*{chapitre~VI, section 8, num\'{e}ro~6, proposition~6}) and, in the case of $I$-complete arc covers, $I$-adically completing (which preserves algebraic closedness, see \cite{BGR84}*{Section 3.4, Proposition~3}), we may restrict to those $V'$ of dimension $\le 1$ in  \eqref{test-diagram}  that have an algebraically closed fraction field and, in the case of $I$-complete arc covers, are $I$-adically complete. Similarly, one then loses no generality by assuming that $\Frac(V)$ be algebraically closed. 

For example,
\benuma
\m \label{arc-fflat}
any faithfully flat $A \ra A'$ is an arc cover: to see this, we may assume that $A = V$, lift the specialization of points in $\Spec(V)$ to $\Spec(A')$ (see \cite{EGAIV2}*{proposition 2.3.4~(i)}), and use the maximality of valuation rings with respect to domination; 

\m \label{arc-pi-fflat}
any $A \ra A'$ that is faithfully flat modulo powers of a finitely generated ideal $I \subset A$ is an $I$-complete arc cover: we may assume that $A = V$ for an $I$-adically complete valuation ring $V$ of rank $\le 1$, assume that $A'$ is $I$-adically complete, and use the resulting injectivity $V \hra A'$ to lift the specialization (due to the $I$-adic completeness of $A'$, the closure of the generic $V$-fiber of $\Spec(A')$ meets the closed $V$-fiber, so one applies \cite{SP}*{Lemma~\href{https://stacks.math.columbia.edu/tag/0903}{0903}} to conclude); 

\m \label{arc-integral}
any integral $A \ra A'$ that is surjective on spectra is an arc cover: now one uses going up to lift the specialization (see \cite{SP}*{Lemma~\href{https://stacks.math.columbia.edu/tag/00GU}{00GU}}). 
\eenum
 
As the name suggests, finite collections of ring maps $\{A \ra A'_i \}_{i \in I}$ for which $A \ra \prod_{i \in I} A'_i$ an arc cover (resp.,~an $I$-complete arc cover) are the covering maps for a Grothendieck topology on commutative rings (resp.,~on $\bZ[x_1, \dotsc, x_n]$-algebras where $I = (\varpi_1, \dotsc, \varpi_n)$ with $x_i \mapsto \varpi_i$), the \emph{arc topology} (resp.,~the \emph{$I$-complete arc topology}). On perfectoids this topology is insensitive to tilting as follows. 
\epp

\blem \label{arc-tilt}
Let $A \ra A'$ be a map of perfectoid rings, let $\varpi \in A$ with $\varpi^p \mid p$ be such that $A$ and $A'$ are $\varpi$-adically complete, and let $\varpi^\flat \in A^\flat$ be such that $(\varpi^\flat)^\sharp$ is a unit multiple of $\varpi$ \up{see \uS\uref{perfectoid-def}}. Then $A \ra A'$ is a $\varpi$-complete arc cover if and only if its tilt $A^\flat \ra A'^\flat$ is a $\varpi^\flat$-complete arc cover.
\elem

\bpf
By \S\ref{pp:arc-topology}, the condition of being a $\varpi$-complete (resp.,~$\varpi^\flat$-complete) arc cover may be phrased to only involve maps to $\varpi$-adically (resp.,~$\varpi^\flat$-adically) complete valuation rings of dimension $\le 1$ with algebraically closed fraction fields, and such are perfectoid by \eqref{p-oid-mod-pi}. It then remains to recall from \Cref{tilt-summary} that the tilting equivalence identifies such valuation rings, respects their dimensions, and matches $\varpi$-adic completeness with $\varpi^\flat$-adic completeness. 
\epf

We will exploit the following convenient base of the $\varpi$-complete arc topology. 

\blem \label{lem:arc-base}
Every ring $A$ \up{resp.,~with a $\varpi \in A$} has an arc \up{resp.,~a $\varpi$-complete arc} cover $A \ra \prod_{i \in I} V_i$ whose factors $V_i$ are valuation rings \up{resp.,~$\varpi$-adically complete valuation rings} of dimension $\le 1$ with algebraically closed fraction fields. 
\elem
 
\bpf 
For each prime $\fp \subset A$, choose an algebraic closure $\ov{k(\fp)}$ of the residue field at $\fp$. Let $\fp$ vary and let $I$ be the set of valuation subrings $V_i \subset \ov{k(\fp)}$ of dimension $\le 1$ containing the image of $A$ and with $\Frac(V_i) = \ov{k(\fp)}$. To check that the resulting $A \ra \prod_{i \in I} V_i$ is a desired arc cover, we note that, by the choice of $I$, any map $A \ra V$ to a valuation ring of dimension $\le 1$ with an algebraically closed fraction field factors through some $A \ra V_i$ and use \S\ref{pp:arc-topology}.
For the $\varpi$-complete arc aspect, it suffices to instead take  $A \ra \prod_{i \in I} \wh{V}_i$ where $\wh{V}_i$ is the $\varpi$-adic completion of $V_i$ (see \S\ref{pp:arc-topology}). 
\epf

As we now verify, the arc covers constructed in the previous lemma have no nonsplit \'{e}tale covers. 

\blem \label{lem:ultrafilters}
Let $\{V_i\}_{i \in I}$ be valuation rings. The connected components of $\Spec(\prod_{i \in I} V_i)$ are the 
\[
\tst \Spec(\prod_{\sU} V_i) \qxq{for ultrafilters} \sU \qxq{on} I
\]
\up{where $\prod_\sU \ce \varinjlim_{I' \in \sU} \prod_{i \in I'}$}. In particular, if all the $\Frac(V_i)$ are algebraically closed, then each quasi-compact open  $U \subset \Spec(\prod_{i \in I} V_i)$ has no nonsplit \'{e}tale covers and its connected components are spectra of valuation rings with algebraically closed fraction fields. 
\elem

\bpf
For $I' \subset I$, let $e_{I'} \in \prod_{i \in I} V_i$ be the idempotent whose coordinates at $I'$ (resp.,~at $I \setminus I'$) are $0$ (resp.,~$1$). Since $e_{I' \cap I''} = e_{I'} + e_{I''} - e_{I' \cup I''}$, for any prime ideal $\fp \subset \prod_{i \in I} V_i$, the set $\sU_\fp \ce \{ I' \, \vert\, e_{I'} \in \fp\}$ is an ultrafilter on $I$. The assignment $\fp \mapsto \sU_\fp$ gives a continuous map 
\[
\tst\Spec(\prod_{i \in I} V_i) \ra \gB I
\] 
to the Stone--\v{C}ech compactification of $I$: indeed, the sets 
\[
U_{I'} \ce \{ \sU \, \vert \, I' \in \sU\} \subset \gB I \qxq{with} I' \subset I
\]
are a base of opens for $\gB I$, and the preimage of $U_{I'}$ is the open 
\[
\tst\Spec(\prod_{i \in I'} V_i) \subset \Spec(\prod_{i \in I} V_i).
\] 
Any $\sU \in \gB I$ is the intersection of its neighborhoods $U_{I'}$ for $I' \in \sU$, so the preimage of $\sU$ is precisely $\Spec(\prod_{\sU} V_i)$. 

The preceding paragraph works for any rings $\{V_i\}_{i \in I}$; however, if the $V_i$ are valuation rings, then $\prod_{\sU} V_i$ is a valuation ring with the fraction field $\prod_{\sU} \Frac(V_i)$, and the latter is algebraically closed whenever so are all the $\Frac(V_i)$. Thus, since $\gB I$ is totally disconnected and every quasi-compact open of some $\Spec(\prod_{\sU} V_i)$ is the spectrum of a valuation ring (see \cite{SP}*{Lemma~\href{https://stacks.math.columbia.edu/tag/088Y}{088Y}}), the connected component aspects of the claim follow. Moreover, if the $\Frac(V_i)$ are all algebraically closed, then, by the above, the local rings of $U$ are strictly Henselian. A limit argument then shows that every \'{e}tale cover of $U$ may be refined by a Zariski cover, and, thanks to \cite{SP}*{Lemma~\href{https://stacks.math.columbia.edu/tag/0968}{0968}}, the latter has a section. 
\epf

Our approach to tilting \'{e}tale cohomology builds on the following arc descent result of Bhatt--Mathew. 

\begin{theorem}\label{arc-descent}
Let $A$ be a ring \up{resp.,~with a finitely generated ideal $I \subset A$} and let $\sF$ be a torsion sheaf on $A_\et$. On the category of $A$-algebras $A'$, the functor
\[
\tst A' \mapsto R\Gamma_\et(A', \sF)
\]
{\upshape(}resp., 
\[
A' \mapsto R\Gamma_\et(\Spec(\wh{A'})\setminus V(I), \sF), \qx{where the completion is $I$-adic})
\]
satisfies hyperdescent in the arc \up{resp.,~$I$-complete arc} topology, and the functor 
\[
A' \mapsto R\Gamma_\et(A'^h, \sF), \qx{where $(-)^h$ denotes the $I$-Henselization,}
\]
satisfies hyperdescent in the $I$-complete arc topology.
\end{theorem}

\begin{proof} 
Since the functors in question are bounded below, descent for them implies hyperdescent, so we focus on arguing descent. 
By \cite{Gab94a}*{Theorem~1}, we have 
\[
R\Gamma_\et(A', \sF) \cong R\Gamma_\et(A'/IA', \sF)
\]
for $I$-Henselian $A'$, so the last assertion follows from the rest and \S\ref{pp:arc-topology}. Moreover, the descent claims were settled in \cite{BM20}*{Theorem~5.4} and, respectively, \cite{Mat20}*{Remark~5.18} with \cite{BM20}*{Corollary~6.17}, except that for the functor 
\[
A' \mapsto R\Gamma_\et(\Spec(\wh{A'})\setminus V(I), \sF)
\]
\emph{loc.~cit.}~used the arc$_I$ topology instead. The latter is the variant of the arc topology in which in \eqref{test-diagram} one requires $I$ to map to nonzero subideals of the maximal ideals of the  valuation rings $V$ and $V'$ of dimension $\le 1$. By replacing such $V$ by its $I$-adic completion, we see that every $I$-complete arc cover is an arc$_I$ cover, so our descent claim follows.
\end{proof}

\brem \label{BM-hens-comp}
By \cite{BM20}*{Theorem~6.11} (or \cite{ILO14}*{expos\'{e} XX, section 4.4}), we have 
\[
\tst 
R\Gamma_\et(\Spec(A'^h)\setminus V(I),\sF) \isomto R\Gamma_\et(\Spec(\widehat{A'}) \setminus V(I),\sF),
\]
where $(-)^h$ denotes the $I$-Henselization, so the following functor also satisfies $I$-complete arc hyperdescent on $A$-algebras $A'$:
\[
A' \mapsto R\Gamma_\et(\Spec(A'^h)\setminus V(I), \sF).
\]
\erem

We are ready for the promised algebraic approach to tilting \'{e}tale cohomology of perfectoids. 

\bthm \label{baby-tilting}
Let $p$ be a prime, let $A$ be a ring, let $\varpi \in A$ with $\varpi^p \mid p$ be such that $A$ is $\varpi$-Henselian, has bounded $\varpi^\infty$-torsion, and its $\varpi$-adic completion is perfectoid, and let
\[
\tst \Spec(A[\f 1\varpi]) \subset U \subset \Spec(A) \qxq{and} \Spec(A^\flat[\f{1}{\varpi^\flat}]) \subset U^\flat \subset \Spec(A^\flat)
\]
be opens  whose complements agree via \eqref{same-mod-pi} with $A^\flat \ce \varprojlim_{a \mapsto a^p} (A/\varpi A)$. There are identifications of sets of idempotents
\[
\Idem(U) \cong \Idem(U^\flat) \qx{compatibly with orthogonality,}
\]
and of \'{e}tale cohomology 
\[
R\Gamma_\et(U, G) \cong R\Gamma_\et(U^\flat, G) \qx{for every torsion abelian group $G$,}
\]
functorial in $A$, $U$, and $G$. In particular, for a closed $Z \subset \Spec(A/\varpi A)$ and a torsion abelian~group~$G$,
\be \label{eqn:baby-tilting-supports}
R\Gamma_Z(A, G) \cong R\Gamma_Z(A^\flat, G).
\ee
\ethm

\bpf
The claim about \eqref{eqn:baby-tilting-supports} follows from the rest and the cohomology with supports triangle. 

By, for instance, \cite{BC19}*{Theorems 2.3.1 and 2.3.4}, base change to the $\varpi$-adic completion of $A$ changes neither $\Idem(U)$ nor $R\Gamma_\et(U, G)$, so we assume that $A$ is $\varpi$-adically complete and, in particular, perfectoid. The $p$-power map of any ring induces the identity map on the set of idempotents of that ring, so the claim about $\Idem(U)$ when $U$ is either $\Spec(A)$ or $\Spec(A[\f 1\varpi])$ follows from the functorial, compatible, multiplicative isomorphisms \eqref{monoid-id} and \eqref{monoid-id-pi}, namely, from
\[
\tst \varprojlim_{a \mapsto a^p} A \cong A^\flat \qxq{and} \varprojlim_{a \mapsto a^p} (A[\f 1\varpi]) \cong A^\flat[\f 1 {\varpi^\flat}].
\]
For a general $U$, by glueing and limit arguments, giving an idempotent on $U$ amounts to giving an idempotent $e$ on $A[\f 1\varpi]$ together with a compatible under pullback system of idempotents $e_B$ on the localizations $\wt{B}$ of $A$ along variable principal affine opens $\Spec(B) \subset U_{A/\varpi} \subset \Spec(A/\varpi)$ subject to the condition that after inverting $\varpi$ each $e_B$ agrees with the pullback of $e$. Moreover, by Beauville--Laszlo glueing \cite{SP}*{Lemma~\href{https://stacks.math.columbia.edu/tag/0BNR}{0BNR}},\footnote{In this proof, one may avoid the Beauville--Laszlo glueing by replacing $\wt{B}$ by its $\varpi$-Henselization and using \cite{BC19}*{Theorem~2.3.1} again (resp.,~\cite{BC19}*{Theorem~2.3.4} for $R\Gamma_\et$ in place of $\Idem$), but this comes at the expense of having to consider self-intersections in the limit arguments. The use of the Beauville--Laszlo technique was suggested by Arnab Kundu.} in this description we may replace $\wt{B}$ by its $\varpi$-adic completion. By \Cref{cor:ind-etale}, this completion is perfectoid and, by \eqref{same-mod-pi}, its tilt is the $\varpi^\flat$-adic completion of 
the localization of $A^\flat$ along $\Spec(B)$. Thus, the analogous description of the idempotents on $U^\flat$ and the settled cases $U = \Spec(A)$ and $U = \Spec(A[\f 1\varpi])$ 
give the desired functorial identification $\Idem(U) \cong \Idem(U^\flat)$ that is compatible with orthogonality.  

The analogous glueing (or descent) argument carried out with $R\Gamma_\et$ in place of $\Idem$, 
which this time uses formal glueing for \'{e}tale cohomology in place of Beauville--Laszlo glueing to pass to completions, so, concretely, it uses \Cref{arc-descent} and \cite{BM20}*{Theorem~6.4},
 reduces us to exhibiting compatible identifications 
\[
R\Gamma_\et(U, G) \cong R\Gamma_\et(U^\flat, G)
\]
in the cases when $U = \Spec(A)$ or $U = \Spec(A[\f 1\varpi])$ (functorially in $A$ and $G$). For this, we first treat the case when $A = \prod_{i \in I} V_i$ for $\varpi$-adically complete valuation rings $V_i$ over $A$ with algebraically closed fraction fields  (
such $V_i$ are perfectoid by \eqref{p-oid-mod-pi} and hence, by \Cref{recognize}~\ref{R-c}, so is $\prod_{i \in I} V_i$). For such $A$, we have $A^\flat \cong \prod_{i \in I} V_i^\flat$. Thus, \Cref{lem:ultrafilters} implies that $A$ and $A^\flat$, as well as $A[\f 1\varpi]$ and $A^\flat[\f 1{\varpi^\flat}]$, have no nonsplit \'{e}tale covers. In particular, both $R\Gamma_\et(U, G)$ and $R\Gamma_\et(U^\flat, G)$ are concentrated in degree zero where they are given by locally constant $G$-valued functions on $U$ and $U^\flat$, respectively. Due to the functorial identification $\Idem(U) \cong \Idem(U^\flat)$, the clopens of $U$ are in a functorial bijection with those of $U^\flat$, compatibly with the relation of  disjointness (which amounts to orthogonality of the corresponding idempotents). Thus, the spaces of locally constant $G$-valued functions on $U$ and $U^\flat$ are functorially identified, which settles the case when $A = \prod_{i \in I} V_i$ as above. 

By \Cref{lem:arc-base}, the $A$ that are products $\prod_{i \in I} V_i$ as above with each $V_i$ of rank $\le 1$ form a base of the $\varpi$-complete arc topology of $A$. By \Cref{tilt-summary} and \Cref{arc-tilt}, tilting matches this base with its analogue for the $\varpi^\flat$-complete arc topology of $A^\flat$. Thus, to deduce the remaining case of a general $A$, it remains to combine the already functorially settled case $A = \prod_{i \in I} V_i$ with the $\varpi$-complete arc descent supplied by \Cref{arc-descent}.
\epf


\csub[The ind-syntomic generalization of Andr\'{e}'s lemma] \label{And-lem}

Andr\'{e}'s lemma, which originated in \cite{And18a,And18b}, says that up to passing to a perfectoid cover elements of a perfectoid admit compatible $p$-power roots. This is useful for constructing perfectoids above a Noetherian local ring $(R, \fm)$ with $\Char(R/\fm) = p$ beyond regular $R$: one writes $\wh{R}$ as a quotient of a regular ring, chooses a faithfully flat perfectoid cover of the latter (as in \Cref{lem:towers} below), uses Andr\'{e}'s lemma to ensure that the equations cutting out $R$ have compatible $p$-power roots, and then kills these roots (the relevance of such roots 
is seen already in \Cref{recognize}~\ref{R-b}). This mechanism is how we will use Andr\'{e}'s lemma in the proof of \Cref{main}.

The goal of this section is to present a generalization of Andr\'{e}'s lemma stated in \Cref{Andre-lem} below. 
More precisely, in Andr\'{e}'s work the refining perfectoid cover was almost faithfully flat modulo powers of $p$ (see, for instance, \cite{Bha18a}*{Theorem~1.5}), which was improved to actual faithful flatness by Gabber--Ramero in \cite{GR18}*{Theorem~16.9.17}  at the cost of ``decompleting.'' We follow their method to improve further to ind-syntomicity and to eliminate torsion freeness assumptions. Ind-syntomicity modulo powers of $p$ was achieved in \cite{BS19}*{Theorem~7.14, Remark~7.15} by a different 
argument and, as we explain in the proof of \Cref{main-pf}, 
suffices for our purposes, so a pragmatic reader could skip this section. 




We begin with the following ``integral'' variant of the approximation lemma \cite{Sch12}*{Corollary~6.7~(i)}. 

\blem \label{approx-lem}
Let $A$ be a perfectoid ring, let $\varpi \in A$ with $\varpi^p \mid p$ be such that $A$ is $\varpi$-adically complete, let $a \in A$, and let $m \ge 0$. There is an $a' \in A^\flat$ 
such that for every continuous valuation $\abs{\cdot}$ on  $A$  with $\abs{A} \le 1$ we have
\be \label{eqn:AL-prel}
|a| \le | \varpi^{pm}| \qxq{if and only if}  |a'^\sharp| \le |\varpi^{pm}|,
\ee
more generally, such that for every $\abs{ \cdot}$ above we have
\be \label{AL-eq}
\tst \abs{a - a'^\sharp} \le \abs{p} \cdot \max(|a'^\sharp|, \abs{\varpi}^{pm}).
\ee
\elem

Here continuity means that $\abs{\varpi^n}$ for $n \ge 0$ becomes smaller than any element of the value group. 

\bpf
For completeness, we give a proof; see \cite{KL15}*{Corollary~3.6.7} and \cite{GR18}*{Corollary~16.6.26} for other variants. We loosely follow the argument from \cite{KL15} whose main inputs are \cite{Ked13}*{Lemmas~5.5 and 5.16}.

We focus on \eqref{AL-eq} because it implies \eqref{eqn:AL-prel} by the nonarchimedean triangle inequality. Also, we assume that $\varpi^p$ is a unit multiple of $p$ (see \S\ref{perfectoid-def}): this change of $\varpi$ does not increase $\abs{\varpi}^p$ and only enlarges the collection of valuations in question. In addition, we 
choose a generator $\xi$ of $\Ker(\theta \colon W(A^\flat) \surjects A)$, so that 
$\f{\xi - [\ov{\xi}]}{p} \in W(A^\flat)^\times$ where $\ov{\bullet}$ denotes the residue class modulo $p$ (see \S\ref{perfectoid-def}). Let $z_0 \in W(A^\flat)$ be a fixed lift of $a \in A$, 
recursively define further lifts
\[
\tst z_{n + 1} \ce z_n - \xi  \p{\f{\xi - [\ov{\xi}]}{p}}\i \!\! \p{\f{z_n - [\ov{z}_n]}{p}} = [\ov{z}_n] - [\ov{\xi}] \p{\f{\xi - [\ov{\xi}]}{p}}\i\!\! \p{\f{z_n - [\ov{z}_n]}{p}} \in W(A^\flat), 
\]
and set 
\[
a' \ce \ov{z}_m.
\]
To check that $a'$ satisfies \eqref{AL-eq}, we begin by noting that a continuous valuation $\abs{\cdot}$ on $A$ defines a $\varpi^\flat$-adically continuous valuation $\abs{\cdot}_{\flat}$ on $A^\flat$ by\footnote{The triangle inequality follows from the continuity of $\abs{\cdot}$ and the formula 
\[
(x + x')^\sharp = \lim_{n \ra \infty} ((x^\sharp)^{\frac{1}{p^n}} + (x'^\sharp)^{\frac{1}{p^n}})^{p^n}
\]
that one deduces from \eqref{monoid-id} and the fact that, by induction, $b^{p^{n - 1}} \bmod p^nA$ for $b \in A$ only depends on $b \bmod pA$.} $x \mapsto \abs{x^\sharp}$ (see \eqref{monoid-id}). For $z \in W(A^\flat)$, we set 
\[
\tst \abs{z}_{\sup} \ce \max_{j \ge 0} (|z_{(j)}|_\flat) 
\]
in terms of the unique expansion
\[
\tst  z = \sum_{j \ge 0}\, [z_{(j)}]\cdot p^j \in W(A^\flat);
\]
we will only use this in inequalities ``$\le$'' to abbreviate ``every $|z_{(j)}|_\flat$ is $\le$'' (so the attainment of the $\max$ need not concern us). Since $z = (z_{(j)}^{p^j})_{j \ge 0}$ is the Witt vector expansion, the nature of Witt vector addition and multiplication  \cite{BouAC}*{chapitre IX, section 1, num\'{e}ro 3, a) et b)} ensures that the map given by $z \mapsto \abs{z}_{\sup}$ satisfies the nonarchimedean triangle inequality and  is submultiplicative.  Consequently, since
 \[
\tst a - a'^\sharp =  \theta(z_m - [\ov{z}_m]) = \sum_{j \ge 1} ((z_m - [\ov{z}_m])_{(j)})^\sharp\cdot p^j, 
\]
it suffices to show that
\be \label{AL-eq-re}
\abs{z_m - [\ov{z}_m]}_{\sup} \le \max(|\ov{z}_m|_\flat, |\varpi^\flat|_\flat^{pm}).
\ee
By the definition of $z_{n + 1}$ and the fact that $\ov{\xi}$ is a unit multiple of $(\varpi^\flat)^p$ (see \S\ref{perfectoid-def}), we have
\[
\tst 
\abs{ z_{n + 1} - [\ov{z}_n]}_{\sup} \le |\varpi^\flat|_\flat^p \cdot |z_n|_{\sup}.
\]
Thus, for the least $0 \le N \le \infty$ with $\abs{\ov{z}_N}_\flat > |\varpi^\flat|_\flat^{p(N + 1)}$ (so $N$ depends on $\abs{\cdot}$), induction on $n$ gives
\[
\abs{z_n}_{\sup} \le |\varpi^\flat|_\flat^{pn} \q \x{for} \q n \le N,
\]
which settles \eqref{AL-eq-re} when $m \le N$. 
In the remaining case $m > N$, the preceding displays still give $\abs{ z_{N + 1} - [\ov{z}_N]}_{{\sup}} \le |\varpi^\flat|_\flat^{p(N + 1)}$, so the choice of $N$ 
and the triangle inequality give $\abs{z_{N + 1}}_{\sup} = \abs{\ov{z}_N}_\flat$ and $\abs{\ov{z}_{N + 1}}_\flat = \abs{\ov{z}_N}_\flat$. Thus, by repeating with $N + 1$ in place of $N$ we get $\abs{z_{N + 2} - [\ov{z}_{N + 1}]}_{\sup} < \abs{\ov{z}_{N + 1}}_\flat$, so also $\abs{z_{N + 2}}_{\sup} = \abs{\ov{z}_{N + 1}}_\flat$ and $\abs{\ov{z}_{N + 2}}_\flat = \abs{\ov{z}_{N + 1}}_\flat$. Iteration 
gives the sufficient 
\[
\abs{z_m}_{\sup} = \abs{\ov{z}_m}_\flat.\qedhere
\]
\epf


We will use the approximation lemma in conjunction with the following standard fact.

\blem[Special case of \cite{GR18}*{Corollary~15.4.27~(ii)}] \label{nilp-lem}
Let $A$ be a ring equipped with the $\varpi$-adic topology for a nonzerodivisor $\varpi \in A$. An element $a \in A[\f{1}{\varpi}]$ is topologically nilpotent \up{that is, $a^n \in \varpi A$ for large $n$} if and only if $\abs{a} < 1$ for any continuous valuation $\abs{\cdot}$ on $A[\f{1}{\varpi}]$ with $\abs{A} \le 1$.
\elem

\bpf
The `only if' is clear: if $a \in A[\f{1}{\varpi}]$ is topologically nilpotent and $\abs{\cdot}$ is continuous, then any $n > 0$ with $a^n \in \varpi A$ satisfies $|a|^n = |a^n| \le |\varpi| < 1$, so $|a| < 1$. For the `if,' we first use \cite{Hub93a}*{Lemma~3.3~(i)} to see that $a$ lies in the integral closure of $A$ in $A[\f{1}{\varpi}]$, so its powers are bounded in $A[\f{1}{\varpi}]$. We let $A^{\circ\circ} \subset A$ be the ideal of topologically nilpotent elements 
and consider the $A$-subalgebra 
$A[\f{1}{a}] \subset (A[\f{1}{\varpi}])[\f{1}{a}]$ generated by $\f{1}{a}$. 
If $A^{\circ\circ}\cdot A[\f{1}{a}]$ is the unit ideal of $A[\f{1}{a}]$, then $a$ satisfies an equation 
\[
\tst a^N + \sum_{i = 0}^{N - 1} a_i \cdot a^i = 0 \q \x{in} \q A[\f{1}{\varpi}] \q \x{with} \q a_i \in A^{\circ\circ}.
\]
In this case, $a$ is topologically nilpotent because so are the $a_i \cdot a^i$ by the boundedness of $\{a^i\}_{i \ge 0}$. Thus, we are left with the case when $A^{\circ\circ} \cdot A[\f{1}{a}]$ lies in a maximal ideal $\fm \subset A[\f{1}{a}]$. In turn, $\fm$ contains a minimal prime $\fp$ of $A[\f{1}{a}]$, which extends to a minimal prime $\wt{\fp}$ of $(A[\f{1}{\varpi}])[\f{1}{a}]$ (see \cite{SP}*{Lemmas~\href{https://stacks.math.columbia.edu/tag/00E0}{00E0} and~\href{https://stacks.math.columbia.edu/tag/00FK}{00FK}}). The target of the injection 
\[
\tst A[\f{1}{a}]/\fp \hra (A[\f{1}{\varpi}])[\f{1}{a}]/\wt{\fp}
\] 
is a domain, so it has a valuation subring $V$ that dominates $(A[\f{1}{a}]/\fp)_\fm$ (see \cite{SP}*{Lemma~\href{https://stacks.math.columbia.edu/tag/00IA}{00IA}}). The ideal $\bigcap_{n \ge 0} \varpi^n V \subset V$ is prime, so 
\[
\tst\overline{V} \ce V/(\bigcap_{n \ge 0} \varpi^n V)
\]
is a valuation ring in which the powers of $\varpi$ get arbitrarily close to $0$. Thus, the map $A[\f{1}{\varpi}] \ra \ov{V}[\f{1}{\varpi}]$ gives rise to a continuous valuation $\abs{\cdot}$ with $\abs{A} \le 1$ and $|\f{1}{a} | \le 1$. The latter contradicts $\abs{a} < 1$.
\epf

As a final preparation for the promised variant of Andr\'{e}'s lemma, we review ind-syntomic ring maps.

\bpp[Ind-fppf and ind-syntomic ring maps] \label{ind-fppf-maps}
A ring map $A \ra A'$ is \emph{ind-fppf} (resp.,~\emph{ind-syntomic}) if $A'$ is a filtered direct limit of faithfully flat, finitely presented (resp.,~syntomic\footnote{
A ring map $A \ra A'$ \emph{syntomic} if $\Spec(A')$ is covered by spectra of $A$-algebras of the form $A[x_1, \dotsc, x_n]/(f_1, \dotsc, f_c)$ with each $A[x_1, \dotsc, x_n]/(f_1, \dotsc, f_i)$ flat over $A$ and $f_1, \dotsc, f_c$ a regular sequence in $A[x_1, \dotsc, x_n]_\fp$ for every prime $\fp \supset (f_1, \dotsc, f_c)$ (by \cite{SP}*{Lemmas~\href{https://stacks.math.columbia.edu/tag/00SY}{00SY} and~\href{https://stacks.math.columbia.edu/tag/00SV}{00SV}}, this definition agrees with its counterpart \cite{SP}*{Definition~\href{https://stacks.math.columbia.edu/tag/00SL}{00SL}} used in \emph{op.~cit.}). 
}) $A$-algebras.\footnote{The distinction between ind-fppf and merely faithfully flat maps is subtle: for instance, if $R$ is a Noetherian local ring and $\wh{R}$ is its completion, then $R \ra \wh{R}$ is flat but, by \cite{Gab96}*{Proposition~1}, not ind-fppf when it has a nonreduced fiber and $R$ is a $\bQ$-algebra (as happens 
in \cite{FR70}*{proposition~3.1}). For further examples of maps that are faithfully flat but not ind-fppf, see \cite{SP}*{Section~\href{https://stacks.math.columbia.edu/tag/0ATE}{0ATE}}.}
It is useful to note that $A \ra A'$ is ind-fppf if and only if it is faithfully flat and $A'$ is a filtered direct limit of flat, finitely presented $A$-algebras. Concretely, $A \ra A'$ is ind-fppf (resp.,~ind-syntomic) if and only if every $A$-algebra map $B \ra A'$ with $B$ finitely presented over $A$ factors as $B \ra S \ra A'$ with $S$ faithfully flat, finitely presented (resp.,~syntomic) over $A$ (see \cite{SP}*{Lemma~\href{https://stacks.math.columbia.edu/tag/07C3}{07C3}}). In particular, ind-fppf and ind-syntomic maps are stable under composition\footnote{For instance, to show that the composition of ind-fppf maps $A \ra A'$ and $A' \ra A''$ is ind-fppf, 
for a test $B \ra A''$ over $A$ we factorize $A' \ra A' \tensor_A B \ra S \ra A''$ with $A' \ra S$ faithfully flat and finitely presented, express $A'$ as a filtered direct limit of faithfully flat, finitely presented $A$-algebras, and then descend the factorization using limit formalism, faithful flatness of $A \ra A''$, and \cite{EGAIV3}*{corollaire 11.2.6.1} (in the syntomic version, we use \cite{SP}*{Lemma~\href{https://stacks.math.columbia.edu/tag/0C33}{0C33}} instead).} 
and base change. A finite product or a  filtered direct limit of ind-fppf (resp.,~ind-syntomic) $A$-algebras is ind-fppf (resp.,~ind-syntomic). Certainly, faithfully flat ind-syntomic maps are ind-fppf.
\epp

\bthm \label{Andre-lem} 
Let $A$ be a ring, let $\varpi \in A$ with $\varpi^p \mid p$ be such that it has compatible $p$-power roots $\varpi^{1/p^n} \in A$, and suppose that either 
\benumr
\m \lab{AL-i}
$A$ is a $\varpi$-adically complete perfectoid\uscolon or

\m \lab{AL-ii}
$A$ is $\varpi$-Henselian, its $\varpi$-adic completion is perfectoid, and $\varpi$ is a nonzerodivisor in $A$.
\eenum
There are a faithfully flat, ind-syntomic, 
$\varpi$-Henselian $A$-algebra $A'$ whose $\varpi$-adic completion $\wh{A'}$ is perfectoid and a $\varpi$-divisible ideal $I' \subset A'$ with $A'/I'$ faithfully flat over $A$ such that every monic $P \in A'[T]$ has a root $\gA_P \in A'/I'$ with compatible $\gA_P^{1/p^n} \in A'/I'$, 
in particular, the $\gA_P^{1/p^n}$ exist in 
$\wh{A'}$. 
\ethm

\brem
The perfectoid $\wh{A'}$ contains compatible $p$-power roots of every $a \in A$, and $\wh{A'}/(\varpi)$ is faithfully flat over $A/(\varpi)$. Thus, the preceding theorem  recovers the original lemma of Andr\'{e} \cite{Bha18a}*{Theorem~1.5}, in which one only required $\wh{A'}/(\varpi)$ to be \emph{almost} faithfully flat over $A/(\varpi)$.
\erem

\bpp[]\bpf[Proof of Theorem \uref{Andre-lem}]
The final assertion follows from the rest because $\wh{A'} \cong \wh{A'/I'}$ by the $\varpi$-divisibility of $I'$.
 By \Cref{lem:perfectoid-quotient-of-torsion-free}, the perfectoid $A$ in \ref{AL-i} is a quotient of a perfectoid $\wt{A}$ that is $\wt{\varpi}$-torsion free and $\wt{\varpi}$-adically complete  for a lift $\wt{\varpi} \in \wt{A}$ of $\varpi$ with $\wt{\varpi}^p \mid p$. Once some $\wt{A}'$ with a $\wt{\varpi}$-divisible ideal $\wt{I}' \subset \wt{A}'$ works for $\wt{A}$ with respect to $\wt{\varpi}$, its quotient $A' \ce \wt{A}' \tensor_{\wt{A}} A$ with the image $I' \subset A'$ of $\wt{I}'$ works for $A$ (
see \Cref{recognize}~\ref{R-a}). 
This reduces \ref{AL-i} to \ref{AL-ii}.

For the rest of the proof, we assume \ref{AL-ii} and build on the argument of \cite{GR18}*{Theorem~16.9.17}, which established a similar result without the ind-syntomic aspect. We may then restrict to those $P$ that belong to the set $\cP$ of all the monic polynomials in $A[T]$: indeed, since $A'$ inherits the assumption \ref{AL-ii}, we may \emph{a posteriori} iterate the construction countably many times to build a tower
\[
A \equalscolon A_0' \ra A_1' \ra A_2' \ra \dotsc \qxq{and $\varpi$-divisible ideals} I_n' \subset A_n' \qxq{for} n > 0
\]
such that $I_n' \subset A_n'$ satisfy the requirements with respect to the monic polynomials in $A_{n - 1}'[T]$; since 
\[
\tst A_n'/(\sum_{1 \le i \le n} I_i'A_n') \qxq{is faithfully flat over} A_{n - 1}'/(\sum_{1\le i \le n - 1} I_i'A_{n -1}')
\]
and hence, by induction, also over $A$, the ring $A_\infty' \ce \varinjlim_{n \ge 0} A_n'$ with its $\varpi$-divisible ideal $I_\infty' \ce \sum_{i \ge 1} I_i' A_\infty'$ then satisfies the requirements with respect to the monic polynomials in $A_\infty'[T]$ (see \eqref{p-oid-mod-pi} and the end of \S\ref{perfectoid-def}).  With $\cP$ fixed, we may drop the  requirement that $A'$ be $\varpi$-Henselian---indeed, we may acquire this \emph{a posteriori} by replacing $A'$ by its $\varpi$-Henselization: since $\varpi$ lies in the maximal ideals of $A$, this does not lose faithful flatness (see \cite{SP}*{Lemma~\href{https://stacks.math.columbia.edu/tag/00HP}{00HP}}). Thus, dropping $\varpi$-Henselianity and restricting to $P \in \cP$, we first define $A'$ and then, in the rest of the proof, check that it meets the requirements. We let $A_\infty$ be the subring
\[
\tst  \p{A\left[T_P^{1/p^n} \, \big\vert\, P \in \cP,\, n \ge 0\right]}\tst [\f{P(T_P)}{\varpi^m} ]_{ P \in \cP,\, m \ge 0} \subset \p{A\left[T_P^{1/p^n} \, \big\vert\, P \in \cP ,\, n \ge 0  \right]}[\tst\f{1}{\varpi}],
\]
so that
\[
A_\infty[\tst \f{1}{\varpi}] \cong (A[T_P^{1/p^n} \, \big\vert\, P \in \cP ,\, n \ge 0])[\tst\f{1}{\varpi}]
\]
and define a $\varpi$-divisible ideal $I_\infty \subset A_\infty$ by 
\[
\tst I_\infty \ce (\f{P(T_P)}{\varpi^m} \,  \big\vert\, P \in \cP,\, m \ge 0) \subset A_\infty.
\]
Our candidate $A'$ and a $\varpi$-divisible ideal $I' \subset A'$ are (see \S\ref{p-int-cl})
\[
\tst A' \ce (\x{$p$-integral closure of $A_\infty$ in $A_\infty[\tst \f{1}{\varpi}]$})
\]
and
\[
\tst I' \ce (\f{P(T_P)}{\varpi^m} \,  \big\vert\, P \in \cP,\, m \ge 0) \subset A'.
\]
Since each $P(T_P)$ vanishes in $A'/I'$, the class of $T_P$ is a desired root $\gA_P$. Moreover, $I_\infty$ is $\varpi$-divisible,~so
\[
A_\infty/(\varpi^p) \q \x{is a quotient of} \q (A[T_P^{1/p^n} \, \vert\, P \in \cP, n \ge 0])/(\varpi^p),
\]
and hence every element of $A_\infty/(\varpi^p)$ is a $p$-th power (the same holds for $A$ in place of $A_\infty$, see \S\ref{perfectoid-def}~\ref{PD-ii}), to the effect that $\wh{A'}$ is perfectoid by \Cref{p-oid-rec}. Due to the $\varpi$-divisibility of $I_\infty$ and $I'$, the quotients $A_\infty/I_\infty$ and $A'/I'$ are $\varpi$-torsion free, so we have
\[
\tst A_\infty/I_\infty \subset A'/I' \subset (A_\infty/I_\infty)[\f{1}{\varpi}] 
\]
and
\[
\tst (A_\infty/I_\infty)[\f{1}{\varpi}] \cong \p{A\left[T_P^{1/p^n} \, \big\vert\, P \in \cP,\, n \ge 0\right]/(P(T_P)\,  |\,  P \in \cP)}[\f{1}{\varpi}].
\]
The $\varpi$-divisibility of $I'$ and the observation about \eqref{p-int-cl-lem} with \eqref{p-oid-mod-pi} then imply that $A'/I'$ is the $p$-integral closure of $A_\infty/I_\infty$ in $(A_\infty/I_\infty)[\f{1}{\varpi}]$. We may describe $A_\infty/I_\infty$ explicitly as follows: each $P$ is monic, so 
\[
\tst A[T_P^{1/p^n}]_{P,\, n}/(P(T_P))_P \subset \p{A[T_P^{1/p^n} ]_{P,\, n}/(P(T_P))_P}[\f{1}{\varpi}]
\]
and, since  elements of this subring lift to $A_\infty$ (even to $A[T_P^{1/p^n}]_{P,\, n}$), we have 
\[
\tst A_\infty/I_\infty = A[T_P^{1/p^n}]_{P,\, n}/(P(T_P))_P \qxq{inside} \p{A[T_P^{1/p^n} ]_{P,\, n}/(P(T_P))_P}[\f{1}{\varpi}].
\]
In particular, $A'/I'$ is integral over $A$ 
and $(A'/I')[\f{1}{\varpi}]$ is even ind-(finite, module-free) over $A[\f{1}{\varpi}]$. Thus, since  $\varpi \in A$ is a nonzerodivisor, the closed morphism $\Spec(A'/I') \ra \Spec(A)$ is surjective. Moreover, by glueing of flatness \cite{RG71}*{seconde partie, lemme 1.4.2.1}, the desired $A$-flatness of $A'/I'$ will follow from the $A/(\varpi)$-flatness of $(A'/I')/(\varpi) \cong A'/(\varpi)$. In conclusion, it remains to argue that $A'$ is $A$-ind-syntomic.

For the remaining ind-syntomicity of $A'$ over $A$, due the closedness of ind-syntomic maps under filtered direct limits (see \S\ref{ind-fppf-maps}), we may replace $\cP$ by its variable finite subset. Then, since for finite $\cP$ the $A$-algebra $A'$ can equivalently be built iteratively, we may replace $\cP$ by a singleton $ \{ P \}$. To reduce further, we simplify the notation by setting $T \ce T_P$ and for $m \ge 0$ set
\[
A_{m} \ce \p{A[T^{1/p^n} \, \big\vert\, n \ge 0]}\tst [\f{P(T)}{\varpi^m}] \subset \p{A\left[T^{1/p^n} \, \big\vert\, n \ge 0 \right]}[\tst\f{1}{\varpi}] \cong A_{m}[\f{1}{\varpi}] \cong A_{\infty}[\f{1}{\varpi}]
\]
and
\[
\tst A'_m \ce (\x{$p$-integral closure of $A_{m}$ in $A_{m}[\tst \f{1}{\varpi}]$}).
\]
By another passage to a limit, it suffices to show that each $A_m'$ with $m > 0$ is ind-syntomic over $A$. To argue this, we will use the perfectoid nature of $A_0 \cong A\left[T^{1/p^n} \, \big\vert\, n \ge 0 \right]$ and the fact that $\varpi^t, P(T)$ is an $A_0$-regular sequence for any $t \in \bZ[\f{1}{p}]_{\ge 0}$ (since $P$ is monic),  to describe $A_m'$ explicitly. The $A_0$-regularity of $\varpi^m, P$ already implies an explicit description of $A_m$ (see \cite{SP}*{Lemma~\href{https://stacks.math.columbia.edu/tag/0BIQ}{0BIQ}}
):
\be \label{exp-blowup}
\tst A_m \cong A_0[\f{P}{\varpi^m}] \cong A_0[X]/(\varpi^mX - P)
\ee 
and, likewise,
\[
\tst\wh{A_0}[\f{P}{\varpi^m}] \cong \wh{A_0}[X]/(\varpi^m X - P),
\]
where $\wh{A_0}$ is the $\varpi$-adic completion. 

Since $\wh{A_0}$ is perfectoid (see \S\ref{perfectoid-def}), by \Cref{approx-lem}, there is a $Q \in \wh{A_0}$ that admits compatible $p$-power roots $Q^{1/p^j} \in \wh{A_0}$ such that
\be \label{PQ-ineq}
\tst \abs{P - Q} < \max(|Q|, | \varpi^m |)
\ee 
for every continuous valuation $\abs{\cdot}$ on $\wh{A_0}[\f{1}{\varpi}]$ with $| \wh{A_0} | \le 1$. 
Letting 
$\wh{A_0}{}^+$ be the integral closure of 
$\wh{A_0}$ in 
$\wh{A_0}[\f{1}{\varpi}]$, we then have
\[
\tst \{ \abs{\f{P}{\varpi^m}} \le 1 \} = \{ |\f{Q}{\varpi^m}| \le 1 \} \q \x{in} \q \Spa(\wh{A_0}[\f{1}{\varpi}], \wh{A_0}{}^+). 
\]
This agreement implies that if we endow $\wh{A_0}[\f{P}{\varpi^m}]$ and $\wh{A_0}[\f{Q}{\varpi^m}]$ with their $\varpi$-adic topologies, then the continuous valuations $\abs{\cdot}$ on $(\wh{A_0}[\f{P}{\varpi^m}])[\f{1}{\varpi}]$ with $ | \wh{A_0}[\f{P}{\varpi^m}] | \le 1$ are identified with the continuous valuations on $(\wh{A_0}[\f{Q}{\varpi^m}])[\f{1}{\varpi}]$ with $ | \wh{A_0}[\f{Q}{\varpi^m}] | \le 1$. Moreover, \eqref{PQ-ineq} implies that every such valuation satisfies $| \f{P}{\varpi^m} - \f{Q}{\varpi^m} | < 1$. Consequently, by \Cref{nilp-lem}, every
\be \label{PQ-power}
\tst \x{large power of} \ \ \f{P}{\varpi^m} - \f{Q}{\varpi^m} \ \ \x{lies in} \ \ \varpi(\wh{A_0}[\f{P}{\varpi^m}]) \ \  \x{and} \ \ \varpi(\wh{A_0}[\f{Q}{\varpi^m}]).
\ee
\Cref{nilp-lem} and \eqref{PQ-ineq} also imply that $(P - Q)^{p^t} \in \varpi\wh{A_0}$ for some $t \in \bZ_{\ge 0}$, so that, by \eqref{p-oid-mod-pin}, we have $P - Q \in \varpi^{1/p^t}\wh{A_0}$. In particular, 
\be \lab{Qpj-monic}
Q^{1/p^j} \qxq{is monic in} \wh{A_0}/\varpi^{1/p^{j + t}}
\ee
(see \eqref{p-oid-mod-pin}), so the sequence 
\be \label{piQ-regseq}
\varpi^{m/p^{j}}, Q^{1/p^j} \q \x{is $\wh{A_0}$-regular for every $j \ge 0$}
\ee
(see \cite{SP}*{Lemma~\href{https://stacks.math.columbia.edu/tag/07DV}{07DV}}). Consequently, analogously to \eqref{exp-blowup}, we have 
\be \label{exp-blowup-2}
\tst \wh{A_0}[\f{Q^{1/p^j}}{\varpi^{m/p^{j}}}] \cong \wh{A_0}[X^{1/p^j}]/(\varpi^{m/p^j} X^{1/p^j} - Q^{1/p^j}),
\ee
where we chose the label `$X^{1/p^j}$' for the polynomial variable to make evident the resulting identification
\be \label{hat-id}
\tst \wh{A_0}[\f{Q^{1/p^j}}{\varpi^{m/p^{j}}}\, |\, j \ge 0] \cong \wh{A_0}[X^{1/p^j}\, | \, j \ge 0]/(\varpi^{m/p^j} X^{1/p^j} - Q^{1/p^j})_{j \ge 0}.
\ee
It then follows from \eqref{p-oid-mod-pi} that the $\varpi$-adic completion of the subalgebra 
\[
\tst \wh{A_0}[\f{Q^{1/p^j}}{\varpi^{m/p^{j}}}\, |\, j \ge 0] \subset \wh{A_0}[\f{1}{\varpi}]
\]
is perfectoid, and hence, from the observation about \eqref{p-int-cl-lem}, that this subalgebra is $p$-integrally closed. Due to \eqref{PQ-power}, the $p$-integral closure of $\wh{A_0}[\f{P}{\varpi^m}]$ in $\wh{A_0}[\f{1}{\varpi}]$ contains $\f{Q}{\varpi^m}$ and the $p$-integral closure of $\wh{A_0}[\f{Q}{\varpi^m}]$ in $\wh{A_0}[\f{1}{\varpi}]$ contains $\f{P}{\varpi^m}$, so it follows that these two closures agree and both are equal to 
\[
\tst \wh{A_0}[\f{Q^{1/p^j}}{\varpi^{m/p^{j}}}\, |\, j \ge 0].
\]
To describe the sought $p$-integral closure $A_m'$ of $A_0[\f{P}{\varpi^m}]$ in $A_0[\f{1}{\varpi}]$ for $m > 0$, we now take advantage of the preceding analysis over $\wh{A_0}$. We use \eqref{PQ-power} to fix a $d > 0$ such that 
\be \label{PQ-power-2}
\tst (P - Q)^{p^d} \in \varpi^{mp^d}(\wh{A_0}[\f{P}{\varpi^m}]) \q \x{and} \q (P - Q)^{p^d} \in \varpi^{mp^d}(\wh{A_0}[\f{Q}{\varpi^m}]).
\ee
We then fix a 
\[
q \in A_0 \qxq{congruent to}  Q \in \wh{A_0} \q \x{modulo $\varpi^{mp^d}$,}
\]
so that the image of $q$ in $A_0/\varpi^{1/p^t}$ is monic and $\varpi^m, q$ is an $A_0$-regular sequence (compare with \eqref{piQ-regseq}). Consequently, as in \eqref{exp-blowup}, we have 
\[
\tst A_0[\f{q}{\varpi^m}] \cong A_0[X]/(\varpi^mX - q).
\]
By combining this with \eqref{exp-blowup} and \eqref{exp-blowup-2}, we see that both maps
\[
\tst A_0[\f{P}{\varpi^m}] \ra \wh{A_0}[\f{P}{\varpi^m}]
\]
and
\[
\tst A_0[\f{q}{\varpi^m}] \cong A_0[X]/(\varpi^mX - q) \xra{X \mapsto X + \f{q - Q}{\varpi^m}} \wh{A_0}[X]/(\varpi^mX - Q) \cong \wh{A_0}[\f{Q}{\varpi^m}]
\]
induce isomorphisms on $\varpi$-adic completions. Thus, since these maps are compatible with the $\varpi$-adic completion map $A_0 \ra \wh{A_0}$, we get from \eqref{PQ-power-2} that 
\[
\tst (P - q)^{p^d} \in \varpi^{mp^d}(A_0[\f{P}{\varpi^m}]) \q \x{and} \q (P - q)^{p^d} \in \varpi^{mp^d}(A_0[\f{q}{\varpi^m}]).
\]
Consequently, the $p$-integral closures of $A_0[\f{P}{\varpi^m}]$ and $A_0[\f{q}{\varpi^m}]$ in $A_0[\f{1}{\varpi}]$ agree, and hence equal $A_m'$. 

To proceed, we fix $q_j \in A_0$ for $j \ge 0$ such that $q_0 \ce q$ and
\be \label{qj-def}
q_{j} \equiv Q^{1/p^j} \bmod \varpi^{mp} A_0 \q \x{for} \q j > 0.
\ee
Since $q_{j + 1}^p \equiv q_j \bmod \varpi^{mp} A_0$, we have $(\f{q_{j + 1}}{\varpi^{m/p^{j + 1}}})^p - \f{q_j}{\varpi^{m/p^j}} \in A_0$ for every $j \ge 0$, so the subalgebras
\be \label{include-a-lot}
\tst A_0[\f{q}{\varpi^{m}}] \subset \dotsc \subset A_0[\f{q_j}{\varpi^{m/p^j}}] \subset A_0[\f{q_{j + 1}}{\varpi^{m/p^{j + 1}}}] \subset \dotsc \q \x{in} \q A_0[\f{1}{\varpi}]
\ee
are contained in the $p$-integral closure $A_m'$ of $A_0[\f{q}{\varpi^{m}}]$ in $A_0[\f{1}{\varpi}]$. In fact, their union is $p$-this integral closure: to show this, we first note that, due to \eqref{piQ-regseq} and \eqref{qj-def}, the sequence $\varpi^{m/p^j}, q_j$ is $A_0$-regular, and hence, analogously to \eqref{exp-blowup-2}, that
\[
\tst A_0[\f{q_j}{\varpi^{m/p^j}}] \cong A_0[X_j]/(\varpi^{m/p^j} X_j - q_j).
\]
In terms of these identifications, the inclusions \eqref{include-a-lot} become the maps
\[
A_0[X_j]/(\varpi^{m/p^j} X_j - q_j) \xra{X_j\, \mapsto\, X_{j + 1}^p  + \f{q_j - q_{j + 1}^p}{\varpi^{m/p^j}}} A_0[X_{j + 1}]/(\varpi^{m/p^{j + 1}} X_{j + 1} - q_{j + 1})
\]
Since $\varpi^{pm} \mid q_j - q_{j + 1}^p$ and $pm - \f{m}{p^j} \ge 1$, we see from \eqref{hat-id} that the direct limit of these maps modulo $\varpi$ is identified with $(\wh{A_0}[\f{Q^{1/p^j}}{\varpi^{m/p^{j}}}\, |\, j \ge 0])/\varpi$. Since $\wh{A_0}[\f{Q^{1/p^j}}{\varpi^{m/p^{j}}}\, |\, j \ge 0]$ is $p$-integrally closed in $\wh{A_0}[\f{1}{\varpi}]$, it follows from the observation about \eqref{p-int-cl-lem} (applied with $\varpi$ there replaced by $\varpi^{1/p}$) that $A_0[\f{q_j}{\varpi^{m/p^{j}}}\, |\, j \ge 0]$ is $p$-integrally closed in $A_0[\f{1}{\varpi}]$, and hence that it equals $A_m'$.

Thanks to this explicit description of $A_m'$ and the stability of ind-syntomic algebras under filtered direct limits, all that  remains is to show that each $A_0[X_j]/(\varpi^{m/p^j} X_j - q_j)$ is ind-syntomic over $A$. However, $q_j$ comes from  $A[T^{1/p^n}]$ for every large enough $n$ and its image in $(A/\varpi^{1/p^{j + t}})[T^{1/p^n}]$ is monic (see \eqref{qj-def} and \eqref{Qpj-monic}). Thus, the ($(A[T^{1/p^n}])[X_j]$)-regular element $\varpi^{m/p^j} X_j - q_j$ stays regular on every $A$-fiber of $(A[T^{1/p^n}])[X_j]$. Consequently, each $(A[T^{1/p^n}])[X_j]/(\varpi^{m/p^j} X_j - q_j)$ is a syntomic $A$-algebra (see \cite{SP}*{Lemma~\href{https://stacks.math.columbia.edu/tag/00SW}{00SW}}), and it remains to form the direct limit in $n$.
\epf
\epp


The following consequence of Andr\'{e}'s lemma gives convenient ``semiperfectoid'' covers of $\bZ_{(p)}$-algebras. 

\bcor \label{cor:semiperfectoid-cover}
Every ring $A$ that is $p$-Zariski in the sense that $1 + pA \subset A^\times$ admits a faithfully flat map $A \ra A_\infty$ such that the $p$-adic completion of $A_\infty^\red$ is perfectoid, $A_\infty$ is a quotient of a $p$-torsion free, $p$-Henselian ring $\wt{A}_\infty$ whose $p$-adic completion is perfectoid, and every monic polynomial in $\wt{A}_\infty[T]$ has a root in $\wt{A}_\infty$ \up{so the same also holds with $A_\infty$ or $A_\infty^\red$ in place of $\wt{A}_\infty$}.
\ecor

\bpf
The $p$-Zariski condition amounts to $p$ lying in the Jacobson radical, equivalently, in every maximal ideal, of $A$. We recall that the \emph{$p$-Zariskization} of a ring $B$ is the localization $B_{1 + pB}$. By replacing $A$ by the $p$-Zariskization of the countable iteration of the construction 
\[
A \mapsto A[X_a^{1/p^\infty}\, \vert\, a \in A]/(X_a - a \, \vert \, a \in A),
\]
we lose no generality by assuming that every $a\in A$ admits compatible $p$-power roots $a^{1/p^n}$ in $A$. In turn, such an $A$ is then a quotient of the $p$-torsion free $\bZ_{(p)}$-algebra $\wt{A}$ that is the $p$-Zariskization of 
\[
\bZ[p^{1/p^\infty}][X_a^{1/p^{\infty}}\, \vert\, a \in A].
\]
By \Cref{recognize}~\ref{R-0}, the $p$-adic completion of $\wt{A}$ is perfectoid, so we apply \Cref{Andre-lem} to the $p$-Henselization of $\wt{A}$ to build a faithfully flat, $p$-Henselian $\wt{A}$-algebra $\wt{A}_\infty$ whose $p$-adic completion is perfectoid such that every monic polynomial in $\wt{A}_\infty[T]$ has a root in $\wt{A}_\infty$. The quotient $A_\infty \ce \wt{A}_\infty \tensor_{\wt{A}} A$ of $\wt{A}_\infty$ is faithfully flat over $A$ and every element of its nilradical admits a $p$-th root. \Cref{recognize}~\ref{R-c} then ensures that the $p$-adic completion of $A_\infty^\red$ is perfectoid.
\epf


\section{The prime to the characteristic aspects of the main result}

For arguing our purity results, the first task is to dispose of the cases when the order of the coefficients is invertible. For this, we first give a new, perfectoid-based proof of the Gabber--Thomason purity for \'{e}tale cohomology of regular rings in \S\ref{Gabber-Thomason}. We then use it in \S\ref{et-depth-geom-depth} to deduce purity for \'{e}tale cohomology in the general singular case via a local Lefschetz style theorem. In \S\ref{nonabelian}, we present a nonabelian analogue of the results of \S\ref{et-depth-geom-depth}: a generalization of the Zariski--Nagata purity theorem.



\csub[The absolute cohomological purity of Gabber--Thomason] \lab{Gabber-Thomason}

Purity for \'{e}tale cohomology of regular rings, stated precisely in \Cref{abs-coh-pur} (see also footnote~\ref{foot:semipurity}), was conjectured by Grothendieck and settled by Gabber in \cite{Fuj02}, who built on the strategy initiated by Thomason in \cite{Tho84}. Gabber's alternative later proof given in \cite{ILO14}*{expos\'{e} XVI} eliminated the use of algebraic $K$-theory. We present a proof that uses perfectoids, specifically, \Cref{baby-tilting}, to reduce to the positive characteristic case that had been settled by M. Artin already in \cite{SGA4III}*{expos\'{e}~XVI}. The following standard lemmas facilitate the passage to perfectoids.

\blem \label{lem:towers}
Let $(R, \fm)$ be a complete, regular, local ring with residue field $k$.
\benum
\m \label{tower-a}
There is a filtered direct system $\{(R_i, \fm_i)\}_{i \in I}$ of regular, local, finite, flat $R$-algebras that are unramified if so is $R$ \up{see \uS\uref{conv}} such that $\fm_i = \fm R_i$ 
and $(\varinjlim_i R_i, \varinjlim_i \fm_i)$ is a regular local ring whose residue field is an algebraic closure $\ov{k}$ of $k$. 

\m \label{tower-b}
If $R$ is of mixed characteristic $(0, p)$ and $k$ is perfect, then there is a tower $\{R_n\}_{n \ge 0}$ of regular, local, finite, flat $R$-algebras of $p$-power rank over $R$ such that the $p$-adic completion of $R_\infty \ce \varinjlim_{n \ge 0} R_n$ is perfectoid\ucolon explicitly, by the Cohen structure theorem, we have
\[
\qq R \simeq W(k)\llb x_1, \dotsc, x_d\rrb/(p - f), 
\]
where either $f = x_1$ or $f \in (p, x_1, \dotsc, x_d)^2$ \up{the two cases correspond to whether or not $R$ is unramified}, 
and one may choose
\[
\qq \tst R_n \ce W(k)\llb x_1^{1/p^n}, \dotsc, x_d^{1/p^n} \rrb/(p - f)
\]
with
\[
\tst \qq R_\infty \simeq W(k)\llb x_1^{1/p^\infty}, \dotsc, x_d^{1/p^\infty} \rrb/(p - f).
\]
\eenum
\elem

For our somewhat nonstandard use of the notation $\llb \cdot \rrb$ in part \ref{tower-b} above, see \eqref{nst-terminology}.


\bpf
In essence, the claims are restatements of \cite{brauer-purity}*{Lemmas 5.1 and 5.2}: part \ref{tower-a} follows from \cite{brauer-purity}*{Lemma 5.1} and its proof, whereas part \ref{tower-b} follows from \cite{brauer-purity}*{Lemma 5.2} and its proof. For a prismatic point of view on the construction of $R_\infty$, see \cite{BS19}*{Remark 3.11}.
\epf

\blem \label{lem:et-to-hat}
Let $A$ be a ring, let $a \in A$ be such that $A$ is $a$-Henselian and has bounded $a^\infty$-torsion, 
let $\wh{A}$ be the $a$-adic completion of $A$, and let $\Spec(A[\f 1a]) \subset U \subset \Spec(A)$ be an open. We have
\[
R\Gamma_\et(U, \sF) \isomto R\Gamma_\et(U_{\wh{A}}, \sF)
\]
for every torsion abelian sheaf $\sF$ on $U_\et$. In particular, for every closed subset $Z \subset \Spec(A/aA)$ and every torsion abelian sheaf $\sF$ on $A_\et$, we have
\[
R\Gamma_Z(A, \sF) \isomto R\Gamma_Z(\wh{A}, \sF),
\]
and for a Noetherian ring $R$, an ideal $I \subset R$ such that $R$ is  $I$-Henselian,  and the $I$-adic completion~$\wh{R}$,
\[
R\Gamma_I(R, \sF) \isomto R\Gamma_I(\wh{R}, \sF) 
\]
for every torsion abelian sheaf $\sF$ on $R_\et$.
\elem

\bpf
The claims are special cases of \cite{BC19}*{Theorem 2.3.4, Corollary 2.3.5 (e)}, although we could also use earlier references \cite{Fuj95}*{Corollary~6.6.4} or \cite{ILO14}*{expos\'{e} XX, section 4.4}; 
see also \Cref{cor:excision} below.
\epf

\bthm \lab{abs-coh-pur}
For a regular local ring $(R, \fm)$ and a commutative, finite, \'{e}tale $R$-group $G$ whose order is invertible in $R$, 
\[
H^i_\fm(R, G) \cong 0 \q \text{for} \q i < 2 \dim(R).
\]
\ethm

\bpf
We use the local-to-global spectral sequence \cite{SGA4II}*{expos\'{e}~V, proposition~6.4} to assume that $R$ is strictly Henselian and then that $G \simeq \bZ/\ell\bZ$ for a prime $\ell$. We then use \Cref{lem:et-to-hat} to assume 
that $R$ is also complete. Thus, by the Cohen structure theorem \cite{Mat89}*{Theorem~29.7}, if $R$ is equicharacteristic, then $R \simeq k\llb x_1, \dotsc, x_d\rrb$ for a field $k$ and, by \Cref{lem:et-to-hat} again, we may assume that $R$ is the Henselization of $\bA^d_k$ at the origin. For this $R$ the claim was settled already in \cite{SGA4III}*{expos\'{e}~XVI, th\'{e}or\`{e}me~3.7}, so from now on we assume that our complete, strictly Henselian $R$ is of mixed characteristic $(0, p)$.

Since multiplication by $p$ is an automorphism of $\bZ/\ell\bZ$, the trace map \cite{SGA4III}*{expos\'{e}~XVII, sections~6.3.13 et 6.3.14, proposition 6.3.15~(iv)} allows us to replace $R$ by any module-finite, flat $R$-algebra $R'$ of $p$-power rank over $R$ such that $R'$ is a regular local ring. Thus, by \Cref{lem:towers}~\ref{tower-a} and a limit argument, we may pass to a tower to reduce to the case  when the residue field $k$ of $R$ is algebraically closed (we use \Cref{lem:et-to-hat} to complete $R$ again). We then likewise use \Cref{lem:towers}~\ref{tower-b} to  reduce to showing that
\[
H^{i}_{(p,\, x_1,\, \dotsc,\, x_d)}(R_\infty, \bZ/\ell\bZ) \cong 0 \qxq{for}  i < 2d 
\]
with
\[
R_\infty \cong W(k)\llb x_1^{1/p^\infty}, \dotsc, x_d^{1/p^\infty} \rrb/(p - f), \ \  f\in \fm^2 \cup \{ x_1 \},
\]
knowing that the $p$-adic completion $\wh{R}_\infty$ of $R_\infty$ is perfectoid. The tilt $\wh{R}_\infty^\flat$ of $\wh{R}_\infty$ reviewed in \eqref{monoid-id} is the $\ov{f}$-adic completion of $k\llb (x_1^{\flat})^{1/p^\infty}, \dotsc, (x_d^{\flat})^{1/p^\infty} \rrb$ for some $\ov{f} \in k\llb (x_1^{\flat})^{1/p^\infty}, \dotsc, (x_d^{\flat})^{1/p^\infty} \rrb$: explicitly, 
\[
\wh{R}_\infty^\flat \cong \varprojlim\limits_{z \mapsto z^p} \p{k\llb x_1^{1/p^\infty}, \dotsc, x_d^{1/p^\infty}\rrb/(f)},
\]
which via the map $(z_n)_{n \ge 0} \mapsto (z_n^{p^n})_{n \ge 0}$ is identified with
\[
\tst\invlim\limits_n \p{k\llb (x_1^{\flat})^{1/p^\infty}, \dotsc, (x_d^{\flat})^{1/p^\infty}\rrb /(\ov{f}^{p^n})},  
\]
where $x_i^\flat$ corresponds to the $p$-power compatible system $(x_i^{1/p^n})_{n \ge 0}$. Thus, by \eqref{eqn:baby-tilting-supports} with \Cref{lem:et-to-hat} (the latter removes the $\ov{f}$-adic completion), we are reduced to showing 
\[
H^{i}_{(x_1^\flat, \dotsc, x_d^\flat)}(k\llb (x_1^{\flat})^{1/p^\infty}, \dotsc, (x_d^{\flat})^{1/p^\infty}\rrb, \bZ/\ell\bZ) \cong 0 \qxq{for} i < 2d.
\]
By the perfection-invariance of \'{e}tale cohomology, we may replace 
\[
k\llb (x_1^{\flat})^{1/p^\infty}, \dotsc, (x_d^{\flat})^{1/p^\infty}\rrb \qxq{by} k\llb x_1^{\flat}, \dotsc, x_d^{\flat}\rrb,
\]
which brings us to the already discussed equal characteristic case.
\epf

\brems
\remi
As Teruhisa Koshikawa pointed out to us, the above argument also reduces the full absolute cohomological purity for \'{e}tale cohomology, namely, the statement that for all $n \in \bZ_{> 0}$ invertible in $R$ the \'{e}tale-sheafification $\cH^i_\fm$ of the cohomology with supports $H^i_\fm$ satisfies
\[
\qq \cH^i_\fm(-, \bZ/n\bZ) \cong \begin{cases} 0, & \x{for $i \neq 2\dim(R)$,} \\ \underline{\bZ/n\bZ}(-\dim(R)), & \x{for $i = 2\dim(R)$,} \end{cases}
\]
to positive characteristic. Indeed, the isomorphism in degree $i = 2\dim(R)$ is induced by the cycle class map 
\[
\qqq \underline{\bZ/n\bZ}(-\dim(R)) \ra \cH^{2\dim(R)}_\fm(-, \bZ/n\bZ),
\]
which one first argues to be injective as in \cite{Fuj02}*{Lemma 2.3.1}. The bijectivity then becomes the matter of bounding the nonzero stalk of the target, which may be done after passing to $\wh{R}^\flat_\infty$. The vanishing in degrees $i \neq 2\dim(R)$ reduces to positive characteristic as in the proof of \Cref{abs-coh-pur}.

\remi
Another way to pass to the tilt, without using \Cref{baby-tilting}, is to use diamonds developed in \cite{Sch17}. Namely, we consider the ``punctured adic spectrum'' of $R_\infty$ defined as
\[
\tst \qqq U^\ad_{R_\infty} \ce \Spa(R_\infty, R_\infty) \setminus \{ x_1 = \dotsc = x_d = 0 \}
\]
and isomorphic to
\[
\tst \qqq \bigcup_{i = 1}^{d} \Spa\p{R_\infty\langle\f{x_1,\, \dotsc,\, x_d}{x_i}\rangle, R_\infty\langle\f{x_1,\, \dotsc,\, x_d}{x_i}\rangle^+}\!,
\] 
where $R_\infty$ is endowed with its $(x_1, \dotsc, x_d)$-adic topology, so that $U^\ad_{R_\infty}$ is an analytic adic space over $\bZ_p$ (to simplify we ignore the issue of showing that the appearing Huber pairs are sheafy).\footnote{In fact, by \cite{Nak20}*{Example~2.1.4}, each $\Spa$ appearing in the union on the right side is even an affinoid perfectoid.} Likewise, we endow $\wh{R}_\infty^\flat$ with the $(x_1^\flat, \dotsc, x_d^\flat)$-adic topology and consider 
\[
\tst\qqq U^\ad_{\wh{R}_\infty^\flat} \ce \Spa(\wh{R}_\infty^\flat, \wh{R}_\infty^\flat) \setminus \{ x_1^\flat = \dotsc = x_d^\flat = 0 \},
\] 
which is isomorphic to
\[
\tst \qqq \bigcup_{i = 1}^{d} \Spa\p{\wh{R}_\infty^\flat\langle\f{x_1^\flat,\, \dotsc,\, x_d^\flat}{x_i^\flat}\rangle, \wh{R}_\infty^\flat\langle\f{x_1^\flat,\, \dotsc,\, x_d^\flat}{x_i^\flat}\rangle^+}\!
\]
and is a perfectoid space because the coordinate rings of the appearing affinoids inherit perfectness from $\wh{R}_\infty^\flat$ (see \cite{Sch17}*{Proposition 3.5}). By tilting 
(see \Cref{tilt-summary} and, for compatibility of definitions, \cite{BMS18}*{Lemma~3.20}), the universal property of adic localization and \eqref{monoid-id-pi} show that giving a map from a perfectoid space to 
\[
\tst\qqq \Spa\p{R_\infty\langle\f{x_1,\, \dotsc,\, x_d}{x_i}\rangle, R_\infty\langle\f{x_1,\, \dotsc,\, x_d}{x_i}\rangle^+}
\]
amounts to giving a map from its tilt to the perfectoid space 
\[
\tst\qqq\Spa\p{\wh{R}_\infty^\flat\langle\f{x_1^\flat,\, \dotsc,\, x_d^\flat}{x_i^\flat}\rangle, \wh{R}_\infty^\flat\langle\f{x_1^\flat,\, \dotsc,\, x_d^\flat}{x_i^\flat}\rangle^+},
\]
compatibly with overlaps of such rational subsets. Thus, the perfectoid space $U^\ad_{\wh{R}_\infty^\flat}$ represents the diamond associated  by \cite{Sch17}*{Definition~15.5} to the analytic adic space $U^\ad_{R_\infty}$.  Consequently, since the functor that sends an analytic adic space over $\bZ_p$ to its associated diamond induces an equivalence of \'{e}tale sites \cite{Sch17}*{Lemma~15.6}, we obtain the key
\[
\qqq H^i_\et(U_{R_\infty}^\ad, \bZ/\ell\bZ) \cong H^i_\et(U_{\wh{R}_\infty^\flat}^\ad, \bZ/\ell\bZ). 
\]
It remains to set 
\[
\qqq U_{R_\infty} \ce \Spec(R_\infty) \setminus \{ x_1 = \dotsc = x_d = 0\}
\]
and
\[
\qqq U_{R^\flat_\infty} \ce \Spec(R^\flat_\infty) \setminus \{ x^\flat_1 = \dotsc = x^\flat_d = 0\}
\]
and apply \cite{Hub96}*{Theorem~3.2.10}\footnote{Due to the blanket Noetherianity assumption  \cite{Hub96}*{paragraph~1.1.1}, the citation does not apply directly, so one performs a slightly tedious limit argument, similar to the one used in the proof of \cite{brauer-purity}*{Theorem 4.10}. Another way around this is to pass to the adic spectra at the finite levels of the tower and then use \cite{Sch17}*{Proposition 14.9} to pass to the limit of adic spaces.} (for the flanking isomorphisms) to deduce the passage to the tilt:
\[
\qq\ \ H^i_\et(U_{R_\infty}, \bZ/\ell\bZ) \cong  H^i_\et(U_{R_\infty}^\ad, \bZ/\ell\bZ) \cong H^i_\et(U_{\wh{R}_\infty^\flat}^\ad, \bZ/\ell\bZ) \cong H^i_\et(U_{R_\infty^\flat}, \bZ/\ell\bZ).
\]
\erems




\csub[The \'{e}tale depth is at least the virtual dimension] \lab{et-depth-geom-depth}

Purity for \'{e}tale cohomology of possibly singular Noetherian local rings $R$ was settled in the case when $R$ is an excellent $\bQ$-algebra in \cite{SGA2new}*{expos\'{e} XIV, th\'{e}or\`{e}me 5.6} and in the case when $R$ is a complete intersection in \cite{Gab04}*{Theorem~3} (by reduction to \cite{Ill03}*{th\'{e}or\`{e}me~2.6}). In \Cref{ell-semipurity}, we deduce the general case from \Cref{abs-coh-pur}. We begin with the definition of the virtual dimension, which is a numerical invariant of $R$ that has already appeared in the context of purity in \cite{SGA2new}*{expos\'{e} XIV, d\'{e}finition 5.3}.

\bpp[The virtual dimension of a Noetherian local ring] \lab{geom-depth-def}
For a Noetherian local ring $(R, \fm)$, by Cohen's theorem \cite{EGAIV1}*{chapitre 0, th\'{e}or\`{e}me 19.8.8 (i)}, the $\fm$-adic completion $\wh{R}$ is of the form
\[
\wh{R} \cong \wt{R}/I \qq \text{for a complete regular local ring} \q \wt{R} \q \text{and an ideal} \q I \subset \wt{R}.
\]
The \emph{virtual dimension} of $R$ is 
\be \label{eqn:vdim-def}
\vdim(R) \ce \dim(\wt{R}) - (\text{minimal number of generators for $I$})
\ee
and, by \cite{SGA2new}*{expos\'{e} XIV, proposition 5.2}, does not depend on the presentation $\wt{R}/I$. 
By \cite{SGA2new}*{expos\'{e} XIV, proposition 5.4}, 
\be \label{vdim-ineq}
\vdim(R) \le \dim(R) \q \text{with equality iff $R$ is a complete intersection.}
\ee
By definition, $\vdim(R) = \vdim(\wh{R})$, so also $\vdim(R) = \vdim(R^h)$; 
more generally, by \cite{Avr77}*{Proposition 3.6 and equations (3.2.1) and (3.2.2)}, for any flat local homomorphism $R \ra R'$ of Noetherian local rings, we have
\be \label{vdim-add}
\vdim(R') = \vdim(R) + \vdim(R'/\fm R'), 
\ee
so, in particular,
\[
\vdim(R) = \vdim(R^{\sh})
\]
(\emph{loc.~cit.}~proves this for the \emph{complete intersection defect} defined as $\dim(*) - \vdim(*)$ but the dimension $\dim(*)$ is likewise additive, see \cite{EGAIV2}*{corollaire~6.1.2}). 
\epp

\brem
Despite the name ``geometric depth'' used for $\vdim(R)$ in \cite{SGA2new}*{expos\'{e}~XIV, d\'{e}finition 5.3}, in general there is no inequality between $\depth_\fm(R)$ and $\vdim(R)$: a Cohen--Macaulay $R$ that is not a complete intersection has $\depth_\fm(R) > \vdim(R) $, whereas, due to \cite{Bur68} (or \cite{Koh72}*{Theorem~A}) and the Auslander--Buchsbaum formula, any regular local ring $\wt{R}$ has an ideal $I$ generated by three elements such that the quotient 
\[
R \ce \wt{R}/I \qxq{has} \depth_\fm(R) = 0
\]
(and $\vdim(R) \ge \dim(\wt{R}) - 3$). 
\erem

To deduce \Cref{ell-semipurity} from \Cref{abs-coh-pur}, we will use the following Lefschetz hyperplane theorem in local \'{e}tale cohomology. This strategy is close in spirit to the one used by Mich\`{e}le Raynaud in \cite{SGA2new}*{expos\'{e} XIV, th\'{e}or\`{e}me 5.6} in the case when $R$ is an excellent $\bQ$-algebra. 



\blem \label{lem:local-hyperplane}
For a regular local ring $(R, \fm)$, an $f \in \fm$, and an invertible in $R$ prime $\ell$, the map
\[
\q H^i_\fm(R, \bZ/\ell\bZ) \ra H^i_\fm(R/(f), \bZ/\ell\bZ) \qxq{is} \begin{cases} \x{bijective for} \q i < \dim(R) - 1, \\ \x{injective for} \q i = \dim(R) - 1. \end{cases}
\]
\elem

\bpf
Letting $j \colon \Spec(R[\f{1}{f}]) \hra \Spec(R)$ be the indicated open immersion, we need to show that
\be \label{eqn:j!-vanishes}
H^i_{\fm}(R, j_!(\bZ/\ell \bZ)) \cong 0 \q \x{for} \q i < \dim(R).
\ee
Moreover, by the local-to-global spectral sequence \cite{SGA4II}*{expos\'{e}~V, proposition~6.4}, we may assume that $R$ is strictly Henselian and, by \Cref{lem:et-to-hat}, that $R$ is also complete. We will derive \eqref{eqn:j!-vanishes} from Gabber's affine Lefschetz theorem \cite{ILO14}*{expos\'{e} XV, corollaire 1.2.4}, which gives
\be \lab{aff-Lef}
\tst H^i(R[\f{1}{f}], \bZ/\ell\bZ) \cong 0 \q \x{for} \q i > \dim(R).
\ee
Namely, by 
\cite{ILO14}*{expos\'{e} XVII, th\'{e}or\`{e}me 0.2}, the complex of \'{e}tale sheaves $(\mu_\ell^{\tensor \dim(R)})[2\dim(R)]$  on $R$ is dualizing and its $j^!$-pullback is dualizing on $R[\f{1}{f}]$. This pullback is $(\mu_\ell^{\tensor \dim(R)})[2\dim(R)]$ (see \cite{SGA4III}*{expos\'{e}~XVIII, proposition~3.1.8~(iii)}), so, since 
\[
\bR j_* \circ \bD \cong \bD \circ j_!
\]
(see \cite{SGA5}*{expos\'{e}~I, proposition~1.12~(a)}), the vanishing \eqref{aff-Lef} amounts~to
\[
\tst \p{H^{i}\p{\bR\sH om\p{j_!(\bZ/\ell \bZ), (\mu_\ell^{\tensor \dim(R)})[2\dim(R)]}}}_{\fm} \cong 0 \q \text{for} \q i > -\dim(R),
\]
where $(-)_{\fm}$ indicates the stalk. To obtain \eqref{eqn:j!-vanishes} it remains to use local duality in \'{e}tale cohomology \cite{SGA5}*{expos\'{e}~I, \'{e}quation~(4.2.2)} (our dualizing complex is normalized 
as there by \cite{ILO14}*{expos\'{e}~XVI, th\'{e}or\`{e}me~3.1.1}). 
\epf




\bthm \lab{ell-semipurity}
For a Noetherian local ring $(R, \fm)$ and a commutative, finite, \'{e}tale $R$-group $G$ whose order is invertible in $R$,
\[
H^i_\fm(R, G) \cong 0 \q \text{for} \q i < \vdim(R).
\]
\ethm


We will remove the assumption on the order of $G$ in \Cref{ell-semipurity-p} below.

\bpf
\Cref{abs-coh-pur} and \Cref{lem:local-hyperplane} settle the case when $R$ is a hypersurface, that is, a quotient of a regular ring by a principal ideal. Thus, it suffices to show how to reduce from a general $R$ to a hypersurface. This reduction works for any commutative, finite, \'{e}tale $G$, so, to be able to reuse it in the proof of \Cref{ell-semipurity-p}, we drop the assumption on the order of $G$.

By the spectral sequence \cite{SGA4II}*{expos\'{e}~V, proposition~6.4}, we may assume that $R$ is strictly Henselian, and then, by d\'{e}vissage, that $G \simeq \bZ/p\bZ$ for a prime $p$. \Cref{lem:et-to-hat} reduces further to complete $R$, so that %
\[
R \simeq \wt{R}/(f_1, \dotsc, f_n)
\]
for a complete regular local ring $(\wt{R}, \wt{\fm})$ and $f_1, \dotsc, f_n \in \wt{\fm}$ chosen so that $f_1, \dotsc, f_n$ is a minimal generating set for the ideal $(f_1, \dotsc, f_n) \subset \wt{R}$ (see \S\ref{conv}). Suppose that $n > 1$ and consider the rings 
\[
R_1 \ce \wt{R}/(f_1, \dotsc, f_{n - 1}) \qxq{and} R_2 \ce \wt{R}/(f_n),
\]
as well as 
\[
R_{12} \ce \wt{R}/(f_1f_n, \dotsc, f_{n - 1}f_n).
\]
Set theoretically, in $\Spec(\wt{R})$ we have 
\[
\Spec(R_1) \cap \Spec(R_2) = \Spec(R) \q \x{and} \q \Spec(R_1) \cup \Spec(R_2) = \Spec(R_{12}).
\]
Thus, since the \'{e}tale site is insensitive to nilpotents, we obtain the exact sequence
\[
0 \ra \underline{\bZ/p\bZ}_{R_{12}} \xra{z\, \mapsto\, (z|_{R_1},\, z|_{R_2})} \underline{\bZ/p\bZ}_{R_1} \oplus \underline{\bZ/p\bZ}_{R_2} \xra{(x,\, y)\, \mapsto\, x|_R - y|_{R}} \underline{\bZ/p\bZ}_{R} \ra 0
\]
of sheaves on $\Spec(\wt{R})_\et$. Since $\vdim(R_*) \ge \vdim(R) + 1$ for $* \in \{1, 2, 12\}$, the associated long exact sequence of cohomology with supports reduces the desired vanishing to its counterpart for the rings $R_*$. This allows us to decrease $n$, so  we arrive at $n = 1$, that is, at $R$ being a hypersurface $\wt{R}/(f)$. 
\epf



\csub[Zariski--Nagata purity for rings of virtual dimension $\ge 3$] \lab{nonabelian}

The nonabelian analogue of \Cref{ell-semipurity} is the following generalization of the Grothendieck--Zariski--Nagata purity theorem \cite{SGA2new}*{expos\'{e} X, th\'{e}or\`{e}me 3.4} (and of the main result of \cite{Cut95}): for extending finite \'{e}tale covers over a  closed subscheme of a Noetherian scheme, it suffices to assume that the total space have virtual dimension $\ge 3$ at the missing points (instead of even being a complete intersection of dimension $\ge 3$ at these points).
 We learned from de Jong that his generalizations contained in \cite{SP} of the algebraization theorems from \cite{SGA2new} and \cite{Ray75} could be used to prove this---indeed, as the reader will notice, they are the main inputs to the proof of \Cref{non-ab-et}. This section is purely classical and does not use perfectoid inputs or other sections.

\bthm \lab{non-ab-et}
Let $(R, \fm)$ be a Noetherian local ring, set $U_R \ce \Spec(R) \setminus \{ \fm\}$, and let $\Spec(R)_\fet$ \up{resp.,~$(U_{R})_\fet$} denote the category of finite \'{e}tale $R$-schemes \up{resp.,~of finite \'{e}tale $U_{R}$-schemes}.
\benum
\item \lab{NAE-a}
If $\vdim(R) \ge 2$, then the pullback $\Spec(R)_\fet \ra (U_{R})_\fet$ is fully faithful, $U_R$ is connected, and 
\[
\q \pi_1^\et(U_R) \surjects \pi_1^\et(R).
\]

\item \lab{NAE-b}
If $\vdim(R) \ge 3$, then the pullback $\Spec(R)_\fet \ra (U_{R})_\fet$ is an equivalence of categories and
\[
\q \pi_1^\et(U_R) \isomto \pi_1^\et(R).
\]
\eenum
\ethm

\bpf
In \ref{NAE-a}, granted the full faithfulness, the connectedness of $U_R$ follows from that of $\Spec(R)$  by considering sections both of the finite \'{e}tale map $\Spec(R) \bigsqcup \Spec(R) \ra \Spec(R)$ and of its base change to $U_R$. Moreover, by \cite{SGA1new}*{V, 6.9, 6.10}, 
the conclusions about the fundamental groups follow from the claims about the functors. For the latter, patching \cite{FR70}*{proposition~4.2} and flat descent allow us to replace $R$ by its $\fm$-adic completion. Then we may write 
\[
R \simeq \wt{R}/(f_1, \dotsc, f_n)
\]
for a complete regular local ring $(\wt{R}, \wt{\fm})$ and $f_1, \dotsc, f_n \in \wt{\fm}$ chosen so that $f_1, \dotsc, f_n$ is a minimal generating set for the ideal $\ff \ce (f_1, \dotsc, f_n) \subset \wt{R}$ (see \S\ref{conv}). Let $\fU \subset \fX$ be the formal schemes obtained from 
\[
U_{\wt{R}} \ce \Spec(\wt{R}) \setminus \{ \wt{\fm}\} \subset  \Spec(\wt{R})
\]
 by formal $\ff$-adic completion. 
Since  \'{e}tale sites are insensitive to nilpotents, pullback gives equivalences of categories
\[
\fX_\fet \isomto \Spec(R)_\fet \qq \x{and} \qq \fU_\fet \isomto (U_{R})_\fet,
\]
so we are reduced to considering the pullback functor $\fX_\fet \ra \fU_\fet$. By \cite{SP}*{Lemma~\href{https://stacks.math.columbia.edu/tag/09ZL}{09ZL}}, 
the pullback $\Spec(\wt{R})_\fet \isomto \fX_\fet$ is an equivalence, so we only need to consider the pullback $\Spec(\wt{R})_\fet \ra \fU_\fet$.

\benum
\item
The assumption $\vdim(R) \ge 2$ allows us to apply the algebraization of formal sections \cite{SP}*{Lemma~\href{https://stacks.math.columbia.edu/tag/0DXR}{0DXR} (with Lemma~\href{https://stacks.math.columbia.edu/tag/0DX9}{0DX9})} 
(or \cite{Ray75}*{chapitre~I, th\'{e}or\`{e}me~3.9}) to conclude that any morphism between the $\fU$-pullbacks of finite \'{e}tale $\wt{R}$-schemes $Y_1$ and $Y_2$ algebraizes to a morphism between their $U$-pullbacks for some open $U_R \subset U \subset U_{\wt{R}}$, and this algebraization is unique up to shrinking $U$. The complement 
\[
\qq Z \ce \Spec(\wt{R}) \setminus U
\] 
is at most $n$-dimensional because the $f_i|_Z$ cut out the closed point of $Z$. Thus, since $n \le \dim(\wt{R}) - 2$, the codimension of $Z$ in $\Spec(\wt{R})$ is $\ge 2$, to the effect that the algebraized morphism extends uniquely to an $\wt{R}$-morphism between $Y_1$ and $Y_2$ (see \cite{EGAIV2}*{th\'{e}or\`{e}me 5.10.5}). Consequently, the pullback $\Spec(\wt{R})_\fet \ra \fU_\fet$ is fully faithful.

\item
The assumption $\vdim(Rd) \ge 3$ allows us to apply the algebraization of coherent formal modules \cite{SP}*{Lemma~\href{https://stacks.math.columbia.edu/tag/0EJP}{0EJP}} to conclude that any finite \'{e}tale $\fU$-scheme algebraizes to a finite \'{e}tale $U$-scheme $Y$ for some open $U_R \subset U \subset U_{\wt{R}}$ (to algebraize the algebra structure maps, we use \cite{SP}*{Lemma~\href{https://stacks.math.columbia.edu/tag/0DXR}{0DXR}} as in \ref{NAE-a}). 
The complement of $U$ in $\Spec(\wt{R})$ is now of codimension $\ge 3$, 
so, by the Zariski--Nagata purity for regular rings \cite{SGA2new}*{expos\'{e} X, th\'{e}or\`{e}me 3.4 (i)}, we may extend $Y$ to a finite \'{e}tale $\wt{R}$-scheme. Consequently, the pullback $\Spec(\wt{R})_\fet \ra \fU_\fet$ is essentially surjective.
\qedhere
\eenum
\epf

\brem
The connectedness of $U_R$ holds more generally, see, for instance, \cite{SP}*{Lemma~\href{https://stacks.math.columbia.edu/tag/0ECR}{0ECR}} or \cite{Var09}*{Theorem~1.6}.
\erem



\section{Inputs from crystalline and prismatic Dieudonn\'{e} theory} \label{chapter:Dieudonne-theory}

Our eventual source of the characteristic-primary aspects of purity for flat cohomology is a relation to coherent cohomology and the vanishing of the latter in presence of enough depth. To exhibit this relation, we use crystalline and prismatic Dieudonn\'{e} theories that classify commutative, finite, locally free groups of $p$-power order over perfect and perfectoid rings in terms of Dieudonn\'{e} modules. We review the crystalline classification in \S\ref{section:Dieudonne-positive-characteristic} and its prismatic generalization in \S\ref{section:Dieudonne-mixed-characteristic}. 




\csub[Finite, locally free groups of $p$-power order over perfect rings] \label{section:Dieudonne-positive-characteristic}

The positive characteristic case of the key formula \eqref{eqn:key-formula} is a perfect ring variant of Kato--Trihan's \cite{KT03}*{Proposition~5.10}. To establish it in \Cref{KT-input} we build on  Gabber's suggestion to use the pro-fppf site (see \Cref{thm:KTG-input}; alternatively, one could adapt the arguments of \emph{op.~cit.}). A key input 
is the crystalline classification of commutative, finite, locally free groups of $p$-power order over perfect $\bF_p$-schemes due to Berthelot, Gabber, and Lau, which we now review.



\bpp[Crystalline Dieudonn\'{e} modules over perfect bases] \label{Gab-eq}
For a perfect $\bF_p$-scheme $S$, by an unpublished result of Gabber that built on \cite{Ber80}*{th\'{e}or\`{e}me 3.4.1} and was reproved by Lau \cite{Lau13}*{Corollary 6.5} by a different method, there is a covariant equivalence of categories
\[
G \mapsto \bM(G)
\]
from the category of commutative, finite, locally free $S$-groups $G$ that are locally on $S$ of $p$-power order to the category of quasi-coherent 
$W(\sO_S)$-modules $\bM$ equipped with Frobenius (resp.,~inverse Frobenius) semilinear maps $F\colon \bM \ra \bM$ (resp.,~$V\colon \bM \ra \bM$) with $FV = VF = p$ such that for every  affine open $\Spec(R) \subset S$ the $W(R)$-module $\Gamma(R, \bM)$ is finitely presented, killed by a power of $p$, and of projective dimension $\le 1$. The functor is defined by Zariski-local glueing as follows (see \cite{Lau13}*{proof of Corollary 6.5}): Zariski locally on $S$ one finds $p$-divisible groups $H_0$, $H_1$ that fit into an exact sequence
\[
0 \ra G \ra H_0 \ra H_1 \ra 0 \qq \x{and sets} \qq \bM(G) \ce \Coker(\bM(H_0) \hra \bM(H_1)),
\]
where $\bM(H_i) \ce \Gamma((S/\bZ_p)_\cris, \bD(H_i))$ is the evaluation of the covariant Dieudonn\'{e} crystal 
\[
\bD(H_i) \ce \sE xt^1_{(S/\bZ_p)_\cris}(H_i^*, \sO_{(S/\bZ_p)_\cris}), 
\]
which, by \cite{BBM82}*{section 5.3}, is identified with,
 \[
(\sE xt^1_{(S/\bZ_p)_\cris}(H_i, \sO_{(S/\bZ_p)_\cris}))^* 
\]
(the dual of the locally free crystal of $H_i$ defined in \cite{BBM82}*{d\'{e}finition~3.3.6, th\'{e}or\`{e}me~3.3.10}) at the terminal ind-object $\{(S, W_n(\sO_S))\}_{n \ge 0}$ of the crystalline site $(S/\bZ_p)_\cris$. 
Since \emph{op.~cit.}~uses big crystalline sites, the formation of $\bM(G)$ commutes with base change to any perfect $S$-scheme. 
\epp

\beg \label{Gab-eq-ex}
By \cite{BBM82}*{exemple~4.2.16~(i)}, the $W(\sO_S)$-module that underlies $\bb1_n \ce \bM(\bZ/p^n\bZ)$ is
\[
W(\sO_S)/p^n \qxq{with} F = p\cdot\Frob(-), \q V = \Frob\i(-).
\]
Likewise, the $W(\sO_S)$-module that underlies $\bM(\mu_{p^n})$ is  
\[
W(\sO_S)/p^n \qxq{with} F = \Frob(-), \q V = p\cdot\Frob\i(-).
\] 
\eeg

\blem \label{MG-exact}
In \uS\uref{Gab-eq}, the functor $G \mapsto \bM(G)$ and its inverse preserve short exact sequences. 
\elem

\bpf
By \cite{BBM82}*{proposition~1.1.7, th\'{e}or\`{e}me~3.3.3}, the functor $H_i \mapsto \bD(H_i)$ preserves short exact sequences, hence so does $H_i \mapsto \bM(H_i)$ (compare with \cite{BBM82}*{proposition~1.1.19, d\'{e}finition~1.2.1}). 
To deduce the same for the functor $G \mapsto \bM(G)$, thanks to the snake lemma, it suffices to Zariski locally on $S$ embed a short exact sequence $0 \ra G \ra G' \xra{\pi} G'' \ra 0$ into a short exact sequence $0 \ra H_0 \ra H_0' \ra H_0'' \ra 0$ of $p$-divisible groups (see \cite{BBM82}*{lemme~3.3.12}). For this, we choose Zariski local embeddings into $p$-divisible groups $\iota'\colon G' \hra H_0$ and $\iota'' \colon G'' \hra H_0''$ and replace $\iota'$ by 
\[
(\iota',\, \iota''\circ \pi)\colon G' \hra  H_0 \times H_0'' \equalscolon H_0'.
\] 
For the remaining exactness of the inverse, granted that 
\[
0 \ra \bM(G_1) \ra \bM(G_2) \ra \bM(G_3) \ra 0
\]
is a short exact sequence, we need to show that the complex $G_1 \ra G_2 \ra G_3$ is also a short exact sequence. In the case when $S$ is geometric point we may decompose this complex into short exact sequences of finite flat group schemes and conclude by the exactness of $G \mapsto \bM(G)$. Thus, in general we check on   geometric $S$-fibers that $G_1 \isomto \Ker(G_2 \ra G_3)$ (see \cite{EGAIV4}*{corollaire 17.9.5}). It then remains to note that $G_2/G_1 \hra G_3$ becomes an isomorphism after applying $\bM(-)$, so is an isomorphism.
\epf

The equivalence $G \mapsto \bM(G)$ leads to the following description of the low degree cohomology of $G$. 

\bprop \label{lem:syntomic-low-coh-id}
Let $S$ be a perfect $\bF_p$-scheme and let $G$ be a commutative, finite, locally free $S$-group killed by $p^n$. 
 We have the following functorial in $G$ identifications of sheaves on $S_\et$\ucolon
\[\ba \label{H0-cris-id} 
G\cong \,\sH om_S(\bZ/p^n\bZ, G) \overset{\ref{Gab-eq}}{\cong}  \sH om_{W_n(\sO_S),\,F,\, V}(\bb1_n, \bM(G)) \cong \bM(G)^{V = 1}, \\
\sE xt^1_S(\bZ/p^n\bZ, G) \overset{\ref{Gab-eq},\, \ref{MG-exact}}{\cong}  \sE xt^1_{W_n(\sO_S),\, F,\, V}(\bb1_n, \bM(G)) \cong \bM(G)/(V - 1)(\bM(G)),
\ea \]
where $\sE xt^1_S$ denotes the \'{e}tale sheafification of the functor of extensions of fppf $\bZ/p^n\bZ$-module sheaves.
\eprop

\bpf
The full faithfulness of $G \mapsto \bM(G)$ gives the first line of the displayed identification: the map to $\bM(G)^{V = 1}$ is the evaluation $(f \colon \bb1_n\, \ra\, \bM(G))\, \mapsto\, f(1)$. For the second line, we define the last identification as follows. To a local section $m$ of $\bM(G)$, we associate the extension $\bM(G) \oplus \bb1_n$ for which the Verschiebung is determined by $(0, 1) \mapsto (m, 1)$ and the Frobenius is then necessarily determined by $(0, 1) \mapsto (-F(m), p)$ (we write $F$ and $V$ for those of $\bM(G)$). Such extensions for $m$ and $m'$ are isomorphic if and only if the isomorphism of $W_n(\sO_S)$-module extensions determined by $(0, 1) \mapsto (a, 1)$ for some $a \in \bM(G)$ is $V$- and $F$-equivariant. The $V$-equivariance amounts to $(m + a, 1) = (m' + V(a), 1)$, that is, to $m - m' \in (V - 1)(\bM(G))$, and the $F$-equivariance amounts to $(-F(m) + pa, p) = (-F(m') + F(a), p)$ and follows from $V$-equivariance. Since any extension of $\bb1_n$ by $\bM(G)$ is \'{e}tale locally split as an extension of $W_n(\sO_S)$-modules, the claim follows.
\epf

In \Cref{KT-input}, we will upgrade the identifications of \Cref{lem:syntomic-low-coh-id} to a formula that expresses the flat cohomology $R\Gamma(S, G)$ in terms of the quasi-coherent cohomology $R\Gamma(W_n(S), \bM(G))$. For this, we first show in \Cref{thm:KTG-input} that the map $V - 1 \colon \bM(G) \ra \bM(G)$ is pro-fppf locally surjective.


\bpp[The pro-fppf site] \label{X-profppf}
A scheme map $X' \ra X$ is \emph{pro-ppf} (resp.,~\emph{pro-fppf}) if $X'$ may be covered by opens $\Spec(A') \subset X'$ for which there is a factorization 
\[
\Spec(A') \ra \Spec(A) \subset X
\]
in which the $A$-algebra $A'$ is a filtered direct limit of flat, finitely presented $A$-algebras (compare with \S\ref{ind-fppf-maps}) and $\Spec(A) \subset X$ is an open immersion (resp.,~and $X' \ra X$ is faithfully flat). Pro-ppf maps 
$\{X_i' \ra X\}_{i \in I}$ 
form a \emph{pro-fppf cover} if each quasi-compact open of $X$ is a finite union of images of quasi-compact opens of $\bigsqcup_{i \in I} X_i'$.  By \S\ref{ind-fppf-maps}, the pro-ppf maps are stable under base change and composition, so the category of $X$-schemes with pro-fppf covers as coverings defines the \emph{pro-fppf site} of $X$. A pro-fppf cover is also an fpqc cover, so an fpqc sheaf is also a pro-fppf sheaf.
\epp

\blem \label{lem:perf-still-pro-fppf}
Let $S$ be a perfect $\bF_p$-scheme and let $X^\perf$ denote the perfection of an $\bF_p$-scheme $X$. If $\{X_i \ra S\}_{i \in I}$ is an fpqc \up{resp.,~pro-fppf} cover of $S$, then so is $\{X_i^\perf \ra S\}_{i \in I}$.
\elem

\bpf
The composite 
\[
X_i \xra{\Frob^n} X_i \ra S \qxq{factors as} \xymatrix@C=10pt{X_i \ar[r] &S \ar[rr]^-{\Frob^n}_-{\sim} &&S},
\]
so it is fpqc (resp.,~pro-fppf). Thus, the inverse limit $X_i^\perf$ of 
\[
\dotsc \xra{\Frob} X_i \xra{\Frob} X_i
\]
is fpqc (resp.,~pro-fppf) over $S$ (see \S\ref{ind-fppf-maps}).
\epf

\bprop \label{thm:KTG-input}
Let $S$ be a perfect $\bF_p$-scheme and let $G$ be a commutative, finite, locally free $S$-group that is locally on $S$ of $p$-power order. There is a functorial in $G$ short exact sequence 
\be \label{eqn:KTG-input-SES}
0 \ra G \xra{\ref{lem:syntomic-low-coh-id}} \bM(G) \xra{V - 1} \bM(G) \ra 0
\ee
of sheaves on the category of perfect $S$-schemes endowed with the pro-fppf topology. 
\eprop

\bpf
The left exactness follows from \Cref{lem:syntomic-low-coh-id}. For the remaining surjectivity we may work \'{e}tale locally on $S$, so, by \Cref{lem:syntomic-low-coh-id}, we need to show that a given extension 
\[
0 \ra G \ra G' \xra{\pi} \underline{\bZ/p^n\bZ}_S \ra 0
\]
of $\bZ/p^n\bZ$-module sheaves splits over a pro-fppf cover of $S$. The extension splits over the fppf cover $\pi\i(1) \surjects S$, which is a $G$-torsor, 
so \Cref{lem:perf-still-pro-fppf} supplies a desired pro-fppf cover. (The point of using \Cref{lem:perf-still-pro-fppf}, so, relatedly, the pro-fppf rather than simply the fppf topology, is to stay in the realm of perfect base schemes in order to be able to consider the functor $\bM(-)$.)
\epf


\bthm \label{KT-input}
Let $S$ be a perfect $\bF_p$-scheme and let $G$ be a commutative, finite, locally free $S$-group that locally on $S$ is of $p$-power order. The map  of sites $\eps \colon S_\fppf \ra S_{\et}$ gives a functorial triangle
\be \label{eqn:KT-triangle}
R\eps_*(G) \ra \bM(G) \xra{V - 1} \bM(G) \ra (R\eps_*(G))[1] \qxq{on} S_\et.
\ee
In particular, if $G$ is killed by $p^n$, then, 
for any closed subset $Z \subset S$, we have
\be \label{crystalline-cor}
R\Gamma_Z(S, G) \cong R\Gamma_Z(W_n(S), \bM(G))^{V = 1}  
\ee
functorially in $S$, $Z$, and $G$.
\ethm


\bpf
The identification \eqref{crystalline-cor} follows from \eqref{eqn:KT-triangle} by applying $R\Gamma_Z(S, -)$. For the latter, we fix a suitable auxiliary 
cutoff cardinal $\kappa$ with 
$\kappa > \abs{S}$ (see \S\ref{conv}), consider the small pro-fppf site $S_{\profppf,\, \kappa}$ bounded by $\kappa$, its subsite $S_{\profppf,\, \perf,\, \kappa}$ of perfect schemes, and the morphisms
\[
\xymatrix@C=15pt@R=12pt{
S_{\profppf,\, \kappa} \ar[r]_-b \ar[d]_-a & S_{\profppf,\, \perf,\, \kappa} \ar[d]_-c \\
S_\fppf \ar[r]^-{\eps} & S_\et.
}
\]
By limit arguments, 
$R^{\ge 1}a_*(G) \cong 0$, so 
\[
R\eps_*(G) \cong R(\eps \circ a)_*(G) \cong R(c \circ b)_*(G) \cong Rc_*(Rb_*(G)).
\]
By \Cref{lem:perf-still-pro-fppf}, the functor $b_*$ is exact, so $R\eps_*(G) \cong Rc_*(G)$. Moreover, by 
faithfully flat descent for modules (see \cite{SP}*{Lemma~\href{https://stacks.math.columbia.edu/tag/023M}{023M}}), we have 
\[
R^{\ge 1}c_*(\bM(G)) \cong 0.
\]
Thus, by applying $Rc_*(-)$ to the short exact sequence \eqref{eqn:KTG-input-SES} we get the desired \eqref{eqn:KT-triangle} (independently of the choice of $\kappa$). 
\epf

\bcor \label{cor:perfect-vanish}
For a perfect $\bF_p$-algebra $A$ and a commutative, finite, locally free $A$-group $G$ of $p$-power order, 
\[
H^i(A, G) \cong 0 \qxq{for} i \ge 2.
\]
\ecor

\bpf
Since affines have no higher quasi-coherent cohomology, \eqref{crystalline-cor} with $Z = S$ suffices. 
\epf

\Cref{KT-input} implies much of the positive characteristic case of purity for flat cohomology. We record this in \Cref{thm:main-pos-char} because the intervening auxiliary lemmas are also important in the general case. The following example illustrates why positive characteristic is significantly simpler.

\beg \label{eg:depth-example}
For a Noetherian local ring $(R, \fm)$, by the cohomological characterization of depth, $H^i_\fm(R, \bG_a) \cong 0$ for $i < \depth(R)$ (compare with \cite{SGA2new}*{expos\'{e} III, proposition~3.3 (iv)}). Thus, if $R$ is also an $\bF_p$-algebra, then the Frobenius-kernel and the Artin--Schreier sequences give 
\be \label{eqn:depth-vanishing}
H^i_\fm(R, \gA_p) \cong 0 \qxq{and} H^i_\fm(R, \bZ/p\bZ) \cong 0 \q \x{for} \q i < \depth(R).
\ee
Complete intersection $R$ have $\depth(R) = \dim(R)$, so for them \eqref{eqn:depth-vanishing} gives cases of \Cref{main}. 
\eeg

The complete intersection assumption fully manifests itself in the following lemma.

\blem \label{mod-out-by-more}
Let $A$ be a ring, let $f_1, \dotsc, f_m \in A$ be a regular sequence, and let $G$ be a commutative $(A/(f_1, \dotsc, f_m))$-group scheme that is either smooth or finite locally free. For all $n_1, \dotsc, n_m > 0$ and $f_1^{1/n_1}, \dotsc, f_m^{1/n_m} \in A$ and all ideals $I \subset A/(f_1, \dotsc, f_m)$ containing some $(A/(f_1, \dotsc, f_m))$-regular sequence $a_1, \dotsc, a_d$ of length $d$, the map
\[
H^i_I(A/(f_1, \dotsc, f_m), G) \ra H^i_I(A/(f_1^{1/n_1}, \dotsc, f_m^{1/n_m}), G) 
\]  
is injective for $i < d$ and bijective for $i < d - 1$.
\elem

\bpf
We use the B\'{e}gueri sequence \eqref{Begueri-0} and the five lemma to assume that $G$ is smooth.
By \cite{SP}*{Lemma~\href{https://stacks.math.columbia.edu/tag/07DV}{07DV}}, the sequence $f_1^{i_1/n_1}, \dotsc, f_m^{i_m/n_m}, a_1, \dotsc, a_d$ is regular for all $i_1, \dotsc, i_m \ge 1$. Thus, if $n_1 > 1$, then, by induction on $m$,  the square-zero ideal $J$ that cuts out the closed immersion 
\[
j \colon \Spec(A/(f_1^{(n_1 - 1)/n_1}, f_2, \dotsc, f_m)) \hra \Spec(A/(f_1, f_2, \dotsc, f_m))
\]
is free as an $(A/(f^{1/n_1}_1, f_2, \dotsc, f_m))$-module.
 Moreover, deformation theory supplies the short exact sequence of \'{e}tale sheaves
\[
0 \ra \sH om (e^*(\Omega^1_{G/(A/(f_1,\, \dotsc,\, f_m))}), J) \ra G \ra j_*(G_{A/(f^{(n_1 - 1)/n_1}_1,\, f_2,\, \dotsc,\, f_m)}) \ra 0,
\]
where the $\sH om$ is of $A/(f_1,\, \dotsc,\, f_m )$-module sheaves and $e$ is the unit section of $G$ (see, for instance, \cite{topology-torsors}*{Lemma~B.14 and its proof}). Since $G$ is smooth, the $(A/(f^{1/n_1}_1,  f_2, \dotsc, f_m))$-module 
\[
\Hom_{A/(f_1,\, \dotsc,\, f_m )} (e^*(\Omega^1_{G/(A/(f_1,\, \dotsc,\, f_m))}), J)
\]
is finite projective. Consequently, by decreasing induction on $t$ and the regularity of the sequence $a_1, \dotsc, a_d$, for $i < d - t$ we have
\[
\tst H^i_I\p{A/(f_1,\, \dotsc,\, f_m ), \p{\Hom (e^*(\Omega^1_{G/(A/(f_1,\, \dotsc,\, f_m))}), J)}\!/(a_1, \dotsc, a_t)} \cong 0.
\]
The $t = 0$ case of this vanishing and the preceding short exact sequence then imply that
\[
H^i_I(A/(f_1,\, \dotsc,\, f_m ), G) \ra H^i_I(A/(f^{(n_1 - 1)/n_1}_1,\, f_2,\, \dotsc,\, f_m ), G) 
\]
is injective for $i < d$ and bijective for $i < d - 1$. By repeating this argument for 
\[
\Spec(A/(f^{(n_1 - 2)/n_1}_1, f_2,\, \dotsc,\, f_m)) \hra \Spec(A/(f^{(n_1 - 1)/n_1}_1, f_2,\, \dotsc,\, f_m ))
\]
and so on, we eventually replace $f_1$ by $f_1^{1/n_1}$ and then conclude by induction on $m$.
\epf

Another useful reduction is the following passage to a cover (compare with \cite{brauer-purity}*{Remark 2.7}). 

\blem \lab{finite-cover} 
Let $(R, \fm) \ra (R', \fm')$ be a finite, flat, local map of Noetherian local rings that are 
complete intersections, let $n < \dim(R)$ be an integer, and let $G$ be a commutative, finite, flat $R$-group. If for each local ring $(A, \fq)$ at some maximal ideal of some finite nonempty self-product $R' \tensor_R \dotsc \tensor_R R'$ Theorem~\uref{main} holds for $A$ and $G$ in cohomological degrees $< n$ in the sense that
\[
H^i_\fq(A, G) \cong 0 \qxq{for} i < n, 
\]
then 
\[
H^i_\fm(R, G) \cong 0 \qxq{for} i < n \qxq{and} H^n_\fm(R, G) \hra H^n_{\fm'}(R', G).
\]
\elem

\bpf
The maximal ideals $\fq$ in question are the primes above $\fm$. Each $A$ is of dimension $\dim(R)$ and, by \eqref{vdim-ineq}--\eqref{vdim-add} (or by \cite{SP}*{Lemma~\href{https://stacks.math.columbia.edu/tag/09Q7}{09Q7}
}), is a complete intersection. Thus, the assumption gives 
\[
H^i_\fm(R' \tensor_R \dotsc \tensor_R R', G) \cong 0 \qxq{for} i < n.
\]
It then suffices to use the spectral sequence 
\[
E_2^{ij} = H^i(R'/R, H^j_\fm(-, G)) \Ra H^{i + j}_\fm(R, G)
\]
that results from fppf descent for $R'\mapsto R\Gamma_\fm(R',G)$ (and that could also be derived by choosing an injective resolution of $G$ and considering the \v{C}ech complexes of its terms with respect to $R'/R$).
\epf


\bprop \label{thm:main-pos-char}
For a complete, Noetherian, local $\bF_p$-algebra $(R, \fm)$ that is a complete intersection with a perfect residue field $k$ and a commutative, finite, flat $R$-group $G$ of $p$-power order,
\[
H^i_\fm(R, G) \cong 0 \qxq{for} i < \dim(R).
\]
\eprop

\bpf
We use induction on $i$ simultaneously for all $R$. By the Cohen structure theorem (see \S\ref{conv}), 
\[
R \simeq k\llb t_1, \dotsc, t_N \rrb/(f_1, \dotsc, f_m)
\]
for a $k\llb t_1, \dotsc, t_N \rrb$-regular sequence $f_1, \dotsc, f_m$. By \Cref{finite-cover}, 
we may pass to the limit of the rings $k\llb t_1^{1/p^j}, \dotsc, t_N^{1/p^j} \rrb/(f_1, \dotsc, f_m)$ 
to reduce to
\[
H^i_{(t_1,\, \dotsc,\, t_N)}(k\llb t_1^{1/p^{\infty}}, \dotsc, t_N^{1/p^{\infty}} \rrb/(f_1, \dotsc, f_m), G) \overset{?}{\cong} 0 \qxq{for} i < N - m.
\]
The ring $k\llb t_1^{1/p^\infty}, \dotsc, t_N^{1/p^\infty} \rrb$ is perfect, so 
\Cref{mod-out-by-more} and a limit argument reduce to 
\[
H^i_{(t_1,\, \dotsc,\, t_N)}(k\llb t_1^{1/p^{\infty}}, \dotsc, t_N^{1/p^{\infty}} \rrb/(f_1^{1/p^{\infty}}, \dotsc, f_m^{1/p^{\infty}}), G) \overset{?}{\cong} 0 \qxq{for} i < N - m.
\]
Since 
\[
R_\infty \ce k\llb t_1^{1/p^{\infty}}, \dotsc, t_N^{1/p^{\infty}} \rrb/(f_1^{1/p^{\infty}}, \dotsc, f_m^{1/p^{\infty}})
\]
is a perfect $\bF_p$-algebra, \eqref{crystalline-cor} reduces us to
\[
H^i_{(p,\, t_1,\, \dotsc,\, t_N)}(W(R_\infty), \bM(G)) \overset{?}{\cong} 0 \qxq{for} i < N - m. 
\] 
Since the $W(R_\infty)$-module $\bM(G)$ is of projective dimension $\le 1$ (see \S\ref{Gab-eq}), this, in turn, reduces to 
\[
H^i_{(p,\, t_1,\, \dotsc,\, t_N)}(W(R_\infty), W(R_\infty)) \overset{?}{\cong} 0 \qxq{for} i < N - m + 1. 
\] 
By \cite{SP}*{Lemma~\href{https://stacks.math.columbia.edu/tag/07DV}{07DV}}, an $R$-regular sequence $a_1, \dotsc, a_{N - m} \in \fm$ is $R_\infty$-regular, so the sequence $a_0 \ce p$, $a_1$, $\dotsc$, $a_{\dim(R)}$ is $W(R_\infty)$-regular. Decreasing induction on $j$ then gives the sufficient
\[
H^i_{(p,\, t_1,\, \dotsc,\, t_N)}(W(R_\infty), W(R_\infty)/(a_0, \dotsc, a_j)) \cong 0 
\]
for $i < N - m - j$ and $-1 \le j \le N - m$.
\epf




\csub[Finite, locally free groups of $p$-power order over perfectoids]\label{section:Dieudonne-mixed-characteristic}

The classification of commutative, finite, locally free groups of $p$-power order over perfect $\bF_p$-algebras 
was extended to perfectoid $\bZ_p$-algebras by Lau \cite{Lau18} in the case $p > 2$ and by the second-named author \cite{SW19}*{Appendix to Lecture~17} in general by using ideas from integral $p$-adic Hodge theory. In \cite{ALB20}, Ansch\"utz--Le Bras drew a parallel to the crystalline theory 
by relating this classification to the prismatic point of view. We will use these results for formulating (and proving) the general case of the key formula \eqref{eqn:key-formula}, so we now review them and include some relevant for us complements.

\bpp[Prismatic Dieudonn\'{e} modules over perfectoid rings] \label{pp:prismatic-review}\label{thm:dieudonne}
For a perfectoid $\bZ_p$-algebra $A$ and a fixed generator $\xi$ of $\Ker(\bA_\Inf(A) \surjects A)$, by \cite{ALB20}*{Theorem 5.4} (which builds on \cite{SW19}*{Theorem~17.5.2}), there is a covariant, compatible with base change 
equivalence of categories
\[
G \mapsto \bM(G) \ce \Ext^1_{A_\Prism}(G^*, \mathcal O_\Prism)
\]
from the category of commutative, finite, locally free $A$-groups $G$ of $p$-power order to the category of finitely presented, $p$-power torsion $\bA_\Inf(A)$-modules $\bM$ of projective dimension $\le 1$ equipped with Frobenius (resp.,~inverse Frobenius) semilinear maps 
\[
F \colon \bM \ra \bM \ \ \x{(resp.,}\ \ V \colon \bM \ra \bM) \ \ \x{that satisfy} \ \ FV = \Frob(\xi) \ \ \x{and}\ \ VF = \xi.
\]
The $\Ext^1$ above is in the absolute prismatic site of $A$ 
and, by \cite{BS19}*{Lemma~4.8}, is an $\bA_\Inf(A)$-module. The Frobenius $F$ is induced by the Frobenius of the prismatic structure sheaf $\cO_{\Prism}$, and the only role of $\xi$ is to define the Verschiebung $V$. By \cite{ALB20}*{Theorem 4.44 and the proof of Theorem~5.4}, in the case when $A$ is a perfect $\bF_p$-algebra, 
the functor $G \mapsto \bM(G)$ may be identified with its crystalline counterpart discussed in \S\ref{Gab-eq}. 
By construction, in the case when $G$ is the $p^n$-torsion of a $p$-divisible group, $\bM(G)$ is a finite projective $(\bA_\Inf(A)/p^n\bA_\Inf(A))$-module.
\epp

\beg \label{explicit-M(G)}
By \cite{ALB20}*{Lemma 4.75}, the $\bA_\Inf(A)$-module that underlies $\bM(\mu_{p^n})$ is 
\[
\bA_\Inf(A)/p^n\bA_\Inf(A) \qxq{with} F = \Frob(-),\q V = \xi\cdot\Frob\i(-).
\]
Likewise, 
the $\bA_\Inf(A)$-module that underlies ${\mathbf{1}}_n \ce \bM(\bZ/p^n\bZ)$ is  
\[
\bA_\Inf(A)/p^n\bA_\Inf(A) \qxq{with} F = \Frob(\xi)\cdot\Frob(-), \q V = \Frob\i(-).
\] 
\eeg

\blem \label{MG-exact-2}
In \uS\uref{pp:prismatic-review}, the functor $G \mapsto \bM(G)$ and its inverse preserve short exact sequences. 
\elem

\bpf
Cartier duality is exact, so the claim is part of \cite{ALB20}*{Theorem~5.4}. 
\epf


Similarly to \Cref{lem:syntomic-low-coh-id}, we obtain the following description of low degree cohomology of $G$.
%

\bprop \label{Hom-Ext-describe}
For a perfectoid $\bZ_p$-algebra $A$ and a commutative, finite, locally free $A$-group $G$ of $p$-power order, we have  functorial in $A$ and $G$ identifications
\begin{equation} \label{MG-fixed-pts}
\tst   G(A) \cong \bM(G)^{V = 1} \qxq{and} H^1(A, G) \cong \bM(G)/(V - 1)\bM(G). 
\end{equation}
\eprop

\bpf
For any $p^n$ that kills $G$, we have 
\[
H^1(A, G) \cong \Ext^1(\underline{\bZ/p^n\bZ}_A, G)
\]
(extensions of fppf $\bZ/p^n\bZ$-module sheaves), so the same argument as for \Cref{lem:syntomic-low-coh-id} (with \S\ref{pp:prismatic-review}, \Cref{explicit-M(G)}, and \Cref{MG-exact-2} in place of \S\ref{Gab-eq}, \Cref{Gab-eq-ex}, and \Cref{MG-exact}) gives the claim.
\epf

\brem
In \Cref{Hom-Ext-describe}, the $H^{\ge 2}(A, G)$ vanish: we will deduce this in \Cref{no-higher-coho} from its positive characteristic case of \Cref{cor:perfect-vanish} via the $p$-adic continuity formula of \S\ref{section:continuity-formula}.
\erem

We turn to analyzing the prismatic side of the key formula \eqref{eqn:key-formula}: we show that it satisfies $p$-complete arc descent in \Cref{M(G)-coho-descent} and then arc locally relate it to flat cohomology in \Cref{M(G)-local-analysis}. A key input to our 
arc descent results is the following lemma that in essence restates \cite{BS19}*{Proposition~8.10}.

\blem \label{arc-sheaves}
The following functors satisfy hyperdescent for those $p$-complete arc hypercovers whose terms are perfectoid $\bZ_p$-algebras\ucolon
\begin{equation*} \label{}
A \mapsto \bA_\Inf(A), \q A \mapsto A, \qxq{and} A \mapsto A^\flat.
\end{equation*}
\elem

\bpf
Since $\Ker(\theta\colon \bA_{\Inf}(A) \surjects 
A)$ and $\Ker(\bA_\Inf(A) \surjects A^\flat)$ are principal, generated by nonzerodivisors $\xi$ and $p$, and $\bA_\Inf(A)$ is $p$-adically complete, it suffices to treat $A \mapsto A^\flat$. Thus, fixing $A$ and letting $\wh{(-)}$ denote derived $p^\flat$-adic completion, \Cref{arc-tilt} (with \Cref{recognize}~\ref{R-a}) 
reduces us to showing that the functor $S \mapsto \wh{S}$ satisfies $p^\flat$-complete arc hyperdescent on the category of perfect $A^\flat$-algebras (by \S\ref{perfectoid-def} and \eqref{no-big-tor}, the tilts of perfectoid $A$-algebras are derived $p^\flat$-adically complete, so on them this functor agrees with $S \mapsto S$). Since the functor is bounded below (even concentrated in degree $0$), hyperdescent for it is equivalent to descent.  Moreover, it suffices to show arc descent: indeed, if $S \ra S'$ is a $p^\flat$-complete arc cover, then 
\[
\tst S \ra S' \times S[\frac 1{p^\flat}]
\]
is an arc cover and the functor 
has identical values on the two \v{C}ech nerves. We then instead consider the functor $S\mapsto \wh{S_\perf}$ defined on \emph{all} $A^\flat$-algebras, where 
\[
\tst S_\perf \ce \varinjlim_{s \mapsto s^p} S,
\]
and then reduce further to showing arc descent for $S \mapsto S_\perf$ on the category of all $A^\flat$-algebras $S$. By \cite{BS17}*{Theorem~4.1~(i), Proposition~4.5}, 
the functor $S \mapsto S_\perf$ satisfies $v$-descent, so, by  \cite{BM20}*{Proposition~4.8 and its proof}, 
it also satisfies arc descent. 
\epf

\bprop \label{M(G)-coho-descent}
Let $A$ be a perfectoid $\bZ_p$-algebra,  let   $Z \subset \Spec(A/pA)$ be a closed subset, let $a \in \bA_\Inf(A)$, and 
let $G$ be a commutative, finite, locally free $A$-group of $p$-power order. 
 The functors
\[
\tst A' \mapsto \bM(G_{A'})[\f 1a] \qxq{and} A' \mapsto R\Gamma_Z(\bA_\Inf(A'), \bM(G_{A'})) 
\]
satisfy hyperdescent for those $p$-complete arc hypercovers whose terms are perfectoid $A$-algebras.
\eprop

\bpf
For any cosimplicial abelian group $M^\bullet$, the associated complex 
\[
M^0 \ra M^1 \ra \dotsc
\]
represents $R\lim_\Delta(M^\bullet)$, 
so, since localization is exact, we may assume that $a = 1$. 
By  \S\ref{pp:prismatic-review}, we have 
\[
\bM(G_{A'}) \cong \bM(G) \tensor_{\bA_\Inf(A)} \bA_\Inf(A')
\]
and there are finite projective $\bA_\Inf(A)$-modules $M_i$ that fit into an exact sequence 
\[
0 \ra M_0 \ra M_1 \ra \bM(G) \ra 0.
\]
Since $\bM(G)$ is $p$-power torsion and $\bA_\Inf(A')$ is $p$-torsion free, this sequence stays exact after base change to $\bA_\Inf(A')$. Thus, for the claim about $A' \mapsto \bM(G_{A'})$ it suffices 
note that, by \Cref{arc-sheaves}, the functors 
\[
A' \mapsto M_i \tensor_{\bA_\Inf(A)} \bA_\Inf(A')
\]
satisfy hyperdescent for those $p$-complete arc hypercovers whose terms are perfectoid $A$-algebras. For the claim about the functor
\[
A' \mapsto R\Gamma_Z(\bA_\Inf(A'), \bM(G_{A'})),
\]
we use the functorial triangle 
\be
\ba
\tst R\Gamma_Z(\bA_\Inf(A'), \bM(G_{A'})) \ra &R\Gamma(\bA_\Inf(A'), \bM(G_{A'})) \\ &\tst\ra R\varprojlim_z \p{R\Gamma(\bA_\Inf(A')[\f{1}{z}], \bM(G_{A'})[\f 1z])}
\ea
\ee
where $z$ ranges over the elements of $\bA_\Inf(A)$ that vanish on $Z$. This reduces us to the settled claim about $A' \mapsto \bM(G_{A'})[\f 1a]$ because 
\[
\tst R\Gamma(\bA_\Inf(A')[\f{1}{z}], \bM(G_{A'})[\f 1z]) \cong \bM(G_{A'})[\f 1z]. \qedhere
\]
\epf

The promised arc-local analysis uses the following 
general lemma about modules on infinite products.

\blem \label{product-lemma}
For rings $\{R_i\}_{i \in I}$ and a finitely presented $(\prod_{i \in I} R_i)$-module $M$, we have 
\begin{equation*}
\tst M \isomto \prod_{i \in I} (M \tensor_R R_i).
\end{equation*}
\elem

\bpf
Set $R\ce \prod_{i\in I} R_i$ and choose a resolution $R^n\to R^m\to M\to 0$. Both $-\otimes_R R_i$ and 
infinite products are exact, so the claim for $M$ reduces to the evident case of a finite 
free $R$-module. 
\epf


\bprop \label{M(G)-products}
Let $\{A_i \}_{i \in I}$ be perfectoid $\bZ_p$-algebras, let $A \ce \prod_{i \in I} A_i$ \up{which is perfectoid by Proposition \uref{recognize}~\uref{R-c}}, and let $G$ be a commutative, finite, locally free $A$-group of $p$-power order. We have
\[
\tst \bM(G) \isomto \prod_{i \in I} \bM(G_{A_i})
\]
over 
\[
\tst \bA_\Inf(A) \isomto \prod_{i \in I} \bA_\Inf(A_i)
\]
compatibly with $F$ and $V$, granted that we choose compatible $\xi$ in $\bA_\Inf(A)$ and $\bA_\Inf(A_i)$ to define $V$ \up{see \uS\uref{pp:prismatic-review}}.
\eprop

\bpf
Since $G \mapsto \bM(G)$ commutes with base change (see \S\ref{pp:prismatic-review}), compatibility with $F$ and $V$ is automatic and \Cref{product-lemma} gives the claim.
\epf

\bprop \label{M(G)-local-analysis}
Let $\{A_i\}_{i \in I}$ be $p$-adically complete valuation rings  
of rank $\le 1$ such that  the fraction field $K_i$ of $A_i$ is algebraically closed, let $A \ce \prod_{i \in I} A_i$, let $K \ce \prod_{i \in I} K_i$, let $a \in A$, and let $G$ be a commutative, finite, locally free $A$-group of $p$-power order. The map
\[
\tst V - 1 \colon \bM(G)[\f 1{[a^\flat]}] \ra \bM(G)[\f 1{[a^\flat]}]
\]
is surjective \up{for any $\xi \in \bA_\Inf(A)$ used to define $V$ of $\bM(G)$, see \uS\uref{pp:prismatic-review}}, 
\begin{equation*}
\tst H^j(A[\f 1a], G) \cong 0 \qxq{for} j \ge 1,
\end{equation*}
and there is a unique commutative~square
\[
\tst \xymatrix@C=40pt{
G(A) \ar@{->>}[d]\ar[r]_-{\eqref{MG-fixed-pts}}^-{\sim} & \bM(G)^{V = 1} \ar@{->>}[d] \\
G(A[\f 1{a}]) \ar[r]^-{\sim} & (\bM(G)[\f 1{[a^\flat]}])^{V = 1},
}
\]
in which the vertical arrows are bijective whenever $a$ is a nonzerodivisor. 
\eprop

\emph{Proof.}
In the statement, $A$ is perfectoid by \eqref{p-oid-mod-pi} and \Cref{recognize}~\ref{R-c}, the element $a^\flat \in A^\flat \cong \prod_{i \in I} A_i^\flat$ is a system of compatible $p$-power roots of $a \in A$ (see \eqref{monoid-id}), and the map  
\[
\tst V \colon \bM(G)[\f 1{[a^\flat]}] \ra \bM(G)[\f 1{[a^\flat]}]
\]
is defined by 
\[
\tst V(\f{m}{[(a^\flat)^n]}) \ce \f{V(m)}{[(a^\flat)^{n/p}]}.
\]
By decomposing $A$ into subproducts, we may replace $I$ by parts of a finite partition of $I$. Thus, since the case $a = 0$ is clear, we assume that $a$ is a nonzerodivisor. In this case, since, by \Cref{tilt-summary}, each $A_i^\flat$ is a valuation ring, $a^\flat$ is also a nonzerodivisor, and the uniqueness aspect 
will follow from the bijectivity of the vertical maps. Moreover, $A \hra A[\f 1a] \hra K$, so, since $G(A) \isomto G(K)$ by the valuative criterion of properness, 
also $G(A) \isomto G(A[\f 1a])$.
Due to the local structure of $A[\f 1a]$ described in \Cref{lem:ultrafilters}, the valuative of properness applied locally on $\Spec(A[\f{1}{a}])$ also ensures that $X(A[\f 1a]) \isomto X(K)$
for every $G_{A[\f{1}{a}]}$-torsor $X$. Since $K$ is a product of algebraically closed fields, 
this shows that 
\[
\tst H^1(A[\f 1a], G) \cong 0.
\]
In contrast, \Cref{lem:ultrafilters} and the B\'{e}gueri sequence \eqref{Begueri-0} show that 
\[
\tst H^j(A[\f 1a], G) \cong 0 \qxq{for} j \ge 2.
\]
The remaining claims concern $\bM(G)$, and we use \cite{BBM82}*{th\'{e}or\`{e}me~3.1.1} to cover $\Spec(A)$ by finitely many affine opens $\Spec(A')$ over which $G$ embeds into a $p$-divisible group. Since each $A_i$ is local, the map $A \ra A_i$ factors through some $A \ra A'$, so we subdivide $I$ to assume that there is an exact sequence 
\begin{equation*} \label{}
0 \ra G \ra G' \ra G'' \ra 0
\end{equation*}
in which $G'$ is a truncated $p$-divisible group. Since 
\[
G(A) \isomto \bM(G)^{V = 1} \qxq{and} H^1(A, G) \cong 0,
\]
the functor 
\[
\tst G \mapsto \bM(G)^{V = 1} \qx{is exact}
\]
and, by \Cref{MG-exact-2}, the functor 
\[
\tst G \mapsto (\bM(G)[\f 1{[a^\flat]}])^{V = 1} \qx{is left exact.}
\]
 Thus, by snake lemma, for the injectivity and then also the surjectivity of the right vertical map in the diagram, 
we may replace $G$ by $G'$ to assume that $G$ is a truncated $p$-divisible group. By then $p$-adically filtering $G$ and again using the snake lemma, we may assume that $G$ is also killed by $p$. In this case, $\bM(G)$ is a projective $A^\flat$-module (see \S\ref{pp:prismatic-review}) and the right vertical map is injective because $a^\flat$ is a nonzerodivisor in $A^\flat$. 
For its surjectivity, fix an 
\[
\tst m \in (\bM(G)[\f 1{a^\flat}])^{V = 1}.
\]
By the $\Frob\i$-semilinearity of $V$, if $bm \in \bM(G)$ for a $b \in A^\flat$, then also $b^{1/p}m \in \bM(G)$. Such elements of $A^\flat[\f 1{a^\flat}]$ lie in $A^\flat$, so, since $\bM(G)$ is a direct summand of a finite free $A^\flat$-module, $m \in \bM(G)$. 

For the remaining surjectivity of 
\[
\tst \tst V - 1 \colon \bM(G)[\f 1{[a^\flat]}] \surjects \bM(G)[\f 1{[a^\flat]}],
\]
by \Cref{MG-exact-2}, the bijectivity of the right vertical arrow, 
and the snake lemma, the functor 
\[
\tst G \mapsto \bM(G)[\f 1{[a^\flat]}]/(V - 1)(\bM(G)[\f 1{[a^\flat]}]) \qx{is exact.}
\]
Thus, as before, we may first assume that $G$ is a truncated $p$-divisible group and then that it is also killed by $p$, so that $\bM(G)$ is a finite projective $A^\flat$-module. We then first claim that $V - 1$ is surjective after passing to the limit over all the nonzerodivisors   $a^\flat$.

\bcl \label{first-claim}
Letting $K_i^\flat$ be the fraction field of $A_i^\flat$, the map $V - 1$ is surjective on
\[
\tst \bM(G) \tensor_{A^\flat} (\prod_{i \in I} K_i^\flat) \overset{\ref{product-lemma}}{\cong}  \prod_{i \in I} (\bM(G_{A_i}) \tensor_{A_i^\flat} K_i^\flat).
\] 
\ecl

\bpf
We may assume that $I$ is a singleton $\{ i \}$. If $A_i$ is an $\bF_p$-algebra, then 
\[
\bM(G) \tensor_{A_i^\flat} K_i^\flat \cong \bM(G_{K_i})
\]
and $V - 1$ is surjective by \eqref{MG-fixed-pts}. 
Otherwise $\xi$ is a nonzerodivisor in $A^\flat$, so both $F$ and $V$ are bijective on 
\[
M \ce \bM(G_{A_i}) \tensor_{A_i^\flat} K_i^\flat.
\]
Thus, by \cite{Kat73}*{Proposition 4.1.1 and its proof}, $M^{F = 1}$ is a finite $\bF_p$-module and $M \cong M^{F = 1} \tensor_{\bF_p} K_i^\flat$ compatibly with the Frobenius. By choosing an $\bF_p$-basis for $M^{F = 1}$, the desired surjectivity of $V - 1$ on $M$ then amounts to the surjectivity of $\xi \cdot \Frob\i(*) - *$ on $K_i^\flat$, equivalently, of $\xi \cdot * - \Frob(*)$. The latter translates to the solvability in $K_i^\flat$ of any equation $x^p - \xi x + b = 0$  with $b \in K_i^\flat$, which follows from the algebraic closedness of $K_i^\flat$ (see \Cref{tilt-summary}).
\epf

Thanks to \Cref{first-claim}, 
it remains to show that if an $m \in \bM(G)[\f 1{a^\flat}]$ is of the form $V(x) - x$ for some $x \in \bM(G) \tensor_{A^\flat} (\prod_{i \in I} K_i^\flat)$, then $x \in \bM(G)[\f 1{a^\flat}]$. For this, it suffices to show that $x$ lies in each stalk of $\bM(G)$ at a variable point of $\Spec(A[\f 1{a^\flat}])$. 
By \Cref{lem:ultrafilters}, each local ring of $\Spec(A[\f 1{a^\flat}])$ is a valuation ring whose fraction field is a localization of $K$, so we are reduced to the following claim.

\bcl
For a perfect $\bF_p$-algebra $W$ that is a valuation ring with fraction field $L$, a finite free $W$-module $M$, and a $\Frob\i$-semilinear $V \colon M \ra M$, any $x \in M_L$ with $V(x) - x \in M$ lies in $M$.
\ecl

\emph{Proof. }
In the statement, the map $V$ on $M_L$ is defined by the same formula as in the beginning of the proof of \Cref{M(G)-local-analysis}. Suppose that $w \in W$ is such that $wx \in M$, so that $x = \f{m}{w}$ for some $m \in M$. Then, since $V(x) - x = \f{V(m)}{w^{1/p}} - x$ lies in $M$, we get that also $w^{1/p}x \in M$. However, since $M$ is finite free and $w$ is arbitrary subject to $wx \in M$, this means that $x \in M$.
\QEDD









\section{Flat cohomology with finite, locally free coefficients} \label{reduction}

Our argument for 
purity for flat cohomology uses several new properties of fppf cohomology with coefficients in commutative, finite, locally free group schemes. We establish these properties in this chapter by combining deformation theory discussed in \S\ref{sec:defthy} with crystalline Dieudonn\'{e} theory discussed in \S\ref{section:Dieudonne-positive-characteristic}. It is convenient to extend the statements to the setting of fppf cohomology of animated rings: even for usual rings this removes unnecessary assumptions and makes the proofs possible because our reductions involve derived $p$-adic completions and derived base changes that leave the realm of usual rings. The ultimate goal of these reductions is to end up with perfect $\bF_p$-algebras, 
which may then be treated by using the key formula \eqref{crystalline-cor} in positive characteristic. 


\csub[Deformation theory over animated rings] \label{sec:defthy}

A crucial tool in our reductions is deformation theory, carried out in the setting of simplicial rings. We will, however, not work with the latter in the strict sense of simplicial objects in the category of (commutative, as always) rings: instead, we will consider the $\infty$-category obtained by inverting the weak equivalences, that is, the maps that induce weak equivalences of the underlying simplicial sets, equivalently, the maps that induce quasi-isomorphisms of the underlying simplicial abelian groups considered as connective chain complexes via the Dold--Kan correspondence. The passage from the category of commutative rings to this $\infty$-category is an instance of a general procedure discussed in \S\ref{animation} that Dustin Clausen, inspired by Beilinson's \cite{Bei07},\footnote{In \cite{Bei07}, Beilinson lifts certain equalities in $K_0$ to actual homotopies in the $K$-theory space (that is, in the $K$-theory anima in the terminology we use), which he refers to as ``animations'' of that equality.} suggested to term ``animation.'' For instance, the animation of the category of sets---the $\infty$-category of ``animated sets'' or, briefly, of ``anima''---is simply the $\infty$-category of ``spaces'' in the sense of Lurie: it is the $\infty$-category obtained from the category of simplicial sets (or topological spaces) by inverting weak equivalences.


To put animation into context, we begin with the following general category-theoretic review.

\bpp[Free generation by $1$-sifted colimits] \label{sfp-fun}
For a category $\sC$ that has filtered colimits, we let $\sC^{\mathrm{fp}} \subset \sC$ be the full subcategory of those $X \in \sC$ that are \emph{of finite presentation} (also called \emph{compact}) in the sense that $\Hom_\sC(X, -)$ commutes with filtered colimits. Finite colimits in $\sC$ of objects of $\sC^{\mathrm{fp}}$ lie in $\sC^{\mathrm{fp}}$, and we have fully faithful embeddings\footnote{\label{cat-foot} To see that the functor 
\[
f\colon \Ind(\sC^{\mathrm{fp}}) \ra \sC
\]
supplied by the universal property of $\Ind(\sC^{\mathrm{fp}})$ is an embedding, we use the argument of \cite{HTT}*{Proposition 5.5.8.22} as follows. For a fixed $X \in \sC^{\mathrm{fp}}$, the full subcategory of the $Y \in \Ind(\sC^{\mathrm{fp}})$ with 
\[
\Hom_{\Ind(\sC^{\mathrm{fp}})}(X, Y) \isomto \Hom_\sC(X, f(Y))
\]
contains $\sC^{\mathrm{fp}}$ and is stable under filtered colimits, so it is all of $\Ind(\sC^{\mathrm{fp}})$. Thus, the full subcategory of the $X \in \Ind(\sC^{\mathrm{fp}})$ such that 
\[
\Hom_{\Ind(\sC^{\mathrm{fp}})}(X, Y) \isomto \Hom_\sC(f(X), f(Y)) \qxq{for all} Y \in \Ind(\sC^{\mathrm{fp}})
\]
contains $\sC^{\mathrm{fp}}$; since it is also stable under filtered colimits, it must be all of $\Ind(\sC^{\mathrm{fp}})$.} 
\[
\sC^{\mathrm{fp}} \hra \Ind(\sC^{\mathrm{fp}}) \hra \sC.
\]

For a category $\sC$ that has $1$-sifted\footnote{\label{foot-sifted}A small category $\sD$ is \emph{$1$-sifted}---usually simply called ``sifted'' in traditional category theory, but we want to reserve the term ``sifted'' for the $\infty$-categorical concept---if $\sD$-indexed colimits commute with finite products in the category $\Set$ of sets (see \cite{ARV10}*{Remark~1.1~(i)} for a concrete description; for context, we recall that $\sD$ is filtered if and only if $\sD$-indexed colimits commute with finite limits in $\Set$). 
For example, the category $\Delta^\op$ that indexes simplicial objects is $1$-sifted: $\Delta^\op$-indexed colimits, that is, \emph{geometric realizations}, are computed after restricting to $\Delta^\op_{\le 1}$, 
which is $1$-sifted \cite{ARV10}*{Example~1.2}; the $\Delta^\op_{\le 1}$-indexed colimits are \emph{reflexive coequalizers}. If the domain of a functor $F \colon \sC \ra \sC'$  has finite colimits, then $F$ commutes with $1$-sifted colimits if and only if it commutes with filtered colimits and reflexive coequalizers (see~\cite{ARV10}*{Theorem~2.1}).
} colimits, 
we let $\sC^{\mathrm{sfp}} \subset \sC$ be the full subcategory of those $X \in \sC$ that are \emph{strongly of finite presentation} (also called \emph{compact projective} when $\sC$ has all colimits) in the sense that $\Hom_{\sC}(X, -)$ commutes with $1$-sifted colimits. Finite coproducts in $\sC$ of objects of $\sC^\sfp$ lie in $\sC^\sfp$, and, letting $\mathrm{sInd}$ denote the $1$-sifted ind-category (the subcategory of $\Fun((\sC^\sfp)^\op,\Set)$ generated under $1$-sifted colimits by the Yoneda image of $\sC^\sfp$), we have fully faithful embeddings 
\[
\sC^\sfp \hra \mathrm{sInd}(\sC^\sfp) \hra \sC 
\]
(compare with \cref{cat-foot}). 
If $\sC$ is cocomplete and \emph{generated under colimits} by $\sC^\sfp$ in the sense that $\sC$ has no proper cocomplete subcategory containing $\sC^\sfp$, then\footnote{\label{sInd-coco}Indeed, $\mathrm{sInd}(\sC^\sfp)$ inherits cocompleteness from $\sC$: since a product of $1$-sifted diagrams is $1$-sifted, 
it inherits the existence of finite coproducts from $\sC^\sfp$, so, by taking filtered limits, it has arbitrary coproducts, and it remains to recall that any colimit is a reflexive coequalizer of coproducts.} 
\be \label{sInd-all}
\mathrm{sInd}(\sC^\sfp) \isomto \sC.
\ee
Consequently, $1$-sifted-colimit-preserving functors $F$ from such a $\sC$ correspond 
to functors from $\sC^{\sfp}$, and $F$ commutes with all colimits if and only if $F|_{\sC^{\sfp}}$ commutes with finite coproducts. By \cite{Mac98}*{Chapter V, Section 8, Corollary} and the following proposition, 
if $\sC^\sfp$ is small, then $\sC^{\mathrm{fp}}$ is also small and the category of functors $F \colon \sC^\op \ra \Set$ that bring colimits in $\sC$ to limits in $\Set$ is nothing else than the essential image of the Yoneda embedding of $\sC$; equivalently, $\sC$ is the category of functors
\[
(\sC^\sfp)^\op\to \Set
\]
that bring finite coproducts in $\sC^\sfp$ to products in $\Set$. 
\epp

\begin{proposition} 
Let $\sC$ be a cocomplete category generated under colimits by $\sC^\sfp$. The finitely presented objects of $\sC$ \up{that is, the objects of $\sC^{\mathrm{fp}}$} are precisely the coequalizers \up{equivalently, the reflexive coequalizers} of objects in $\sC^\sfp$ and 
$\sC$ is generated under filtered colimits by $\sC^{\mathrm{fp}}$. In~particular, 
\[
\Ind(\sC^{\mathrm{fp}})\isomto \sC.
\]
\end{proposition}

\begin{proof} 
The coequalizer $X$ of parallel arrows $Y\rightrightarrows Z$ agrees with the (reflexive) coequalizer of $Y\sqcup Z\rightrightarrows Z$, so the parenthetical claim follows. 
Moreover, the equalizer of $\Hom(Z,-)\rightrightarrows \Hom(Y,-)$ is $\Hom(X,-)$, so, since equalizers commute with filtered colimits, if $Y, Z \in \sC^\sfp$, then	 $X\in \sC^{\mathrm{fp}}$.

A colimit is a coequalizer of coproducts, so any $X\in \sC$ is a coequalizer of some 
\[
\tst \bigsqcup_{i\in I} Y_i\rightrightarrows \bigsqcup_{j\in J} Z_j \qxq{with} Y_i,Z_j\in \sC^\sfp.
\]
Since the $Y_i$ are finitely presented, for every finite subset $I'\subset I$ there is a finite subset $J'\subset J$ such that one has a subdiagram $\bigsqcup_{i\in I'} Y_i\rightrightarrows \bigsqcup_{j\in J'} Z_j$. Its coequalizer $X_{I',\,J'}$ is finitely presented by the above, so, by taking the filtered colimit over all such choices of $I', J'$, we express $X$ as the filtered colimit of the $X_{I',\, J'}$, so that $\sC$ is generated under filtered colimits by $\sC^{\mathrm{fp}}$.

It remains to see that every $X\in \sC^{\mathrm{fp}}$ is a coequalizer of objects of $\sC^\sfp$. The preceding arguments imply that $X$ is a retract of some coequalizer $X'$ of a diagram $Y\rightrightarrows Z$ in $\sC^\sfp$; let $f: X'\to X'$ be the corresponding idempotent endomorphism of $X'$. Then $X$ is the coequalizer of $X'\rightrightarrows X'$, where the two maps are the identity and $f$. Since $Z\in \sC^\sfp$, the map $f: X'\to X'$ can be lifted to a map $\wt{f}: Z\to Z$, and then $X$ is also the coequalizer of $Z\sqcup Y\rightrightarrows Z$ where the two maps are the given ones on $Y$ and the identity (resp.,~$\wt{f}$) on $Z$.
\end{proof}

\beg \label{sfp-egs}
The following cocomplete categories $\sC$ are generated under colimits by $\sC^\sfp$: 
\benuma
\m \label{sfp-egs-Set}
$\Set$ of sets: $\Set^{\mathrm{sfp}}$ consists of the finite sets;

\m
$\Gp$ of groups: $\Gp^{\mathrm{sfp}}$ consists of the free groups on finite sets;

\m 
$\Ab$ of abelian groups: $\Ab^{\mathrm{sfp}}$ consists of the finitely generated, free abelian groups;

\m
$\Ring$ of (unital, commutative) rings: $\Ring^{\mathrm{sfp}}$ consists of the retracts of finite type, polynomial $\bZ$-algebras, in other words, of $\bZ$-algebras $R$ that are quotients $\pi \colon \bZ[x_1, \dotsc, x_n] \surjects R$ such that there is a $\bZ$-algebra map $\iota \colon R \ra \bZ[x_1, \dotsc, x_n]$ with	 $\pi \circ \iota = \id_R$. 
\eenum
The claimed descriptions of the subcategories $\sC^\sfp$ follow from the case $\sC = \Set$ and \cite{HA}*{Corollary~4.7.3.18}, which in each of the cases characterizes $\sC^\sfp$ as the full subcategory consisting of the retracts of the ``finite free'' objects (one also uses the Nielsen--Schreier theorem, according to which a subgroup of a free group is free, so that any retract of a finitely generated, free group inherits these properties). \emph{Loc.~cit.}~applies because the forgetful functor from the respective category to sets commutes with $1$-sifted colimits, that is, with filtered colimits and reflexive coequalizers: the former is clear and for the latter we note that the set-theoretic equivalence relation generated by the parallel arrows of reflexive equalizers preserves the algebraic structures (thanks to the built in simultaneous splitting).
\eeg

\bpp[The animation of a category] \label{animation}
 For a cocomplete 
 category $\sC$ generated under colimits by $\sC^{\sfp}$ (so $\sC \cong \mathrm{sInd}(\sC^\sfp)$ by \eqref{sInd-all}), 
 the \emph{animation} of $\sC$ is the $\infty$-category $\mathrm{Ani}(\sC)$ freely  generated under sifted colimits by $\sC^{\sfp}$, that is, $\Ani(\sC)$ has sifted colimits\footnote{For siftedness in the $\infty$-categorical context, see \cite{HTT}*{Definition 5.5.8.1 and what follows}; prototypical examples are filtered colimits and geometric realizations (that is, $\Delta^\op$-indexed colimits), and in some sense all sifted colimits are generated from these.} and a functor
 \[
 \sC^\sfp \ra \Ani(\sC)
 \]
such that 
\[
\Fun_{\mathrm{sifted}}(\Ani(\sC), \sA) \isomto \Fun(\sC^\sfp, \sA)
\]
for any $\infty$-category $\sA$ that has sifted colimits, where $(-)_{\mathrm{sifted}}$ indicates the full subcategory of functors that commute with sifted colimits (equivalently, with filtered colimits and geometric realizations). This characterization determines $\Ani(\sC)$ uniquely, 
whereas \cite{HTT}*{Definition 5.5.8.8, Proposition 5.5.8.10 (4), Proposition 5.5.8.15} ensure its existence. By \cite{HTT}*{Theorem 5.5.1.1, Proposition 5.5.8.10 (1)}, the $\infty$-category $\Ani(\sC)$ is presentable, so, by \cite{HTT}*{Definition~5.5.0.1, Corollary 5.5.2.4}, it is complete and cocomplete. If the objects of $\sC$ are widgets, then those of $\mathrm{Ani}(\sC)$ are \emph{animated widgets}, except we abbreviate $\Ani(\Set)$ to $\Ani$ and the term `animated set' to \emph{anima} (plural: \emph{anima}). For a comparison between $\Ani(\sC)$ and constructions going back to Quillen's \cite{Qui67}, see \cite{HTT}*{Section 5.5.9, especially, Corollary~5.5.9.3}.

By \Cref{sfp-egs}~\ref{sfp-egs-Set} and \cite{HTT}*{Definition 5.5.8.8}, the $\infty$-category $\Ani$ of anima is obtained from the category of simplicial sets by inverting weak equivalences, and for a general $\sC$ as above, for which $\sC^\sfp$ is small, $\Ani(\sC)$ is the $\infty$-category of functors $(\sC^{\sfp})^\mathrm{op}\to \mathrm{Ani}$ that take finite coproducts in $\sC^{\sfp}$ to products in $\Ani$. 
By \cite[Lemma 5.5.8.14]{HTT}, any such functor can be lifted to a 
functor that is representable by a simplicial object of $\sC$ (even of $\Ind(\sC^{\sfp}) \subset \mathrm{sInd}(\sC^\sfp) \cong \sC$).
In fact, \cite[Corollary 5.5.9.3]{HTT} (with the final paragraph of \S\ref{sfp-fun} above) describes $\mathrm{Ani}(\sC)$ as the $\infty$-category obtained from the category of simplicial objects of $\sC$ by inverting weak equivalences with respect to a suitable model structure induced by the Quillen model structure on the category $\mathrm{sSet}$ of simplicial sets. 

By \cite{HTT}*{Remark 5.5.8.26, Proposition 5.5.6.18}, composition of a functor
\[
(\sC^\sfp)^\op \ra \Ani
\]
with the truncation $\tau_{\le n} \colon \Ani \ra \Ani$ induces a truncation functor 
\[
\tau_{\le n} \colon \Ani(\sC) \ra \Ani(\sC)
\]
that is left adjoint to the inclusion of the full subcategory of $n$-truncated objects of $\Ani(\sC)$ (in the sense of \cite[Definition~5.5.6.1]{HTT}). In particular, by the last aspect of \S\ref{sfp-fun}, there is a fully faithful inclusion $\sC\hookrightarrow \mathrm{Ani}(\sC)$ that identifies $\sC$ with the full subcategory of the $0$-truncated objects of $\mathrm{Ani}(\sC)$; the functor $\pi_0 \ce \tau_{\le 0}$ is left adjoint to the inclusion and is given by composition with the  connected component functor $\pi_0\colon \mathrm{Ani}\to \mathrm{Set}$. 

In particular, for a functor $F: \sC\to \sD$ between cocomplete categories as above, if $F$ preserves $1$-sifted colimits, then it induces a unique functor 
\[
\Ani(F)\colon \Ani(\sC)\to \Ani(\sD),
\]
the \emph{animation} of $F$, that preserves sifted colimits, that restricts to 
\[
F\colon \sC^\sfp\to \sD\subset \Ani(\sD)
\]
on $\sC^\sfp\subset \Ani(\sC)$, and that is such that $\pi_0\circ \Ani(F) = F \circ \pi_0$. In general this operation is not compatible with composition; this is akin to the formation of derived functors that only compose well under certain assumptions.
\epp

\begin{proposition}\label{prop:composeanimation} 
Let $F\colon \sC\to \sD$ and $G\colon \sD\to \sE$ be $1$-sifted-colimit-preserving functors between cocomplete categories generated under colimits by their strongly finitely presented objects.
\benum
\item \lab{CA-a}
There is a natural transformation from the composite $\Ani(G)\circ \Ani(F)$ to $\Ani(G\circ F)$.

\item \lab{CA-b}
If either $F(\sC^\sfp) \subset \Ind(\sD^\sfp)$ in $\sD$ or  $(\Ani(G))(F(\sC^\sfp)) \subset \sE$ in $\Ani(\sE)$, then the natural transformation 
\[
\qq \Ani(G)\circ \Ani(F)\to \Ani(G\circ F)
\]
of \ref{CA-a} is an equivalence.
\eenum
\end{proposition}


\begin{proof} 
Both $\Ani(G) \circ \Ani(F)$ and $\Ani(G \circ F)$ are  functors $\Ani(\sC)\to \Ani(\sE)$ that preserve sifted colimits, so it suffices to compare their restrictions to $\sC^\sfp$. Such restriction of the first functor is $X\mapsto \Ani(G)(F(X))$, while that of the second one is $X\mapsto G(F(X))$. However, $\pi_0 \circ \Ani(G) = G\circ \pi_0$ and $F(X)$ is $0$-truncated, so we have the desired natural transformation
\[
\Ani(G)(F(X))\to \pi_0 (\Ani(G)(F(X))) \cong G(\pi_0 (F(X))) \cong G(F(X)).
\]
For the second part, we need to see that this is an equivalence if $F(X)$ is a filtered colimit of objects of $\sD^\sfp$ or if $(\Ani(G))(F(X))$ is $0$-truncated. The latter is clear and for the former we note that the class of $Y \in \sD$ such that $\Ani(G)(Y)\isomto G(Y)$ contains $\sD^\sfp$ and is stable under filtered colimits.
\end{proof}

\beg \label{ani-egs}
The animations of the categories $\Gp$, $\Ab$, and $\Ring$ may be described as follows.
\benuma
\m 
The $\infty$-category $\mathrm{Ani}(\mathrm{Gp})$ of animated groups is obtained from the category of simplicial groups by inverting weak equivalences and, by a classical theorem (see \cite{HA}*{Theorem~5.2.6.10, Corollary~5.2.6.21}), is identified with the $\infty$-category of $\mathbb E_1$-groups (also known as associative groups) in $\mathrm{Ani}$.

\m 
The $\infty$-category $\mathrm{Ani}(\mathrm{Ab})$
of animated abelian groups is obtained from the category of simplicial abelian groups by inverting weak equivalences and, by the Dold--Kan correspondence (see \cite{HA}*{Theorem~1.2.3.7}), 
is identified with the connective part\footnote{As pointed out to us by Hesselholt, it is pleasing to recover $\cD^{\le 0}(\bZ)$ in this way because the cochain complex requirement $d^2 = 0$ becomes part of a solution to a universal problem instead of a construction. } $\mathcal{D}^{\leq 0}(\mathbb Z)\subset \mathcal D(\mathbb Z)$ of the derived $\infty$-category of $\mathbb Z$ 
(however, $\Ani(\mathrm{Ab})$ is \emph{not} equivalent to what might be called commutative groups in $\mathrm{Ani}$, namely, it is not equivalent to the $\infty$-category of  $\mathbb E_\infty$-groups in $\mathrm{Ani}$). 
The $\infty$-category $\Ani(\Ab)$ is also identified with the $\infty$-category of abelian group objects in $\Ani$.\footnote{Recall that an \emph{abelian group object} (or a \emph{$\mathbb Z$-module object}) in an $\infty$-category $\sC$ that has finite products is a contravariant functor from the category of finite free $\mathbb Z$-modules to $\sC$ that commutes with finite products.} 

\m 
The $\infty$-category $\mathrm{Ani}(\mathrm{Ring})$ 
of animated rings is obtained from the category of simplicial commutative rings by inverting weak equivalences. 
\eenum
Since the forgetful functors $\Ring \ra \Ab \ra \Set$ commute with $1$-sifted colimits (see \Cref{sfp-egs}), they induce functors $\Ani(\Ring) \ra \Ani(\Ab) \ra \Ani$. In this case, Proposition~\ref{prop:composeanimation} ensures that the functors compose well. Moreover, these functors admit left adjoints, given by the animations of the usual left adjoints; in particular, these forgetful functors commute with all limits.
\eeg

\bpp[Animated rings and modules] \label{ani-rings}
For an animated ring $A$, we write $a\in A$ for a map $\ast\xra{a} A$ of anima 
(equivalently, a map $\mathbb Z[X]\ra A$ of animated rings), call $a$ an element of $A$, and set 
\[
\tst A[\f 1a] \ce A \tensor^\bL_{\bZ[X]} \bZ[X, \f 1X] \qxq{and} A/^\bL a \ce A \tensor_{\bZ[X],\, X\mapsto 0}^\bL \bZ.
\] 
Up to equivalence, the datum of an $a\in A$ amounts to that of an element of $\pi_0(A)$. More generally, elements $a_1, \dotsc, a_n \in A$ correspond to a map 
\[
\bZ[X_1, \dotsc, X_n] \xra{X_i\, \mapsto\, a_i} A
\]
of animated rings, and we set 
\[
A/^\bL (a_1, \dotsc, a_n) \ce A \tensor_{\bZ[X_1, \dotsc, X_n],\, X_i\mapsto 0}^\bL \bZ,
\]
so that
\[
A/^\bL (a_1, \dotsc, a_n) \cong ((A/^\bL a_1)/^\bL \dotsc )/^\bL a_n.
\]
Thanks to  \Cref{ani-egs}, every animated ring $A$ has its associated graded ring of homotopy groups 
\[
\tst \pi_*(A) \ce \bigoplus_{n \ge 0} \pi_n(A);
\]
the $m$-truncation functor of \S\ref{animation} gives the universal map $A \ra \tau_{\le m}(A)$ to an animated ring with vanishing homotopy $\pi_{> m}(-)$. To work with animated algebras over a base ring $R$, one either starts with the category of $R$-algebras and animates it or considers animated rings equipped with a map from $R$---the two perspectives are equivalent (compare with \cite[Corollary~25.1.4.3]{SAG}).

For an animated ring $A$, the $\infty$-category $\Mod(A)$ of $A$-modules is, by definition, the $\infty$-category of modules over the underlying $\bE_1$-ring of $A$, compare with \cite[Notation~25.2.1.1]{SAG}. The $\infty$-category of \emph{animated $A$-modules} is its subcategory 
\[
\tst \Mod^{\le 0}(A) \subset \Mod(A)
\]
of connective objects. When $A$ is a usual ring, $\Mod(A)$ is nothing else but the derived $\infty$-category of $A$ (see \cite{HA}*{Theorem~7.1.2.13}) and $\Mod^{\le 0}(A)$ agrees with the animation of the category of $A$-modules (so there is no clash of terminology). For a general animated ring $A$, the $\infty$-category $\Mod^{\le 0}(A)$ is identified with the $\infty$-category of modules in animated abelian groups over $A$ (regarded as an $\bE_1$-algebra in animated abelian groups), which may reasonably be called animated $A$-modules.

Equivalently, one may define the $\infty$-category of animated rings $A$ equipped with an animated $A$-module $M$ by animating the category of rings equipped with modules, compare with \cite[Proposition~25.2.1.2]{SAG}. One can then define various forms of derived tensor products between animated rings or animated modules by animating the usual functors. In particular, for a diagram $B\leftarrow A\to C$ of animated rings, one may define the animated ring $B\tensor^\bL_A C$, by animating the usual functor on rings.
\epp

\bpp[The cotangent complex of an animated ring] \label{pp:cot-cx}
For an animated ring $A$ and an animated $A$-module $M$, we define an animated ring $A\oplus M$, the prototypical example of a ``square-zero extension,'' by animating the corresponding functor defined on usual rings equipped with modules, compare with \cite[Construction~25.3.1.1]{SAG}. The animated ring $A \oplus M$ comes with maps from and to $A$ and, as can be checked on underlying anima, the functor 
\[
(A,M)\mapsto A\oplus M
\]
commutes with 
limits. 

A \emph{derivation} of an animated ring $A$ with values in an animated $A$-module $M$ is a map $A\to A\oplus M$ of animated rings splitting the projection $A\oplus M\to A$. We follow \cite[Definition~25.3.1.4]{SAG} in writing $\mathrm{Der}(A,M)$ for the anima of derivations of $A$ with values in $M$. By \cite[Proposition~25.3.1.5]{SAG} or, rather, by the theorem on corepresentable functors \cite[Proposition~5.5.2.7]{HTT}, there is a universal derivation: for an animated ring $A$, the \emph{cotangent complex} $L_{A/\mathbb Z}$ is the universal animated $A$-module equipped with a derivation of $A$ with values in $L_{A/\mathbb Z}$, that is, such that postcomposition induces an equivalence of anima 
\[
\Hom_A(L_{A/\mathbb Z}, M)\cong \mathrm{Der}(A,M) \qx{for all animated $A$-modules $M$.}
\]
When $A$ is $0$-truncated, this $L_{A/\mathbb Z}$ agrees with the usual cotangent complex, see \cite{SAG}*{Example~25.3.1.8}.

More generally, one defines the cotangent complex of a map of animated rings $f\colon A'\to A$ by repeating the above definitions verbatim, defining $A'$-derivations of $A$ with values in $M$ as maps of animated $A'$-algebras $A\to A\oplus M$ splitting the projection. By \cite[Remark~25.3.2.4]{SAG}, this agrees with the definition of $L_{A/A'}$ 
as the cofiber of the map 
\[
L_{A'/\mathbb Z}\otimes^{\mathbb L}_{A'} A\to L_{A/\mathbb Z},
\]
so there is a transitivity triangle
\[
L_{A'/\mathbb Z}\otimes^{\mathbb L}_{A'} A\to L_{A/\mathbb Z}\to L_{A/A'}\to (L_{A'/\mathbb Z}\otimes^{\mathbb L}_{A'} A)[1],
\]
and for any morphism $A' \ra B'$ of animated rings with $B \ce A \tensor^\bL_{A'} B'$, we have
\be \label{cot-cx-bc}
L_{A/A'} \tensor^\bL_A B \cong L_{B/B'}. 
\ee
\epp

\bpp[Square-zero extensions of animated rings] \label{pp:sq-zero}
A \emph{square-zero extension} of animated rings is the datum of a map of animated rings $f\colon A'\to A$, an animated $A$-module $M$ (the \emph{ideal} of the square-zero extension), and a pullback square of animated rings
\be\label{sq-zero}
\ba\xymatrix{
A'\ar[r]^-f\ar[d] & A\ar[d]^-{i}\\
A\ar[r]^-{s} & A\oplus (M[1]),
}\ea\ee
where $i$ is the inclusion and $s$ is a derivation of $A$ with values in $M[1]$ (that is, $s$ is a map  of animated rings that splits the projection $A\oplus (M[1]) \to A$, see \S\ref{pp:cot-cx}). Equivalently, the $\infty$-category of square-zero extensions of animated rings is the $\infty$-category of pairs $(A,M)$ of an animated ring $A$ and an animated $A$-module $M$ equipped with a derivation $s$ of $A$ with values in $M[1]$.
Indeed, this defines $A'$ as the equalizer of 
\[
A\rightrightarrows A\oplus (M[1]),
\]
where the two maps are $s$ and the inclusion.

To define a commutative diagram as in the definition it suffices to define a map $L_{A/A'}\to M[1]$: indeed, by \S\ref{pp:cot-cx}, this gives a derivation, that is, a map between the $A'$-algebras $A$ and $A\oplus (M[1])$ splitting the projection (even as $A'$-algebras; one forgets that part of the information).
\epp

\begin{example} \label{sq-zero-eg}
Let us give several examples of square-zero extensions.
\benuma
\item Taking any pair $(A,M)$ and $s$ to be the inclusion, by the commutation of $(A,M)\mapsto A\oplus M$ with limits in $M$, we recover the trivial square-zero extension $A' \cong A\oplus M$. 

\item \label{sq-0-2}
Let $A'\surjects A$ be a square-zero extension of usual rings with 
\[
\qq M \ce \Ker(A'\surjects A),
\]
which is an $A$-module. To find a map $s\colon A\to A\oplus (M[1])$ that gives $A'\to A$ the structure of a square-zero extension of animated rings, we need to exhibit  a map $L_{A/A'}\to M[1]$, and for this it suffices to recall from \cite{SP}*{Lemmas~\href{https://stacks.math.columbia.edu/tag/08US}{08US} and \href{https://stacks.math.columbia.edu/tag/07BP}{07BP}} that $\tau_{\leq 1} (L_{A/A'}) \cong M[1]$.

\item \label{sq-0-3}
Let $A'$ be an $(m+1)$-truncated animated ring, set $A \ce \tau_{\leq m} (A')$, and consider the $\pi_0( A)$-module $M\ce \pi_{m+1}( A')$ as an animated $A$-module. There is a map $s\colon A\to A\oplus (M[m+2])$ that realizes $A'$ as a square-zero extension of $A$: to define the corresponding map $L_{A/A'} \ra M[m + 2]$, we recall from \cite[Proposition~25.3.6.1]{SAG} that $\tau_{\leq m+2} (L_{A/A'})\simeq M[m+2]$.
\end{enumerate}
\end{example}

We now apply these ideas to deformation theory, in particular, we derive the crucial  \Cref{cor:deftheorygroup}.

\bpp[Deformation-theoretic setup] \label{def-setup}
For a ring $R$, a commutative, flat, affine $R$-group scheme $G$ is automatically a $\mathbb Z$-module object in the $\infty$-category opposite to that of animated rings over $R$. 
It follows that for animated $R$-algebras $A$, the anima $G(A)$ of $A$-valued points has a functorial $\mathbb Z$-module structure, and thus becomes an animated abelian group. We are interested in the behavior of $G(A)$ under square-zero extensions, so we  consider the functor from the $\infty$-category of animated $A$-modules $M$ to that of animated abelian groups defined by
\[
M\mapsto T(M) \ce \mathrm{Fib}(G(A\oplus M)\to G(A)).
\]
Since $G$ is affine, the functor $A\mapsto G(A)$ commutes with 
 limits (as can be checked on underlying anima), so the functor $T$ commutes with 
 limits. 
 It 
 then  follows from \cite[Proposition 5.5.2.7]{HTT} that $T$ it is corepresentable by some $\mathbb Z$-module $L_{G_A/A}$ in animated $A$-modules.

Let $e \colon \Spec(R)\ra G$ be the unit section. If $G$ is of finite presentation, then $e^*(L_{G/R})$ has perfect amplitude $[-1, 0]$ (see \cite{Ill72}*{chapitre~VII, \'{e}quation~(3.1.1.3)}); if $G$ is even smooth, then $e^*(L_{G/R})$ even has perfect amplitude $[0, 0]$. In particular, in these cases the $\mathbb Z$-module structure on $e^*(L_{G/R})$ is the trivial one (obtained from the animated $R$-module structure by restriction of scalars): indeed, a priori $e^*(L_{G/R})$ is a module over the $\bE_\infty$-ring $\mathbb Z\otimes_{\mathbb S}^\bL R$ but, being $1$-truncated, it is a module over $\tau_{\leq 1}(\mathbb Z\otimes_{\mathbb S}^\bL R)\cong R$.
\epp

\begin{proposition} \label{L-underlie}
In \uS\uref{def-setup}, the animated $A$-module that underlies $L_{G_A/A}$ is $e^\ast(  L_{G/R})\otimes_R^{\mathbb L} A$. In particular, the formation of $L_{G_A/A}$ commutes with base change. 
\end{proposition}

\begin{proof} Let $G=\Spec(S)$. Then $T$, as a functor to anima, sends $M$ to the anima of $R$-algebra maps $S\to A\oplus M$ whose projection to $A$ is identified with the composite $S\to R\to A$ where the first map corresponds to the unit section of $G$, in other words, to that of $R$-algebra maps $S\to R\oplus M$ whose first component is the unit section. But this functor is also corepresented by $e^\ast(L_{G/R})\otimes_R^{\mathbb L} A$.
\end{proof}


\bthm\label{cor:deftheorygroup}
Let $R$ be a ring, let $G$ be a flat, affine, commutative $R$-group of finite presentation, and let $e$ be the unit section of $G$. For a square-zero extension of animated $R$-algebras $A'\surjects A$ with ideal $M$,
 there is the following functorial fiber sequence of animated abelian groups\uc
\be \label{eqn:def-thy-basic}
G(A')\to G(A)\to \Hom_R(e^*(L_{G/R}), M[1]).
\ee
\ethm

\begin{proof} 
Since the functor $A\mapsto G(A)$ commutes with limits, by applying $G(-)$ to the Cartesian square \eqref{sq-zero} that is part of the structure of a square-zero extension gives a Cartesian square
\[
\xymatrix{
G(A')\ar[r]\ar[d] & G(A)\ar[d]^{G(i)}\\
G(A)\ar[r] & G(A\oplus (M[1]))
}\]
of animated abelian groups. The map $G(i)$ is split by the projection, so its cofiber is identified with the fiber of $G(A\oplus (M[1]))\to G(A)$, which is $\Hom_A(L_{G_A/A},M[1])$ by the definition of $L_{G_A/A}$. The conclusion now follows from \Cref{L-underlie}. 
\end{proof}

\brem \label{rem-proj-amp}
When thinking of animated abelian groups as connective objects of $\cD(\bZ)$, the last term of \eqref{eqn:def-thy-basic} is the connective cover $\tau_{\geq 0}$ of the $R\Hom$. In practice, $G$ is of finite presentation, so $e^*(L_{G/R})$ has perfect Tor-amplitude in $[-1,0]$ (see \S\ref{def-setup}) and the truncation is not necessary. However, in the latter case the fibre sequence \eqref{eqn:def-thy-basic} in animated abelian groups may fail to be a fibre sequence in $\mathcal D(\mathbb Z)$ because the last map may not be surjective on $\pi_0$. On fppf cohomology, this issue disappears by \Cref{prop:defthy} below.
\erem


\csub[Flat cohomology of animated rings] \label{sec:fppf-coho-simplicial}

Flat cohomology in the setting of animated rings is at the heart of our approach to exhibiting new properties of flat cohomology of usual rings. We define the former in this section and record its basic features, for instance, hyperdescent and convergence of Postnikov towers in \Cref{prop:smooth-fppf-et-comp} and a key deformation-theoretic cohomology triangle in \Cref{prop:defthy}. We begin by discussing the basic properties of flatness in the setting of animated rings. 

\bpp[Flat and \'{e}tale maps of animated rings] \label{pp:flat-modules}
An animated module $M$ over an animated ring $A$ is \emph{flat} (resp.,~\emph{faithfully flat}) if $\pi_0(M)$ is a flat (resp.,~faithfully flat) $\pi_0(A)$-module and 
\[
\pi_i(A) \otimes_{\pi_0(A)}\pi_0(M) \isomto \pi_i(M) \qxq{for all $i$}
\]
or, more succinctly, 
\[
\pi_*(A) \otimes_{\pi_0(A)}\pi_0(M) \isomto \pi_*(M),
\]
so that the graded $\pi_*(A)$-module $\pi_*(M)$ is flat (resp.,~faithfully flat). If $M$ is flat, then for any animated $A$-module $M'$ we have
\[
\Tor^{\pi_*(A)}_i(\pi_*(M'), \pi_*(M)) \cong \begin{cases} \pi_*(M') \tensor_{\pi_0(A)} \pi_0(M)  &\x{for $i = 0$,} \\ 0 &\x{for $i > 0$,} \end{cases}
\]
so the spectral sequence \cite{Qui67}*{Chapter~II, Section~6, Theorem~6 (b)} gives 
\[
\pi_*(M' \tensor_A^\bL M) \cong \pi_i(M') \tensor_{\pi_0(A)} \pi_0(M).
\]
In particular, (resp.,~faithful) flatness is stable under base change along a map $A \ra A'$ of animated~rings.
\epp

In the animated setting, flatness is insensitive to infinitesimal thickenings as follows.

\blem \label{lem:reduced-flatness}\lab{support-lem}
Let $A \ra A'$ be a map of animated rings that induces a surjective map $\pi_0(A) \surjects \pi_0(A')$ with nilpotent kernel. An animated $A$-module $M$ is \up{resp.,~faithfully} flat if and only if so is the animated $A'$-module $A' \tensor^\bL_A M$. In particular, $M \simeq 0$ if and only if $A' \tensor^\bL_A M  \simeq 0$.
\elem

\bpf
A flat $M$ vanishes if and only if $\pi_0(M)$ vanishes, which may be tested modulo any nilpotent ideal of $\pi_0(A)$, so the `in particular' follows from the main assertion. Moreover, in the latter the `only if' is clear and for the `if' we may focus on flatness because the support of $\pi_0(M)$ is insensitive to base change to $\pi_0(A')$. For the flatness, we first consider the special case when $A' = \pi_0(A)$, that is, we first claim that $M$ is $A$-flat if and only if $\pi_0(A) \tensor^\bL_A M$ is $\pi_0(A)$-flat.

For this, since base change to $\tau_{\le m}(A)$ does not affect the $\pi_i(M)$ with $i \le m$, we lose no generality by assuming that $A$ is $m$-truncated for some $m > 0$ and, by induction, need to show that $M$ is $A$-flat if $\tau_{\le m - 1}(A) \tensor^\bL_A M$ is $\tau_{\le m -1}(A)$-flat. However, the latter assumption and \S\ref{pp:flat-modules} give
\[\ba
\pi_m(A)[m] \tensor^\bL_A M &\cong \pi_m(A)[m] \tensor^\bL_{\tau_{\le m - 1}(A)} (\tau_{\le m - 1}(A)  \tensor^\bL_A M) \\
&\cong (\pi_m(A) \tensor_{\pi_0(A)} \pi_0(M))[m],
\ea\]
and the exact triangle 
\[
(\pi_m(A)[m]) \tensor^\bL_A M \ra M \ra \tau_{\le m - 1}(A)  \tensor^\bL_A M
\]
 then shows that $M$ is $A$-flat.

The settled case when $A'$ is $0$-truncated allows us to replace $A$ and $A'$ by $\pi_0(A)$ and $\pi_0(A')$, and hence assume that $A$ and $A'$ are $0$-truncated. Moreover, induction on the order of nilpotence of the ideal $I \ce \Ker(A \surjects A')$ allows us to assume that $I^2 = 0$. In this case, \S\ref{pp:flat-modules} gives
\[
I \tensor^\bL_A M \cong I \tensor^\bL_{A'} (A' \tensor^\bL_A M) \cong I \tensor_A M,
\]
so the exact triangle $I \tensor^\bL_A M \ra M \ra A' \tensor^\bL_A M$ shows that $M$ is $0$-truncated. Thus, since $A' \tensor^\bL_A M$ is $0$-truncated, we have $\Tor_1^A(A', M) \cong 0$, so that $M$ is $A$-flat by the  flatness criterion \cite{SP}*{Lemma~\href{https://stacks.math.columbia.edu/tag/051C}{051C}}.
\epf

\bpp[Grothendieck topologies on the $\infty$-category of animated rings] \label{pp:topologies}
A map $f\colon A\to A'$ of animated rings is \emph{flat} (resp.,~\emph{faithfully flat}) if $A'$ is flat (resp.,~faithfully flat) as an animated $A$-module (see \S\ref{pp:flat-modules}), concretely, if 
\[
\pi_0(f)\colon \pi_0(A) \to \pi_0(A')
\]
has the same property and  
\[
\pi_i(A) \otimes_{\pi_0(A)}\pi_0(A') \isomto \pi_i(A') \q \x{for all $i$}.
\]
A map $f$ of animated rings is \emph{\'{e}tale} if it is flat and $\pi_0(f)$ is \'{e}tale. A flat map $f$ is \emph{of finite presentation} (resp.,~\emph{finite}) if so is $\pi_0(f)$. All of these properties are inherited by the map $\pi_*(f)$ of graded rings. Moreover, by \S\ref{pp:flat-modules}, they are stable under composition and base change.

A map $f\colon A\to A'$ of animated rings is an \emph{fpqc cover} (resp.,~\emph{fppf cover}; resp.,~\emph{\'{e}tale cover}) if it is faithfully flat (resp.,~faithfully flat and of finite presentation; resp.,~faithfully flat and \'{e}tale). The stability properties above imply that such are covering maps for a Grothendieck topology on the $\infty$-category of animated rings (see \cite{HTT}*{Definition 6.2.2.1, Remark 6.2.2.3}). 
Of course, if $A$ is $0$-truncated, then so is $A'$, to the effect that one does not obtain new covers of $0$-truncated rings.
\epp

The \'{e}tale site of an animated ring is insensitive to derived structure as follows. 



\bprop \label{rem:no-new-etale}
For an animated ring $A$, the $\pi_0(-)$ \up{or base change} functor from \'etale \up{resp.,~finite \'{e}tale} $A$-algebras to \'etale \up{resp.,~finite \'{e}tale} $\pi_0(A)$-algebras is an equivalence of $\infty$-categories. 
\eprop

\bpf 
The two functors agree because $A'\otimes_A^{\mathbb L} \pi_0(A) \cong \pi_0(A')$ for any $A$-\'{e}tale (or even $A$-flat) $A'$ (see \S\ref{pp:topologies}). In particular, by \eqref{cot-cx-bc} and \Cref{support-lem}, we have $L_{A'/A}\cong0$ for any $A$-\'{e}tale $A'$. 

To prove the full faithfulness, it suffices to argue that for any $A$-\'{e}tale $A'$ and any animated $A$-algebra $B$, the following map is an equivalence of anima:
\[
\Hom_A(A', B)\to \Hom_A(A',\pi_0(B))
\]
Since $B \isomto R\lim_n(\tau_{\leq n} (B))$, by induction it suffices to show that 
\[
\Hom_A(A',\tau_{\leq n} (B))\to \Hom_A(A',\tau_{\leq n-1} (B)) \qx{is an equivalence of anima.}
\]
Since $\tau_{\leq n} (B)\to \tau_{\leq n-1} (B)$ admits the structure of a square-zero extension (see \Cref{sq-zero-eg}~\ref{sq-0-3}) and $\Hom_A(A',-)$ commutes with 
limits, it then suffices to argue that for any trivial square-zero extension $C\oplus M$ of an animated $A$-algebra~$C$,
\[
\Hom_A(A',C)\to \Hom_A(A', C\oplus M) \qx{is an equivalence of anima.}
\]
But this map has an evident section, whose fibers are given by maps $L_{A'/A}\to M$ by the definition of the cotangent complex. Since $L_{A'/A}\cong 0$, the claim follows.

For the remaining essential surjectivity, it is enough to handle the \'etale case (finiteness can be checked on $\pi_0$). Ideally, one should prove the result by deformation theory, using the vanishing of the cotangent complex, but 
we give an ad hoc argument. Namely, Zariski localizations can be lifted (for $f\in A$, one can form $A[\frac 1f]$ by base change from the universal case $\mathbb Z[f]\to \mathbb Z[f, \frac 1f]$), and \'etale algebras can be constructed Zariski locally. However, Zariski locally an \'etale map is standard \'etale, whose explicit description gives a lift to $A$ by lifting the defining elements from $\pi_0(A)$ to $A$.
\epf

\bpp[Cohomology over an animated ring] \label{pp:simplicial-coho-def}
For an animated ring $A$, we define the $\infty$-topoi of \'etale, fppf, or fpqc sheaves  over $A$ (the latter for an implicit sufficiently large cardinal bound $\kappa$ as in \S\ref{conv}) by considering the corresponding $\infty$-category of presheaves, that is, of functors from animated $A$-algebras \'etale/fppf/flat over $A$ to anima, and taking the full subcategory of those functors that satisfy the sheaf condition for the respective notion of covers. The inclusion into all presheaves has a left adjoint, the \emph{sheafification}. The same applies to (pre)sheaves with values in any $\infty$-category, so for a presheaf $F$ with values in $\mathcal D(\mathbb Z)$ we let 
\[
R\Gamma_\et(A, F) \in \mathcal D(\mathbb Z) \qxq{and} R\Gamma_\fppf(A,F)\in \mathcal D(\mathbb Z)
\]
denote the values at $A$ of the \'{e}tale and the fppf sheafifications of $F$ and write $H^i_\et(A, F)$ and $H^i_\fppf(A, F)$ for the resulting cohomology groups. Since a $0$-truncated $A$ does not attain new \'{e}tale or flat covers in the animated setting (see \S\ref{pp:topologies}), for $0$-truncated $A$ this definition reproduces the classical \'{e}tale and fppf cohomology, respectively. In this article, we will get by with cohomology in the affine setting, but for any open $U \subset \Spec(\pi_0(A))$ we also set
\[
\tst R\Gamma_\fppf(U_A, F) \ce R\lim_{A'} (R\Gamma_\fppf(A', F))
\]
 where $A'$ ranges over those animated $A$-algebras fppf over $A$ for which the map
 \[
 \Spec(\pi_0(A')) \ra \Spec(\pi_0(A))
 \]
 factors over $U$; the subscript in $U_A$ reminds us that we are not merely forming the usual flat cohomology of the scheme $U$. We often abbreviate $R\Gamma_\fppf$ and $H^i_\fppf$ to $R\Gamma$ and $H^i$, respectively.
 
 

It is useful to keep in mind that in this setting cohomology need not vanish in negative degrees: for instance, if $A$ is an animated algebra over a commutative ring $R$ and $G$ is a commutative, affine $R$-group,\footnote{Throughout this article we limit ourselves to group schemes defined over a base classical ring $R$ because this suffices for our applications and is expedient. It would be more natural to allow groups $G$ defined over $A$ itself, but it seems unclear how to correctly define commutative finite flat group schemes over animated rings (in particular, so as to admit Cartier duals and B\'{e}gueri resolutions \eqref{Begueri-0}). We expect that with the correct definition, all such commutative finite flat group schemes may arise via base change from classical rings.} then, by fpqc descent (see, for instance, \cite{SAG}*{Remark~D.6.3.6}), 
\be \label{coho-neg-deg}
G(A) \isomto \tau^{\le 0}(R\Gamma_\fppf(A, G)).
\ee
We recall that a sheaf $F$ is a hypersheaf if it satisfies the sheaf condition with respect to hypercovers. This is automatic when $F$ is $n$-truncated for some $n$; for example, if $F$ is a sheaf of coconnective complexes (as in the usual setting of cohomology). Another important example is that quasi-coherent sheaves are hypersheaves in the \'etale, fppf, and fpqc sites: the quasi-coherent sheaf defined by some animated module $M$ is the limit of the sheaves defined by its truncations $\tau_{\leq n}(M)$, all of which are truncated, so the claim follows as limits of hypersheaves are hypersheaves (see also \cite{SAG}*{Corollary~D.6.3.4}).
\epp

Despite possible negative degree cohomology, we have the following hyperdescent and Postnikov convergence result; it also extends Grothendieck's fppf-\'{e}tale comparison \cite{Gro68c}*{th\'{e}or\`{e}me~11.7} to animated~rings.

\bthm \label{prop:smooth-fppf-et-comp} 
Let $R$ be a ring and let $G$ be an affine, commutative $R$-group that is either smooth or finite locally free. The functor $A\mapsto R\Gamma_\fppf(A,G)$ satisfies fppf hyperdescent on animated $R$-algebras $A$ and
\be \label{eqn:fppf-hyperdescent}
\tst R\Gamma_\fppf(A,G)\isomto R\lim_n (R\Gamma_\fppf(\tau_{\leq n} (A),G)).
\ee
If $G$ is smooth, then even the functor $A\mapsto R\Gamma_\et(A,G)$ satisfies fppf hyperdescent and, in particular,
\be \label{et-fppf-agree}
R\Gamma_\et(A,G) \isomto R\Gamma_\fppf(A,G).
\ee
\ethm

\begin{proof} 
For a finite, locally free $G$, the B\'{e}gueri resolution
\be \label{Begueri}
0\to G\to \mathrm{Res}_{G^\ast/R} (\mathbb G_m) \to Q\to 0
\ee
is exact on the fppf site of any animated $R$-algebra $A$ because the map $\mathrm{Res}_{G^\ast/R} (\mathbb G_m) \to Q$ is faithfully flat and finitely presented. Thus, it reduces us to the case when $G$ is smooth. Moreover, since $G$ is affine, for any $A$ we have the Postnikov tower equivalence
\be \label{eqn:G-affine-Postnikov}
\tst G(A)\isomto R\lim_n (G(\tau_{\leq n} (A))).
\ee
By induction on $n$ and \Cref{cor:deftheorygroup}, the fiber $G(\tau_{\leq n}(A))^0$ of the map
\[
G(\tau_{\leq n}(A))\to G(\pi_0(A))
\]
satisfies fppf hyperdescent (see \S\ref{pp:simplicial-coho-def}). 
 By \eqref{eqn:G-affine-Postnikov}, we have the identification
\be \label{eqn:G0-Postnikov}
\tst G(A)^0 \cong R\lim_n \p{G(\tau_{\leq n}(A))^0}, 
\ee
where $G(A)^0$ is the fiber of $G(A)\to G(\pi_0(A))$. Thus, the functor $A \mapsto G(A)^0$ satisfies fppf hyperdescent. By then sheafifying the fiber sequence
\[
G(A)^0 \to G(A)\to G(\pi_0(A))
\]
for the \'etale topology and using \Cref{rem:no-new-etale}, we obtain a fiber sequence 
\be \label{eqn:after-fppf-sheafification}
G(A)^0 \to R\Gamma_\et(A,G)\to R\Gamma_\et(\pi_0(A),G)
\ee
(see \S\ref{pp:simplicial-coho-def}). The base change of an fppf hypercover of $A$ along the map $A \ra \pi_0(A)$ is an fppf hypercover of $\pi_0(A)$ obtained by forming $\pi_0(-)$ levelwise, so the functor $A\mapsto R\Gamma_\et(\pi_0(A), G)$ satisfies fppf hyperdescent by Grothendieck's \cite{Gro68c}*{th\'{e}or\`{e}me~11.7}. Thus, the outer terms of \eqref{eqn:after-fppf-sheafification} satisfy fppf hyperdescent in $A$, and hence so does the middle term $A \mapsto R\Gamma_\fppf(A, G)$. By combining \eqref{eqn:G0-Postnikov} with \eqref{eqn:after-fppf-sheafification} applied with $\tau_{\le n}(A)$ in place of $A$, we obtain \eqref{eqn:fppf-hyperdescent}.
\end{proof}

We will use the following mild strengthening of the Postnikov completeness of $A\mapsto R\Gamma(A, G)$. 

\bcor \label{cor:fppf-hyperdescent}
Let $R$ be a ring, let $G$ be a commutative affine $R$-group that is either smooth or finite locally free, and let $A$ be an animated $R$-algebra. For a tower of maps $\dotsc \ra A_{n + 1} \ra A_n \ra \dotsc \ra A_0$ of animated $A$-algebras such that $\tau_{\le n}(A) \isomto \tau_{\le n}(A_n)$ for all $n$, we have
\be \label{eqn:similar-to-Postnikov-tower}
\tst R\Gamma_\fppf(A,G)\isomto R\lim_n (R\Gamma_\fppf(A_n,G)).
\ee
\ecor

\bpf
We consider the inverse limit diagram $\{ R\Gamma_\fppf(\tau_{\le m}(A_n), G) \}_{m,\, n \ge 0}$. If one first forms $R\lim$ in $m$ and afterwards in $n$, then, by \eqref{eqn:fppf-hyperdescent} with $A_n$ in place of $A$, one obtains the right side of \eqref{eqn:similar-to-Postnikov-tower}. If, on the other hand, one first forms $R\lim$ in $n$ and afterwards in $m$, then, by the assumption on the $A_n$, one obtains $R\lim_m(R\Gamma(\tau_{\le m}(A), G))$, which, by \eqref{eqn:fppf-hyperdescent} again, is $R\Gamma(A, G)$. 
\epf

The following sheafification of the deformation-theoretic \Cref{cor:deftheorygroup} plays a central role below.

\bthm \label{prop:defthy}
Let $R$ be a ring, let $G$ be a commutative, affine $R$-group that is either smooth or finite locally free, and let $e$ be the unit section of $G$. For a square-zero extension of animated $R$-algebras $A'\surjects A$ with ideal $M$, there is the following functorial fiber sequence in $\mathcal D(\mathbb Z)$\ucolon
\[
R\Gamma_\fppf(A',G)\to R\Gamma_\fppf(A,G)\to R\Hom_R(e^*(L_{G/R}), M[1]).
\]
\ethm

\begin{proof} 
For smooth $G$, this follows from \Cref{prop:smooth-fppf-et-comp} and \Cref{cor:deftheorygroup} by \'etale sheafification (the right-most term is $1$-connective, so the fibre sequence in animated abelian groups gives a fibre sequence in $\mathcal D(\mathbb Z)$). For finite, locally free $G$, as in the proof of \Cref{prop:smooth-fppf-et-comp}, the B\'{e}gueri resolution \eqref{Begueri}
reduces us to the smooth case.
\end{proof}

The following description of the positive degree flat cohomology of animated rings with suitable coefficients complements \eqref{coho-neg-deg}, which gave a description of the negative degree cohomology.


\bcor \label{cor:nil-nil-deform-cohomology}
Let $R$ be a ring and let $G$ be a smooth \up{resp.,~finite, locally free}, affine, commutative $R$-group. For an animated $R$-algebra $A$, the map
\[
\tst H^i_\fppf(A,G) \ra H^i_\fppf(\pi_0(A),G)  
\]
is surjective for all $i$ and is bijective for $i\ge 0$ \up{resp.,~for $i \ge 1$}.
\ecor

\bpf
The finite locally free case reduces to the smooth one via the B\'{e}gueri sequence \eqref{Begueri}. In the smooth case, $e^*(L_{G/R})$ in \Cref{prop:defthy} is a projective module concentrated in degree $0$, so 
\[
H^i_\fppf(\tau_{\le n}(A), G) \ra H^i_\fppf(\tau_{\le n - 1}(A), G) \qxq{is} \begin{cases} \xq{surjective for} i \ge -1, \\ \xq{bijective for} \, \, \, i \ge 0. \end{cases}
\]
The Postnikov convergence \eqref{eqn:fppf-hyperdescent} then gives our claim.
\epf

Deformation theory has the following consequence for the insensitivity to nonreduced structure. 

\bcor \label{no-nil-feelings}
For a ring $R$, an ideal $I \subset R$ whose elements are nilpotent, and a smooth \resp{finite, locally free}, affine, commutative $R$-group $G$, the map
\[
H^i_\fppf(R, G) \ra H^i_\fppf(R/I, G) \ \ \x{is}\ \ \begin{cases} \xq{surjective for} i\ge 0 \qx{\resp{for $\q i \ge 1$}}, \\   \xq{bijective for} \ \  i\ge 1 \qx{\resp{for $\q i \ge 2$}.} \end{cases}
\]
\ecor

For finite, locally free $G$, we will extend this result to general Henselian pairs in \Cref{inv-Hens-pair}. For smooth $G$, the same statement does not hold for arbitrary Henselian pairs, but see \cite{BC19}*{Proposition~2.1.4, Theorem 2.1.6, Remarks 2.1.7 and 2.1.8} for positive results in low cohomological degrees and counterexamples to general statements, as well as \cite{SGA3IIInew}*{expos\'{e}~XXIV, lemme~8.1.8, remarque~8.1.9} for positive results for \v{C}ech cohomology.

\bpf
The ring $R$ is a filtered direct union of its finite type $\bZ$-subalgebras $R'$ and $R/I$ is a similar direct union of the $R'/(R' \cap I)$. Thus, limit formalism reduces us to the case when $R$ is Noetherian, so that $I$ is nilpotent and, arguing by induction, even square-zero. In this case, \Cref{prop:defthy} (with \Cref{sq-zero-eg}~\ref{sq-0-2}) supplies the long exact sequence
\[
\dotsc \ra H^{i}(R\Hom_R(e^*(L_{G/R}), I)) \ra H^i(R, G) \ra H^i(R/I, G) \ra
\dotsc
\]
To get a desired vanishing of the flanking terms it now remains to recall from \S\ref{def-setup} that in the smooth \resp{finite, locally free} case, $e^*(L_{G/R})$ has perfect amplitude $[0, 0]$ \resp{$[-1, 0]$}.
\epf


\csub[The $p$-adic continuity formula for flat cohomology] \label{section:continuity-formula}

The ultimate driving force of our analysis of new properties of fppf cohomology is the positive characteristic case of the key formula that we established in \Cref{KT-input}. To deduce mixed characteristic phenomena from this positive characteristic statement, we rely on the $p$-adic continuity formula that we exhibit in \Cref{thm:padiclimit} below (see also \Cref{adic-continuity} for a subsequent extension to adically complete rings). This formula has the flavor of a flat cohomology counterpart of Gabber's affine analog of proper base change for \'{e}tale cohomology \cite{Gab94a}*{Theorem~1} and is new already for $p$-adically complete, $p$-torsion free rings. The proof of this case does not require animated inputs, but for the sake of maximal applicability we directly treat the general case. The technique is to 
 reduce to rings that have no nonsplit fppf covers via the following lemma.

\blem \label{lem:large-cover}
For each ring $R$, there is an ind-fppf $R$-algebra $\wt{R}$ that has no nonsplit fppf covers. For each ideal $I \subset R$ contained in every maximal ideal, there is an $I\wt{R}$-Henselian such $\wt{R}$. 
\elem

\bpf
The first claim is the $I = 0$ case of the second. Fix a set $\sS$ of representatives for isomorphism classes of faithfully flat, finitely presented $R$-algebras $R'$, and consider the ind-fppf $R$-algebra
\[
\tst R_1 \ce \bigotimes_{R' \in \sS} R' \qxq{and its $I$-Henselization} R_1^h.
\]
By iterating this construction (with $R_1^h$ in place of $R$ and so on), we obtain a tower of ind-fppf $R$-algebras $R_1 \ra R_1^h \ra  R_2 \ra R_2^h \ra \dotsc$ that are faithfully flat (see \cite{SP}*{Lemma~\href{https://stacks.math.columbia.edu/tag/00HP}{00HP}}
). By construction,
\[
\tst \wt{R} \ce \varinjlim_{n > 0} R_n \cong \varinjlim_{n > 0} R^h_n
\] 
is $I\wt{R}$-Henselian. By a limit argument, every fppf cover $\wt{f} \colon \wt{R} \ra \wt{S}$ descends to an fppf cover $f\colon R_n \ra S$ for some $n$. There is an $R_n$-morphism $S\ra R_{n + 1} \ra \wt{R}$, so the cover $\wt{f}$ has a section.
\epf

\brem
\Cref{lem:large-cover} continues to hold with an analogous proof if in its formulation ind-fppf/fppf is replaced by ind-\'{e}tale/\'{e}tale, or by ind-smooth/smooth, or by ind-syntomic/syntomic.
\erem




Another input to the $p$-adic continuity formula is the following lemma of Beauville--Laszlo type.

\blem \label{lem:derived-BL-key}
For a map $A \ra A'$ of animated rings and an $a\in A$ such that $A/^\bL a \isomto A'/^\bL a$, we have
\[
\tst 
A \isomto A' \times_{A'[\f{1}{a}]} A[\f{1}{a}] \q\qx{\up{even in the derived $\infty$-category $\cD(\bZ)$}.}
\]
\elem

\bpf
Consider the fiber $M$ of the morphism 
\[
\tst A \ra \mathrm{Fib}\p{A' \oplus A[\f{1}{a}] \ra A'[\f{1}{a}]},
\]
where the second arrow is the difference map. We have $M[\f{1}{a}] \cong 0$ and $M/^\bL a \cong 0$; the former shows that the homotopy groups of $M$ are $a$-power torsion, and the latter that multiplication by $a$ is an isomorphism on them. In conclusion, $M \cong 0$, so also $A \isomto A'\times_{A'[\f{1}{a}]} A[\f{1}{a}]$. 
\epf

The final input is the following lemma that will be useful for us on multiple occasions.

\blem \label{lem:sq-0-trick}
For an animated ring $A$, an element $a \in A$, and the derived $a$-adic completion $\wh{A} \ce R\lim_{n > 0}(A/^\bL a^n)$, the ring $\pi_0(\wh{A})$ is a square-zero extension of the $a$-adic completion of $\pi_0(A)$. 
\elem

\bpf
To analyze $\pi_0(\wh{A})$ we use the exact sequence 
\[
\tst 0 \ra \varprojlim_{n > 0}^1 (\pi_1(A/^\bL a^n)) \ra \pi_0(\wh{A}) \ra \varprojlim_{n > 0} (\pi_0(A)/(a^n)) \ra 0.
\]
It suffices to note that, since $\varprojlim_{n > 0}^2 \cong 0$, the ideal 
\[
\tst \varprojlim_{n > 0}^1 (\pi_1(A/^\bL a^n)) \subset \pi_0(\wh{A})
\]
is square-zero because the limit filtration on 
\[
\tst \pi_0(\wh{A}) \cong \pi_0(R\lim_{n > 0}(A/^\bL a^n))
\]
 is multiplicative.
\epf

\begin{theorem}\label{thm:padiclimit} 
For a prime $p$, a ring $R$, a commutative, finite, locally free $R$-group $G$ of $p$-power order, and an animated $R$-algebra $A$ such that $\pi_0(A)$ is $p$-Henselian, 
\be \label{eqn:continuity-formula}
\tst R\Gamma(A,G)\isomto R\lim_n (R\Gamma(A/^\bL p^n,G)).
\ee
\end{theorem}

\begin{proof}
For an initial reduction to $0$-truncated $A$, we begin by considering the system
\[
\{ R\Gamma(\tau_{\le m}(A)/^\bL p^n, G) \}_{m,\, n \ge 0}.
\]
The map $A/^\bL p^n \ra \tau_{\le m}(A)/^\bL p^n$ induces an isomorphism on truncations $\tau_{\le m}$, so \eqref{eqn:similar-to-Postnikov-tower} shows that forming $R\lim_m$ followed by $R\lim_n$ gives the right side of \eqref{eqn:continuity-formula}. On the other hand, by \eqref{eqn:fppf-hyperdescent}, forming $R\lim_n$ followed by $R\lim_m$ gives the left side granted that we know \eqref{eqn:continuity-formula} for truncated $A$. To reduce from the latter to $0$-truncated $A$ by induction on the truncation level, we use \Cref{prop:defthy}: for any square-zero extension $A'\to A$ with kernel $M$ we have
\be \label{def-triangle-eq}
R\Gamma(A', G) \ra R\Gamma(A, G) \ra R\Hom_R(e^*(L_{G/R}), M[1]),
\ee
so it suffices to argue that $R\Hom_R(e^*(L_{G/R}), M[1])$ is insensitive to replacing $M$ by its derived $p$-adic completion $\wh{M}$. For this, it suffices to observe that the cofiber of $M\to \wh{M}$ is a $\mathbb Z[\f 1p]$-module, whereas $e^*(L_{G/R})$ vanishes after inverting $p$. In conclusion, without losing generality, $A$ is $0$-truncated.

\bcl \label{claim:descent}
On $A$-algebras, the functor $A' \mapsto R\Gamma((A')^h_p, G)$ (where $(-)^h_p$ is the $p$-Henselization) 
satisfies hyperdescent in the topology whose covers are the filtered direct limits of flat, finitely presented maps that are faithfully flat modulo $p$.
\ecl

\bpf
Any filtered direct limit of flat, finitely presented maps that is faithfully flat modulo $p$ is a $p$-complete arc cover (see \S\ref{pp:arc-topology}~\ref{arc-pi-fflat}). Thus, by \Cref{BM-hens-comp}, the functor $A' \mapsto R\Gamma((A')^h_p[\f 1p], G)$ satisfies the desired hyperdescent. The cohomology with supports sequence then reduces us to arguing the same for the functor 
\be \label{elt-excision}
A' \mapsto R\Gamma_{\{ p = 0\}}((A')^h_p, G), \ \ \x{which is simply}\ \ A' \mapsto R\Gamma_{\{ p = 0\}}(A', G)
\ee
thanks to excision (for this last identification, see \cite{Mil80}*{Chapter~III, Proposition~1.27}\footnote{\label{foot:Milne}The proof of \cite{Mil80}*{Chapter~III, Proposition~1.27}, written for \'{e}tale cohomology, also works for flat cohomology: after shifting degrees as there, one reduces to the case $i = 0$, which is a claim about the restriction of an fppf sheaf to the small \'{e}tale site.}). Moreover, since this functor takes   coconnective values, descent (as opposed to hyperdescent) would suffice. The functor 
also vanishes on $A[\f 1p]$-algebras, so the descent assertion is insensitive to replacing the cover $A' \ra A''$ by $A' \ra A'' \times A'[\f 1p]$. Thus, we may assume that our cover is faithfully flat. By limit arguments, both functors $A' \mapsto R\Gamma(A', G)$ and $A' \mapsto R\Gamma(A'[\f 1p], G)$ satisfy  descent with respect to ind-fppf maps. Consequently, so does the functor $A' \mapsto R\Gamma_{\{ p = 0\}}(A', G)$, and the claim follows.
\epf


By \Cref{lem:large-cover}, we may build a hypercover $A^\bullet$ of $A$ in the topology whose covers are the filtered direct limits of flat, finitely presented maps that are faithfully flat modulo $p$ in such a way that each $A^i$ is $p$-Henselian and admits no nonsplit fppf covers. By \Cref{claim:descent}, the left side of \eqref{eqn:continuity-formula} satisfies hyperdescent with respect to this hypercover; by also using the deformation triangle \eqref{def-triangle-eq} and faithfully flat descent for modules, so does the right side. 
Effectively, we may replace $A$ by $A^i$ to reduce to the case when our ($0$-truncated) $p$-Henselian ring $A$ has no nonsplit fppf covers, a case in which
\be \label{eqn:dream-cohomology}
G(A) \isomto R\Gamma(A,G)
 \ee
(see \S\ref{pp:simplicial-coho-def}). We claim that for $n > 0$ we also have
\be \label{eqn:mod-pn-desire}
G(A/^\bL p^n) \isomto R\Gamma(A/^\bL p^n,G)
\ee
or, equivalently, that
\[
H^i(A/^\bL p^n,G)\cong0 \qxq{for} i>0.
\]
For this, \Cref{cor:nil-nil-deform-cohomology} gives $H^i(A/^\bL p^n,G) \isomto H^i(A/(p^n),G)$ for $i > 0$, so the B\'{e}gueri resolution \eqref{Begueri-0} reduces us to showing that for a smooth, affine $R$-group $Q$ we have 
\[
Q(A) \surjects Q(A/(p^n)) \qxq{and} H^i(A/(p^n), Q) \cong 0 \qxq{for} i > 0.
\]
The surjectivity follows from \cite{SP}*{Proposition~\href{https://stacks.math.columbia.edu/tag/07M8}{07M8}}. On the other hand, by 
\cite{SP}*{Lemma~\href{https://stacks.math.columbia.edu/tag/04D1}{04D1}}, every \'etale $(A/(p^n))$-algebra lifts to an \'{e}tale $A$-algebra, and hence is split, so the vanishing follows, too.

In the view of \eqref{eqn:dream-cohomology}--\eqref{eqn:mod-pn-desire}, it remains to show that for our $p$-Henselian, $0$-truncated $A$,
\be \label{eqn:pi0-and-ML-desire}
\tst G(A) \isomto R\lim_{n > 0}( G(A/^\bL p^n)) \q\x{(in the derived $\infty$-category $\mathcal D(\mathbb Z)$).}
\ee
For the derived $p$-adic completion $\widehat{A} \ce R\lim_n (A/^\bL p^n)$ of $A$, by  \Cref{lem:sq-0-trick}, the ring $\pi_0(\wh{A})$ 
is a square-zero extension of the $p$-adic completion of $\pi_0(A)$, so \Cref{BM-hens-comp} and \Cref{rem:no-new-etale} 
imply that 
\[
\tst G(A[\f{1}{p}]) \isomto G(\widehat{A}[\f{1}{p}]).
\]
\Cref{lem:derived-BL-key} and the affineness of $G$ then give
\[
\tst G(A) \isomto G(\widehat{A})\times_{G(\widehat{A}[\tfrac 1p])} G(A[\tfrac 1p]) \cong G(\wh{A}) \cong  \lim_{n > 0}( G(A/^\bL p^n))
\]
(limit in the $\infty$-category $\Ani(\Ab)$). Thus, to deduce \eqref{eqn:pi0-and-ML-desire}, it suffices to show that the system $\{\pi_0 (G(A/^\bL p^n))\}_{n > 0}$ is Mittag--Leffler. Let $p^k$ be a power that kills $e^*(L_{G/R})$, and consider an $n > k$. By \Cref{cor:deftheorygroup}, the obstruction to lift an $x \in \pi_0(G(A/^\bL p^n))$ to $\pi_0(G(A/^\bL p^{2n}))$ is an
\[
\tst \gA \in \pi_0\!\p{\Hom_R\!\p{e^*(L_{G/R}), A\tensor^\bL_\bZ(\frac{p^n\bZ}{p^{2n}\bZ})[1]}\!},
\]
where this last $\pi_0$ injects into
\[
\tst\pi_0\!\p{\Hom_R\!\p{e^*(L_{G/R}), A\tensor^\bL_\bZ(\frac{p^n\bZ}{p^{2n}\bZ} \tensor^\bL_\bZ \frac{\bZ}{p^k\bZ})[1]}\!}\!.
\]
Since the object $e^*(L_{G/R})$ is of projective amplitude $[-1, 0]$ (see \S\ref{def-setup}), the truncation triangle 
\[
\tst (\frac{p^{2n - k}\bZ}{p^{2n}\bZ})[1] \ra \frac{p^n\bZ}{p^{2n}\bZ} \tensor^\bL_\bZ \frac{\bZ}{p^k\bZ} \ra \frac{p^n\bZ}{p^{n + k}\bZ} \ra \p{\frac{p^{2n - k}\bZ}{p^{2n}\bZ}}[2]
\]
is exact and shows that the last displayed group injects into
\[
\tst \pi_0\!\p{\Hom_R\!\p{e^*(L_{G/R}), A\tensor^\bL_\bZ (\frac{p^n\bZ}{p^{n + k}\bZ})[1]}\!}\!.
\]
Consequently, by functoriality of \Cref{cor:deftheorygroup}, the obstruction to lifting $x$ to $\pi_0(G(A/^\bL p^{n + k}))$ is also $\gA$. In other words, if $x$ is in the image of $\pi_0(G(A/^\bL p^{n + k}))$, then it is also in the image of $\pi_0(G(A/^\bL p^{2n}))$ and, by replacing $n$ by $2n - k$ and iterating, we see that $x$ is in the image of $\pi_0(G(A/^\bL p^N))$ for every $N \ge n$, so that the system $\{\pi_0 (G(A/^\bL p^n))\}_{n > 0}$ is indeed Mittag--Leffler.
\end{proof}

\brem
In the final part of the proof, instead of checking the Mittag--Leffler condition we could also apply \Cref{lem:lame-case} (that does not use \Cref{thm:padiclimit} or the subsequent parts of \S\S\ref{section:continuity-formula}--\ref{sec:descent}).
\erem

The following concrete consequence of \Cref{thm:padiclimit} extends \Cref{cor:perfect-vanish} beyond perfect $\bF_p$-algebras. Its case when $A$ is a $0$-truncated $\bF_p$-algebra and $G = \mu_{p}$ was settled in \cite{Tre80}*{Theorem}.

\bcor \label{no-higher-coho}
Let $R$ be a ring, let $p$ be a prime, and let $G$ be a commutative, finite, locally free $R$-group of $p$-power order. For each animated $R$-algebra $A$ such that $\pi_0(A)$ is $p$-Henselian, 
\[
\tst H^i(A, G) \cong 0 \qxq{for} i \ge 3 
\]
{\upshape(}resp., 
\[
H^i(A, G) \cong 0 \qxq{for} i \ge 2 \quad \x{if $(\pi_0(A)/(p))^\red$ is perfect{\upshape)}.}
\]
\ecor

\bpf
\Cref{cor:nil-nil-deform-cohomology} reduces us to $0$-truncated $A$. By \Cref{no-nil-feelings}, 
for $n \ge 2$, the map
\[
H^i(A/(p^n), G) \ra H^i(A/(p^{n - 1}), G) \qxq{is} \begin{cases}\xq{surjective for} &i \ge 1, \\  \xq{bijective for} &i \ge 2. \end{cases}
\]
Thus,	 \Cref{thm:padiclimit} (with \Cref{cor:nil-nil-deform-cohomology} again) reduces us to the case when $A$ is an $\bF_p$-algebra. By \Cref{no-nil-feelings} again,
 we may then replace $A$ by $A^\red$ to assume that $A$ is reduced. If the resulting $A$ is perfect, then the Dieudonn\'{e}-theoretic \Cref{cor:perfect-vanish} gives the claim. Otherwise, we may assume that $A$ is Noetherian and consider the morphism of sites $\eps\colon A_\fppf \ra A_\et$ with its spectral sequence
\[
H^i_\et(A, R^j\eps_*(G)) \Rightarrow H^{i +j}_\fppf(A, G).
\]
By \cite{Gro68c}*{th\'{e}or\`{e}me~11.7} and the B\'{e}gueri sequence \eqref{Begueri-0}, we have 
\[
R^{\ge 2}\eps_*(G) \cong 0.
\] 
The \'{e}tale cohomological $p$-dimension of $A$ is $\le 1$ (see \cite{SGA4III}*{expos\'{e}~X, th\'{e}or\`{e}me~5.1}), so we obtain the desired $H^{\ge 3}_\fppf(A, G) \cong 0$. 
\epf



%



\Cref{thm:padiclimit} has the following consequence for the passage to the derived $p$-adic completion  for flat cohomology. At least for $0$-truncated $A$ and under additional bounded $p^\infty$-torsion assumptions, this could also be argued directly by an argument similar to the one used for \cite{BC19}*{Theorem~2.3.3~(d)}.

\bcor \label{heavy-handed}
Let $R$ be a ring, let $p$ be a prime, let $G$ be a commutative, finite, locally free $R$-group of $p$-power order, and let $A \ra A'$ be a map of animated $R$-algebras such that both $\pi_0(A)$ and $\pi_0(A')$ are $p$-Henselian and $A/^\bL p\isomto A'/^\bL p$. For each open 
\[
\tst \Spec(\pi_0(A)[\f{1}{p}]) \subset U \subset \Spec(\pi_0(A))
\]
with complement $Z \ce \Spec(\pi_0(A)) \setminus U$, we have \up{with definitions as in \uS\uref{pp:simplicial-coho-def}}
\[
\tst R\Gamma(U_A, G) \isomto R\Gamma(U_{A'}, G), \qxq{so also} R\Gamma_Z(A, G) \isomto R\Gamma_Z(A', G).
\]
\ecor

We will complement this last isomorphism with a general excision result of this sort in \Cref{cor:excision}. For a version of \Cref{heavy-handed} beyond the $p$-adic setting, see \Cref{cor:more-algebrization} below.

\bpf
\Cref{lem:sq-0-trick} ensures that $\pi_0(\wh{A})$ is a square-zero extension of the $p$-adic completion of $\pi_0(A)$, so is $p$-Henselian. Thus, we may replace $A'$ by $\wh{A}$ and use \Cref{BM-hens-comp} (with \Cref{rem:no-new-etale}) to obtain the case $U = \Spec(A[\f{1}{p}])$:
\be \lab{U-1p-case}
\tst R\Gamma(A[\f 1p], G) \isomto R\Gamma(A'[\f 1p], G).
\ee
In general, it remains to see that 
\[
R\Gamma_{\{p = 0\}}(U, G) \isomto R\Gamma_{\{p = 0\}}(U_{A'}, G).
\]
 For this, we may work locally on $U$, so, by passing to an affine cover and forming $p$-Henselizations (which do not change the $R\Gamma_{\{ p = 0\}}$; see \Cref{prop:excision-prep} below for a much more general result, in our case \eqref{eqn:exicision-prep-square} is Cartesian by descent), we reduce to $U = \Spec(\pi_0(A))$. Then  the continuity formula \eqref{eqn:continuity-formula} gives 
 \[
 R\Gamma(A, G) \isomto R\Gamma(A', G),
 \]
 so, due to \eqref{U-1p-case}, also the following desired isomorphism:
 \[
 R\Gamma_{\{p = 0\}}(A, G) \isomto R\Gamma_{\{p = 0\}}(A', G). \qedhere
 \]
\epf



\csub[Excision for flat cohomology and reduction to complete rings] \label{sec:excision}

To reduce purity for flat cohomology to the case of complete rings, we exhibit a general excision property of flat cohomology of (animated) rings in \Cref{cor:excision}, which vastly extends its special cases that appear in the literature: for instance, \cite{DH19}*{Lemma~2.6} proves it for excellent, Henselian discrete valuation rings (see also \cite{Maz72}*{Lemma~5.1} for an earlier special case). 
 The argument uses animated deformation theory of \S\S\ref{sec:defthy}--\ref{sec:fppf-coho-simplicial} and the $p$-adic continuity formula \eqref{eqn:continuity-formula} to eventually reduce to the positive characteristic case \eqref{crystalline-cor} of the key formula. 
 The bulk of it is captured by \Cref{prop:excision-prep}, which itself uses the following auxiliary lemma.

\blem \label{lem:auxiliary-perfect-semiperfect-control}
Let $A$ be a perfect $\bF_p$-algebra, let $A'$ be a semiperfect $A$-algebra, let $a \in A$ be such that $A/^\bL a\isomto A'/^\bL a$, and let 
\[
\tst A'_\perf \ce \varinjlim_{x \mapsto x^p} A'
\]
be the perfection of $A'$.   
The module 
\[
\Ker(A' \surjects A'_\perf) \qxq{is uniquely $a$-divisible and} A'/^\bL a\isomto A'_\perf/^\bL a.
\]
\elem

\bpf
Since $A$ is perfect, its $a^\infty$-torsion is bounded. Thus, since 
\[
A/^\bL a^n \isomto A'/^\bL a^n \qxq{and hence also} A\langle a^n\rangle \isomto A'\langle a^n\rangle,
\] 
the $a^\infty$-torsion of $A'$ is also bounded. Consequently, the $a$-adic completion of $A'$ agrees with the derived $a$-adic completion. The latter agrees with the $a$-adic completion $\wh{A}$ of $A$, and $\wh{A}$ is perfect, so the completion map $A' \ra \wh{A}$ factors through $A' \surjects A'_\perf$ and then even through the $a$-adic completion of $A'_\perf$. We obtain a factorization 
\[
A'/(a) \surjects A'_\perf/(a) \ra \wh{A}/(a) \cong A/(a)
\]
in which, by assumption, the composition is an isomorphism. Thus, both maps in this composition are isomorphisms. By repeating the same argument with $a^n$ in place of $a$, we get that $A' \ra A'_\perf$ induces an isomorphism on $a$-adic completions. Due to bounded $a^\infty$-torsion, these completions agree with their derived counterparts, so the resulting isomorphism $A'/^\bL a\isomto A'_\perf/^\bL a$ and the snake lemma give the unique $a$-divisibility of  $\Ker(A' \surjects A'_\perf)$.
\epf

\blem \label{prop:excision-prep} 
For a ring $R$, a commutative, finite, locally free $R$-group $G$, a map $A\to A^\prime$ of animated $R$-algebras, and an $a\in A$ such that $A/^\bL a \isomto A^\prime/^\bL a$, 
we have
\be \label{eqn:excision-prep-supp}
R\Gamma_{\{a = 0\}}(A, G) \isomto R\Gamma_{\{a = 0\}}(A', G),
\ee
equivalently,  the following square is Cartesian\ucolon
\be \label{eqn:exicision-prep-square} \ba\xymatrix{
R\Gamma(A,G)\ar[r]\ar[d] & R\Gamma(A[\tfrac 1a],G)\ar[d]\\
R\Gamma(A^\prime,G)\ar[r] & R\Gamma(A^\prime[\tfrac 1a],G).
}
\ea\ee
\elem

\begin{proof} 
The two formulations are equivalent because the fibers of the maps 
\[
\tst R\Gamma_{\{a = 0\}}(A, G) \ra R\Gamma_{\{a = 0\}}(A', G)
\]
and
\[
\tst R\Gamma(A,G) \ra R\Gamma(A^\prime,G) \times_{R\Gamma(A^\prime[\frac 1a],\,G)} R\Gamma(A[\frac 1a],G)
\]
are isomorphic, 
so we will focus on the Cartesian square statement. By decomposing $R$ into direct factors, we may assume that the order of $G$ is constant. By then expressing $G$ as the direct product of its primary factors, we may assume that this order is a power of a prime $p$. 

\begin{claim} \label{claim:def-theoretic-reduction}
For a square-zero extension $A\surjects B$ of animated $R$-algebras and the base changed square-zero extension $A' \surjects B' \ce A' \tensor_A^\bL B$, the square
\[\xymatrix{
R\Gamma(A,G)\ar[r]\ar[d] & R\Gamma(A[\tfrac 1a],G) \ar[d] \\
R\Gamma(A^\prime,G)\ar[r] & R\Gamma(A^\prime[\tfrac 1a],G)
}
\]
is Cartesian if and only if so is the square
\[\xymatrix{
  R\Gamma(B,G)\ar[r]\ar[d] & R\Gamma(B[\tfrac 1a],G)\ar[d] \\
  R\Gamma(B',G)\ar[r] & R\Gamma(B'[\tfrac 1a],G).
}
\]
\end{claim}

\bpf
For the ideal $M$ of $A\surjects B$, \Cref{lem:derived-BL-key} gives 
\[
\tst M\isomto (A' \tensor^\bL_A M) \times_{(A' \tensor^\bL_A M)[\f 1a]} M[\f 1a].
\]
This identification persists after applying $R\Hom_R(e^*(L_{G/R}), -)$, so the  triangle
\be \label{eqn:def-theoretic-triangle-middle-of-pf}
R\Gamma(A,G)\to R\Gamma(B, G)\to R\Hom_R(e^*(L_{G/R}), M[1])
\ee
of \Cref{prop:defthy} and its analogues after base change to $A'$, $A[\f 1a]$, and $A'[\f 1a]$ show that the natural map between the fibers of the map 
\[
\tst R\Gamma(A,G) \ra R\Gamma(A^\prime,G) \times_{R\Gamma(A^\prime[\frac 1a],\,G)} R\Gamma(A[\tfrac 1a],G)
\]
and of its analogue after base change to $B$ is an isomorphism, 
and the claim follows. 

For later use, we note that, due to the isomorphism between the fibers, we actually obtain a sharper variant: for instance, if $A$ is $0$-truncated, then, by repeating the argument inductively and forming a filtered direct limit, we see that the claim holds for any surjection $A\surjects B$ whose kernel is nil.
\epf

\begin{claim} \label{claim:def-theoretic-reduction-2}
For a square-zero extension of animated $R$-algebras $A\surjects B$  whose ideal $M$ is an $A[\f{1}{a}]$-module, the following square is Cartesian:
\[\xymatrix{
R\Gamma(A,G)\ar[r]\ar[d] & R\Gamma(B,G)\ar[d]   \\
R\Gamma(A[\tfrac 1a],G) \ar[r] & R\Gamma(B[\tfrac 1a],G).
}
\]
For  $0$-truncated $A$, the same holds for any surjection $A\surjects B$ whose kernel is nil and an $A[\f 1a]$-module. 
\end{claim}

\bpf
The assumption on $M$ implies that the term $R\Hom_R(e^*(L_{G/R}), M[1])$ in the deformation-theoretic triangle \eqref{eqn:def-theoretic-triangle-middle-of-pf} does not change once we replace $A$ by $A[\f{1}{a}]$. Thus, the same argument that gave the equivalence of \eqref{eqn:excision-prep-supp} and \eqref{eqn:exicision-prep-square} implies the first part of the claim.  For $0$-truncated $A$, by iteration, the claim holds for any surjection $A\surjects B$ whose kernel is nilpotent and an $A[\f 1a]$-module. By forming filtered direct limits, we may then weaken nilpotence to being nil. 
\epf

The main stages of the subsequent argument are:
\benuma
\item
to reduce to the case when $R$ is an $\bF_p$-algebra (so that $A$ and $A'$ are animated $\bF_p$-algebras); 

\m
to reduce to the case when $A$ is a $0$-truncated, perfect $\bF_p$-algebra;

\item
to reduce to the case when $A$ and $A'$ are $0$-truncated $\bF_p$-algebras with $A$ perfect;

\item
when $A$ and $A'$ are perfect $\bF_p$-algebras, to conclude using Dieudonn\'{e}-theoretic results of \S\ref{section:Dieudonne-positive-characteristic};

\item
to use the preceding step to reduce to the case when $A$ and $A'$ are perfect $\bF_p$-algebras.
\eenum
\emph{(1) Reduction to the case when $R$ is an $\bF_p$-algebra.} For any animated $R$-algebra $S$, letting  $S_{(p)}^h$ denote its $p$-Henselization (defined via \Cref{rem:no-new-etale}), descent supplies a functorial Cartesian~square
\[\xymatrix{
R\Gamma(S,G)\ar[r]\ar[d] & R\Gamma(S_{(p)}^h,G)\ar[d]\\
R\Gamma(S[\tfrac 1p],G)\ar[r] & R\Gamma(S_{(p)}^h[\tfrac 1p],G).
}\]
By applying this with $S$ replaced by, in turn, $A$, $A'$, $A[\f{1}{a}]$, and $A'[\f{1}{a}]$, we see that each term of \eqref{eqn:exicision-prep-square} is a glueing of its version after inverting $p$ with the version after $p$-Henselizing along the version where we first $p$-Henselize and then invert $p$. Therefore, since limits commute, it suffices to show that the following analogues of the square \eqref{eqn:exicision-prep-square} are all Cartesian:
\[\xymatrix{
R\Gamma(A[\tfrac 1p],G)\ar[r]\ar[d] & R\Gamma(A[\tfrac 1{pa}],G)\ar[d] & R\Gamma(A_{(p)}^h,G)\ar[r]\ar[d] & R\Gamma((A[\tfrac 1a])_{(p)}^h,G)\ar[d]\\
R\Gamma(A^\prime[\tfrac 1p],G)\ar[r] & R\Gamma(A^\prime[\tfrac 1{pa}],G), & R\Gamma(A_{(p)}^{\prime h},G)\ar[r] & R\Gamma((A^\prime[\tfrac 1a])_{(p)}^h,G),
}\]
\[\xymatrix{
R\Gamma(A_{(p)}^h[\tfrac 1p],G)\ar[r]\ar[d] & R\Gamma((A[\tfrac 1a])_{(p)}^h[\tfrac 1p],G)\ar[d]\\
R\Gamma(A_{(p)}^{\prime h}[\tfrac 1p],G)\ar[r] & R\Gamma((A^\prime[\tfrac 1a])_{(p)}^h[\tfrac 1p],G).
}\]
By \Cref{rem:no-new-etale}, the first and the last of these three squares depend only on the $\pi_0(-)$ of the animated rings involved. In addition, since $A/^\bL a^n \isomto A'/^\bL a^n$, we have $\pi_0(A)/(a^n) \isomto \pi_0(A')/(a^n)$, and, by \Cref{arc-descent} and \Cref{BM-hens-comp} (with \S\ref{pp:arc-topology}),  the functors 
\[
\tst R'\mapsto R\Gamma(R'[\frac 1p],G) \qxq{and} R'\mapsto R\Gamma(R'^h_{(p)}[\frac 1p],G)
\]
are arc sheaves on the category of $R$-algebras $R'$. Consequently, since arc descent implies formal gluing squares \cite{BM20}*{Theorem~6.4}, 
the first and the last squares are Cartesian. 
By the $p$-adic continuity formula \eqref{eqn:continuity-formula}, the remaining Cartesianness of the second square reduces us to the case when $A$ is over $\bZ/p^n\bZ$ for some $n > 0$, a case in which this square is nothing else than \eqref{eqn:exicision-prep-square}. Moreover,  \Cref{claim:def-theoretic-reduction} reduces us further to $n = 1$, in other words, we have achieved the promised overall reduction to the case when $R$ is an $\bF_p$-algebra.

\emph{(2) Reduction to the case when $A$ is a $0$-truncated, perfect $\bF_p$-algebra.}  
The map 
\[
A' \ra A' \tensor^\bL_A \tau_{\le n}(A)
\]
is an isomorphism after applying $\tau_{\le n}(-)$, so, by \eqref{eqn:fppf-hyperdescent} and \eqref{eqn:similar-to-Postnikov-tower}, 
\[
\tst R\Gamma(A, G) \isomto R\lim_n\p{R\Gamma(\tau_{\le n}(A), G)}
\]
and 
\[
R\Gamma(A', G) \isomto R\lim_n\p{R\Gamma(A' \tensor_A^\bL \tau_{\le n}(A), G)}\!,
\]
and likewise after inverting $a$. Thus, we may assume that $A$ is $n$-truncated for some $n > 0$. \Cref{claim:def-theoretic-reduction} then reduces to $n = 0$, that is, to the animated $\bF_p$-algebra $A$ being a usual $\bF_p$-algebra. The $A$-algebra 
\[
\tst A_\infty \ce A[X_a^{1/p^\infty}\, \vert\, a \in A]/(X_a - a\, \vert\, a \in A)
\]
is ind-fppf and semiperfect, and the same holds for its tensor self-products over $A$. Thus, ind-fppf descent for fppf cohomology allows us to replace $A$ by such a tensor self-product and $A'$ by its base change to reduce to the case when $A$ is semiperfect. The last paragraph of the proof of \Cref{claim:def-theoretic-reduction} then allows us to divide out the nil-ideal $\Ker(A\surjects A_\perf)$ to reduce to the case when $A$ is perfect (here and below $A_\perf \ce \varinjlim_{x \mapsto x^p} A$). 

\emph{(3) Reduction to the case when $A$ and $A'$ are $0$-truncated $\bF_p$-algebras with $A$ perfect.}  
Our $0$-truncated, perfect $\bF_p$-algebra $A$ has bounded $a^\infty$-torsion, and $A/^\bL a^n \isomto A'/^\bL a^n$. Thus, each $\pi_{i}(A')$ with $i \ge 2$ is uniquely $a$-divisible and there is an exact sequence 
\[
0 \ra \pi_1(A') \xra{a^n} \pi_1(A') \ra A\langle a^n\rangle \ra \pi_0(A') \xra{a^n} \pi_0(A'),
\]
which, by letting $n$ grow, shows that $\pi_0(A')$ has bounded $a^\infty$-torsion and that $\pi_1(A')$ is also uniquely $a$-divisible. Then even 
\[
A\langle a^\infty\rangle \isomto (\pi_0(A))\langle a^\infty\rangle, \qxq{so also} A/^\bL a \isomto \pi_0(A')/^\bL a.
\]
Moreover, since each $\pi_i(A')$ with $i \ge 1$ is an $A[\f{1}{a}]$-module, the deformation-theoretic \eqref{eqn:def-theoretic-triangle-middle-of-pf} implies that for each $n > 0$,
\[\xymatrix{
R\Gamma(\tau_{\le n}(A'),G)\ar[r]\ar[d] & R\Gamma(\tau_{\le n - 1}(A'),G) \ar[d] \\
R\Gamma(\tau_{\le n}(A')[\tfrac 1a],G)\ar[r] & R\Gamma(\tau_{\le n - 1}(A')[\tfrac 1a],G)
}
\]
is a Cartesian square. Thus, the same square with $\pi_0(A')$ in place of $\tau_{\le n - 1}(A')$ is also Cartesian. By passing to the inverse limit in $n$ and using \eqref{eqn:fppf-hyperdescent}, we then find that the right square in
\be \label{eqn:two-Cartesian-squares-argument} \ba\xymatrix{
R\Gamma(A, G) \ar[d] \ar[r] & R\Gamma(A',G)\ar[r]\ar[d] & R\Gamma(\pi_0(A'),G) \ar[d] \\
R\Gamma(A[\f 1a], G) \ar[r] & R\Gamma(A'[\tfrac 1a],G)\ar[r] & R\Gamma(\pi_0(A')[\tfrac 1a],G)
}
\ea \ee
is Cartesian. Thus, the sought Cartesianness of the left square reduces to that of the outer one, 
so we may replace $A'$ by $\pi_0(A')$ to reduce to the case when $A$ and $A'$ are $\bF_p$-algebras with $A$ perfect.

\emph{(4) Conclusion in the case when $A$ and $A'$ are perfect  $\bF_p$-algebras.}  
We now treat the case when the $\bF_p$-algebra $A'$ is also perfect. Then 
\Cref{KT-input} identifies the map
\[
R\Gamma_{\{a = 0\}}(A, G) \ra R\Gamma_{\{a = 0\}}(A', G)
\]
with
\[
 R\Gamma_{(p,\, a)}(W(A), \bM(G))^{V - 1} \ra R\Gamma_{(p,\,  a)}(W(A'), \bM(G_{A'}))^{V - 1}.
\]
For showing that the latter is an isomorphism, by d\'{e}vissage, we may remove $(-)^{V - 1}$, replace $\bM(G)$ by a projective $W(A)$-module (see \S\ref{Gab-eq}), 
and then even replace it by a free $W(A)$-module. In other words, we are reduced to showing that the following map  is an isomorphism:
\[
R\Gamma_{(p,\,a)}(W(A), W(A)) \ra R\Gamma_{(p,\, a)}(W(A'), W(A')).
\]
The cohomology groups of the fiber of this map are $p$-power torsion, so we need to show that multiplication by $p$ is an automorphism on them. Thus, the sequences 
\[
0 \ra W(A) \xra{p} W(A) \ra A \ra 0 \qxq{and} 0 \ra W(A') \xra{p} W(A') \ra A' \ra 0
\]
reduce us to showing that the fiber of the map
\[
R\Gamma_{(p,\, a)}(W(A), A) \ra R\Gamma_{(p,\, a)}(W(A'), A'), 
\]
that is, of the map
\[
R\Gamma_{\{ a = 0\}}(A, A) \ra R\Gamma_{\{ a = 0\}}(A', A'),
\]
vanishes. The last fiber agrees with that of the map
\[
\tst R\Gamma(A, A) \ra R\Gamma(A', A') \times_{R\Gamma(A'[\f 1a],\, A'[\f 1a])} R\Gamma(A[\f 1a], A[\f 1a]),
\]
that is, with the fiber of the map $A\ra A' \times_{A'[\f 1a]}A[\f 1a]$, which is an isomorphism by \Cref{lem:derived-BL-key}.

\emph{(5) Conclusion in the general case.}  Having established the case when both of our $\bF_p$-algebras $A$ and $A'$ are perfect, we return to the situation in which only $A$ is perfect and let $\wh{A}$ be the $a$-adic completion of $A$. Since $A$ has bounded $a^\infty$-torsion, $\wh{A}$ agrees with its derived counterpart, so that we have a map $A' \ra \wh{A}$ with $A'/^\bL a \isomto \wh{A}/^\bL a$. We consider the analogue of \eqref{eqn:two-Cartesian-squares-argument} with the perfect $\bF_p$-algebra $\wh{A}$ in place of $\pi_0(A')$. By the preceding step of the overall argument, in this analogue the outer square is Cartesian, so we are reduced to showing that the right one is, too. In other words, we may replace $A$ and $A'$ by $A'$ and $\wh{A}$, respectively, to reduce to the case when the $\bF_p$-algebra $A$ is arbitrary but $A'$ is perfect. In this situation, we repeat the reduction to perfect $A$ and note that it transforms our perfect $A'$ into an animated $\bF_p$-algebra for which $\pi_0(A')$ is semiperfect (the passage to semiperfect $A$ leaves $A'$ semiperfect, and the subsequent derived base change along $A\surjects A_\perf$ may introduce higher homotopy). After this, we repeat the reduction that uses the boundedness of the $a^\infty$-torsion of $A$ to replace $A'$ by $\pi_0(A')$ and are left with the case in which $A$ is perfect and $A'$ is semiperfect. \Cref{lem:auxiliary-perfect-semiperfect-control} then ensures that the nil-ideal $\Ker(A' \surjects A'_\perf)$ is an $A[\f 1a]$-module. Thus, \Cref{claim:def-theoretic-reduction-2} shows the Cartesianness of the square 
\[\xymatrix{
R\Gamma(A',G)\ar[r]\ar[d] & R\Gamma(A'_\perf,G)\ar[d]   \\
R\Gamma(A'[\tfrac 1a],G) \ar[r] & R\Gamma(A'_\perf[\tfrac 1a], G).
}
\]
By \Cref{lem:auxiliary-perfect-semiperfect-control}, we have $A'/^\bL a \isomto A'_\perf/^\bL a$, so the same argument as for \eqref{eqn:two-Cartesian-squares-argument} allows us to replace $A'$ by $A'_\perf$. However, then both $A$ and $A'$ are perfect, a case that we already settled. 
\end{proof}

Before deducing the sought \Cref{cor:excision}, we clarify its excision condition. 

\blem \label{comp-explain}
For a map $f\colon A \ra A'$ of animated rings, an ideal  
\[
I = (a_1, \dotsc, a_r) \subset \pi_0(A) \qxq{satisfies} \pi_0(A)/I \isomto (\pi_0(A)/I) \tensor^\bL_A A'
\]
if and only if $f$ induces an isomorphism after iteratively forming derived $a_i$-adic completions for $i = 1,\dotsc, r$, equivalently, if and only if
\[
A/^\bL(a_1^n, \dotsc, a_r^n) \isomto A'/^\bL(a_1^n, \dotsc, a_r^n) \qxq{for every} n > 0;
\] 
in particular, all of these equivalent conditions depend only on the ideal $I$ and not on the $a_i$.
\elem

\bpf
By \S\ref{ani-rings}, the iterated derived $a_i$-adic completion of $A$ is identified with 
\[
\tst \varprojlim_{n_1,\, \dotsc,\, n_r \ge 0} (A/^\bL(a_1^{n_1}, \dotsc, a_r^{n_r})),
\]
and likewise for $A'$. Thus, since the inverse subsystem where all the $n_i$ are equal is final, $f$ induces an isomorphism on iterated derived $a_i$-adic completions if and only if 
\be \label{A-Apr-map}
A/^\bL(a_1^n, \dotsc, a_r^n) \isomto A'/^\bL(a_1^n, \dotsc, a_r^n) \qxq{for every} n > 0.
\ee
If this holds for $n = 1$, then, since 
\[
\tst (\pi_0(A)/I) \tensor^\bL_A A' \cong (\pi_0(A)/I) \tensor^\bL_{A/^\bL(a_1,\, \dotsc,\, a_r)} A'/^\bL(a_1, \dotsc, a_r),
\] 
also $\pi_0(A)/I \isomto (\pi_0(A)/I) \tensor^\bL_A A'$. Conversely, if this last map is an isomorphism, then
the cofiber of the map \eqref{A-Apr-map} is an animated $(A/^\bL(a_1^n, \dotsc, a_r^n))$-module that vanishes after applying $(\pi_0(A)/I) \tensor^\bL_{A/^\bL(a_1^n,\, \dotsc,\, a_r^n)} -$, and hence  itself vanishes by \Cref{support-lem}.
\epf

\bthm\label{cor:excision} 
Let $R$ be a ring and let $G$ be a commutative, finite, locally free $R$-group. For each map $A\to A^\prime$ of animated $R$-algebras and each finitely generated ideal $I\subset \pi_0(A)$ such that 
\[
\pi_0(A)/I \isomto (\pi_0(A)/I)\otimes^\bL_A A'
\]
\up{see Lemma~\uref{comp-explain}}, we have
\[
R\Gamma_I(A, G) \isomto R\Gamma_{I}(A^\prime,G). 
\]
\ethm

\begin{proof} 
By \Cref{comp-explain}, we may write $I = (a_1, \dotsc, a_r)$ and assume that $A'$ is the iterated derived $a_i$-adic completion of $A$ for $i = 1, \dotsc, r$, 
and then, by arguing inductively, assume instead that $A'$ is the derived $a$-adic completion of $A$ for some $a \in I$. 
There is a functorial fiber sequence
\[
\tst R\Gamma_I(A,G)\to R\Gamma_{\{a = 0\}}(A,G)\to R\lim_{f \in I} \p{R\Gamma_{\{a = 0\}}(A[\f 1f], G)}, 
\]
and likewise for $A'$, so the claim follows from \Cref{prop:excision-prep} applied to $A\to A^\prime$ and to its localizations.
\end{proof}

We are ready to show that the validity of \Cref{main} depends only on the completion $\wh{R}$. 

\bcor \label{reduce-to-complete}
For a Noetherian local ring $(R, \fm)$ and a commutative, finite, flat $R$-group $G$,
\[
H^i_\fm(R, G) \isomto H^i_{\fm}(\wh{R}, G) \qxq{for every} i \in \bZ.
\]
In particular, Theorem \uref{main} reduces to its case when the complete intersection $R$ is $\fm$-adically complete.
\ecor

\bpf
Indeed, $\wh{R}$ is $R$-flat with $R/\fm \isomto \wh{R}/\fm \wh{R}$, 
so \Cref{cor:excision} gives the claim.
\epf

\brem \label{rem:easy-complete}
Under additional assumptions on $R$ or $G$, previous results suffice to reduce \Cref{main} to complete $R$. 
For instance, if $R$ is excellent (as one can sometimes reduce to using \cite{Pop19}*{Corollary 3.10}), 
then \cite{BC19}*{Lemma~2.1.3} (with elementary excision \cite{Mil80}*{Chapter~III, Proposition~1.27}) 
gives 
\[
H^i_\fm(R, G) \hra H^i_{\fm}(\wh{R}, G).
\]
If instead $G$ is \'{e}tale, then \Cref{lem:et-to-hat} (or already \cite{Fuj95}*{Corollary~6.6.4}) suffices. 
\erem


\csub[$p$-complete arc descent over perfectoids and fpqc descent] \label{sec:descent}

As we saw in \Cref{M(G)-coho-descent}, the Dieudonn\'{e} module side of the key formula \eqref{eqn:key-formula} satisfies $p$-complete arc hyperdescent. The main goal of this section is to establish the same for the flat cohomology side in \Cref{thm:strongdescent-perfectoid}, which will lead us to showing the key formula in \S\ref{perfectoid-version}. The main inputs are the $p$-adic continuity formula \eqref{eqn:continuity-formula}, deformation theory, and the $p$-complete arc hyperdescent for the structure presheaf on perfectoid rings \cite{BS19}*{Proposition~8.10} (see \Cref{arc-sheaves}).


\bthm \label{thm:strongdescent-perfectoid}
For a prime $p$, a perfectoid ring $A$, a commutative, finite, locally free $A$-group $G$, and a closed $Z \subset \Spec(A/pA)$, 
both functors
\[
A' \mapsto R\Gamma(A', G) \qxq{and} A' \mapsto R\Gamma_Z(A', G) 
\] 
satisfy hyperdescent for those $p$-complete arc hypercovers whose terms are perfectoid $A$-algebras.
\ethm

\bpf
For any $a \in A$, by \Cref{cor:ind-etale}, the $p$-adic completion $\wh{A[\f1a]}$ of $A[\f1a]$ is perfectoid. Moreover, by excision \Cref{cor:excision}, we have
\[
\tst R\Gamma_{\{p=0\}}(A[\f 1a],G) \isomto R\Gamma_{\{p=0\}}(\wh{A[\f 1a]},G),
\]
as one can also deduce from the simpler \Cref{heavy-handed} (with \eqref{elt-excision} and \eqref{no-big-tor}) for the $p$-primary part of $G$ and from the consequence \cite{BM20}*{Theorem 6.4} of arc descent \Cref{arc-descent} for the prime to $p$ part of $G$).
 Letting $a \in A$ range over the elements vanishing on $Z$, the functorial fiber sequence
\be \label{eqn:p-for-Z}
\tst R\Gamma_Z(A,G)\to R\Gamma_{\{p=0\}}(A,G)\to R\lim_a \p{R\Gamma_{\{p=0\}}(A[\f 1a],G)}
\ee
therefore reduces us to the case when $Z = \Spec(A/pA)$. Moreover, by \Cref{arc-descent}, the functor $A' \mapsto R\Gamma(A'[\f 1p], G)$ satisfies hyperdescent for those $p$-complete arc hypercovers whose terms are perfectoid $A$-algebras, so the cohomology with supports triangle reduces us to the functor 
\[
A' \mapsto R\Gamma(A', G).
\]
For the latter, the arc descent aspect of \Cref{arc-descent} and invariance of \'{e}tale cohomology under Henselian pairs \cite{Gab94a}*{Theorem 1} take care of the prime to $p$-factor of $G$. Thus, we may assume that $G$ is of $p$-power order and use the $p$-adic continuity formula \eqref{eqn:continuity-formula} to reduce to the functor 
\[
A' \mapsto R\Gamma(A'/^\bL p^n, G).
\] 
By \Cref{arc-sheaves}, the functor $A' \mapsto A'/^\bL p$ satisfies hyperdescent for those $p$-complete arc hypercovers whose terms are perfectoid $A$-algebras. Thus, the deformation theoretic \Cref{prop:defthy} 
allows us to decrease $n$ to reduce to showing $p$-complete arc hyperdescent for the functor $A' \mapsto R\Gamma(A'/^\bL p, G)$. To reduce further, we choose a system $\pi^{1/p^n} \in A$ of compatible $p$-power roots with $\pi$ a unit multiple of $p$ in $A$ (see \S\ref{perfectoid-def}). As above, by \Cref{arc-sheaves}, each functor $A' \mapsto A'/^\bL \pi^{1/p^n}$ satisfies hyperdescent for those $p$-complete arc hypercovers whose terms are perfectoid $A$-algebras $A'$. Thus, by iteratively applying \Cref{prop:defthy} we see that for any $p$-complete arc hypercover $A' \ra A'^\bullet$ whose terms are perfectoid $A$-algebras, the fiber of the hyperdescent comparison map
\[
\tst R\Gamma(A'/^\bL \pi^{1/p^n}, G) \ra R\lim_\Delta( R\Gamma(A'^\bullet/^\bL \pi^{1/p^n}, G))
\]
maps isomorphically to the corresponding fiber in which $n$ replaced by $n + 1$. Thus, by passing to the filtered direct limit over all $n \ge 0$, we even reduce to considering the functor 
\[
\tst A' \mapsto R\Gamma(\varinjlim_n (A'/^\bL \pi^{1/p^n}), G).
\]
Since $A'/ (\pi^{1/p^\infty}) \cong (A'/pA')^\red$ (see \eqref{p-oid-mod-pi}) and $A'\langle \pi\rangle = A'\langle \pi^{1/p^\infty} \rangle$ (see \eqref{no-big-tor}), we have 
\[
\tst \varinjlim_n (A'/^\bL \pi^{1/p^n}) \isomto (A'/pA')^\red.
\]
Thus, we may replace $A$ by the perfect $\bF_p$-algebra $(A/pA)^\red$ to reduce to showing that the functor $A' \mapsto R\Gamma(A', G)$ satisfies hyperdescent for those arc hypercovers whose terms are perfect $(A/pA)^\red$-algebras. By \eqref{crystalline-cor}, this functor is nothing else than 
\[
A' \mapsto R\Gamma(W(A'), \bM(G_{A'}))^{V = 1},
\]
so \Cref{M(G)-coho-descent} gives the claim.
\epf

The method also leads to the following agreement of fppf cohomology with fpqc cohomology.

\begin{theorem}\label{thm:strongdescent} 
For a ring $R$ and a commutative, finite, locally free $R$-group $G$, the functor 
\[
A\mapsto R\Gamma_\fppf(A,G)
\]
satisfies fpqc hyperdescent on animated $R$-algebras $A$. 
\end{theorem}

\begin{proof} 
We may restrict to $G$ of $p$-power order for a prime $p$ and have a functorial Cartesian square
\[\xymatrix{
R\Gamma_\fppf(A,G)\ar[r]\ar[d] & R\Gamma_\fppf(A_{(p)}^h,G)\ar[d]\\
R\Gamma_\fppf(A[\tfrac 1p],G)\ar[r] & R\Gamma_\fppf(A_{(p)}^h[\tfrac 1p],G),
}\]
where $(-)_{(p)}^h$ denotes $p$-Henselization (see \Cref{rem:no-new-etale}). By \eqref{et-fppf-agree} and \Cref{rem:no-new-etale}, 
\[
\tst R\Gamma_\fppf(A[\frac 1p], G) \cong R\Gamma_\et(\pi_0(A)[\frac 1p], G)
\]
and 
\[
\tst R\Gamma_\fppf(A_{(p)}^h[\frac 1p],G) \cong R\Gamma_\et(\pi_0(A)_{(p)}^h[\frac 1p], G),
\]
so arc descent, that is, \Cref{arc-descent} and \Cref{BM-hens-comp} (with \S\ref{pp:arc-topology}~\ref{arc-fflat}), reduces us to showing fpqc hyperdescent for the functor
\[
A\mapsto R\Gamma_\fppf(A_{(p)}^h,G).
\]
\Cref{thm:padiclimit} then allows us to instead consider the functor 
\[
A\mapsto R\Gamma_\fppf(A/^\bL p^n,G).
\]
By fpqc hyperdescent for modules (see \S\ref{pp:simplicial-coho-def}), the deformation-theoretic \Cref{prop:defthy} reduces us further to $n = 1$. In other words, we have reduced to $R$ being an $\bF_p$-algebra. Postnikov-convergence of \Cref{cor:fppf-hyperdescent} then allows us to assume that $A$ is $n$-truncated, and we apply \Cref{prop:defthy} again  to assume further that $A$ is even $0$-truncated. 

The $A$-flat $A'$ are then also $0$-truncated, so $R\Gamma_\fppf(-, G)$ takes coconnective values on them. In particular, it remains to show fpqc descent as opposed to hyperdescent. 
Moreover, by ind-fppf descent as in the proof of \Cref{prop:excision-prep}, we may assume that $A$ is semiperfect. By applying the deformation-theoretic \Cref{prop:defthy} and passing to the direct limit over the nilpotent ideals of $A$, we may even replace $A$ by its perfection (compare with the proof of \Cref{prop:excision-prep}).  For perfect $A$, we strengthen the sought claim: we will show that any faithfully flat $A\ra A'$ with  \v{C}ech nerve $A'^\bullet$ is of universal descent for $R\Gamma_\fppf(-, G)$ in the sense that for any $A$-algebra $B$, 
\[
\tst R\Gamma_\fppf(B,G) \isomto  R\lim_\Delta \p{R\Gamma_\fppf(A^{\prime\bullet} \tensor_A B,G)}\!.
\]
The advantage of allowing any $B$ is that then, by \cite{LZ17}*{Lemma~3.1.2 (3)}, we can replace $A\to A'$ by any refinement. Thus, by \Cref{lem:perf-still-pro-fppf}, we may assume that $A'$ is also perfect. Then we repeat the same reductions as above to first assume by ind-fppf descent that $B$ is semiperfect and then, by deformation theory, perfect. Once $A$, $A'$, and $B$ are perfect, so is $B'=B\otimes_A A'$, to the effect that we have reduced the original claim to the case when both $A$ and $A'$ are perfect $\bF_p$-algebras and $G$ is of $p$-power order. \Cref{KT-input} then reduces us to showing that
\[
\tst R\Gamma(W(A),\bM(G_A))^{V = 1} \isomto  R\lim_\Delta \p{R\Gamma(W(A^{\prime\bullet}), \bM(G_{A^{\prime\bullet}}))^{V = 1}}\!.
\]
For this, by considering fiber sequences, we may first drop the superscripts ``$V = 1$,'' and then, by resolving $\bM(G_A)$ by finite projective $W(A)$-modules (see \S\ref{Gab-eq}) and expressing the latter as direct summands of finite free $W(A)$-modules, we reduce to showing that 
\[
\tst W(A) \isomto  R\lim_\Delta \p{W(A^{\prime\bullet})}\!.
\]
Since the Witt vectors of a perfect ring is a derived inverse limit of its reductions modulo powers of $p$, we may now drop $W(-)$ from both sides and conclude by fpqc descent.
\end{proof}

We conclude the section with the following consequence for $p$-completely faithfully flat descent. 

\begin{corollary}\label{cor:descentsupport} 
Let $R$ be a ring, let $p$ be a prime, let $Z \subset \Spec(R/pR)$ be a closed subset, and let $G$ be a commutative, finite, locally free $R$-group of $p$-power order. On animated $R$-algebras $A$, both
\[
A\mapsto R\Gamma(A^h_{(p)}, G) \qxq{and} A\mapsto R\Gamma_Z(A,G)
\]
satisfy hyperdescent for cosimplicial algebras $A \ra A^\bullet$ such that $A/^\bL p \ra A^\bullet/^\bL p$ is an fpqc hypercover.
\end{corollary}

\begin{proof} 
The $p$-Henselization $A^h_{(p)}$ has no effect on $A$ modulo powers of $p$ and, by \Cref{lem:reduced-flatness}, the map $A \ra A^\bullet$ is a faithfully flat hypercover modulo powers of $p$. 
Thus, the claim about the functor $A\mapsto R\Gamma(A^h_{(p)}, G)$ follows from the $p$-adic continuity formula \eqref{eqn:continuity-formula} and \Cref{thm:strongdescent}. 
For the claim about $A\mapsto R\Gamma_Z(A, G)$, the fiber sequence \eqref{eqn:p-for-Z} reduces us to the case when $Z = \Spec(R/pR)$. In this case, by excision \eqref{elt-excision}, there is a functorial fiber sequence
\[
\tst R\Gamma_{\{p=0\}}(A,G)\to R\Gamma(A^h_{(p)},G)\to R\Gamma(A^h_{(p)}[\frac 1p],G), 
\]
which, due to the settled case of $A\mapsto R\Gamma(A^h_{(p)}, G)$, reduces us to considering $A\mapsto R\Gamma(A^h_{(p)}[\frac 1p],G)$. By \eqref{et-fppf-agree} and \Cref{rem:no-new-etale}, this functor agrees with $A\mapsto R\Gamma_\et(\pi_0(A)^h_{(p)}[\frac 1p],G)$, which, by \Cref{BM-hens-comp} (with \S\ref{pp:arc-topology}~\ref{arc-pi-fflat}), satisfies hyperdescent for $p$-completely faithfully flat covers. 
\end{proof}



\csub[The continuity formula for flat cohomology] \label{sec:genl-continuity}


We wish to supplement the $p$-adic continuity formula of \Cref{thm:padiclimit} with a general continuity formula for flat cohomology that we present in \Cref{adic-continuity} below. We will not use this general formula other than in \S\ref{sec:genl-continuity} (whose results will not be used elsewhere in this article), but we give some of its consequences in \Cref{inv-Hens-pair,cor:more-algebrization,cor:Gab-ann}. The continuity formula concerns those animated rings that are derived $I$-adically complete as follows.

\bpp[Derived $I$-adic completeness] \label{pp:derived-complete}
An animated ring $A$ is \emph{derived $I$-adically complete} for an ideal $I \subset \pi_0(A)$ if $A$ is derived $a$-adically complete for every $a \in I$, in other words, if
\[
\tst A \isomto R\lim_{n > 0} (A/^\bL a^n) \qxq{for} a \in I.
\]
Similarly to \S\ref{ani-rings}, we may consider $A$ as an object of $D(\bZ[X])$ via the map $\bZ[X] \ra A$ given by $X\mapsto a$. Thus, by \cite{BS15}*{Proposition~3.4.4}, derived $I$-adic completeness of $A$ is equivalent to each $\pi_i(A)$ being derived $a$-adically complete for every $a \in I$, and \cite{BS15}*{Lemma~3.4.12} allows us to restrict to those $a$ that lie in a fixed generating set of $I$. Thus, derived $I$-adic completeness is stable under truncations and if $I = (a_1, \dotsc, a_n)$ is finitely generated, then, arguing as in the proof of \Cref{comp-explain}, we see that $A$ is derived $I$-adically complete if and only if
\be \label{eqn:continuity-condition}
\tst A \isomto R\lim_{n > 0} (A/^\bL (a_1^n, \dotsc, a_r^n)).
\ee
Derived $I$-adic completeness of $A$ is weaker than $I$-adic completeness of homotopy groups: if each $\pi_i(A)$ is (classically) $I$-adically complete, then $A$ is derived $I$-adically complete, see \cite{BS15}*{Lemma~3.4.13}.
\epp

Derived $I$-adic completeness of $A$ implies $I$-Henselianity of $\pi_0(A)$ as follows. 


\blem \label{lem:DC-comp-Hens}
For an animated ring $A$ and an ideal $I \subset \pi_0(A)$ such that $A$ is derived $I$-adically complete, the ring $\pi_0(A)$ is $I$-Henselian. 
\elem

\bpf
By \S\ref{pp:derived-complete}, 
the animated aspect of the statement is illusory: we may replace $A$ by $\pi_0(A)$ to assume that $A$ is $0$-truncated. In this case, we use \cite{SP}*{Lemma~\href{https://stacks.math.columbia.edu/tag/0G1S}{0G1S}} to reduce further to the case when $I$ is principal, generated by an $a \in I$, so that, by \Cref{lem:sq-0-trick}, the ring $A$ an extension of its $a$-adic completion by a square-zero ideal $J$. 
Thus, since the $a$-adic completion 
is $a$-Henselian, \cite{SP}*{Lemma~\href{https://stacks.math.columbia.edu/tag/0DYD}{0DYD}} ensures that $A$ is $(J + (a))$-Henselian, and then also $a$-Henselian. 
\epf

The following special case of \Cref{inv-Hens-pair} is an input to the overall proof of the continuity formula.


\blem \lab{baby-inv-Hens-pair}
For a Henselian pair $(R, I)$ and a commutative, finite, locally free $R$-group $G$, 
\be \lab{BIHP-map}
H^2(R, G) \isomto H^2(R/I, G).
\ee
\elem

\bpf
By \cite{BC19}*{Theorem~2.1.6}, the map \eqref{BIHP-map} is injective and, for every smooth, affine $R$-group~$Q$, we have
\be \label{eqn:BC-input}
H^1(R, Q) \isomto H^1(R/I, Q) 
\ee
and
\[
\Ker(H^2(R, Q) \ra H^2(R/I, Q)) = \{*\}.
\]
Thus, the B\'{e}gueri sequence \eqref{Begueri-0} and the five lemma reduce the surjectivity of \eqref{BIHP-map} to that of the map 
\[
H^2(R', \bG_m)_\tors \ra H^2(R'/IR', \bG_m)_\tors \qxq{where} R' \ce \Gamma(G^*, \sO_{G^*}).
\]
For the latter, by \cite{Gab81}*{Chapter II, Theorem 1}, every element of the target $H^2(R'/IR', \bG_m)_\tors$ comes from some $H^1(R'/IR', \PGL_n)$, so it suffices to apply \eqref{eqn:BC-input} to $\PGL_{n}$ over the Henselian pair $(R', IR')$.
\epf

We use these lemmas to show the following case of \Cref{adic-continuity} that will serve as an input to the general case. Its vanishing condition holds when $A$ is of characteristic $p$ and $G$ is of $p$-power order (see \Cref{no-higher-coho}), so this case essentially includes the continuity formula in positive characteristic.

\blem \label{lem:lame-case}
Let $R$ be a ring, let $G$ be a commutative, finite, locally free $R$-group, let $A$ be an animated $R$-algebra, and let $a \in A$ be an element such that $A$ is derived $a$-adically complete and
\[
H^{i}(A, G) \cong H^i(A/^\bL a, G) \cong 0 \qxq{for} i \ge 3
\]
{\upshape(}equivalently, 
\[
H^{i}(\pi_0(A), G) \cong H^i(\pi_0(A)/(a), G) \cong 0 \qxq{for} i \ge 3).
\]
Then the continuity formula holds\ucolon
\be \label{eqn:lame-map}
\tst R\Gamma(A, G) \isomto R\lim_{n > 0}( R\Gamma(A/^\bL a^n, G)).
\ee
\elem

\bpf
The parenthetical reformulation of the vanishing condition follows from \Cref{cor:nil-nil-deform-cohomology}. This condition and \Cref{cor:nil-nil-deform-cohomology,no-nil-feelings} ensure that the map \eqref{eqn:lame-map} is an isomorphism in cohomological degrees $i \ge 3$. Likewise, by also using \Cref{baby-inv-Hens-pair,lem:DC-comp-Hens} and the surjectivity aspect of \Cref{no-nil-feelings} in cohomological degree $i = 1$ we see that the map \eqref{eqn:lame-map} is an isomorphism in cohomological degree $i = 2$. As for nonpositive degrees, it suffices to observe the following identification that results from $G$ being affine and $A$ being derived $a$-adically complete:
\[
\tst G(A) \isomto R\lim_{n > 0} G(A/^\bL a^n) \qx{(limit in the $\infty$-category $\Ani(\Ab)$).}
\]
To proceed with the remaining cohomological degree $1$, one possible approach (suggested by the referee) is to carry out a discussion of the cohomological classification of $G$-torsors over animated rings and then deduce the claim from the equivalence $\Vect(A) \simeq \lim_n \Vect(A/^\bL a^n)$ for $\infty$-categories of vector bundles. For the sake of expedience, we opt for a more pedestrian approach that is based on the following general reduction to the case when $A$ is $0$-truncated and $a$-adically complete.

By \cite{BS15}*{Proposition~3.4.4}, the derived $a$-adic completeness of $A$ amounts to the derived $X$-adic completeness of each $\pi_i(A)$ viewed as a $\bZ[X]$-module via the map $\bZ[X] \mapsto \pi_0(A)$ given by $X \mapsto a$. In particular, the truncations $\tau_{\le m}(A)$ are again derived $a$-adically complete. Thus, by applying the strengthened version of Postnikov completeness presented in \Cref{cor:fppf-hyperdescent}, we may replace $A$ by the $\tau_{\le m}(A)$ to assume that $A$ is $m$-truncated for some $m \ge 0$. If $m > 0$, then, by \Cref{sq-zero-eg}~\ref{sq-0-3}, such an $A$ is a square-zero extension of $\tau_{\le m - 1}(A)$ by $\pi_m(A)[m]$, so that the deformation-theoretic triangle
\[
R\Gamma(A, G) \ra R\Gamma(\tau_{\le m - 1}(A), G) \ra R\Hom_R(e^*(L_{G/R}), \pi_m(A)[m + 1])
\]
supplied by \Cref{prop:defthy}, its analogues after derived reduction modulo $a^n$, and the derived $a$-adic completeness of $\pi_m(A)$ allow us to replace $A$ by $\tau_{\le m - 1}(A)$. By decreasing $m$ in this way, we reduce to the case when our animated $R$-algebra $A$ that satisfies the condition on the vanishing of cohomology is $0$-truncated. Moreover, since $A$ is derived $a$-adically complete, by \Cref{lem:sq-0-trick}, it is an extension of its $a$-adic completion by an ideal of square-zero. By \cite{SP}*{Lemma~\href{https://stacks.math.columbia.edu/tag/05GG}{05GG}, Proposition~\href{https://stacks.math.columbia.edu/tag/091T}{091T}}, the $a$-adic completion of $A$ is automatically derived $a$-adically complete, so the same holds for the square-zero ideal in question. Thus, by repeating the deformation-theoretic reduction once more and using \Cref{no-nil-feelings} to retain the vanishing condition, we reduce to the case when our $R$-algebra $A$ is $a$-adically complete.

In this case, for the remaining claim about cohomological degree $i = 1$, we first show that the map 
\[
\tst R\lim_{n > 0}( R\Gamma(A/^\bL a^n, G)) \ra R\lim_{n > 0}( R\Gamma(A/(a^n), G))
\]
is an isomorphism in cohomological degree $1$. Due to \Cref{cor:nil-nil-deform-cohomology}, for this it suffices to show that
\be \label{eqn:lim1-iso}
\tst \varprojlim_{n > 0}^1 H^0(A/^\bL a^n, G) \isomto \varprojlim_{n > 0}^1 G(A/(a^n)).
\ee
\Cref{sq-zero-eg}~\ref{sq-0-3} ensures that $A/^\bL a^n$ is a square-zero extension of $A/(a^n)$ by $(A\langle a^n\rangle)[1]$, so \Cref{prop:defthy} gives the exact sequences
\[
0 \ra H^1(R\Hom_R(e^*(L_{G/R}), A\langle a^n\rangle)) \ra H^0(A/^\bL a^n, G) \ra G(A/(a^n)) \ra 0. 
\]
Since $A$ is $a$-adically and so also derived $a$-adically complete, we have 
\[
R\lim(\{A\langle a^n\rangle\}_{n > 0}) \cong 0,
\]
so
\[
R\lim(\{\Hom_A(P, A\langle a^n\rangle)\}_{n > 0}) \cong 0
\]
for every finite projective $A$-module $P$. Since $e^*(L_{G/R})$ has perfect amplitude $[-1, 0]$ (see \S\ref{def-setup}), we conclude that the systems 
\[
\tst \{H^i(R\Hom_R(e^*(L_{G/R}), A\langle a^n\rangle))\}_{n > 0}
\]
have vanishing $\varprojlim_{n\ge 0}$ and $\varprojlim_{n > 0}^1$ (to argue this concretely, one uses the $\varprojlim^1$ sequences \cite{BK72}*{Chapter IX, Propositions 2.3, Remark~2.6}). Thus, the displayed short exact sequences above give
\[
\tst \varprojlim_{n > 0} H^0(A/^\bL a^n, G) \cong \varprojlim_{n> 0} G(A/(a^n)) 
\]
and
\[
\tst \varprojlim_{n > 0}^1 H^0(A/^\bL a^n, G) \cong \varprojlim_{n> 0}^1 G(A/(a^n)),
\]
so that, in particular, \eqref{eqn:lim1-iso} is an isomorphism, as desired. 

All that remains to argue for our $a$-adically complete $R$-algebra $A$ is that the map 
\[
\tst R\Gamma(A, G) \ra R\lim_{n > 0}( R\Gamma(A/(a^n), G))
\]
is an isomorphism in cohomological degree $i = 1$, which amounts to the exactness of the sequence
\be \label{eqn:H1-seq}
\tst 0 \ra \varprojlim^1_{n > 0} G(A/(a^n)) \ra H^1(A, G) \ra \varprojlim_{n > 0} H^1(A/(a^n), G) \ra 0.
\ee
Since $A$ is $a$-adically complete and $G$ is finite, locally free, a $G$-torsor $X$ amounts to a compatible sequence of $G_{A/(a^n)}$-torsors 
\[
\tst (X_n, \iota_n \colon X_{n + 1}|_{A/(a^n)} \isomto X_n)_{n\ge 0}
\]
with specified torsor isomorphisms $\iota_n$. In particular, the third arrow is surjective and its kernel consists of the isomorphism classes of those systems for which each $X_n$ is trivial. By \cite{BK72}*{Chapter IX, Section 2.1}, the elements of $\varprojlim_{n > 0}^1 G(A/(a^n))$ are orbits of the sequences $(x_n) \in \prod_{n > 0} G(A/(a^n))$ under the action of $\prod_{n > 0} G(A/(a^n))$ given by 
\be \label{eqn:lim1-action}
(g_n) \cdot (x_n) = (g_n + x_n - g_{n + 1}|_{A/(a^n)}).
\ee
The sequence $(x_n)$ amounts to a sequence of torsor isomorphisms $\iota_n$ as above with each $X_n$ being a trivial torsor. A sequence $(g_n)$ amounts to a change of trivializations of the $X_n$, and the effect of this change on the $\iota_n$ amounts precisely to the formula \eqref{eqn:lim1-action}. From this optic, the first map of \eqref{eqn:H1-seq} is indeed the inclusion of the classes of those $(X_n, \iota_n)$ for which each $X_n$ is trivial.
\epf

The general continuity formula \eqref{eqn:genl-continuity} that we are pursuing and  fpqc hyperdescent of \Cref{thm:strongdescent} imply that fppf cohomology with commutative, finite, locally free coefficients should satisfy $a$-completely faithfully flat hyperdescent on derived $a$-adically complete animated rings. We now establish this hyperdescent directly to later use it as an input to the proof of the continuity formula.

\blem \label{lem:a-hyperdescent}
Let $R$ be a ring, let $G$ be a commutative, finite, locally free $R$-group, and let $A$ be an animated $R$-algebra. For each cosimplicial animated $A$-algebra $A^\bullet$ and each $a \in A$ such that $A$ and each $A^i$ are derived $a$-adically complete and the map $A \ra A^\bullet$ is a faithfully flat hypercover modulo $a$,
\be \label{eqn:a-hyperdescent}
\tst R\Gamma(A, G) \isomto R\lim_\Delta(R\Gamma(A^\bullet, G)). 
\ee
\elem

\bpf
Thanks to \Cref{lem:reduced-flatness}, the map $A \ra A^\bullet$ is a faithfully flat hypercover modulo every $a^n$. Moreover, by decomposing $G$ into primary factors, we assume that it is of $p$-power order for a prime $p$. By \Cref{lem:sq-0-trick} and \cite{SP}*{Lemmas~\href{https://stacks.math.columbia.edu/tag/0DYD}{0DYD} and~\href{https://stacks.math.columbia.edu/tag/09XI}{09XI}}, for an animated $A$-algebra $A'$, the $\pi_0(-)$ of the derived $a$-adic completion $(A'^h_{(p)})\wh{\ }$ of the $p$-Henselization of $A'$ (defined using \Cref{rem:no-new-etale}) is $p$-Henselian. Thus, the $p$-adic continuity formula \eqref{eqn:continuity-formula} applies and, together with the fact that the formation of the derived $a$-adic completion commutes with reduction modulo $p$, gives the identification of functors 
\[
\tst R\Gamma(((-)^h_{(p)})\wh{\ }, G) \cong R\lim_{n > 0} R\Gamma(-/^\bL p^n, G)
\]
on derived $a$-adically complete animated $A$-algebras. Each functor $R\Gamma(-/^\bL p^n, G)$ satisfies an analogue of the desired hyperdescent \eqref{eqn:a-hyperdescent} thanks to the fpqc descent of \Cref{thm:strongdescent} and to \Cref{lem:lame-case} (we recall that the vanishing condition of the latter holds for inputs whose $\pi_0(-)$ is $p$-Henselian, for instance, whose $\pi_0(-)$ is even killed by some $p^n$, see \Cref{no-higher-coho}). Thus, the functor $R\Gamma(((-)^h_{(p)})\wh{\ }, G)$ also satisfies this analogue, to the effect that it suffices to show  the analogue of \eqref{eqn:a-hyperdescent} for the following functor on derived $a$-adically complete animated $A$-algebras:
\be \label{eqn:functor-fib}
\Fib\p{R\Gamma(-, G) \ra R\Gamma(((-)^h_{(p)})\wh{\ }, G) }\!.
\ee
For this, we will use an excision trick to replace $G$ by $j_!(G)$, where 
\[
\tst j \colon \Spec(R[\f 1p]) \ra \Spec(R)
\]
is the indicated open immersion and $j_!(G)$ is taken in the \'{e}tale topology. Namely, for $a$-adically complete animated $A$-algebras $A'$, the map $A' \ra (A'^h_{(p)})\wh{\ }$ is an isomorphism modulo $p$, so the excision \Cref{prop:excision-prep} and its counterpart for \'{e}tale cohomology supplied by \cite{BM20}*{Theorems~1.15 and~5.4} imply that the following commutative squares of functors are Cartesian when evaluated on such $A'$:
\[
\xymatrix{
R\Gamma(-, G) \ar[d] \ar[r] & R\Gamma(((-)^h_{(p)})\wh{\ }, G) \ar[d] \\
R\Gamma((-)[\f 1p], G) \ar[r] & R\Gamma(((-)^h_{(p)})\wh{\ }\,[\f 1p], G),
}
\]
\[
\xymatrix{
  R\Gamma(-, j_!(G)) \ar[d] \ar[r] & R\Gamma(((-)^h_{(p)})\wh{\ }, j_!(G)) \ar[d] \\
 R\Gamma((-)[\f 1p], G) \ar[r] & R\Gamma(((-)^h_{(p)})\wh{\ }\,[\f 1p], G),
}
\]
where in the second square the cohomology is taken in the \'{e}tale topology. Thus, the functor \eqref{eqn:functor-fib} agrees with the functor
\[
\tst \Fib\p{R\Gamma(-, j_!(G)) \ra R\Gamma(((-)^h_{(p)})\wh{\ }, j_!(G)) }\!.
\]
We conclude that it suffices to show the analogue of \eqref{eqn:a-hyperdescent} for each of the functors
\[
R\Gamma(-, j_!(G)) \qxq{and} R\Gamma(((-)^h_{(p)})\wh{\ }, j_!(G))
\]
where the cohomology is \'{e}tale. By \Cref{rem:no-new-etale}, \'{e}tale cohomology depends only on the $\pi_0(-)$ of the animated ring in question. Moreover, by \Cref{lem:DC-comp-Hens} and the invariance of \'{e}tale cohomology with torsion coefficients under Henselian pairs \cite{Gab94a}*{Theorem~1}, on derived $a$-adically complete animated $A$-algebras these functors agree with the functors 
\[
R\Gamma((-)/^\bL a, j_!(G)) \qxq{and} R\Gamma(((-)/^\bL p)/^\bL a, j_!(G)) \cong R\Gamma(((-)/^\bL a)/^\bL p, j_!(G)).
\]
Thus, the desired hyperdescent for them with respect to $A \ra A^\bullet$ follows from faithfully flat hyperdescent for \'{e}tale cohomology with torsion coefficients, which itself is a special case of \Cref{arc-descent} (and the fact that faithfully flat maps are arc covers, as reviewed in \S\ref{pp:arc-topology}~\ref{arc-fflat}).
\epf

We are ready for the promised general continuity formula for flat cohomology. We thank Akhil Matthew for pointing our attention to such a statement in analogy with \cite{DM17}*{Theorem~1.5}.

\bthm \label{adic-continuity}
Let $R$ be a ring and let  $G$ be a commutative, finite, locally free $R$-group. For an animated $R$-algebra $A$ and an ideal $I = (a_1, \dotsc, a_r) \subset \pi_0(A)$ such that $A$ is derived $I$-adically complete, we have the following continuity formula\ucolon
\be \label{eqn:genl-continuity}
\tst R\Gamma(A, G) \isomto R\lim_{n > 0} (R\Gamma(A/^\bL (a_1^n, \dotsc, a_r^n), G)). 
\ee
\ethm

\bpf
The case $r = 0$ is clear, so we assume that $r > 0$. By \S\ref{pp:derived-complete} and the proof of \Cref{comp-explain}, derived $I$-adic completeness amounts to $A$ being equal to its iterated derived $a_i$-adic completion for $i = 1, \dotsc, r$. 
Moreover, each $A/^\bL (a_1^n, \dotsc, a_{r - 1}^n)$ inherits derived $a_r$-adic completeness from $A$ and 
\[
\tst R\lim_{n > 0} (A/^\bL (a_1^n, \dotsc, a_r^n)) \cong R\lim_{n > 0} (R\lim_{m \ge 0} (A/^\bL (a_1^n, \dotsc, a_{r - 1}^n, a_r^m))),
\]
and likewise after first applying $R\Gamma(-, G)$. Thus, since \eqref{eqn:continuity-condition} also holds with $a_r$ omitted, we induct on $r$ to reduce  to the case $r = 1$. From now on we place ourselves in this case and set $a \ce a_1$. 

By \Cref{thm:strongdescent} and \Cref{lem:a-hyperdescent}, both sides of \eqref{eqn:genl-continuity} satisfy hyperdescent with respect to cosimplicial animated $A$-algebras $A^\bullet$ such that each $A^i$ is derived $a$-adically complete and the map $A \ra A^\bullet$ is a faithfully flat hypercover modulo every $a^n$. In particular, \eqref{eqn:genl-continuity} holds for $A$ once it holds for each $A^i$. We use \Cref{rem:no-new-etale} and \Cref{lem:large-cover} to construct such an $A^\bullet$ for which each $A^i$ is the derived $a$-adic completion of an animated $A$-algebra $A'^i$ that has no nonsplit \'{e}tale covers. \Cref{lem:sq-0-trick} ensures that $\pi_0(A^i)$ is a square-zero extension of the $a$-adic completion of $\pi_0(A'^i)$, so \Cref{rem:no-new-etale} and \cite{SP}*{Lemmas~\href{https://stacks.math.columbia.edu/tag/09XI}{09XI} and~\href{https://stacks.math.columbia.edu/tag/04D1}{04D1}} imply that $\pi_0(A^i)$ also has no nonsplit \'{e}tale covers. Thus, by applying \Cref{rem:no-new-etale} one more time and renaming $A^i$ to $A$, we are reduced to the case when our $a$-adically complete animated $R$-algebra $A$, equivalently, $\pi_0(A)$, has no nonsplit \'{e}tale covers. In this case, the B\'{e}gueri sequence \eqref{Begueri-0} shows that $H^{j}(A, G) \cong 0$ for $j \ge 2$. By \Cref{rem:no-new-etale} and \cite{SP}*{Lemma~\href{https://stacks.math.columbia.edu/tag/04D1}{04D1}}, each $A/^\bL a^n$ also has no nonsplit \'{e}tale covers, so the B\'{e}gueri sequence and \eqref{et-fppf-agree} also show that $H^{j}(A/^\bL a^n, G) \cong 0$ for $j \ge 2$. Thus, applying \Cref{lem:lame-case} gives the conclusion.
\epf

\beg \label{eg:continuity}
Let $R$ be a ring that is derived complete (for example, complete) with respect to a finitely generated ideal $I = (a_1, \dotsc, a_r) \subset R$ and let $G$ be a commutative, finite, locally free $R$-group scheme. Since $R$ is derived $a_i$-adically complete for $i = 1, \dotsc, r$, one argues as in the beginning of the proof of \Cref{adic-continuity} that the condition \eqref{eqn:continuity-condition} holds. Thus, \eqref{eqn:genl-continuity} and \Cref{cor:nil-nil-deform-cohomology,no-nil-feelings} show that for $i \ge 2$ the map
\be \label{eqn:derived-prelim-case}
H^i(R, G) \ra H^i(R/I, G) \qxq{is} \begin{cases}\xq{surjective for} i \ge 1, \\ \xq{bijective for}\ \, i \ge 2.  \end{cases}
\ee
Similarly, for $i = 1$ they give a short exact sequence
\[
\tst 0 \ra \varprojlim^1_{n > 0} H^{0}(R/^\bL(a_1^n, \dotsc, a_r^n), G) \ra H^1(R, G) \ra \varprojlim_{n > 0} H^1(R/I^n, G) \ra 0.
\]
The animated rings $R/^\bL(a_1^n, \dotsc, a_r^n)$ are $r$-truncated, so if $R$ is $I$-adically complete, then we may argue as in the proof of \Cref{lem:lame-case} with $A\langle a^n\rangle$ replaced by the positive homotopy groups of the $R/^\bL(a_1^n, \dotsc, a_r^n)$ to inductively replace $R/^\bL(a_1^n, \dotsc, a_r^n)$ by its $j$-truncation for $j = r - 1, \dotsc, 0$ in the short exact sequence displayed above. In effect, if $R$ is $I$-adically complete, then the sequence above takes the more concrete form
\[
\tst 0 \ra \varprojlim^1_{n > 0} G(R/I^n) \ra H^1(R, G) \ra \varprojlim_{n > 0} H^1(R/I^n, G) \ra 0.
\]
\eeg


\brem
Contrary to the $p$-adic continuity formula \eqref{eqn:continuity-formula}, even in the case when the ideal $I$ in \Cref{adic-continuity} is principal, the general continuity formula \eqref{eqn:genl-continuity} does not hold if $A$ is merely $I$-Henselian (in the sense that $\pi_0(A)$ is $I$-Henselian). Indeed, if the Henselian version held, then together with the complete version \eqref{eqn:genl-continuity} it would imply that the cohomology groups $H^1(\bF_p\{t\}, \mu_p)$ and $H^1(\bF_p\llb t \rrb, \mu_p)$ are isomorphic, where $\bF_p\{t\}$ is the $t$-Henselization of $\bF_p[t]$. However, the Kummer sequence shows that $H^1(\bF_p\{t\}, \mu_p)$ is countable, whereas $H^1(\bF_p\llb t \rrb, \mu_p)$ is not.
\erem

As the following corollary shows, the insufficiency of Henselianity is a low degree phenomenon.
This complements \cite{BC19}*{Theorem~2.1.6}, which showed that for a smooth, quasi-affine group scheme $Q$, the functor $H^1(-, Q)$ is invariant under Henselian pairs. 

\bcor \lab{inv-Hens-pair}
For a Henselian pair $(R, I)$ and a commutative, finite, locally free $R$-group $G$,
\[
H^i(R, G) \ra H^i(R/I, G) \qxq{is} \begin{cases}\xq{surjective for} i \ge 1, \\ \xq{bijective for}\ \, i \ge 2.  \end{cases}
\]
\ecor

The case $G = \mu_n$ of this corollary amounts to an unpublished result of Gabber. Moreover, when $R$ is Henselian local, it continues to hold for any commutative, flat, finitely presented $R$-group algebraic space $G$, see \cite{topology-torsors}*{Proposition B.13} which essentially restates \cite{Toe11}*{corollaire 3.4}, but beyond local $R$ there are counterexamples even when $G = \bG_m$, see \cite{BC19}*{Remark 2.1.8}. 

\bpf
The surjectivity for $i = 1$ follows from \cite{BC19}*{Theorem~2.1.6 (b)} (and does not require $G$ to be commutative), so we assume that $i \ge 2$. Moreover, we use limit formalism for flat cohomology to assume that $(R, I)$ is the Henselization of a finite type $\bZ$-algebra along some ideal, so that, by \cite{SP}*{Lemma~\href{https://stacks.math.columbia.edu/tag/0AGV}{0AGV}}, the ring $R$ is Noetherian and, by \cite{SP}*{Lemmas~\href{https://stacks.math.columbia.edu/tag/0AH3}{0AH3} and~\href{https://stacks.math.columbia.edu/tag/0AH2}{0AH2}}, the fibers of the map $R \ra \wh{R}$ to the $I$-adic completion are geometrically regular. Thus, \cite{BC19}*{Lemma~2.1.3} (which is based on Popescu's theorem) allows us to assume that $R$ is even complete. This case follows from \eqref{eqn:derived-prelim-case}.
\epf

This invariance under Henselian pairs leads to the following algebraization statement for flat cohomology, which complements analogous algebraization for \'{e}tale cohomology \cite{BC19}*{Corollary~2.1.20} and for torsors under smooth, quasi-affine groups \cite{BC19}*{Corollary~2.1.22}. 

\bcor \label{cor:algebraization}
Let $B$ be a topological ring that has an open nonunital subring $B' \subset B$ such that $B'$ is Henselian as a nonunital ring and has an open neighborhood base of zero consisting of ideals of $B'$, and let $\wh{B}$ be the completion of $B$. For each commutative, finite, locally free $B$-group $G$, 
\[
H^i(B, G) \isomto H^i(\wh{B}, G) \qxq{for} i \ge 2. 
\]
\ecor

\bpf
By \Cref{inv-Hens-pair}, the axiomatic criterion \cite{BC19}*{Theorem~2.1.15} applies and gives the claim.
\epf

\beg
Letting $R\{t\}$ denote the $t$-Henselization of $R[t]$ for a ring $R$, one may choose $B \ce R\{t\}[\f 1t]$ and $B' \ce tR\{t\}$ with $B'$ equipped with its $t$-adic topology, so that $\wh{B} \cong R\llp t \rrp$. Another example is that of a Henselian pair $(A, I)$: one may choose $B \ce A$ and $B' \ce I$ with $B'$ equipped with the coarse topology, so that $\wh{B} \cong A/I$ (thus, \Cref{cor:algebraization} recovers the $i \ge 2$ case of \Cref{inv-Hens-pair}). For further examples of possible $B$ and $B'$, see \cite{BC19}*{Example~2.1.18}.
\eeg

The continuity formula also allows us to extend \Cref{heavy-handed} beyond a $p$-adic case as follows.

\bcor \label{cor:more-algebrization}
Let $R$ be a ring, let $G$ be a commutative, finite, locally free $R$-group, let $A \ra A'$ be a map  of animated $R$-algebras, and let $I \subset \pi_0(A)$ be a finitely generated ideal such that $\pi_0(A)$ is $I$-Henselian, $\pi_0(A')$ is $I(\pi_0(A'))$-Henselian, and $\pi_0(A)/I \isomto (\pi_0(A)/I) \tensor^\bL_A A'$. Letting 
\[
U \ce \Spec(\pi_0(A)) \setminus \Spec(\pi_0(A)/I)
\]
be the complement of the vanishing locus of $I$, we have
\[
H^i(R\Gamma(U_A, G)) \isomto H^i(R\Gamma(U_{A'}, G)) \qxq{for} i \ge 2. 
\]
\ecor

\bpf
The excision of \Cref{cor:excision}, the cohomology with supports sequence, and the five lemma reduce us to showing that 
\[
H^i(A, G) \isomto H^i(A', G) \qxq{for} i \ge 2.
\]
By \Cref{cor:nil-nil-deform-cohomology,inv-Hens-pair}, this map is identified with 
\[
H^i(\pi_0(A)/I, G) \ra H^i(\pi_0(A')/I\pi_0(A'), G),
\]
which is an isomorphism because, by our assumptions, even 
\[
\pi_0(A)/I \isomto \pi_0(A')/I\pi_0(A').\qedhere
\] 
\epf

We conclude this section with an algebraization result whose special case with $G = \mu_n$ was announced in \cite{Gab93}*{Theorem~2.8~(ii)}. For an argument for \cite{Gab93}*{Theorem~2.8~(i)}, see \cite{BC19}*{Corollary~2.3.5~(b)--(c)}.

\bcor\label{cor:Gab-ann}
Let $R \ra R'$ be a map  of Noetherian rings, let $G$ be a commutative, finite, locally free $R$-group, and let $I \subset R$ be an ideal such that $R$ is $I$-Henselian, $R'$ is $IR'$-Henselian, and $R/I^n \isomto R'/I^nR'$ for $n > 0$. For each open
\[
\Spec(R)\setminus \Spec(R/I) \subset U \subset \Spec(R),
\]
we have
\be \label{eqn:Gab-ann}
H^i(U, G) \isomto H^i(U_{R'}, G) \qxq{for} i \ge 2.
\ee
\ecor

\bpf
Since the map $R \ra R'$ is an isomorphism on $I$-adic completions, we lose no generality by assuming that $R'$ is the $I$-adic completion of $R$, so that the map $R \ra R'$ is flat. Due to this flatness, \Cref{cor:more-algebrization} applies and settles the case $U = \Spec(R) \setminus V(I)$. In the general case, we first note that, by \cite{BC19}*{Corollary~2.3.5~(c)}, the map \eqref{eqn:Gab-ann} is injective for $i = 2$. Thus, setting $Z \ce U \cap \Spec(R/I)$, we use the cohomology with supports sequence and the five lemma to reduce to showing that 
\be \label{eqn:Gab-ann-temp}
R\Gamma_Z(U, G) \isomto R\Gamma_Z(U_{R'}, G). 
\ee
Both sides of \eqref{eqn:Gab-ann-temp} are Zariski sheaves on $U$, so we may argue locally on $U$ to reduce to the case $U = \Spec(R)$ and then use the excision of   \Cref{cor:excision} to conclude. 
\epf


\csub[Adically faithfully flat descent for cohomology with supports] \label{sec:genl-continuity}

The continuity formula of \Cref{adic-continuity} and the fpqc descent of \Cref{thm:strongdescent} imply that flat cohomology of $I$-adically complete animated rings satisfies $I$-completely faithfully flat hyperdescent. We complement this by establishing the same for flat cohomology with supports in $I$, see \Cref{thm:adic-fpqc}. The results of this section will not be used elsewhere in this article.

\blem \label{lem:M-discrete}
Let $A$ be a perfectoid ring that is derived $I$-adically complete for a finitely generated ideal $I = (a_1, \dotsc, a_r) \subset A$ and let $M$ be a derived $I$-adically complete animated $A$-module. If each $M/^\bL(a_1^n, \dotsc, a^n_r)$ is flat as an animated $A/^\bL(a_1^n, \dotsc, a_r^n)$-module, then $M$ is $0$-truncated.
\elem

\bpf
By \Cref{lem:arc-base}, there is a $p$-complete arc hypercover $A \ra A^\bullet$ whose terms $A^i$ are perfectoid rings that are products of $p$-adically complete valuation rings of dimension $\le 1$ with algebraically closed fraction fields. By \Cref{arc-sheaves}, we have $A \isomto R\lim_\Delta A^\bullet$, so, for $n \ge 0$, also
\[
\tst A/^\bL(a_1^n, \dotsc, a^n_r) \isomto R\lim_\Delta (A^\bullet/^\bL(a_1^n, \dotsc, a^n_r)).
\]
Thus, if each $M/^\bL(a_1^n, \dotsc, a^n_r)$ is $A/^\bL(a_1^n, \dotsc, a_r^n)$-flat, as we assume from now on, then we also have
\[
\tst M/^\bL(a_1^n, \dotsc, a^n_r) \isomto R\lim_\Delta ((M \tensor_A^\bL A^\bullet)/^\bL(a_1^n, \dotsc, a^n_r)), 
\]
and so also
\[
\tst M \isomto R\lim_\Delta (M \wh{\tensor}_A^\bL A^\bullet),
\]
where $\wh{\tensor}_A^\bL$ denotes the iterated derived $a_j$-adic completion of the derived tensor product. Since $M$ is an animated $A$-module, it is connective. Thus, due to the last isomorphism, it is enough to show that each $M \wh{\tensor}_A^\bL A^i$ is $0$-truncated. The latter is an animated module over the iterated derived $a_j$-adic completion $\wh{A^i}$ of $A^i$. By the explicit nature of $A^i$, this $\wh{A^i}$ is the product of a subset of the valuation rings comprising $A^i$, in particular, $\wh{A^i}$ is also perfectoid. In conclusion, we may replace $A$ by $\wh{A^i}$ and $M$ by $M \wh{\tensor}_A^\bL A^i$, respectively, and subdivide into further subproducts if needed, to reduce to the case when the ideal $I \subset A$ is principal, generated by some $a \in A$ that has compatible $p$-power roots in $A$.

In this case, we set $M' \ce \varprojlim_{n > 0} (M/a^n M)$ and we claim that the canonical map $M \ra M'$ is an isomorphism---this will imply that $M$ is $0$-truncated, as desired. Both $M$ and $M'$ are derived $a$-adically complete (see \S\ref{pp:derived-complete}), so it suffices to check that $M/^\bL a \isomto M'/^\bL a$. \Cref{lem:reduced-flatness} and the assumption on $M$ imply that $M/^\bL a$ is $A/^\bL a$-flat, and \cite{Yek18}*{Theorem~2.8} implies that the map in question is an isomorphism on $\pi_0(-)$. Thus, by \S\ref{pp:flat-modules}, it suffices to check that $M'/^\bL a$ is $A/^\bL a$-flat or, by \Cref{lem:reduced-flatness} again, that $M' \tensor^\bL_A A/a$ is $0$-truncated and $A/a$-flat. This, however, is a special case of \cite{Yek18}*{Theorem~6.9 (with Theorem~4.3)} (to apply \emph{loc.~cit.},~we note that, in the terminology there, the ideal $(a) \subset A$ is weakly proregular because $A\langle a \rangle = A\langle a^\infty \rangle$ by \eqref{no-big-tor} above).
\epf

\bthm \label{thm:adic-fpqc}
Let $R$ be a ring, let $G$ be a commutative, finite, locally free $R$-group, let $A$ be an animated $R$-algebra, and let $I = (a_1, \dotsc, a_r) \subset \pi_0(A)$ be  an ideal. For each cosimplicial animated $A$-algebra $A^\bullet$ such that 
the map $A/^\bL(a_1, \dotsc, a_r) \ra A^\bullet/^\bL(a_1, \dotsc, a_r)$ is a faithfully flat hypercover, 
\be \label{eqn:adic-fpqc}
\tst R\Gamma_I(A, G) \isomto R\lim_\Delta(R\Gamma_I(A^\bullet, G)).
\ee
If $A$ and the terms of $A^\bullet$ are derived $I$-adically complete, then, letting the complement of the vanishing locus of $I$ be $U \subset \Spec(\pi_0(A))$, equivalently, we have
\be \label{eqn:adic-U-fpqc}
\tst R\Gamma(A, G) \isomto R\lim_\Delta(R\Gamma(A^\bullet, G))
\ee
and
\be \label{eqn:adic-U-fpqc-2} 
\tst R\Gamma(U_A, G) \isomto R\lim_\Delta(R\Gamma(U_{A^\bullet}, G)).
\ee
\ethm

\bpf
By the continuity formula of \Cref{adic-continuity} and the fpqc descent of \Cref{thm:strongdescent}, if $A$ and the terms of $A^\bullet$ are derived $I$-adically complete, then
\[
\tst R\Gamma(A, G) \isomto R\lim_\Delta(R\Gamma(A^\bullet, G)).
\]
In particular, the cohomology with supports triangle then implies that \eqref{eqn:adic-fpqc} and \eqref{eqn:adic-U-fpqc-2} are equivalent. In general, by excision \Cref{cor:excision}, the claimed \eqref{eqn:adic-fpqc} is insensitive to replacing $A$ and the terms of $A^\bullet$ by their iterated $a_i$-adic completions for $i = 1, \dotsc, r$, so we lose no generality by assuming that $A$ and the terms of $A^\bullet$ are all derived $I$-adically complete.

By decomposing into primary factors, we may assume that $G$ is of $p$-power order for a prime $p$. If $p$ is invertible in $\pi_0(A)$, then $G$ is \'{e}tale over $R$, so, by \eqref{et-fppf-agree} and \Cref{rem:no-new-etale}, we have
\[
R\Gamma_I(A, G) \cong R\Gamma_I(\pi_0(A), G) \qxq{and} R\Gamma_I(A^\bullet, G) \cong R\Gamma_I(\pi_0(A^\bullet), G). 
\]
By \Cref{lem:DC-comp-Hens}, the rings $\pi_0(A)$ and $\pi_0(A^\bullet)$ are $I$-Henselian. Moreover, by \S\ref{pp:arc-topology}~\ref{arc-pi-fflat} and \Cref{lem:reduced-flatness}, the map $\pi_0(A) \ra \pi_0(A^\bullet)$ is an $I$-complete arc hypercover. Thus, in the case when $p$ is invertible in $\pi_0(A)$, the claim follows from $I$-complete hyperdescent for \'{e}tale cohomology, more precisely, from \Cref{BM-hens-comp}. In general, this case shows that the third term of the triangle 
\[
\tst R\Gamma_{I + (p)}(A, G) \ra R\Gamma_I(A, G) \ra R\Gamma_I(A[\f1p], G)
\]
satisfies the analogue of \eqref{eqn:adic-fpqc}. Thus, we may replace $I$ by $I + (p)$ to henceforth assume that~$p \in I$.

Our next goal is to reduce to the case when $A$ is an animated $\bF_p$-algebra, and for this, as in the proof of \Cref{lem:a-hyperdescent}, we will use the excision trick of replacing $G$ by the extension by zero $j_!(G)$ taken in the \'{e}tale topology, where 
\[
\tst j \colon \Spec(R[\f 1p]) \ra \Spec(R)
\]
is the indicated open immersion. Namely, letting $A'$ range over those animated $A$-algebras fppf over $A$ for which $\Spec(\pi_0(A')) \ra \Spec(\pi_0(A))$ factors over $U$ (compare with \S\ref{pp:simplicial-coho-def}), by the excision \Cref{prop:excision-prep} and its counterpart for \'{e}tale cohomology supplied by \cite{BM20}*{Theorems~1.15 and~5.4}, we have the Cartesian squares
\[
\xymatrix{
R\Gamma(U_A, G) \ar[d] \ar[r] & R\lim_{A'}R\Gamma(A^{\prime h}_{(p)}, G) \ar[d]   \\
R\Gamma(A[\f 1p], G) \ar[r] & R\lim_{A'}R\Gamma(A^{\prime h}_{(p)}[\f 1p], G), 
}
\]
\[
\xymatrix{
R\Gamma(U_A, j_!(G)) \ar[d] \ar[r] & R\lim_{A'}R\Gamma(A^{\prime h}_{(p)}, j_!(G)) \ar[d] \\
R\Gamma(A[\f 1p], G) \ar[r] & R\lim_{A'}R\Gamma(A^{\prime h}_{(p)}[\f 1p], G),
}
\]
where $(-)^h_{(p)}$ denotes the $p$-Henselization and in the second square the cohomology is taken in the \'{e}tale topology. Of course, we also have the analogous squares for the terms of $A^\bullet$, and we claim that the common fiber of the horizontal maps in the two squares above satisfies hyperdescent with respect to $A \ra A^\bullet$, that is, that it satisfies the analogues of \eqref{eqn:adic-fpqc}, \eqref{eqn:adic-U-fpqc}, and \eqref{eqn:adic-U-fpqc}. For this, it suffices to show the same claim for the terms $R\Gamma(U_A, j_!(G))$ and $R\lim_{A'}R\Gamma(A^{\prime h}_{(p)}, j_!(G))$. For $R\Gamma(U_A, j_!(G))$, this again follows from $I$-complete hyperdescent for \'{e}tale cohomology, that is, from \Cref{BM-hens-comp}. On the other hand, by invariance of \'{e}tale cohomology under Henselian pairs \cite{Gab94a}*{Theorem~1},
\[
\tst R\lim_{A'}R\Gamma(A^{\prime h}_{(p)}, j_!(G)) \cong R\Gamma(U_{A'}/^\bL p, j_!(G)) \cong R\Gamma(U_{A'/^\bL p}, j_!(G)).
\]
Thus, $I$-complete arc hyperdescent for \'{e}tale cohomology discussed in \Cref{BM-hens-comp} also handles the terms $R\lim_{A'}R\Gamma(A^{\prime h}_{(p)}, j_!(G))$. In conclusion, the common fibers of the horizontal maps do indeed satisfy the analogues of \eqref{eqn:adic-fpqc}, \eqref{eqn:adic-U-fpqc}, and \eqref{eqn:adic-U-fpqc-2}. By inspecting the top horizontal map of the left square, this means that our overall desired conclusion reduces to its analogue for the term $R\lim_{A'}R\Gamma(A^{\prime h}_{(p)}, G)$. By the $p$-adic continuity formula of \Cref{thm:padiclimit}, this term is identified with 
\[
\tst R\lim_{n > 0} R\Gamma(U_{A/^\bL p^n}, G).
\]
In effect, it suffices to show that $R\Gamma(U_{A/^\bL p^n}, G) \isomto R\Gamma(U_{A^\bullet/^\bL p^n}, G)$ for every $n > 0$, in other words, we are allowed to replace $A$ by $A/^\bL p^n$ in the overall claim we are seeking to prove.

To reduce further to the desired $n = 1$, we now establish insensitivity to square-zero extensions: for a square-zero extension $A \ra \ov{A}$ by an animated $\ov{A}$-module $M$, the desired \eqref{eqn:adic-fpqc} holds if and only if it holds after base change to $\ov{A}$. For this, by excision of \Cref{cor:excision}, the claim is insensitive to replacing $A$ and $\ov{A}$ by their iterated $a_i$-adic completions for $i = 1, \dotsc, r$ and $A^\bullet$ by its corresponding base changes. Thus, since these completions form a square-zero extension by the corresponding completion of $M$ (see \S\ref{pp:sq-zero}), we lose no generality by assuming that $A$, $\ov{A}$, and $M$ are all derived $I$-adically complete. Fpqc hyperdescent for modules, as discussed at the end of \S\ref{pp:simplicial-coho-def}, then gives
\[
\tst M \isomto R\lim_\Delta(M \tensor^\bL_A A^\bullet).
\]
Thus, the deformation-theoretic triangle of \Cref{prop:defthy} and its counterparts after base change to $A^\bullet$ give the claimed insensitivity to square-zero extensions.

In the rest of the proof we focus on the remaining case when $A$ is an animated $\bF_p$-algebra. Moreover, Postnikov completeness of \Cref{cor:fppf-hyperdescent} allows us to replace $A$ and $A^\bullet$ by $\tau_{\le n}(A)$ and $\tau_{\le n} \tensor^\bL_A A^\bullet$, respectively,  for a variable $n$, so 
we may assume that $A$ is $n$-truncated. By then iteratively combining the insensitivity to square-zero extensions with \Cref{sq-zero-eg}~\ref{sq-0-3}, we reduce to $0$-truncated $A$. Once $A$ is a $0$-truncated $\bF_p$-algebra, we consider the ind-fppf, faithfully flat, semiperfect $A$-algebra
\[
A_\infty \ce A[X_a^{1/p^\infty}\, |\, a \in A]/(X_a - a\, |\, a \in A).
\]
The terms of the \v{C}ech nerve of the ind-fppf cover $A \ra A_\infty$ are all semiperfect $\bF_p$-algebras, so, by ind-fppf descent for flat cohomology, we lose no generality by replacing the hypercover $A \ra A^\bullet$ by its bases changes to these terms to reduce to the case when $A$ is a semiperfect $\bF_p$-algebra. By iteratively using the insensitivity to square-zero extensions and passing to a filtered direct limit, we may then replace $A$ by $A^\red$ and $A^\bullet$ by its base change to $A^\red$ to reduce further to a perfect $\bF_p$-algebra $A$. 

By excision \Cref{cor:excision} as before, we may replace $A$ and the $A^j$ by their iterated derived $a_i$-adic completions for $i = 1, \dotsc, r$ to assume that $A$ and the $A^j$ are derived $I$-adically complete: \Cref{recognize}~\ref{R-e} ensures that the resulting $A$ is still a perfect $\bF_p$-algebra. By \Cref{lem:M-discrete}, the $A^j$ are then also $0$-truncated, to the effect that the functor $R\Gamma_I(-, G)$ takes coconnective values on them. Thus, it suffices to show descent instead of hyperdescent, more precisely, we lose no generality by assuming that $A \ra A^\bullet$ is of \v{C}ech type, associated to a map of $\bF_p$-algebras $A \ra A'$ such that $A/^\bL (a_1, \dotsc, a_r) \ra A'/^\bL (a_1, \dotsc, a_r)$ is faithfully flat. This property is preserved by the preceding reductions, so we repeat them once more to again reduce to $A$ being perfect (this time, to preserve the \v{C}ech property, we only $I$-adically complete $A$ and $A'$ and not the other terms of $A^\bullet$). 

Once $A \ra A^\bullet$ is of \v{C}ech type with a perfect $\bF_p$-algebra $A$ and a $0$-truncated $A'$, we claim that the desired descent \eqref{eqn:adic-fpqc} holds even after replacing $A \ra A^\bullet$ by its base change to any animated $A$-algebra $B$. The advantage of this claim is that, by \cite{LZ17}*{Lemma~3.1.2 (3)},  it is insensitive to replacing $A \ra A'$ by a refinement $A \ra A' \ra A''$. Thus, we let $A'' \ce \varinjlim_{a \mapsto a^p} A'$ be the perfection of $A'$: since $A$ is perfect, the map $A/^\bL (a_1, \dotsc, a_r) \ra A'' /^\bL (a_1, \dotsc, a_r)$ is still faithfully flat. In effect, we may assume that both $A$ and $A'$ are perfect at the cost of having to show \eqref{eqn:adic-fpqc} after base change to any animated $A$-algebra $B$. We then repeat the preceding reductions for this base change to $B$ to reduce further to the case when $B$ is also a perfect $\bF_p$-algebra. Since $A$ and $A'$ are also perfect, by \cite{BS17}*{Proposition~11.6}, the $B \tensor^\bL_A A^j$ are then also perfect $\bF_p$-algebras. 

In conclusion, we are left with the case when the hypercover $A \ra A^\bullet$ is of \v{C}ech type and only involves perfect $\bF_p$-algebras. As before, we pass to iterated derived $a_i$-adic completions of its terms and use excision \Cref{cor:excision} to arrange that $A$ and all the $A^j$ are derived $I$-adically complete (and still perfect by \Cref{recognize}~\ref{R-e}). 
\Cref{KT-input} then reduces us to showing that
\[
\tst R\Gamma_{(p,\, I)}(W(A), \bM(G))^{V = 1} \isomto R\lim_{\Delta} (R\Gamma_{(p,\, I)}(W(A^\bullet), \bM(G_{A^\bullet}))^{V = 1}),
\]
where we abusively write $(p, I)$ for the ideal $J \subset W(A)$ generated by $p$ and the Teichm\"{u}llers of the $a_i$. For this, the mapping fiber triangle allows us to remove the superscripts $(-)^{V = 1}$. Similarly, since the crystalline Dieudonn\'{e} module $\bM(G)$ is of projective dimension $\le 1$ over $W(A)$ and its formation commutes with base change (see \S\ref{Gab-eq}), it suffices to show that 
\[
\tst R\Gamma_{J}(W(A), W(A)) \isomto R\lim_{\Delta} (R\Gamma_{J}(W(A^\bullet), W(A^\bullet))).
\]
By \cite{SP}*{Lemma~\href{https://stacks.math.columbia.edu/tag/0954}{0954}}, we have 
\[
\tst R\Gamma_{J}(W(A), W(A)) \cong \varinjlim_{n > 0} R\Hom_{W(A)}(W(A)/J^n, W(A)),
\]
and likewise for $W(A^\bullet)$, so it suffices to show that 
\[
\tst R\Hom(W(A)/J^n, W(A)) \isomto R\lim_{\Delta}(R\Hom(W(A^\bullet)/J^nW(A^\bullet), W(A^\bullet))),
\]
where the $R\Hom$ are over $W(A)$ and $W(A^\bullet)$, respectively, and the maps are induced by base change. Since $p$ is a nonzerodivisor in $W(A)$ and $W(A^j)$, our assumption about faithful flatness modulo $(a_1, \dotsc, a_n)$ and \Cref{lem:reduced-flatness} ensure that 
\[
W(A)/J^n \tensor^\bL_{W(A)} W(A^\bullet) \cong W(A^\bullet)/J^nW(A^\bullet).
\]
Consequently, the target of the last displayed map is identified with
\[
\tst 
R\lim_{\Delta}(R\Hom(W(A)/J^n, W(A^\bullet))) \cong R\Hom(W(A)/J^n, R\lim_{\Delta}W(A^\bullet)),
\]
where the $R\Hom$ are over $W(A)$. Moreover, by faithfully flat descent and derived $J$-adic completeness, we have 
\[
\tst W(A)/^\bL(p, a_1, \dotsc, a_r) \isomto R\lim_{\Delta}(W(A^\bullet)/^\bL(p, a_1, \dotsc, a_r)),
\]
so also
\[
\tst W(A) \isomto R\lim_{\Delta} W(A^\bullet).
\]
By combining this isomorphism with the aforementioned identification of the target, we obtain the conclusion.
\epf






\section{The characteristic-primary aspects of the main result} \label{bad-characteristics}

We have gathered all the ingredients we need to exhibit purity for flat cohomology. 
 In \S\ref{perfectoid-version}, we establish the general case of the key formula \eqref{eqn:key-formula} and obtain the perfectoid version of flat purity. In \S\ref{grand-finale}, we then deduce the remaining ``bad residue characteristic'' cases of our main purity results.



\csub[The key formula and purity for flat cohomology of perfectoids] \label{perfectoid-version}

We are ready for the key formula that relates flat cohomology of a perfectoid ring to quasi-coherent cohomology of its $\bA_\Inf$ with values in prismatic Dieudonn\'{e} modules $\bM(G)$ reviewed in \S\ref{pp:prismatic-review}. 




\begin{theorem}\label{thm:perfectoidsupport}  
For a prime $p$, a perfectoid ring $A$, a commutative, finite, locally free, $A$-group $G$ of $p$-power order, and a closed $Z\subset \Spec (A/pA)$, we have a functorial in $A$, $G$, and $Z$ identification
\be \label{eqn:key-formula-pf}
R\Gamma_Z(A,G)\cong R\Gamma_{Z}(\bA_\Inf(A), \bM(G))^{V = 1}.
\ee
\end{theorem}

Here we choose the same $\xi$ when defining $V$ over perfectoid $A$-algebras, see \uS\uref{pp:prismatic-review} and \uS\uref{perfectoid-def}.

\begin{proof}
We may assume that $A$ is a $\bZ_p$-algebra and, by passing to an inverse limit in the end if necessary, we may assume that $Z$ is the vanishing locus of a finite number of elements of $A$. 

Let us begin with the case when $A = \prod_{i \in I} A_i$ for perfectoid valuation rings $A_i$ of rank $\le 1$ that have algebraically closed fraction fields. For such $A$, our closed subset $Z$ is cut out by a single $a \in A$ with $a \mid p$. 
We choose compatible $p$-power roots $a^\flat \in A^\flat$ of $a$, so that $a^\flat \mid p^\flat$ and $a^\flat$ cuts out $Z \subset \Spec(A^\flat/p^\flat A^\flat)$ (see \S\ref{perfectoid-def}, especially \eqref{monoid-id}--\eqref{same-mod-pi}, as well as \Cref{tilt-summary}). By \Cref{M(G)-local-analysis}, we have
\[
\tst H^i(A,G) \cong H^i(A[\frac 1a],G) \cong 0 \qxq{for} i \ge 1,
\]
the map $V - 1$ is surjective on $\bM(G)$ and $\bM(G)[\f 1{[a^\flat]}]$, 
and there is a unique commutative square
\begin{equation*}
\tst \xymatrix@C=40pt{
G(A) \ar@{->>}[d]\ar[r]^-{\eqref{MG-fixed-pts}}_-{\sim} & \bM(G)^{V = 1} \ar@{->>}[d] \\
G(A[\f 1{a}]) \ar[r]^-{\sim} & (\bM(G)[\f 1{[a^\flat]}])^{V = 1}.
}
\end{equation*}
Thus, by the cohomology with supports sequence, 
\[
R\Gamma_Z(A, G) \qxq{and} R\Gamma_{Z}(\bA_\Inf(A), \bM(G))^{V = 1}
\]
are concentrated in degree $0$ and identified. Due to the functoriality of the isomorphism \eqref{MG-fixed-pts} and the uniqueness of the above diagram, this identification is functorial in $A$, $Z$, and $G$, as desired.

For general $A$, by \Cref{thm:strongdescent-perfectoid}, the left side of \eqref{eqn:key-formula-pf} satisfies hyperdescent for those $p$-complete arc hypercovers whose terms are perfectoid $\bZ_p$-algebras. By \Cref{M(G)-coho-descent}, so does the right side. To then deduce the general case from the already settled case of $\prod_{i \in I} A_i$ as above, it remains to recall from \Cref{lem:arc-base} that such products form a base of the $p$-complete arc topology of $A$.
\end{proof}

With the key formula in hand, a similar argument to the one we used in positive characteristic at the end of the proof of \Cref{thm:main-pos-char} now gives purity for flat cohomology of perfectoid rings.

\bthm\label{thm:perfectoid-purity}
For a prime $p$, a perfectoid ring $A$, a commutative, finite, locally free $A$-group $G$ of $p$-power order, a closed $Z \subset \Spec(A/pA)$, and a regular sequence $a_1, \dotsc, a_d \in A$ that vanishes on~$Z$,
\[
H^i_Z(A, G) \cong 0 \qxq{for} i < d.
\]
\ethm

\bpf
By \eqref{eqn:key-formula-pf} and a long exact cohomology sequence, it suffices to show that 
\be \label{eqn:pp-coho-van}
H^i_Z(\bA_\Inf(A), \bM(G)) \cong 0 \qxq{for} i < d. 
\ee
By \S\ref{pp:prismatic-review}, the $\bA_\Inf(A)$-module $\bM(G)$ is finitely presented of projective dimension $\le 1$, so
\[
H^i_Z(\bA_\Inf(A), \bA_\Inf(A)) \cong 0 \qxq{for} i < d + 1
\]
would suffice. For this, we use the $\bA_\Inf(A)$-regular sequence $a_0 \ce \xi, a_1, \dotsc, a_d$, where $\xi$ is a generator of $\Ker(\theta\colon \bA_\Inf(A) \surjects A)$ (see \S\ref{perfectoid-def}). Namely, since the $a_i$ vanish on $Z$ and $H^i_Z(\bA_\Inf(A), M)$ is supported on $Z$ for every $\bA_\Inf(A)$-module $M$, decreasing induction on $-1 \le j \le d$ gives the sufficient
\[
H^i_Z(\bA_\Inf(A), \bA_\Inf(A)/(a_0, \dotsc, a_j)) \cong 0 \qxq{for} i < d - j. \qedhere
\]
\epf


\csub[Purity for flat cohomology of local complete intersections]\label{grand-finale}


We conclude the proof of purity for flat cohomology by reducing to its settled perfectoid case in \Cref{main-pf}.
The following lemmas help to pass to 
completions in the appearing perfectoid towers.

\blem[\cite{Yek18}*{Theorems 1.3 and 1.5}
] \label{flat-lemma}
Let $I \subset R$ be an ideal in a Noetherian ring $R$. Every $I$-adically complete $R$-module $M$ that is $I$-completely flat \up{meaning that $M\otimes^{\bL}_R R/I$ is concentrated in degree $0$ and $R/I$-flat} is flat. In particular, the $I$-adic completion of a flat $R$-module is flat. \QED
\elem

\blem \label{slip-in-complete}
Let $A$ be a ring, let $a \in A$, and let $\wh{A}$ be the $a$-adic completion. Each $A$-regular sequence $a_1, \dotsc, a_n \in A$ such that $A/(a_1, \dotsc, a_i)$ has bounded $a^\infty$-torsion for $0 \le i \le n$ is $\wh{A}$-regular~and
\[
\wh{A}/(a_1, \dotsc, a_n) \isomto (A/(a_1, \dotsc, a_n))\wh{\ }.
\]
\elem

\bpf
The bounded torsion assumption implies that the derived $a$-adic completions of the short exact sequences 
\[
0 \ra A/(a_1, \dotsc, a_{i - 1}) \xra{a_i}  A/(a_1, \dotsc, a_{i - 1}) \ra  A/(a_1, \dotsc, a_{i}) \ra 0 
\]
for $1 \le i \le n$ agree with their classical $a$-adic completions. In particular, we obtain short exact sequences 
\[
0 \ra (A/(a_1, \dotsc, a_{i - 1}))\wh{\ } \xra{a_i}  (A/(a_1, \dotsc, a_{i - 1}))\wh{\ } \ra  (A/(a_1, \dotsc, a_{i}))\wh{\ } \ra 0,
\]
which show the claim.
\epf

\begin{theorem} \label{main-pf}
For a Noetherian local ring $(R,\mathfrak m)$ that is a complete intersection and a commutative, finite, flat $R$-group $G$,
\[
H^i_{\mathfrak m}(R,G) \cong 0 \qxq{for} i < \dim(R).
\]
\end{theorem}

\begin{proof} 
By decomposing into primary factors, we may assume that $G$ is of $p$-power order for a prime $p$. \Cref{ell-semipurity} settles the case when $p$ is invertible in $R$, so we assume that $p = \Char(R/\fm)$. Moreover, by \Cref{reduce-to-complete}, we may assume that $R$ is $\fm$-adically complete, so that there is an unramified, complete, regular, local ring $(\wt{R}, \wt{\fm})$ and a regular sequence $f_1, \dotsc, f_n \in \wt{\fm}$ such that
\[
R \simeq \wt{R}/(f_1, \dotsc, f_n)
\]
(see \S\ref{conv}
). We will argue the desired vanishing by induction on $i$ for all $R$ at once. 

We use \Cref{lem:towers}~\ref{tower-a} to find a filtered direct system of regular, local, finite, flat $\wt{R}$-algebras $\wt{R}_j$ with $\varinjlim_j \wt{R}_j$ a regular local ring with an algebraically closed residue field. By the inductive assumption, \Cref{finite-cover}, 
and a limit argument, we may replace $R$ by $(\varinjlim_j \wt{R}_j)/(f_1, \dotsc, f_n)$ and then apply \Cref{reduce-to-complete} again to assume that $\wt{R}$ has an algebraically closed residue field. Once this is arranged, the passage to a tower argument carried out with \Cref{lem:towers}~\ref{tower-b} instead supplies a faithfully flat 
$\wt{R}$-algebra $\wt{R}_\infty$ whose $p$-adic completion is perfectoid 
for which we need to show that
\[
H^i_\fm(\wt{R}_\infty/(f_1, \dotsc, f_n), G) \cong 0 \qxq{for} i < \dim(R). 
\]
By 
\Cref{heavy-handed} 
and \Cref{slip-in-complete}, we may replace $\wt{R}_\infty$ by its perfectoid 
$p$-adic completion. Thus, it suffices to show that for any $p$-torsion free perfectoid ring $A$ that is $p$-completely faithfully flat over $\wt{R}$ in the sense that $A/p^nA$ is faithfully flat over $\wt{R}/p^n\wt{R}$ for $n > 0$,
\be \label{ultimate-desire}
H^i_\fm(A/(f_1, \dotsc, f_n), G) \cong 0 \qxq{for} i < \dim(R). 
\ee
By \Cref{flat-lemma}, such an $A$ is even $\wt{R}$-flat, so the sequence $f_1,\ldots,f_n$ is $A$-regular. By Andr\'{e}'s lemma, that is, by \Cref{Andre-lem}, there is an ind-syntomic, faithfully flat $A$-algebra $A'$ whose $p$-adic completion $\wh{A'}$ is perfectoid and contains compatible $p$-power roots $f_i^{1/p^\infty}$ for $i = 1, \dotsc, n$. A limit argument then gives the spectral sequence
\[
E_1^{st} = H^t_\fm((\underbrace{A'\otimes_A \ldots\otimes_A A'}_{s + 1})/(f_1, \dotsc, f_n),G) \Rightarrow H^{s + t}_\fm(A/(f_1, \dotsc, f_n), G).
\]
By \Cref{heavy-handed} and \Cref{slip-in-complete} again (with elementary excision as in \cref{foot:Milne} to $p$-Henselize the tensor products),  we may replace the $A'\otimes_A \ldots\otimes_A A'$ by their $p$-adic completions, which are perfectoid by \Cref{recognize}~\ref{R-a}. Thus, the spectral sequence above allows us to assume that our perfectoid $A$ in \eqref{ultimate-desire} contains compatible $p$-power roots $f_i^{1/p^\infty}$.\footnote{Another way to carry out this reduction is to use the version \cite{BS19}*{Theorem~7.14
} of Andr\'{e}'s lemma. Then $A'$ is only $p$-completely faithfully flat over $A$ but is a perfectoid right away and contains compatible $p$-power roots $f_i^{1/p^\infty}$. One then combines \Cref{cor:descentsupport} and \Cref{cor:excision} (with \Cref{slip-in-complete} again) to obtain the spectral sequence with $p$-adic completions already inside. This avoids \Cref{Andre-lem} at the cost of relying on heavier inputs from Chapter \ref{reduction}.} 

For every $R$-regular sequence $r_1, \dotsc, r_{\dim(R)} \in \wt{\fm}$, the sequence $f_1, \dotsc, f_n$, $r_1, \dotsc, r_{\dim(R)}$ is $A$-regular, so \Cref{mod-out-by-more} together with a limit argument reduces us to showing that 
\be \label{main-pf-eq-4}
H^i_\fm(A/(f_1^{1/p^\infty}, \dotsc, f_n^{1/p^\infty}), G) \cong 0 \q \x{for} \q i < \dim(R).
\ee
By the $\wt{R}$-flatness of $A$ and \cite{SP}*{Lemma~\href{https://stacks.math.columbia.edu/tag/07DV}{07DV}}, for all $m_1, \dotsc, m_n \ge 0$, every permutation of the sequence $f_1^{1/p^{m_1}}, \dotsc, f_n^{1/p^{m_n}}, r_1, \dotsc, r_{\dim(R)}$ is $A$-regular. Thus, by induction on the number of nonzero exponents $m_\ell$, the $A$-module 
\[
\tst A/(f_1^{1/p^{m_1}}, \dotsc, f_n^{1/p^{m_n}}, r_1, \dotsc, r_{j})
\]
is isomorphic to a submodule of $A/(f_1, \dotsc, f_n, r_1, \dotsc, r_j)$ for $j = 0, \dotsc, \dim(R)$. In particular, by the $\wt{R}$-flatness of $A$, its $p^\infty$-torsion is killed by $p^N$ for some fixed $N > 0$ that does not depend on the $m_\ell$ or on $j$. By forming colimits, this $p^N$ then kills every $(A/(f_1^{1/p^\infty}, \dotsc, f_n^{1/p^\infty}, r_1, \dotsc, r_j))\langle p^\infty\rangle$, so \Cref{slip-in-complete} ensures that $r_1, \dotsc, r_{\dim(R)}$ is still a regular sequence in the $p$-adic completion $(A/(f_1^{1/p^\infty}, \dotsc, f_n^{1/p^\infty}))\wh{\ }$. However, by \Cref{recognize}~\ref{R-b}, the latter is perfectoid, so \Cref{thm:perfectoid-purity} gives
\[
H^i_\fm((A/(f_1^{1/p^\infty}, \dotsc, f_n^{1/p^\infty}))\wh{\ }, G) \cong 0 \q \x{for} \q i < \dim(R).
\]
By \Cref{heavy-handed}, 
this vanishing gives the desired \eqref{main-pf-eq-4}.
\end{proof}

For \'{e}tale $G$, the variant of purity that involves the virtual dimension (defined in \S\ref{geom-depth-def}) follows, too.

\bthm \lab{ell-semipurity-p}
For a Noetherian local ring $(R, \fm)$ and a commutative, finite, \'{e}tale $R$-group $G$,
\[
H^i_\fm(R, G) \cong 0 \q \text{for} \q i < \vdim(R).
\]
\ethm

\bpf
\Cref{ell-semipurity} and its proof settle the case when the order of $G$ is invertible in $R$ and reduce the rest to the case when 
\[
G = \bZ/p\bZ \qxq{with} p = \Char(R/\fm) > 0
\]
and $R$ is a quotient of a regular local ring by a principal ideal. Such an $R$ is a complete intersection for which \Cref{main-pf} gives
\[
H^i_\fm(R, \bZ/p\bZ) \cong 0 \qxq{for} i < \vdim(R) \overset{\eqref{vdim-ineq}}{=} \dim(R). \qedhere
\]
\epf

\brem
To avoid repetitiveness, we deduced \Cref{ell-semipurity-p} from \Cref{main-pf}, although the proof of the latter simplifies significantly for $G = \bZ/p\bZ$. For example, to pass to completions in this case, we may replace Corollaries \ref{heavy-handed} and \ref{reduce-to-complete} by the simpler \Cref{lem:et-to-hat}. Moreover, since the \'{e}tale site is insensitive to nilpotents, there is no need to appeal to \Cref{mod-out-by-more} when reducing to \eqref{main-pf-eq-4}. Finally, there is no need to refer to \Cref{thm:perfectoid-purity} in the end: \Cref{baby-tilting} directly reduces to positive characteristic granted that one uses \cite{GR18}*{Proposition 16.4.17} to transfer depth. 
\erem

\brem
In contrast to \Cref{ell-semipurity-p}, in \Cref{main-pf} we cannot drop the complete intersection assumption and replace $\dim(R)$ by $\vdim(R)$: for example, for $R \ce \bF_p\llb x, y, z, t\rrb/(x^2, y^2, xz - yt)$, the nonzero element $xy \in R$ dies on $U_R \ce \Spec(R)\setminus \{ \fm\}$, so it is nonzero in $H^0_\fm(R, \gA_p)$. 
\erem

We close the section with a slight sharpening of \Cref{main-pf} in the case when $R$ is regular. 

\bthm \label{thm:main-pf-regular}
For a regular local ring $(R, \fm)$ that is not a field and a commutative, finite, flat $R$-group $G$,
\[
H^i_\fm(R, G) \cong 0 \qxq{for} i \le \dim(R).
\]
\ethm

\bpf
\Cref{main-pf} gives the vanishing for $i < \dim(R)$, so we may focus on the cohomological degree $i = \dim(R)$. Moreover, by decomposing $G$ into primary pieces and using \Cref{ell-semipurity}, we may assume that $G$ is of $p$-power order with $p = \Char(R/\fm) > 0$. We then use \cite{BBM82}*{th\'{e}or\`{e}me~3.1.1} to embed $G$ into a truncated $p$-divisible group and combine the resulting cohomology sequence with \Cref{main-pf} to reduce to $G$ itself being a truncated $p$-divisible group. The filtration by $p$-power torsion then allows us to assume that, in addition, $G$ is killed by $p$.

As in the proof of \Cref{main-pf}, we may assume that $R$ is $\fm$-adically complete. As there, we then use \Cref{lem:towers}~\ref{tower-a}, \Cref{finite-cover}, and \Cref{main-pf} to assume, in addition, that the residue field $k \ce R/\fm$ is algebraically closed. As in \Cref{lem:towers}~\ref{tower-b}, the Cohen theorem then shows that
\[
R \simeq W(k)\llb x_1, \dotsc, x_d\rrb/(p - f) \qxq{with either} f = x_1 \qxq{or} f \in (p, x_1, \dotsc, x_n)^2,
\]
and, since $R$ is not a field, $d > 0$. We then analogously use \Cref{finite-cover} and \Cref{main-pf} to pass to the tower supplied by \Cref{lem:towers}~\ref{tower-b}, and hence to reduce to showing that
\[
H^d_{\fm_\infty}(R_\infty, G) \cong 0,
\]
where
\[
R_\infty \simeq W(k)\llb x_1^{1/p^\infty}, \dotsc, x_d^{1/p^\infty}\rrb/(p - f), \q \fm_\infty \ce (p, x_1, \dotsc, x_d).
\]
Moreover, for showing this vanishing, \Cref{heavy-handed} allows us to replace $R_\infty$ by its $p$-adic completion $\wh{R}_\infty$, which is perfectoid. By \Cref{flat-lemma}, the sequence $x_1, \dotsc, x_d$ is $\wh{R}_\infty$-regular, so the key formula \eqref{eqn:key-formula-pf} and the vanishing \eqref{eqn:pp-coho-van} reduce us to showing that 
\[
H^d_{Z}(\bA_\Inf(\wh{R}_\infty), \bM(G))^{V = 1} \cong 0, 
\]
where $Z \subset \Spec(\bA_\Inf(\wh{R}_\infty))$ is the closed point. Since $G$ is the $p$-torsion of a $p$-divisible group and $\wh{R}_\infty^\flat$ is local, $\bM(G)$ is a finite free $\wh{R}_\infty^\flat$-module equipped with a $\Frob\i$-semilinear map $V \colon \bM(G) \ra \bM(G)$, so our task is to show that $V$ has no nonzero fixed points on $H^d_{Z}(\wh{R}_\infty^\flat, \bM(G))$. For this, we first describe $\wh{R}_\infty^\flat$-module $H^d_Z(\wh{R}_\infty^\flat, \wh{R}_\infty^\flat)$.

As we saw in the proof of \Cref{abs-coh-pur}, the tilt $\wh{R}_\infty^\flat$ is the $\ov{f}$-adic completion of $k\llb (x_1^\flat)^{1/p^\infty}, \dotsc, (x_d^\flat)^{1/p^\infty} \rrb$ for some $\ov{f} \in (x_1^\flat, \dotsc, x_d^\flat)$. By \Cref{flat-lemma}, the sequence $x_1^\flat, \dotsc, x_d^\flat$ is $\wh{R}_\infty^\flat$-regular. In particular, similarly to the proof of \Cref{thm:perfectoid-purity}, from the exact sequences 
\[
0 \ra \wh{R}_\infty^\flat/(x_1^\flat, \dotsc, x_{j - 1}^\flat) \xra{x_j^\flat} \wh{R}_\infty^\flat/(x_1^\flat, \dotsc, x_{j - 1}^\flat) \ra \wh{R}_\infty^\flat/(x_1^\flat, \dotsc, x_{j}^\flat) \ra 0
\]
we get $H^j_Z(\wh{R}_\infty^\flat, \wh{R}_\infty^\flat) \cong 0$ for $j < d$ and, letting the transition maps be the indicated multiplications, we have
\[
\ba
\tst H^d_Z(\wh{R}_\infty^\flat, \wh{R}_\infty^\flat) &\cong H^d_Z(\wh{R}_\infty^\flat, \wh{R}_\infty^\flat) \langle (x_1^\flat)^\infty \rangle \\ &\tst\cong \varinjlim_{x_1^\flat} H^{d - 1}_Z(\wh{R}_\infty^\flat, \wh{R}_\infty^\flat/(x_1^\flat)^n) \cong H^{d - 1}_Z(\wh{R}_\infty^\flat, \wh{R}_\infty^\flat[\f 1{x_1^\flat}]/\wh{R}_\infty^\flat).
\ea\]
Continuing in this way, since $k\llb (x_1^\flat)^{1/p^\infty}, \dotsc, (x_d^\flat)^{1/p^\infty} \rrb$ and its $\ov{f}$-adic completion agree modulo each $((x_1^\flat)^{n_1}, \dotsc, (x_d^\flat)^{n_d})$, we find that $H^d_Z(\wh{R}_\infty^\flat, \wh{R}_\infty^\flat)$ agrees with its analogue for $k\llb (x_1^\flat)^{1/p^\infty}, \dotsc, (x_d^\flat)^{1/p^\infty} \rrb$ and that, concretely, it is given by the quotient of $(k\llb (x_1^\flat)^{1/p^\infty}, \dotsc, (x_d^\flat)^{1/p^\infty} \rrb)[\f 1{x_1\cdots x_d}]$ by the space of those elements whose monomials have at least one nonnegative exponent. Thus, we may identify $H^d_Z(\wh{R}_\infty^\flat, \wh{R}_\infty^\flat)$ with the $k$-vector space with the basis 
\[
\tst\{(x_1^\flat)^{a_1} \cdots (x_d^\flat)^{a_d}\}_{a_1,\, \dotsc,\, a_d\, \in\, \bZ[\f 1p]_{< 0}}.
\] 
Since $\bM(G)$ is a finite free $\wh{R}_\infty^\flat$-module, $H^d_{Z}(\wh{R}_\infty^\flat, \bM(G))$ is a finite direct sum of such $k$-vector spaces. For any hypothetical nonzero element fixed by $V$, we choose a monomial $(x_1^\flat)^{a_1} \cdots (x_d^\flat)^{a_d}$ appearing in it for which the sum $a_1 + \dotsc + a_d$ is the smallest. The effect of $V$ is described by some matrix with coefficients in $k\llb (x_1^\flat)^{1/p^\infty}, \dotsc, (x_d^\flat)^{1/p^\infty} \rrb$ postcomposed with $\Frob\i$, and latter divides each $a_i$ by $p$, so the sum $a_1 + \dotsc + a_d$ strictly increases after applying $V$---more informally, $V$ is ``contracting'' on $H^d_{Z}(\wh{R}_\infty^\flat, \bM(G))$. Consequently, $V$ has no nonzero fixed points on $H^d_{Z}(\wh{R}_\infty^\flat, \bM(G))$, as desired.
\epf


\section{Global purity consequences and the conjectures of Gabber} \label{global}

Our final goal is to deduce global purity consequences from the local \Cref{main-pf} 
and to settle the conjectures of Gabber, as announced in \Cref{thm:Gabber-conjectures}. In particular, we extend purity for the Brauer group settled for regular schemes in \cite{brauer-purity} to the case of complete intersection singularities. 

\csub[Cohomology with finite flat group scheme coefficients]

We begin with the most straight-forward global consequence of \Cref{main-pf}: in \Cref{P-d-pf} we show that on Noetherian schemes with complete intersection singularities, flat cohomology classes with values in commutative, finite, flat group schemes are insensitive to removing closed subschemes of sufficiently large codimension. The reduction of this statement to its local case uses the following concrete manifestation of \'{e}tale descent for fppf cohomology with supports.  

\blem \label{lem-sup-ss}
For a scheme $X$, a closed subset $Z \subset X$, an abelian fppf sheaf $\sF$ on $X$, and the \'{e}tale sheafification $\cH^j_Z(-, \sF)$ of the functor $X' \mapsto H^j_Z(X', \sF)$, there is a 
spectral sequence
\[
E_2^{ij} = H^i_\et(X, \cH^j_Z(-, \sF)) \Ra H^{i + j}_Z(X, \sF). 
\]
\elem

\bpf
One way to show this is by considering injective resolutions, see, for instance, \cite{BC19}*{Lemma~2.3.2}.
\epf

\bthm \label{P-d-pf}
Let $X$ be a scheme, let $G$ be a commutative, finite, locally free $X$-group, and let $Z \subset X$ be a closed subset such that the immersion $X \setminus Z \hra X$ is quasi-compact and each $\sO_{X,\, z}$ with $z \in Z$ is either a complete intersection of dimension $\ge d$ or regular of dimension $\ge d - 1$. The map
\[
H^i(X, G) \ra H^i(X \setminus Z, G)  \q\x{is} \q  \begin{cases}\x{injective for} \q  i < d, \\ \x{bijective for} \q i < d - 1.\end{cases}
\]
\ethm

\bpf
We will deduce the assertion from the local purity Theorems \ref{main-pf} and \ref{thm:main-pf-regular}. One may wish to compare the method of this deduction to \cite{Gab04a}*{Lemma 3.1 and its proof}.

The assertion amounts to the vanishing $H^i_Z(X, G) \cong 0$ for $i < d$. Thus, the  spectral sequence
\[
E_2^{ij} = H^i_\et(X, \cH^j_Z(-, G)) \Ra H^{i + j}_Z(X, G)
\]
of \Cref{lem-sup-ss} reduces us to the case when $X$ is strictly Henselian and $Z \neq \emptyset$ (the quasi-compactness assumption is used in this step to identify the stalks of $\cH^j_Z(-, G)$ via limit formalism, compare with the proof of \cite{brauer-purity}*{Theorem 6.1}). Then $X$ is Noetherian and we will show how to shrink $Z$ to arrive by Noetherian induction at the case $Z = \{\fm\}$ supplied by Theorems \ref{main-pf} and \ref{thm:main-pf-regular}.

Suppose that $Z \neq \{\fm\}$, fix a generic point $z$ of $Z$, and let $U$ be an open neighborhood of $z$ in $X$. The \v{C}ech-to-derived spectral sequence  
\[
\widecheck{H}^p(\{U, X \setminus Z\}, H^q(-, G)) \Ra H^{p+q}(U \cup (X \setminus Z), G)
\]
(a concrete incarnation of Zariski descent for fppf cohomology) gives the Mayer--Vietoris sequence
\[
\dotsc \ra H^i(U \cup (X \setminus Z), G) \ra H^i(U, G) \oplus H^i(X \setminus Z, G) \ra H^i(U \cap (X \setminus Z), G) \ra \dotsc
\]
As $U$ shrinks, it becomes $\Spec(\sO_{X,\, z})$ and $U \cap (X \setminus Z)$ becomes $\Spec(\sO_{X,\, z}) \setminus \{ z\}$, so, by Theorems \ref{main-pf} and \ref{thm:main-pf-regular}, in this limit the maps 
\[
H^i(U, G) \ra H^i(U \cap (X \setminus Z), G) \qxq{become} \begin{cases} \x{injective for}  &i = d - 1, \\ \x{bijective for} & i < d - 1.
\end{cases}
\]
Consequently, for $i < d$, the sequence implies that any $\gA \in H^i(X, G)$  that dies in $H^i(X \setminus Z, G)$ also dies in $H^i(U \cup (X \setminus Z), G)$ for sufficiently small $U$. Likewise, for $i < d - 1$, any $\gB \in H^i(X \setminus Z, G)$ extends to $H^i(U \cup (X \setminus Z), G)$ for some such $U$. This allows us to apply the inductive hypothesis to $Z' \ce Z \cap (X \setminus U) \subsetneq Z$ to conclude.
\epf

Moret-Bailly has a nonabelian version of \Cref{P-d-pf} in \cite{MB85b}*{lemme~2}. For completeness, we include its very mild generalization whose argument is independent of the rest of this article and builds on the omitted one for \emph{loc.~cit.}~(that was explained in \cite{Mar16}*{Chapter~3}). 

\bthm\label{MB-lem}
Let $X$ be a scheme, let $G$ be a finite, locally free $X$-group, and let $Z \subset X$ be a closed subset such that the immersion $j \colon X \setminus Z \hra X$ is quasi-compact and each $\sO_{X,\, z}$ with $z \in Z$ is regular of dimension $\ge 2$. Pullback is an equivalence from the category of $G$-torsors to that of $G_{X\setminus Z}$-torsors. More generally, for any $G$-gerbe $\sB$, the following pullback is an equivalence of~categories\ucolon
\[
\sB(X) \isomto \sB(X\setminus Z), \qxq{so also} \Ker(H^2(X, G) \ra H^2(X\setminus Z, G)) = \{ *\}.
\]
\ethm

\bpf
The claim about the nonabelian $H^2$ amounts to $\sB(X)$ being nonempty whenever so is $\sB(X \setminus Z)$, so it follows from the claim about $\sB$. As for the latter, the classifying stack $\bbB G$ is smooth, so a general $G$-gerbe $\sB$ becomes isomorphic to $\bbB G$ over an \'{e}tale cover of $X$, to the effect that \'{e}tale descent reduces us to the case $\sB = \bbB G$. Thus, we may and will focus on the claim about $G$-torsors. 

By glueing, the assertion is Zariski-local on $X$. Thus, by localizing at a point of $Z$ and  spreading out (which uses the quasi-compactness of $j$), we may assume that $X$ is local and $Z \neq \emptyset$. Then $X$ is Noetherian, so \cite{EGAIV2}*{th\'{e}or\`{e}me 5.10.5} gives the full faithfulness because $\depth_Z(X) \ge 2$. For the 
essential surjectivity, by Noetherian induction and spreading out, we may localize at a generic point of $Z$ to assume that $Z$ is the closed point of $X$. By \cite{EGAIV2}*{corollaire 5.11.4}, for a $G_{X \setminus Z}$-torsor $Y$, the $\sO_X$-algebra $j_*(\sO_Y)$ is coherent, so all we need to show is its flatness as an $\sO_X$-module: the proof of full faithfulness will then uniquely extend the torsor structure map of $Y$ to that of $\un{\Spec}_X(j_*(\sO_Y))$. 

For the remaining $\sO_X$-flatness of $j_*(\sO_Y)$, by a result of Auslander \cite{Aus62}*{Theorem~1.3} that crucially uses the regularity of $X$, it suffices to show that 
\[
\sH om_{\sO_X}(j_*(\sO_Y), j_*(\sO_Y)) \simeq (j_*(\sO_Y))^{\oplus r} \q \x{as $\sO_X$-modules.}
\]
It suffices to argue this over $X\setminus Z$ (see \cite{Grothendieck-Lefschetz}*{Lemma~2.2}), so, since $\sO_G$ is $\sO_X$-free, it suffices to show that
\[
\sH om_{\sO_{X \setminus Z}}(\sO_{G_{X\setminus Z}}, \sO_Y)  \isomto  \sH om_{\sO_{X \setminus Z}}(\sO_Y, \sO_Y)
\]
via
\[
\tst f \mapsto \p{\sO_Y \xra{a} \sO_{G_{X\setminus Z}} \tensor_{\sO_{X\setminus Z}} \sO_Y \xra{(f,\, \id)}  \sO_Y }\!,
\]
where $a$ is the $G$-action morphism. The explicit inverse of this $\sO_{X \setminus Z}$-module homomorphism is
\[
g \mapsto \p{\sO_{G_{X\setminus Z}} \xra{\id  \tensor\, 1} \sO_{G_{X\setminus Z}} \tensor_{\sO_{X\setminus Z}} \sO_Y \xra{(a,\, \id)\i} \sO_Y \tensor_{\sO_{X\setminus Z}} \sO_Y \xra{(g,\, \id)}  \sO_Y }\!. \qedhere
\]
\epf

\brem\label{nonab-rem}
We expect that \Cref{MB-lem} also holds when each $\sO_{X,\, z}$ with $z \in Z$ is either   a complete intersection of dimension $\ge 3$ or regular of dimension $\ge 2$ (compare with \Cref{P-d-pf}). Unfortunately, the argument given above, especially, \cite{Aus62}*{Theorem~1.3}, is specific to regular~$\sO_{X,\, z}$.  
\erem





\csub[The conjectures of Gabber and purity for the Brauer group] \label{passage-section}

We are ready to settle Gabber's conjecture \cite{Gab04}*{Conjecture~3} on the local Picard group.


\bthm
\label{deg-2-case}
For a Noetherian local ring $(R, \fm)$ that is a complete intersection of dimension~$\ge 3$,
\[
\Pic(U_R)_\tors \cong 0, \qxq{where} U_R \ce \Spec(R) \setminus \{ \fm\}.
\]
If $R$ is either of dimension $\ge 4$ or regular of any dimension, then even 
\[
\Pic(U_R) \cong 0.
\]
\ethm

\bpf
The assertion about the case $\dim(R) \ge 4$ was settled in \cite{SGA2new}*{expos\'{e}~XI, th\'{e}or\`{e}me~3.13~(ii)}. Moreover, a line bundle $\sL$ on $U_R$ is trivial if and only if it extends to a line bundle on $R$. For regular $R$, one constructs such an extension either by considering Weil divisors or by first extending $\sL$ as a coherent module and then taking the determinant of a perfect complex representing this module.

For the remaining assertion about $\Pic(U_R)_\tors$, \Cref{main-pf} implies the bijectivity of the left vertical map in the commutative diagram
\[
\xymatrix@R=15pt@C=15pt{
H^1(R, \mu_n) \ar[d]^-{\sim} \ar[r] & H^1(R, \bG_m) \ar[d] & \!\!\!\!\!\!\!\!\!\!\!\!\cong \Pic(R) \cong 0 \\
H^1(U_R, \mu_n) \ar[r] & H^1(U_R, \bG_m) & \!\!\!\!\!\!\!\!\!\!\!\!\!\!\! \cong \Pic(U_R).
}
\]
Since every element of $\Pic(U_R)_\tors$ comes from $H^1(U_R, \mu_n)$ for some $n \ge 0$, the claim follows.
\epf

\brem \label{no-free-lunch}
\Cref{deg-2-case} (so also \Cref{main}) does not hold if $(R, \fm)$ is merely Cohen--Macaulay. For instance, consider the normal local domain
\[
R \ce (\bC\llb x_1, \dotsc, x_n\rrb)^{\bZ/2\bZ}
\]
where the $\bZ/2\bZ$-action is given by $x_i \mapsto -x_i$ for $1 \le i \le n$. 
The map 
\[
R \ra \bC\llb x_1, \dotsc, x_n\rrb
\]
is the normalization in a quadratic extension of the fraction field, is a nontrivial $\bZ/2\bZ$-torsor away from the maximal ideal, and is an inclusion of an $R$-module direct summand (with the antiinvariants as a complementary summand). In particular, a system of parameters for $R$ is also one for $\bC\llb x_1, \dotsc, x_n\rrb$, and $R$ inherits Cohen--Macaulayness from $\bC\llb x_1, \dotsc, x_n\rrb$. 
Thus, for $n \ge 2$, we have $R^\times \isomto H^0(U_R, \bG_m)$ and the torsor gives a nonzero element of $\Pic(U_R)[2]$.
\erem

\bcor \label{global-ci-cor}
For a field $k$ and a global complete intersection $X \subset \bP^n_k$ of dimension $\ge 2$,  
\[
(\Pic(X)/(\bZ \cdot [\sO(1)]))_\tors \cong 0.
\]
If $X$ is of dimension $\ge 3$, then $\Pic(X)$ is even free, generated by $[\sO(1)]$.
\ecor

In the case when $X$ is smooth this corollary was established by Deligne in \cite{SGA7II}*{expos\'{e}~XI, th\'{e}or\`{e}me~1.8}.

\bpf
The assertion about the case $\dim(X) \ge 3$ was settled in \cite{SGA2new}*{expos\'{e} XII, corollaire~3.7}. 
To deduce the rest from \Cref{deg-2-case} we will pass to the affine cone of $X$. Namely, as in \cite{Grothendieck-Lefschetz}*{proof of Theorem~4.1}, the scheme $X$ is the $\Proj$ of a graded $k$-algebra $R\ce k[x_0, \dotsc, x_n]/(f_1, \dotsc, f_{n - d})$ for homogeneous elements $f_1, \dotsc, f_{n - d} \in k[x_0, \dotsc, x_n]$ that form a $k[x_0, \dotsc, x_n]$-regular sequence, and 
\[
\tst R \cong \bigoplus_{m \ge 0} \Gamma(X, \sO_X(m))
\]
compatibly with the gradings. As there, a line bundle $\sL$ on $X$ defines a finite, graded $R$-module 
\[
\tst M_\sL \ce \bigoplus_{m \ge 0} \Gamma(X, \sL(m))
\] 
whose restriction to $\Spec(R) \setminus \{ \fm\}$ with $\fm \ce (x_0, \dotsc, x_n)$ is a line bundle $\wt{\sL}$. As there, $M_\sL$ is the pushforward of $\wt{\sL}$, and if $\wt{\sL}|_{\Spec(R_\fm)\setminus \{\fm\}}$ is free, then so is $M_\sL$, 
in which case $\sL \simeq \sO_X(m)$. Moreover, by \cite{EGAII}*{th\'{e}or\`{e}me~3.4.4}, the $\sO_X$-module associated to $M_\sL$ is $\sL$, so \cite{EGAII}*{propositions 3.2.4, 3.2.6 et 3.4.3} show that $\wt{\sL} \tensor \wt{\sL'} \isomto (\sL \tensor_{\sO_X} \sL')\wt{\ \,}$. It follows that 
\[
\tst \Pic(X)/(\bZ \cdot [\sO(1)]) \hra \Pic(\Spec(R_\fm) \setminus \{ \fm\}),
\]
so \Cref{deg-2-case} gives the claim. 
\epf

\brem\label{eg:from-Hartshorne}
\Cref{global-ci-cor} (so also \Cref{deg-2-case}) 
is sharp: indeed, one cannot drop $(-)_\tors$ because the complete intersection 
\[
X\ce \Proj(k[ x, y, z, w]/(xw - yz))
\]
of dimension $2$ (the Segre embedding of $\bP^1_k \times_k \bP^1_k$) satisfies $\Pic(X) \simeq \bZ \oplus \bZ$, 
and one cannot weaken then dimension assumption because an elliptic curve $E$  over an algebraically closed field has $\#(\Pic(E)_\tors) = \infty$.
\erem



We turn to Gabber's \cite{Gab04}*{Conjecture~2} on the Brauer group of local complete intersections.

\bthm
\label{deg-3-case}
For a Noetherian local ring $(R, \fm)$ that is either a complete intersection of dimension $\ge 4$ or regular of dimension $\ge 2$,
\[
\Br(R) \isomto \Br(U_R), \qxq{where} U_R \ce \Spec(R) \setminus \{ \fm\}. 
\]
\ethm

\bpf
We recall that the Brauer group $\Br(X)$ of a scheme $X$ is defined using Azumaya algebras:
\[
\tst \Br(X) \ce \bigcup_{n > 0} \im\p{H^1(X, \PGL_n) \ra H^2(X, \bG_m)_\tors}\!.
\]
By a result of Gabber and de Jong \cite{dJ02} (see also \cite{CTS21}*{Section 4.2}), we have 
\[
\Br(X) \cong H^2(X, \bG_m)_\tors
\]
whenever $X$ has an ample line bundle, for instance, whenever $X$ is quasi-affine. In our setting, 
$\Pic(R) \cong \Pic(U_R) \cong 0$ (see \Cref{deg-2-case}
), so 
\[
\Br(R)[n] \cong H^2(R, \mu_n) \qxq{and} \Br(U_R)[n] \cong H^2(U_R, \mu_n) \qxq{for} n \ge 0.
\]
Thus, except for the case when $R$ is regular of dimension $2$, the desired conclusion follows by noting that $H^2(R, \mu_n) \isomto H^2(U_R, \mu_n)$ thanks to Theorems \ref{main-pf} and \ref{thm:main-pf-regular}. The just excluded case is actually the most basic and was treated in \cite{Gro68c}*{th\'{e}or\`{e}me~6.1~b)}: in this case, the desired conclusion follows by considering Azumaya algebras and noting that pullback gives an equivalence between the category of vector bundles (resp.,~Azumaya algebras) on $R$ and those on $U_R$.
\epf


\brem
One cannot weaken the dimension assumption of \Cref{deg-3-case}. Indeed, let $S$ be the local ring at the vertex of the affine cone over an elliptic curve over $\bC$, so that $S$ is a $2$-dimensional, normal, Noetherian, local $\bQ$-algebra that is a complete intersection with $\#(\Pic(U_S)_\tors) = \infty$ (see \Cref{eg:from-Hartshorne}). We have  $\Pic(U_S) \hra \Pic(U_{S^{sh}})$ because for any $\sL$ in the kernel, $\Gamma(U_S, \sL)$ is a free $S$-module. Thus, for the $3$-dimensional, strictly Henselian, complete intersection $R \ce S^{sh}\llb x, y\rrb/(xy)$, since $\Pic(U_{S^\sh\llb t \rrb})$ is finitely generated (see, for instance, \cite{Bou78}*{chapitre~V, corollaire~4.9}), the short exact sequence
\[
0 \ra \bG_m \ra (i_{x = 0})_*(\bG_m) \oplus (i_{y = 0})_*(\bG_m) \xra{(r_1,\, r_2)\, \mapsto\, r_1/r_2} (i_{x = y = 0})_*(\bG_m) \ra 0
\]
on $\Spec(R)_\et$ (compare with \cite{Bou78}*{chapitre~IV, lemme~5.1}) shows that $\#(\Br(U_R)) = \infty$, whereas  $\Br(R) \cong 0$. The same reasoning carried out with the Segre embedding of $\bP^1_\bC \times_\bC \bP^1_\bC$ in place of an elliptic curve shows that in \Cref{deg-3-case} the full $H^2(U_R, \bG_m)$ may contain classes that do not come from $H^2(R, \bG_m)$.
\erem


To establish purity for the Brauer group of local complete intersections, we globalize \Cref{deg-3-case} in \Cref{P-b-pf} below. For this, we use the following version of Hartogs' extension principle.


\blem \label{loc-sections}
Let $X$ be a scheme and let $Z \subset X$ be a closed subset. 
\benum
\item \label{LS-a}
If each $\fm_{X,\, z}$ with $z \in Z$ contains an $\sO_{X,\, z}$-regular element, then 
\[
\qq Y(X) \hra Y(X \setminus Z) \q \x{for every separated $X$-scheme $Y$.}
\]

\item \label{LS-b}
If each $\fm_{X,\, z}$ with $z \in Z$ contains an $\sO_{X,\, z}$-regular sequence of length $2$, then
\[
\qq Y(X) \isomto Y(X \setminus Z) \q \x{for every $X$-affine $X$-scheme $Y$.}
\]
\eenum
\elem

\bpf \hfill
\benum
\item
Since $Y$ is separated, its diagonal is a closed immersion. Thus, it suffices to check  that no nonzero local section $f$ of $\sO_X$ vanishes away from $Z$. 
By shrinking $X$, we assume that $f$ is global and let $X \setminus Z \subset U$ be the maximal open on which it vanishes. If $U \neq X$, then we choose a generic point $z$ of $X \setminus U$ and a nonzerodivisor $m \in \fm_{X,\, z}$ to see that $\sO_{X,\, z} \hra \sO_{X,\, z}[\f{1}{m}]$. Since $f$ vanishes in $\sO_{X,\, z}[\f{1}{m}]$, it also vanishes in $\sO_{X,\, z}$, a contradiction.

\item
By \ref{LS-a}, the map is injective, so it suffices to show that every section of $Y$ over $X \setminus Z$ extends (necessarily uniquely) to a section over $X$. Moreover, by working locally on $X$, we may assume that $X$ is affine. We then embed $Y$ into a (possibly infinite dimensional) affine space over $X$ and use \ref{LS-a} to reduce to the case when $Y = \bA^1_X$. In other words, we have reduced to showing that every global section of $X \setminus Z$ extends (necessarily uniquely) to a global section of $X$. 

By glueing, there is the largest open $X \setminus Z \subset U$ such that the global section of $X \setminus Z$ in question extends to a global section of $U$. To show that the inclusion $U \subset X$ is not strict, we suppose otherwise and fix a generic point $z$ of $X \setminus U$. A limit argument reduces us to showing that
\[
\qq \sO_{X,\, z} \isomto \Gamma(\Spec(\sO_{X,\, z}) \setminus \{ z\}, \sO_X).
\]
For Noetherian $\sO_{X,\, z}$, this follows from \cite{EGAIV2}*{th\'{e}or\`{e}me 5.10.5}, and in general we fix an $\sO_{X,\, z}$-regular sequence $m_1, m_2 \in \fm_{X,\, z}$ and seek to show that the complex
\[
\tst \qq \sO_{X,\, z} \hra \sO_{X,\, z}[\f{1}{m_1}] \oplus \sO_{X,\, z}[\f{1}{m_2}] \xra{(a,\, b)\, \mapsto\, a - b} \sO_{X,\, z}[\f{1}{m_1m_2}]
\]
is exact in the middle. This complex is a filtered direct limit of Koszul complexes $K(m_1^n, m_2^n)$ (see \cite{SP}*{Lemma~\href{https://stacks.math.columbia.edu/tag/0913}{0913}}), so it suffices to show that the latter, considered as chain complexes in degrees between $0$ and $2$, have vanishing homology in degree $1$. The sequence $m_1^n, m_2^n$ inherits $\sO_{X,\, z}$-regularity (see \cite{SP}*{Lemma~\href{https://stacks.math.columbia.edu/tag/07DV}{07DV}}), so this vanishing follows from the fact that if $m_1^n a = m_2^n b$ in $\sO_{X,\, z}$, then $b = m_1^n c$ for some $c \in \sO_{X,\, z}$ for which also $a = m_2^n c$.
 \qedhere
\eenum
\epf

\bthm \label{P-b-pf} \label{Brauer-purity}
Let $X$ be a scheme, let $T$ be a finite type $X$-group of multiplicative type, and let $Z \subset X$ be a closed subset such that the open immersion $j \colon X \setminus Z \hra X$ is quasi-compact.
\benum
\item \label{PbP-a} 
If each local ring $\sO_{X,\, z}$ for $z \in Z$ is Noetherian and geometrically parafactorial,\footnote{We recall from \cite{EGAIV4}*{d\'{e}finition 21.13.7} that a local ring $(R, \fm)$ is \emph{parafactorial} if pullback is an equivalence from the category of line bundles on $\Spec(R)$ to those on $\Spec(R) \setminus \{ \fm\}$. A local ring is \emph{geometrically parafactorial} if its strict Henselization is parafactorial. For example, by \cite{EGAIV4}*{exemple 21.13.9 (ii)}, every Noetherian, local, geometrically factorial (in the sense that the strict Henselization is factorial) ring of dimension $\ge 2$ is geometrically parafactorial and, by \cite{SGA2new}*{expos\'{e} XI, th\'{e}or\`{e}me 3.13 (ii)}, so is every local complete intersection of dimension $\ge 4$.} then 
\[\ba
\qqq H^0(X, T) \isomto &H^0(X \setminus Z, T), \q H^1(X, T) \isomto H^1(X \setminus Z, T), 
\\  &H^2(X, T) \hra H^2(X \setminus Z, T). 
\ea
\]

\item \label{PbP-b}
If each local ring $\sO_{X,\, z}$ for $z \in Z$ is either a complete intersection of dimension $\ge 3$ \up{resp.,~$\ge 4$} or  regular of dimension $\ge 2$, then\footnote{Here $(-)_\tors$ denotes classes killed by a locally constant function. For instance, an $\gA \in H^i(X, T)$ lies in $H^i(X, T)_\tors$ if and only if there is a decomposition $X = \bigsqcup_{n \in \bZ_{> 0}} X_n $ into clopens such that each $\gA|_{X_n}$ is killed by $n$.}
\[
\qqq H^1(X, T)_\tors \isomto H^1(X \setminus Z, T)_\tors 
\]
{\upshape(}resp., 
\[
\qqq H^2(X, T)_\tors \isomto H^2(X \setminus Z, T)_\tors).
\]
\eenum
\ethm

The importance of geometric parafactoriality for the $H^2$ aspect of \ref{PbP-a} was noticed in \cite{Str79}*{Teorema~4}. 

\bpf \hfill
\benum
\item
We need to show that $H^i_Z(X, T) \cong 0$ for $i \le 2$. By \cite{EGAIV4}*{proposition~21.13.8}, for each $z \in Z$  we have $\depth(\sO_{X,\, z}) \ge 2$, so \Cref{loc-sections}~\ref{LS-b} gives the $i \le 1$ part of this vanishing. Moreover, as in the proof of \Cref{P-d-pf}, \Cref{lem-sup-ss} reduces us to the case when $X$ is strictly local and $Z \neq \emptyset$. Then $X$ is Noetherian and, by realizing $T$ as the kernel of a morphism between tori, we may assume that $T = \bG_m$. This turns our task into showing that for every line bundle $\sL$ on $X \setminus Z$, the pushforward $j_*(\sL)$ is also line bundle. For the latter, we argue by Noetherian induction, so, since the formation of $j_*(\sL)$ commutes with flat base change, we replace $X$ by its strict Henselization at a generic point of $Z$ to assume that $Z$ is the closed point. In this case, the parafactoriality assumption shows that $j_*(\sL)$ is a line bundle. 

\item
The injectivity follows from \Cref{loc-sections}~\ref{LS-b} (resp.,~from \ref{PbP-a}), which also shows that the clopens of $X$ and $X \setminus Z$ correspond. Thus, for the surjectivity, we may assume that the cohomology class in question is killed by some $n > 0$ and then that $T$ is $n$-torsion. This removes the subscripts `$\tors$,' so \Cref{lem-sup-ss} (with \Cref{loc-sections}~\ref{LS-a} for the vanishing of $\cH^0_Z$) allows us to assume that $X$ is strictly Henselian and $Z \neq \emptyset$. 
Then $X$ is Noetherian and we need to extend a cohomology class on $X\setminus Z$ to $X$. For this, Noetherian induction, limit arguments, and the Mayer--Vietoris sequence used in the proof of \Cref{P-d-pf} allow us to replace $X$ by its localization at a generic point of $Z$. Then $X$ is local and $Z$ is the closed point, so, except for the case when $X$ is regular of dimension $2$, 
Theorems \ref{main-pf} and \ref{thm:main-pf-regular} give the extension. For $H^1$ (resp.,~$H^2$), the remaining case is supplied by \ref{PbP-a} (resp.,~by \Cref{deg-3-case}).
\qedhere
\eenum
\epf

Under  more restrictive assumptions, \Cref{Brauer-purity} extends to higher degree cohomology as follows.

\bthm \label{P-c-pf}
Let $X$ be a scheme, let $T$ be a finite type $X$-group of multiplicative type, let $d \ge 3$ be an integer, and let $Z \subset X$ be a closed subset  such that the open immersion $j \colon X \setminus Z \hra X$ is quasi-compact and each $\sO_{X,\, z}$ for $z \in Z$ either is a complete intersection of dimension $\ge d$ all of whose strict Henselizations are factorial\footnote{To illustrate the assumption, we recall that every regular local ring is factorial and, from \cite{SGA2new}*{expos\'{e} XI, corollaire 3.14}, that every complete intersection whose local rings of height $\le 3$  are factorial (for instance, regular) is factorial.} or is regular of dimension $\ge d - 1$. The map
\[
H^i(X, T) \ra H^i(X \setminus Z, T)  \q\x{is} \q  \begin{cases}\x{injective for} \q  i < d, \\ \x{bijective for} \q i < d - 1.\end{cases}
\]
\ethm

\bpf
We need to show that $H^i_Z(X, T) \cong 0$ for $i < d$, and \Cref{lem-sup-ss} reduces us to $X$ being strictly local with $Z \neq \emptyset$. 
Then $X$ is Noetherian, integral, and we may assume that $T = \bG_m$. By \cite{Gro68b}*{proposition~1.4}, which uses the factoriality assumption, $H^i(X, \bG_m)$ and $H^i(X\setminus Z, \bG_m)$ are torsion for $i \ge 2$. Thus, since 
\[
H^i(X, \bG_m) \isomto H^i(X\setminus Z, \bG_m) \qxq{for} i \le 1
\]
by \Cref{Brauer-purity}~\ref{PbP-a}, all the $H^{i}_Z(X, \bG_m)$ are also torsion. The vanishing $H^i_Z(X, \mu_n) \cong 0$ for $i < d$ supplied by \Cref{P-d-pf} then suffices.
\epf



\end{document}